\documentstyle{amsart}

\newtheorem{theorem}{Theorem}[section]
\newtheorem{lemma}[theorem]{Lemma}
\newtheorem{proposition}[theorem]{Proposition}

\newtheorem{defini}[theorem]{Definition}
\newenvironment{definition}{\begin{defini}\rm}{\end{defini}}

\newtheorem{rem}[theorem]{\it Remark}
\newenvironment{remark}{\begin{rem}\rm}{\end{rem}}

\numberwithin{equation}{section}
\numberwithin{table}{section}
\numberwithin{figure}{section}

\newenvironment{proof}{\begin{pf}}{\end{pf}}

\renewcommand{\P}{\mathord{\mathbb  P}}
\newcommand{\Q}{\mathord{\mathbb  Q}}
\newcommand{\R}{\mathord{\mathbb R}}
\newcommand{\Z}{\mathord{\mathbb Z}}

\newcommand{\Zp}{\Z\sb p}
\newcommand{\Qp}{\Q\sb p}
\newcommand{\Zt}{\Z\sb 2}
\newcommand{\Qt}{\Q\sb 2}

\newcommand{\st}{\subset}

\newcommand{\utsq}[1]{(#1\sp\times)\sp 2}

\newcommand{\opp}[1]{#1\sp{-}}

\font\got=eufm10
\def\Sym{\hbox{\got S}}
\def\UpperHalfPlane{\hbox{\got H}}

\newcommand{\LLL}{\mathord{\mathcal L}}

\newcommand{\PPP}{\mathord{\mathcal P}}
\newcommand{\RRR}{\mathord{\mathcal R}}

\newcommand{\EEE}{\mathord{\mathcal E}}
\newcommand{\SSS}{\mathord{\mathcal S}}
\newcommand{\VVV}{\mathord{\mathcal V}}
\newcommand{\WWW}{\mathord{\mathcal W}}
\newcommand{\GGG}{\mathord{\mathcal G}}
\newcommand{\TTT}{\mathord{\mathcal T}}

\newcommand{\closure}[1]{\overline #1}

\newcommand{\rank}{\mathop{\rm rank}\nolimits}
\newcommand{\disc}{\mathop{\rm disc}\nolimits}
\newcommand{\reddisc}{\mathop{\rm reddisc}\nolimits}
\newcommand{\pexcess}{\text{\rm $p$-excess}}
\newcommand{\texcess}{\text{\rm $2$-excess}}
\newcommand{\length}{\mathop{\rm length}\nolimits}
\newcommand{\ord}[1]{\mathop{\rm ord}\nolimits \sb{#1}}
\newcommand{\Roots}{\mathop{\rm Roots}\nolimits}
\newcommand{\Aut}{\mathop{\rm Aut}\nolimits}

\newcommand{\dual}[1]{#1 \sp{\vee}}

\newcommand{\bigfrac}[2]{\displaystyle\frac{  #1}{\strut #2}}
\newcommand{\euler}{\mathop{\rm euler}\nolimits}
\newcommand{\NS}{\mathop{\rm NS}\nolimits}
\newcommand{\set}[2]{\{\hskip 3pt #1 \hskip 3pt : \hskip 3pt #2 \hskip 3pt\}}
\newcommand{\inv}{\sp{-1}}

\newcommand{\sprime}{\sp\prime}
\newcommand{\SL}{\mathord{{\rm SL}}}

\newcommand{\isom}{\hskip 2pt \smash{\mathop{\to}\limits\sp{\hskip -1pt \sim}}\hskip 2pt}

\def\listup#1{{\footnotesize\smallskip\par\noindent
#1
\smallskip\par\noindent}}
\long\def\verbatim#1{\def\next{#1}%
{\tt\frenchspacing\expandafter\strip\meaning\next}}
\def\strip#1>{}
\newcommand{\cparbox}[2]{\parbox[t]{#1}{\hfil#2\hfil}}
\newcommand{\elldata}[4]{
\framebox {%
		\hbox {%
			\cparbox{1cm}{#1}% No.
			\vrule
			\cparbox{.9cm}{#2}% rank
   \vrule 
			\cparbox{5.35cm}{#3\hfil} %ADE
			\vrule
			\cparbox{3cm}{#4}% %MWs
		}%
	}%
\vskip -1.1pt
}
\newcommand{\vsr}{\vskip 1cm}
\newcommand{\vsrs}{\vskip 0cm}
\newcommand{\Gvsp}{\vskip 3pt}
\newcommand{\GGvsp}{\vskip 2.5pt}

\begin{document}

\title[Elliptic $K3$ surfaces]{On elliptic $K3$ surfaces}

% author  information
\author{Ichiro Shimada}
\address{
Department of Mathematics,
Faculty of Science,
Hokkaido University,
Sapporo 060-0810,
JAPAN
}
%\curraddr{}
\email{shimada@@math.sci.hokudai.ac.jp
}

%\date{}

% at present the "communicated by" line appears only in ERA and PROC
%\commby{}

\dedicatory{}

\begin{abstract}
We make a complete list of all possible $ADE$-types of singular fibers of 
complex elliptic $K3$ surfaces and the torsion parts of
their Mordell-Weil groups. 
\end{abstract}

\maketitle
\section{Introduction}\label{sec:intro}
By virtue of Torelli theorem for the period map
on the moduli of complex $K3$ surfaces (\cite{BurnsRapoport}, \cite{PS}, \cite{Todorov}),
we can study many aspects of $K3$ surfaces 
from the lattice-theoretic point of view.
In this paper, we determine all possible
$ADE$-types of singular fibers of elliptic $K3$ surfaces
using Nikulin's theory of  discriminant forms
of even integral lattices.
We also determine,
for each $ADE$-type of singular fibers,
all possible 
torsion parts of the Mordell-Weil groups.
Throughout this paper,
we use the term ``an elliptic $K3$ surface"
for ``a complex  elliptic $K3$ surface with a distinguished zero section"
and the term ``an elliptic fibration"
for ``a complex Jacobian elliptic fibration".
\par
\medskip
A finite formal sum of the symbols $A\sb l$ $(l\ge 1)$,
$D\sb m$ $(m\ge 4)$ and $E\sb n$ $(n=6, 7, 8)$
with non-negative integer coefficients
is called an $ADE$-type.
For an $ADE$-type
$$
\Sigma :=\sum a\sb l A\sb l + \sum d\sb m D\sb m + \sum e\sb n E\sb n, 
$$
we denote by $\opp{L (\Sigma)}$ the negative-definite root lattice 
generated by a root system of type $\Sigma$,
and by $\rank \Sigma$ the rank of $\opp{L (\Sigma)}$.
By definition, we have
$\rank \Sigma=\sum a\sb l  l + \sum d\sb m  m + \sum e\sb n n$.
\par
\medskip
Let $f: X\to \P\sp 1$ be an elliptic $K3$ surface,
and $O :\P\sp 1 \to X$ the zero section of $f$.
Let $MW \sb f$ be the Mordell-Weil group of $f$.
The torsion part of $MW\sb f$ is a finite abelian group,
which we shall denote by $G\sb f$.
We put 
$$
R\sb f := \set{ p \in \P\sp 1 }{ \text{\rm $f\sp{-1} (p)$ is reducible} },
$$
and, for each $p\in R\sb f$, we
denote by $f\sp{-1} (p)\sp{\sharp}$
the union of irreducible components
of $f\sp{-1} (p)$ that are disjoint from the zero section.
It is known that the cohomology classes of irreducible components
of $f\sp{-1} (p)\sp{\sharp}$ span a negative-definite root lattice
generated by an indecomposable root system of type $A\sb l$, $D\sb m$ or $E\sb n$.
Let $\tau\sb{f, p}$ be the type.
The type of singular fiber $f\sp{-1} (p)$
in the list of Kodaira's classification~\cite{Kodaira}  
is  related to $\tau\sb{f, p}$
in an almost one-to-one way (cf.\ Table~\ref{table:kodaira}).
We define the $ADE$-type $\Sigma\sb f$ of  $f : X\to \P\sp 1$
by
$$
\Sigma\sb f :=\sum\sb{p\in R\sb f} \tau\sb{f, p}.
$$
The N\'eron-Severi lattice $\NS\sb X$ of $X$
contains the sublattice $S\sb f$ generated by
the cohomology classes of  
the irreducible components
of $\cup \sb{p \in R\sb f} f\sp{-1} (p)\sp{\sharp}$,
which is isomorphic to $\opp{L(\Sigma\sb f)}$.
\par
\medskip
Through computer-aided calculation,
we have made the complete list of  pairs $(\Sigma, G)$
of an $ADE$-type $\Sigma$ and a finite abelian group $G$
that can be realized as the   data $(\Sigma\sb f, G\sb f)$
of an elliptic $K3$ surface
$f : X\to \P\sp 1$.
This list $\PPP$ consists of $3693$  pairs.
In this paper,
we present the list $\PPP$,
deduce some geometric facts from it,   and 
explain the algorithm 
for obtaining it.
\par
\medskip
The list $\PPP$ is too large  to be included here
in a naive way.
Therefore we describe $\PPP$ by 
giving a subset $\SSS$ of $\PPP$ 
and a set of transformation rules of $ADE$-types that
generate $\PPP$ from $\SSS$
(cf.\ \S\ \ref{sec:main}).
The reader can obtain $\PPP$
easily using  this description.
\footnote{The list  can also be  retrieved from   the author's homepage.\par
\verbatim{http://www.math.sci.hokudai.ac.jp/~shimada/K3.html}}
\par
\medskip
An elliptic $K3$ surface $f :X\to \P\sp1$ is said to be extremal if
 the sublattice $S\sb f$ attains the maximal  rank $18$.
After the work of Miranda and  Persson~\cite{MP}, 
supplemented by Artal-Bartolo, Tokunaga and Zhang~\cite{ATZ} and Ye~\cite{Ye},
the $ADE$-types of singular fibers 
of extremal elliptic $K3$ surfaces and their Mordell-Weil groups
were completely determined in \cite{ShimadaZhang}.
The list consists of $336$  pairs.
\par
\medskip
One of the remarkable facts that can be read off from the list $\PPP$ is 
that
an $ADE$-type $\Sigma$ is an $ADE$-type of
 an elliptic $K3$ surface 
with   trivial   Mordell-Weil torsion
if and only if $\Sigma$ is obtained from
an $ADE$-type of  an extremal elliptic $K3$ surface 
with   trivial   Mordell-Weil torsion
by {\it elementary transformation};
that is, by deleting vertices from the corresponding Dynkin graph
(cf.\ Theorem~\ref{thm:triv}).
In order to
describe the list of $ADE$-types 
of elliptic $K3$ surfaces with non-trivial Mordell-Weil torsion,
however, 
we have to
forbid to use some types of   elementary transformation
(cf.\ Theorems~\ref{thm:two}-\ref{thm:twotwo}).
\par
\medskip
By Nishiyama~\cite{Nishiyama} and by Besser~\cite{Besser},
the technique of discriminant forms was
used to find out all possible 
elliptic fibrations on special  $K3$
surfaces. In \cite{Urabe1,  Urabe2},
Urabe investigated possible configurations of singular points on
$K3$ surfaces and suggested an existence of
a set of simple rules that generates all possible configurations.
In \cite{Yang1, Yang2},
Yang made the complete list of all possible configurations of
singularities of $ADE$-type on 
plane sextic curves and quartic surfaces
using the technique of discriminant forms
and a computer.
\par
\medskip
This paper is organized as follows.
In \S\ \ref{sec:main},
we describe  $\PPP$ and state some  facts 
about elliptic $K3$ surfaces that can be derived from the list $\PPP$.
In \S\ \ref{sec:local},
we recall the definition and properties 
of local invariants of lattices over $\Z$
according to Conway and Sloane~\cite[Chapter 15]{ConwaySloane}.
In \S\ \ref{sec:theorydiscfrm} and~\ref{sec:givendiscfrm},
we review Nikulin's theory~\cite{Nikulin} of
discriminant forms of even  lattices over $\Z$.
A criterion whether there exists an even integral lattice 
of a given signature and a discriminant form is described in detail 
in \S\ \ref{sec:givendiscfrm}.
This criterion is slightly different from~\cite[Theorem 1.10.1]{Nikulin},
and is more suited to machine  calculation.
In \S\ \ref{sec:roots},
we recall the properties of root lattices.
In \S\ \ref{sec:existenceK3},
we show that it is possible to determine 
by a purely lattice-theoretic calculation
whether a given  pair $(\Sigma, G)$ can be realized as $(\Sigma\sb f, G\sb f)$
of an elliptic $K3$ surface $f : X\to \P\sp 1$.
Here we use Kondo-Nishiyama's lemma on the N\'eron-Severi lattice
of an elliptic $K3$ surface.
In \S\ \ref{sec:making}, we explain our algorithm.
\par
\medskip
The program for making $\PPP$ was written by Maple V.
The author would like to thank Waterloo Maple Incorporation  for developing the nice software.
The author also would like to thank
the referee for suggesting some improvements 
on the first version of the paper.
\par
\medskip
 A compressed version of this paper
 without the table has appeared in 
 \begin{center}
 Michigan Math. J. 47 (2000), 423--446.
 \end{center}
\section{Main results}\label{sec:main}
All results in this section are obtained
simply by looking at the list $\PPP$.
\subsection{Tortion parts of Mordell-Weil groups}
\begin{theorem}\label{G}
The torsion part of the Mordell-Weil group of an elliptic $K3$ surface is 
isomorphic to one of the following{\rm :}
\begin{equation}\label{Gset}
\vcenter{\hbox{\vbox{
\hbox{$(0),\:\:\: \Z/(2),\:\:\:  \Z/(3),\:\:\:  \Z/(4), \:\:\: \Z/(5),\:\:\:  \Z/(6),\:\:\:  \Z/(7), \:\:\: \Z/(8),$}
\vskip 3pt
\hbox{$\Z/(2)\times \Z/(2),\:\:\: \Z/(4)\times\Z/(2),  \:\:\:\Z/(6)\times \Z/(2),$}
\vskip 3pt
\hbox{$\Z/(3)\times\Z/(3), \:\:\: \Z/(4)\times\Z/(4)$.}}}}
\end{equation}
\end{theorem}
For a group $G$ in \eqref{Gset},
we  denote by $\PPP\sp G$ the set of all $ADE$-types $\Sigma$ such that
there exists an elliptic $K3$ surface $f : X\to \P\sp 1$ with 
$\Sigma\sb f = \Sigma$ and $G\sb f \cong G$.
The cardinalities of $\PPP\sp G$ are  
given in Table~\ref{table:PG}.
\begin{table}
\caption{Cardinalities of $\PPP\sp G$}
\label{table:PG}
\def\spaceheight{height 3pt}
\def\smallhsk{\hskip 3pt}
\def\HS{
\spaceheight & \omit  & \spaceheight &
\omit  & \omit  & \omit  &
\omit  & \omit  & \omit  &
\omit  & \omit  & \omit  &
\omit  & \omit  & \omit  &
\omit  & 
\spaceheight & \omit &\spaceheight
\cr
}
\def\HL{\HS\noalign{\hrule}\HS}
\centerline{
\vbox{\offinterlineskip
\hrule
\halign{
\vrule # \hskip 2pt& \strut\smallhsk \hfil#\smallhsk\hfil & \vrule # \hskip 2pt &
\strut\smallhsk \hfil#\smallhsk\hfil & \strut\smallhsk \hfil#\smallhsk\hfil & \strut\smallhsk \hfil#\smallhsk\hfil &
\strut\smallhsk \hfil#\smallhsk\hfil & \strut\smallhsk \hfil#\smallhsk\hfil & \strut\smallhsk \hfil#\smallhsk\hfil &
\strut\smallhsk \hfil#\smallhsk\hfil & \strut\smallhsk \hfil#\smallhsk\hfil & \strut\smallhsk \hfil#\smallhsk\hfil &
\strut\smallhsk \hfil#\smallhsk\hfil & \strut\smallhsk \hfil#\smallhsk\hfil & \strut\smallhsk \hfil#\smallhsk\hfil &
\strut\smallhsk \hfil#\smallhsk\hfil & 
\vrule # \hskip 2pt & \strut\smallhsk \hfil#\smallhsk\hfil & \vrule #
\cr
\HS
& $G$ & &
$[1]$ & $[2]$ & $[3]$  &
$[4]$  & $[5]$  & $[6]$  &
$[7]$  & $[8]$  & $[2,2]$  &
$[4,2]$  & $[6,2]$  & $[3,3]$  &
$[4,4]$  & 
 & total   & 
\cr
\HL
& $|\PPP\sp G|$ & &
$2746$ & $732$ & $85$  &
$41$  & $6$  & $10$  &
$1$  & $1$  & $61$  &
$5$  & $1$  & $3$  &
$1$  & 
 & $ 3693$  & 
\cr
\HS
}
\hrule
}
}
\end{table}
Here, $[a]$ denotes the cyclic group $\Z/(a)$, and
$[a, b]$ denotes $\Z/(a)\times \Z/(b)$.
In particular, $[1]$ denotes the trivial group.
\par
For a positive integer $r$,
let $\PPP\sp G\sb r$ be the subset of $\PPP\sp G$
that consists of $\Sigma \in \PPP\sp G$
with $\rank \Sigma =r$.
Let $f : X\to \P\sp 1$ be an elliptic $K3$ surface.
Since the N\'eron-Severi lattice $\NS\sb X$ of  $X$
is the orthogonal direct sum of $S\sb f\cong L (\Sigma\sb f )\sp{-}$ and the lattice 
of rank $2$ generated by the cohomology classes of
the zero section and a general fiber,
and since the N\'eron-Severi rank of $X$ is at most $20$, 
we always have
$$
 \rank (\Sigma\sb f) \le 18.
$$
Hence  $\PPP\sp G\sb r$ is empty for $r>18$.
\subsection{$ADE$-types of singular fibers}
Next we describe the list $\PPP\sp G$ for
each abelian group $G$ in \eqref{Gset}.
We carry out this task  by three different methods
according to the size of $\PPP\sp G$.
\par
\medskip
{\bf Case 1.} $G\in \{[1], [2], [3], [4], [2,2]\}$ .
We describe $\PPP\sp G$ by giving a subset
$\SSS\sp G\st \PPP\sp G$ 
and a set of transformation rules on $ADE$-types
that generate the whole $\PPP\sp G$
from the subset $\SSS\sp G$. 
\par
Let $\Gamma (\Sigma)$ be the Dynkin graph of the $ADE$-type $\Sigma$.
If we remove a vertex $P$ of $\Gamma (\Sigma)$ and the edges emitting from $P$,
we obtain
the Dynkin graph $\Gamma (\Sigma\sprime)$ of another 
$ADE$-type $\Sigma\sprime$ with $\rank \Sigma \sprime =\rank \Sigma -1$.
In this case,
we say that $\Sigma\sprime$ is obtained from $\Sigma$ by
deleting a vertex.
In other words, 
an $ADE$-type $\Sigma\sprime$ is obtained
from $\Sigma$ by
deleting a vertex
if and only if $\Sigma\sprime$ is 
obtained  by applying to $\Sigma$ one of the substitutions listed in
Table~\ref{table:subs}.
In this Table, we understand that $A\sb 0:=0$.
\begin{definition}
When we can obtain an $ADE$-type $\Sigma\sprime$
from an $ADE$-type $\Sigma$ by
applying substitutions  in Table~\ref{table:subs}
several times,
we say that $\Sigma\sprime$ is obtained from
$\Sigma$ by {\it elementary transformation}.
\end{definition}
\begin{table}\label{table:subs}
\caption{Substitutions.}
\begin{center}
\fbox{
$
\begin{array}{ccl}
\noalign{\vskip 1pt}
A\sb l &\mapsto & A\sb{l\sprime}+A\sb{l-1-l\sprime} \quad
(0\le  l\sprime \le l/2),
\\ 
\noalign{\vskip 5pt}
D\sb{m} &\mapsto & 
\left\{\vcenter{
\hbox{\vbox{
\hbox{$A\sb{m-1}, \quad 2A\sb 1 + A\sb{m-3}, \quad A\sb 3 + A\sb{m-4},$}
\vskip 3pt 
\hbox{$D\sb{m\sprime} + A\sb{m-1-m\sprime}\quad(4\le m\sprime\le  m-1),$}
}
}
}\right.
\\
\noalign{\vskip 5pt}
E\sb{n} &\mapsto & 
\left\{
\vcenter{
\hbox{\vbox{
\hbox{$A\sb{n-1}, \quad D\sb{n-1}, \quad A\sb{1}+A\sb{n-2}, \quad A\sb{1}+A\sb{2}+A\sb{n-4},
\quad A\sb{4}+A\sb{n-5},$}
\vskip 3pt 
\hbox{$D\sb{5}+A\sb{n-6},\quad E\sb{n\sprime}+A\sb{n-1-n\sprime}
\quad(6\le n\sprime\le n-1).$}
}
}
}\right.
\\
\noalign{\vskip 1pt}
\end{array}
$
}
\end{center}
\end{table}
\begin{theorem}\label{thm:triv}
\text{\rm (1)}
The list $\PPP\sp{[1]}\sb{18}$ consists of $199$ elements listed below.
\listup{%
  $  2E_{{8}}+A_{{2}}$,
  $  2E_{{8}}+2A_{{1}}$,
  $  E_{{8}}+E_{{7}}+A_{{3}}$,
  $  E_{{8}}+E_{{7}}+A_{{2}}+A_{{1}}$,
  $  E_{{8}}+E_{{6}}+D_{{4}}$,
  $  E_{{8}}+E_{{6}}+A_{{4}}$,
  $  E_{{8}}+E_{{6}}+A_{{3}}+A_{{1}}$,
  $  E_{{8}}+D_{{10}}$,
  $  E_{{8}}+D_{{9}}+A_{{1}}$,
  $  E_{{8}}+D_{{7}}+A_{{2}}+A_{{1}}$,
  $  E_{{8}}+D_{{6}}+A_{{4}}$,
  $  E_{{8}}+D_{{6}}+2A_{{2}}$,
  $  E_{{8}}+2D_{{5}}$,
  $  E_{{8}}+D_{{5}}+A_{{5}}$,
  $  E_{{8}}+D_{{5}}+A_{{4}}+A_{{1}}$,
  $  E_{{8}}+A_{{10}}$,
  $  E_{{8}}+A_{{9}}+A_{{1}}$,
  $  E_{{8}}+A_{{8}}+A_{{2}}$,
  $  E_{{8}}+A_{{8}}+2A_{{1}}$,
  $  E_{{8}}+A_{{7}}+A_{{2}}+A_{{1}}$,
  $  E_{{8}}+A_{{6}}+A_{{4}}$,
  $  E_{{8}}+A_{{6}}+A_{{3}}+A_{{1}}$,
  $  E_{{8}}+A_{{6}}+2A_{{2}}$,
  $  E_{{8}}+A_{{6}}+A_{{2}}+2A_{{1}}$,
  $  E_{{8}}+2A_{{5}}$,
  $  E_{{8}}+A_{{5}}+A_{{4}}+A_{{1}}$,
  $  E_{{8}}+A_{{5}}+A_{{3}}+A_{{2}}$,
  $  E_{{8}}+2A_{{4}}+2A_{{1}}$,
  $  E_{{8}}+A_{{4}}+A_{{3}}+A_{{2}}+A_{{1}}$,
  $  E_{{8}}+2A_{{3}}+2A_{{2}}$,
  $  2E_{{7}}+A_{{4}}$,
  $  2E_{{7}}+2A_{{2}}$,
  $  E_{{7}}+E_{{6}}+D_{{5}}$,
  $  E_{{7}}+E_{{6}}+A_{{5}}$,
  $  E_{{7}}+E_{{6}}+A_{{4}}+A_{{1}}$,
  $  E_{{7}}+E_{{6}}+A_{{3}}+A_{{2}}$,
  $  E_{{7}}+D_{{11}}$,
  $  E_{{7}}+D_{{9}}+A_{{2}}$,
  $  E_{{7}}+D_{{7}}+A_{{4}}$,
  $  E_{{7}}+D_{{5}}+A_{{6}}$,
  $  E_{{7}}+D_{{5}}+A_{{4}}+A_{{2}}$,
  $  E_{{7}}+A_{{11}}$,
  $  E_{{7}}+A_{{10}}+A_{{1}}$,
  $  E_{{7}}+A_{{9}}+A_{{2}}$,
  $  E_{{7}}+A_{{8}}+A_{{3}}$,
  $  E_{{7}}+A_{{8}}+A_{{2}}+A_{{1}}$,
  $  E_{{7}}+A_{{7}}+A_{{4}}$,
  $  E_{{7}}+A_{{7}}+2A_{{2}}$,
  $  E_{{7}}+A_{{6}}+A_{{5}}$,
  $  E_{{7}}+A_{{6}}+A_{{4}}+A_{{1}}$,
  $  E_{{7}}+A_{{6}}+A_{{3}}+A_{{2}}$,
  $  E_{{7}}+A_{{6}}+2A_{{2}}+A_{{1}}$,
  $  E_{{7}}+A_{{5}}+A_{{4}}+A_{{2}}$,
  $  E_{{7}}+A_{{4}}+A_{{3}}+2A_{{2}}$,
  $  2E_{{6}}+D_{{6}}$,
  $  2E_{{6}}+A_{{6}}$,
  $  2E_{{6}}+2A_{{3}}$,
  $  E_{{6}}+D_{{12}}$,
  $  E_{{6}}+D_{{11}}+A_{{1}}$,
  $  E_{{6}}+D_{{9}}+A_{{3}}$,
  $  E_{{6}}+D_{{9}}+A_{{2}}+A_{{1}}$,
  $  E_{{6}}+D_{{8}}+A_{{4}}$,
  $  E_{{6}}+D_{{7}}+D_{{5}}$,
  $  E_{{6}}+D_{{7}}+A_{{4}}+A_{{1}}$,
  $  E_{{6}}+D_{{6}}+A_{{6}}$,
  $  E_{{6}}+D_{{6}}+A_{{4}}+A_{{2}}$,
  $  E_{{6}}+D_{{5}}+A_{{7}}$,
  $  E_{{6}}+D_{{5}}+A_{{6}}+A_{{1}}$,
  $  E_{{6}}+D_{{5}}+A_{{4}}+A_{{3}}$,
  $  E_{{6}}+A_{{12}}$,
  $  E_{{6}}+A_{{11}}+A_{{1}}$,
  $  E_{{6}}+A_{{10}}+A_{{2}}$,
  $  E_{{6}}+A_{{10}}+2A_{{1}}$,
  $  E_{{6}}+A_{{9}}+A_{{3}}$,
  $  E_{{6}}+A_{{9}}+A_{{2}}+A_{{1}}$,
  $  E_{{6}}+A_{{8}}+A_{{4}}$,
  $  E_{{6}}+A_{{8}}+A_{{3}}+A_{{1}}$,
  $  E_{{6}}+A_{{7}}+A_{{5}}$,
  $  E_{{6}}+A_{{7}}+A_{{4}}+A_{{1}}$,
  $  E_{{6}}+A_{{6}}+A_{{5}}+A_{{1}}$,
  $  E_{{6}}+A_{{6}}+A_{{4}}+A_{{2}}$,
  $  E_{{6}}+A_{{6}}+A_{{4}}+2A_{{1}}$,
  $  E_{{6}}+A_{{6}}+A_{{3}}+A_{{2}}+A_{{1}}$,
  $  E_{{6}}+A_{{5}}+A_{{4}}+A_{{3}}$,
  $  E_{{6}}+2A_{{4}}+A_{{3}}+A_{{1}}$,
  $  D_{{18}}$,
  $  D_{{17}}+A_{{1}}$,
  $  D_{{15}}+A_{{2}}+A_{{1}}$,
  $  D_{{14}}+A_{{4}}$,
  $  D_{{14}}+2A_{{2}}$,
  $  D_{{13}}+D_{{5}}$,
  $  D_{{13}}+A_{{5}}$,
  $  D_{{13}}+A_{{4}}+A_{{1}}$,
  $  D_{{11}}+A_{{6}}+A_{{1}}$,
  $  D_{{11}}+A_{{5}}+A_{{2}}$,
  $  D_{{11}}+A_{{4}}+A_{{2}}+A_{{1}}$,
  $  D_{{11}}+A_{{3}}+2A_{{2}}$,
  $  D_{{10}}+A_{{8}}$,
  $  D_{{10}}+A_{{6}}+A_{{2}}$,
  $  D_{{10}}+2A_{{4}}$,
  $  2D_{{9}}$,
  $  D_{{9}}+D_{{5}}+A_{{4}}$,
  $  D_{{9}}+A_{{9}}$,
  $  D_{{9}}+A_{{8}}+A_{{1}}$,
  $  D_{{9}}+A_{{6}}+A_{{2}}+A_{{1}}$,
  $  D_{{9}}+A_{{5}}+A_{{4}}$,
  $  D_{{9}}+A_{{4}}+2A_{{2}}+A_{{1}}$,
  $  D_{{8}}+A_{{6}}+2A_{{2}}$,
  $  2D_{{7}}+2A_{{2}}$,
  $  D_{{7}}+A_{{10}}+A_{{1}}$,
  $  D_{{7}}+A_{{9}}+A_{{2}}$,
  $  D_{{7}}+A_{{6}}+A_{{5}}$,
  $  D_{{7}}+A_{{6}}+A_{{4}}+A_{{1}}$,
  $  D_{{7}}+A_{{6}}+A_{{3}}+A_{{2}}$,
  $  D_{{7}}+2A_{{4}}+A_{{2}}+A_{{1}}$,
  $  D_{{6}}+A_{{12}}$,
  $  D_{{6}}+A_{{10}}+A_{{2}}$,
  $  D_{{6}}+A_{{8}}+A_{{4}}$,
  $  D_{{6}}+2A_{{6}}$,
  $  D_{{6}}+A_{{6}}+A_{{4}}+A_{{2}}$,
  $  D_{{6}}+2A_{{4}}+2A_{{2}}$,
  $  2D_{{5}}+A_{{8}}$,
  $  2D_{{5}}+2A_{{4}}$,
  $  D_{{5}}+A_{{13}}$,
  $  D_{{5}}+A_{{12}}+A_{{1}}$,
  $  D_{{5}}+A_{{10}}+A_{{2}}+A_{{1}}$,
  $  D_{{5}}+A_{{9}}+A_{{4}}$,
  $  D_{{5}}+A_{{9}}+2A_{{2}}$,
  $  D_{{5}}+A_{{8}}+A_{{5}}$,
  $  D_{{5}}+A_{{8}}+A_{{4}}+A_{{1}}$,
  $  D_{{5}}+2A_{{6}}+A_{{1}}$,
  $  D_{{5}}+A_{{6}}+A_{{5}}+A_{{2}}$,
  $  D_{{5}}+A_{{6}}+A_{{4}}+A_{{2}}+A_{{1}}$,
  $  D_{{5}}+A_{{6}}+A_{{3}}+2A_{{2}}$,
  $  D_{{5}}+A_{{5}}+2A_{{4}}$,
  $  A_{{18}}$,
  $  A_{{17}}+A_{{1}}$,
  $  A_{{16}}+A_{{2}}$,
  $  A_{{16}}+2A_{{1}}$,
  $  A_{{15}}+A_{{2}}+A_{{1}}$,
  $  A_{{14}}+A_{{4}}$,
  $  A_{{14}}+A_{{3}}+A_{{1}}$,
  $  A_{{14}}+A_{{2}}+2A_{{1}}$,
  $  A_{{13}}+A_{{5}}$,
  $  A_{{13}}+A_{{4}}+A_{{1}}$,
  $  A_{{13}}+A_{{3}}+A_{{2}}$,
  $  A_{{13}}+2A_{{2}}+A_{{1}}$,
  $  A_{{12}}+A_{{6}}$,
  $  A_{{12}}+A_{{5}}+A_{{1}}$,
  $  A_{{12}}+A_{{4}}+A_{{2}}$,
  $  A_{{12}}+A_{{4}}+2A_{{1}}$,
  $  A_{{12}}+A_{{3}}+A_{{2}}+A_{{1}}$,
  $  A_{{12}}+2A_{{2}}+2A_{{1}}$,
  $  A_{{11}}+A_{{6}}+A_{{1}}$,
  $  A_{{11}}+A_{{4}}+A_{{2}}+A_{{1}}$,
  $  A_{{10}}+A_{{8}}$,
  $  A_{{10}}+A_{{7}}+A_{{1}}$,
  $  A_{{10}}+A_{{6}}+A_{{2}}$,
  $  A_{{10}}+A_{{6}}+2A_{{1}}$,
  $  A_{{10}}+A_{{5}}+A_{{3}}$,
  $  A_{{10}}+A_{{5}}+A_{{2}}+A_{{1}}$,
  $  A_{{10}}+2A_{{4}}$,
  $  A_{{10}}+A_{{4}}+A_{{3}}+A_{{1}}$,
  $  A_{{10}}+A_{{4}}+2A_{{2}}$,
  $  A_{{10}}+A_{{4}}+A_{{2}}+2A_{{1}}$,
  $  A_{{10}}+2A_{{3}}+A_{{2}}$,
  $  A_{{10}}+A_{{3}}+2A_{{2}}+A_{{1}}$,
  $  2A_{{9}}$,
  $  A_{{9}}+A_{{8}}+A_{{1}}$,
  $  A_{{9}}+A_{{7}}+A_{{2}}$,
  $  A_{{9}}+A_{{6}}+A_{{3}}$,
  $  A_{{9}}+A_{{6}}+A_{{2}}+A_{{1}}$,
  $  A_{{9}}+A_{{5}}+A_{{4}}$,
  $  2A_{{8}}+2A_{{1}}$,
  $  A_{{8}}+A_{{7}}+A_{{2}}+A_{{1}}$,
  $  A_{{8}}+A_{{6}}+A_{{4}}$,
  $  A_{{8}}+A_{{6}}+A_{{3}}+A_{{1}}$,
  $  A_{{8}}+A_{{6}}+A_{{2}}+2A_{{1}}$,
  $  A_{{8}}+A_{{5}}+A_{{4}}+A_{{1}}$,
  $  A_{{8}}+2A_{{4}}+2A_{{1}}$,
  $  A_{{8}}+A_{{4}}+A_{{3}}+A_{{2}}+A_{{1}}$,
  $  2A_{{7}}+2A_{{2}}$,
  $  A_{{7}}+A_{{6}}+A_{{5}}$,
  $  A_{{7}}+A_{{6}}+A_{{4}}+A_{{1}}$,
  $  A_{{7}}+A_{{6}}+A_{{3}}+A_{{2}}$,
  $  A_{{7}}+A_{{6}}+2A_{{2}}+A_{{1}}$,
  $  A_{{7}}+A_{{5}}+A_{{4}}+A_{{2}}$,
  $  A_{{7}}+A_{{4}}+A_{{3}}+2A_{{2}}$,
  $  2A_{{6}}+A_{{4}}+A_{{2}}$,
  $  2A_{{6}}+2A_{{3}}$,
  $  2A_{{6}}+2A_{{2}}+2A_{{1}}$,
  $  A_{{6}}+A_{{5}}+A_{{4}}+A_{{3}}$,
  $  A_{{6}}+A_{{5}}+A_{{4}}+A_{{2}}+A_{{1}}$,
  $  A_{{6}}+2A_{{4}}+A_{{3}}+A_{{1}}$,
  $  A_{{6}}+2A_{{4}}+A_{{2}}+2A_{{1}}$,
  $  A_{{6}}+A_{{4}}+2A_{{3}}+A_{{2}}$,
  $  A_{{6}}+A_{{4}}+A_{{3}}+2A_{{2}}+A_{{1}}$,
  $  2A_{{5}}+2A_{{4}}$,
  $  2A_{{4}}+2A_{{3}}+2A_{{2}}$.
}%
\text{\rm (2)}
An $ADE$-type $\Sigma$ with $r:=\rank \Sigma <18$ is a member of $\PPP\sp{[1]}\sb{r}$
if and only if 
$\Sigma$ is obtained
from a member of $\PPP\sp{[1]}\sb{18}$ by
elementary transformation.
\end{theorem}
%
% 
%  Here starts the tables of forbidden substitutions.
%
% 
%
\begin{table}
\caption{Forbidden substitutions for $[2]$.}
\label{table:forbiddensubs2}
\begin{center}
\fbox{
$
\begin{array}{ccl}
\noalign{\vskip 1pt}
A\sb l &\mapsto & A\sb{l\sprime}+A\sb{l-1-l\sprime} \quad
\textrm{with $l$ odd and $l\sprime$  even}\quad(0\le l\sprime < l/2), 
\\
\noalign{\vskip 3pt}
D\sb{m} &\mapsto & A\sb{m-1}, 
\\
\noalign{\vskip 3pt}
D\sb{m} &\mapsto & A\sb 3 + A\sb{m-4}\quad\textrm{with $m$ even,}
\\
\noalign{\vskip 3pt}
D\sb{m} &\mapsto & D\sb{m\sprime} + A\sb{m-1-m\sprime}\quad\textrm{with $m$ even and $m\sprime$ odd}
\quad (5\le   m\sprime \le m-1),
\\
\noalign{\vskip 3pt}
E\sb{7} &\mapsto & A\sb 6, \: A\sb 4 + A\sb 2, \: E\sb 6.
\\
\noalign{\vskip 1pt}
\end{array}
$
}
\end{center}
\end{table}
\begin{table}
\caption{Forbidden substitutions for $[3]$.}
\label{table:forbiddensubs3}
\begin{center}
\fbox{
$
\begin{array}{ccl}
\noalign{\vskip 1pt}
A\sb l &\mapsto & A\sb{l\sb 1}+A\sb{l\sb 2} \:\:(l\sb 1+l\sb 2=l-1, \: 0\le l\sb 1 \le l\sb 2)\\
\noalign{\vskip 3pt}
&&
\textrm{with $l \bmod 3=2$ ,\:
$l\sb 1 \bmod 3\ne 2$,\: $l\sb 2\sprime \bmod 3 \ne 2$,}
\\
\noalign{\vskip 3pt}
E\sb{6} &\mapsto & A\sb 4+A\sb 1, \: D\sb{5}.
\\
\noalign{\vskip 1pt}
\end{array}
$
}
\end{center}
\end{table}
\begin{table}
\caption{Forbidden substitutions for $[4]$.}
\label{table:forbiddensubs4}
\begin{center}
\fbox{
$
\begin{array}{ccl}
\noalign{\vskip 1pt}
A\sb 1 &\mapsto & 0
\\
\noalign{\vskip 3pt}
A\sb l &\mapsto & A\sb{l\sb 1}+A\sb{l\sb 2} \:\: (l\sb 1+l\sb 2=l-1, \:  0\le l\sb 1 \le l\sb 2)\\
\noalign{\vskip 3pt}
&&
\textrm{with $l \bmod 4=3$ ,\:
$l\sb 1 \bmod 4\ne 3$, \:$l\sb 2\sprime \bmod 4 \ne 3$,}
\\
\noalign{\vskip 3pt}
D\sb{m} &\mapsto & A\sb{m-1},  \: D\sb{m-1}, \: 2A\sb 1 +A\sb{m-3}\quad
\textrm{with $m$ odd}, 
\\
\noalign{\vskip 3pt}
D\sb{m} &\mapsto &  
 D\sb{m-3}+A\sb 2 \quad
\textrm{with $m$ odd and $>6$}.
\\
\noalign{\vskip 1pt}
\end{array}
$
}
\end{center}
\end{table}
\begin{table}
\caption{Forbidden substitutions for $[2,2]$.}
\label{table:forbiddensubs22}
\begin{center}
\fbox{
$
\begin{array}{ccl}
\noalign{\vskip 1pt}
A\sb l &\mapsto & A\sb{l\sprime}+A\sb{l-1-l\sprime} \quad
\textrm{with $l$ odd and $l\sprime$ even}\quad(0\le l\sprime < l/2), 
\\
\noalign{\vskip 3pt}
D\sb{m} &\mapsto & A\sb{m-1},  \:
A\sb 3 + A\sb{m-4}\quad
\textrm{with $m$ even, }
\\
\noalign{\vskip 3pt}
D\sb{m} &\mapsto &   D\sb{m\sprime}+ A\sb{m-1- m\sprime}\quad
\textrm{with $m$ even and $m\sprime$ odd}\quad(5\le m\sprime \le m-1).
\\
\noalign{\vskip 1pt}
\end{array}
$
}
\end{center}
\end{table}
\begin{theorem}\label{thm:two}
\text{\rm (1)}
The list $\PPP\sp{[2]}\sb{18}$ consists of $84$ elements listed below.
\listup{%
 $2E_{{7}}+D_{{4}}$, 
 $2E_{{7}}+A_{{3}}+A_{{1}}$, 
 $E_{{7}}+D_{{10}}+A_{{1}}$, 
 $E_{{7}}+D_{{8}}+A_{{2}}+A_{{1}}$, 
 $E_{{7}}+D_{{7}}+A_{{3}}+A_{{1}}$, 
 $E_{{7}}+D_{{6}}+D_{{5}}$, 
 $E_{{7}}+D_{{6}}+A_{{5}}$, 
 $E_{{7}}+D_{{6}}+A_{{3}}+A_{{2}}$, 
 $E_{{7}}+D_{{5}}+A_{{5}}+A_{{1}}$, 
 $E_{{7}}+A_{{9}}+A_{{2}}$, 
 $E_{{7}}+A_{{9}}+2A_{{1}}$, 
 $E_{{7}}+A_{{7}}+A_{{3}}+A_{{1}}$, 
 $E_{{7}}+A_{{7}}+A_{{2}}+2A_{{1}}$, 
 $E_{{7}}+A_{{5}}+A_{{4}}+2A_{{1}}$, 
 $E_{{7}}+A_{{5}}+2A_{{3}}$, 
 $E_{{7}}+A_{{5}}+A_{{3}}+A_{{2}}+A_{{1}}$, 
 $E_{{7}}+A_{{4}}+2A_{{3}}+A_{{1}}$, 
 $D_{{16}}+A_{{2}}$, 
 $D_{{16}}+2A_{{1}}$, 
 $D_{{14}}+A_{{3}}+A_{{1}}$, 
 $D_{{14}}+A_{{2}}+2A_{{1}}$, 
 $D_{{12}}+D_{{6}}$, 
 $D_{{12}}+D_{{5}}+A_{{1}}$, 
 $D_{{12}}+A_{{4}}+2A_{{1}}$, 
 $D_{{12}}+A_{{3}}+A_{{2}}+A_{{1}}$, 
 $D_{{12}}+2A_{{2}}+2A_{{1}}$, 
 $D_{{10}}+D_{{7}}+A_{{1}}$, 
 $D_{{10}}+D_{{6}}+A_{{2}}$, 
 $D_{{10}}+D_{{5}}+A_{{2}}+A_{{1}}$, 
 $D_{{10}}+A_{{5}}+A_{{3}}$, 
 $D_{{10}}+A_{{4}}+A_{{3}}+A_{{1}}$, 
 $D_{{9}}+A_{{7}}+2A_{{1}}$, 
 $D_{{9}}+A_{{5}}+A_{{3}}+A_{{1}}$, 
 $D_{{8}}+2D_{{5}}$, 
 $D_{{8}}+A_{{9}}+A_{{1}}$, 
 $D_{{8}}+A_{{7}}+A_{{2}}+A_{{1}}$, 
 $D_{{8}}+2A_{{5}}$, 
 $D_{{8}}+A_{{5}}+A_{{4}}+A_{{1}}$, 
 $D_{{8}}+2A_{{3}}+2A_{{2}}$, 
 $D_{{7}}+D_{{6}}+A_{{5}}$, 
 $D_{{7}}+D_{{5}}+A_{{5}}+A_{{1}}$, 
 $D_{{7}}+A_{{9}}+2A_{{1}}$, 
 $D_{{7}}+A_{{7}}+A_{{2}}+2A_{{1}}$, 
 $D_{{6}}+D_{{5}}+A_{{7}}$, 
 $D_{{6}}+D_{{5}}+A_{{5}}+A_{{2}}$, 
 $D_{{6}}+A_{{11}}+A_{{1}}$, 
 $D_{{6}}+A_{{9}}+A_{{3}}$, 
 $D_{{6}}+A_{{9}}+A_{{2}}+A_{{1}}$, 
 $D_{{6}}+A_{{7}}+A_{{4}}+A_{{1}}$, 
 $D_{{6}}+A_{{7}}+A_{{3}}+A_{{2}}$, 
 $D_{{6}}+A_{{7}}+2A_{{2}}+A_{{1}}$, 
 $D_{{6}}+A_{{5}}+A_{{4}}+A_{{3}}$, 
 $D_{{5}}+A_{{11}}+A_{{2}}$, 
 $D_{{5}}+A_{{9}}+A_{{3}}+A_{{1}}$, 
 $D_{{5}}+A_{{9}}+A_{{2}}+2A_{{1}}$, 
 $D_{{5}}+A_{{7}}+A_{{4}}+2A_{{1}}$, 
 $D_{{5}}+2A_{{5}}+A_{{3}}$, 
 $D_{{5}}+A_{{5}}+A_{{4}}+A_{{3}}+A_{{1}}$, 
 $A_{{15}}+A_{{2}}+A_{{1}}$, 
 $A_{{13}}+A_{{4}}+A_{{1}}$, 
 $A_{{13}}+A_{{3}}+2A_{{1}}$, 
 $A_{{13}}+2A_{{2}}+A_{{1}}$, 
 $A_{{13}}+A_{{2}}+3A_{{1}}$, 
 $A_{{11}}+A_{{5}}+2A_{{1}}$, 
 $A_{{11}}+A_{{4}}+3A_{{1}}$, 
 $A_{{11}}+A_{{3}}+A_{{2}}+2A_{{1}}$, 
 $A_{{9}}+A_{{6}}+3A_{{1}}$, 
 $A_{{9}}+A_{{5}}+A_{{4}}$, 
 $A_{{9}}+A_{{5}}+A_{{3}}+A_{{1}}$, 
 $A_{{9}}+A_{{5}}+A_{{2}}+2A_{{1}}$, 
 $A_{{9}}+A_{{4}}+A_{{3}}+2A_{{1}}$, 
 $A_{{9}}+A_{{4}}+A_{{2}}+3A_{{1}}$, 
 $A_{{9}}+2A_{{3}}+A_{{2}}+A_{{1}}$, 
 $A_{{9}}+A_{{3}}+2A_{{2}}+2A_{{1}}$, 
 $2A_{{7}}+2A_{{2}}$, 
 $A_{{7}}+A_{{6}}+A_{{3}}+2A_{{1}}$, 
 $A_{{7}}+2A_{{5}}+A_{{1}}$, 
 $A_{{7}}+A_{{5}}+A_{{4}}+2A_{{1}}$, 
 $A_{{7}}+A_{{5}}+A_{{3}}+A_{{2}}+A_{{1}}$, 
 $A_{{7}}+A_{{4}}+A_{{3}}+A_{{2}}+2A_{{1}}$, 
 $A_{{6}}+2A_{{5}}+2A_{{1}}$, 
 $A_{{6}}+A_{{5}}+2A_{{3}}+A_{{1}}$, 
 $2A_{{5}}+A_{{4}}+A_{{3}}+A_{{1}}$, 
 $A_{{5}}+A_{{4}}+2A_{{3}}+A_{{2}}+A_{{1}}$.
}%
\text{\rm (2)}
Let $\SSS\sp{[2]}$ be the union of $\PPP\sp{[2]}\sb{18}$ 
and the following list.
\listup{%
$2E_{{7}}+A_{{3}}$, 
$E_{{7}}+D_{{10}}$, 
$E_{{7}}+D_{{5}}+A_{{5}}$, 
$D_{{12}}+D_{{5}}$, 
$2D_{{8}}+A_{{1}}$, 
$D_{{7}}+A_{{9}}+A_{{1}}$, 
$D_{{7}}+2A_{{5}}$, 
$D_{{6}}+A_{{11}}$, 
$2D_{{5}}+A_{{7}}$, 
$A_{{15}}+A_{{2}}$, 
$E_{{7}}+A_{{9}}$, 
$D_{{16}}$, 
$2D_{{8}}$, 
$D_{{5}}+A_{{11}}$, 
$A_{{15}}$.
}%
Then 
an $ADE$-type $\Sigma$  is a member of $\PPP\sp{[2]}$
if and only if 
$\Sigma$ is 
a member of $\SSS\sp{[2]}$ or obtained
from a member of $\SSS\sp{[2]}$ by
applying substitutions listed in 
Table~\ref{table:subs}
 but not in Table~\ref{table:forbiddensubs2}.
\end{theorem}
\begin{theorem}\label{thm:three}
\text{\rm (1)}
The list $\PPP\sp{[3]}\sb{18}$ consists of $19$ elements listed below.
\listup{%
$3E_{{6}}$, 
$2E_{{6}}+A_{{5}}+A_{{1}}$, 
$E_{{6}}+A_{{11}}+A_{{1}}$, 
$E_{{6}}+A_{{8}}+2A_{{2}}$, 
$E_{{6}}+A_{{8}}+A_{{2}}+2A_{{1}}$,
$E_{{6}}+2A_{ {5}}+A_{{2}}$, 
$E_{{6}}+A_{{5}}+A_{{3}}+2A_{{2}}$,
$A_{{17}}+A_{{1}}$, 
$A_{{ 14}}+2A_{{2}}$, 
$A_{{14}}+A_{{2}}+2A_{{1}}$,
$A_{{11}}+A_{{5}}+A_{{2}}$, 
$A _{{11}}+A_{{3}}+2A_{{2}}$,
$A_{{11}}+3A_{{2}}+A_{{1}}$, 
$2A_{{8}}+2A _{{1}}$,
$A_{{8}}+A_{{5}}+A_{{3}}+A_{{2}}$, 
$A_{{8}}+A_{{5}}+2A_{{2}}+A_{{ 1}}$,
$A_{{8}}+A_{{4}}+3A_{{2}}$, 
$A_{{8}}+A_{{3}}+3A_{{2}}+A_{{1}}$, 
$2A_{{5}}+A_{{4}}+2A_{{2}}$.
}%
\text{\rm (2)}
Let $\SSS\sp{[3]}$ be  $\PPP\sp{[3]}\sb{18}$.
Then 
an $ADE$-type $\Sigma$  is a member of $\PPP\sp{[3]}$
if and only if 
$\Sigma$ is a member of $\SSS\sp{[3]}$ or obtained
from a member of $\SSS\sp{[3]}$ by
applying substitutions listed in 
Table~\ref{table:subs}
 but not in Table~\ref{table:forbiddensubs3}.
\end{theorem}
\begin{theorem}\label{thm:four}
\text{\rm (1)}
The list $\PPP\sp{[4]}\sb{18}$ consists of $11$ elements listed below.
\listup{%
$D_{{7}}+A_{{11}}$, 
$D_{{7}}+A_{{7}}+A_{{3}}+A_{{1}}$, 
$D_{{7}}+3A_{{3}}+A_{{2}}$, 
$2D_{{5}}+A_{{7}}+A_{{1}}$,
$D_{{5}}+A_{{11}}+2A_{{1}}$, 
$D_{{5}}+ A_{{7}}+A_{{3}}+A_{{2}}+A_{{1}}$,
$A_{{15}}+A_{{3}}$, 
$A_{{15}}+3A_{{1}}$,
$A _{{11}}+2A_{{3}}+A_{{1}}$,
$A_{{11}}+A_{{3}}+A_{{2}}+2A_{{1}}$, 
$A_{{7}} +3A_{{3}}+A_{{2}}$.
}%
\text{\rm (2)}
Let $\SSS\sp{[4]}$ be the union of $\PPP\sp{[4]}\sb{18}$ 
and the following list.
\listup{%
$2D_{{5}}+A_{{7}}$, $A_{{15}}+2A_{{1}}$.
}%
Then 
an $ADE$-type $\Sigma$  is a member of $\PPP\sp{[4]}$
if and only if 
$\Sigma$ is a member of $\SSS\sp{[4]}$ or 
obtained
from a member of $\SSS\sp{[4]}$ by
applying substitutions listed in 
Table~\ref{table:subs}
 but not in Table~\ref{table:forbiddensubs4}.
\end{theorem}
\begin{theorem}\label{thm:twotwo}
\text{\rm (1)}
The list $\PPP\sp{[2,2]}\sb{18}$ consists of $11$ elements listed below.
\listup{%
$D_{{10}}+A_{{5}}+3A_{{1}}$, 
$D_{{10}}+2A_{{3}}+2A_{{1}}$, 
$2D_{{8}}+2A_{{1}}$, 
$D_{{8}}+D_{{6}}+A_{{3}}+A_{{1}}$,
$D_{{8}}+A_{{5}}+A_{{3}}+2 A_{{1}}$, 
$3D_{{6}}$, 
$2D_{{6}}+2A_{{3}}$,
$D_{{6}}+2A_{{5}}+2A_{{1 }}$, 
$D_{{6}}+A_{{5}}+2A_{{3}}+A_{{1}}$,
$A_{{7}}+A_{{5}}+A_{{3}}+3A_{{1 }}$, 
$2A_{{5}}+2A_{{3}}+2A_{{1}}$.
}%
\text{\rm (2)}
Let $\SSS\sp{[2,2]}$ be the union of $\PPP\sp{[2,2]}\sb{18}$ 
and the list
\listup{%
$4D_{{4}}$.
}%
Then 
an $ADE$-type $\Sigma$  is a member of $\PPP\sp{[2,2]}$
if and only if 
$\Sigma$ is a member of $\SSS\sp{[2,2]}$ or obtained
from a member of $\SSS\sp{[2,2]}$ by
applying substitutions listed in 
Table~\ref{table:subs}
 but not in Table~\ref{table:forbiddensubs22}.
\end{theorem}
By these theorems, we can easily generate the complete list $\PPP\sp G$
for $G=[1]$, $[2]$, $[3]$, $[4]$, $[2, 2]$.
Table~\ref{table:CardPPPGr}
shows the cardinalities of  $\PPP\sp G\sb r$.
\begin{table}
\caption{Cardinalities of $\PPP\sp G \sb r$}
\label{table:CardPPPGr}
\def\spaceheight{height 3pt}
\def\smallhsk{\hskip 1pt}
\def\HS{
\spaceheight & \omit  & \spaceheight &
\omit  & \omit  & \omit  &
\omit  & \omit  & \omit  &
\omit  & \omit  & \omit  &
\omit  & \omit  & \omit  &
\omit  & \omit  & \omit  &
\omit  & \omit  & \omit  &
 \spaceheight & \omit &\spaceheight
\cr
}
\def\HL{\HS\noalign{\hrule}\HS}
\centerline{
\vbox{\offinterlineskip
\hrule
\halign{
\vrule # \hskip 2pt&  \strut\smallhsk \hfil#\smallhsk  & \vrule # \hskip 2pt\hfill  &
\strut\smallhsk \hfil#\smallhsk  & \strut\smallhsk \hfil#\smallhsk  & \strut\smallhsk \hfil#\smallhsk  &
\strut\smallhsk \hfil#\smallhsk  & \strut\smallhsk \hfil#\smallhsk  & \strut\smallhsk \hfil#\smallhsk  &
\strut\smallhsk \hfil#\smallhsk  & \strut\smallhsk \hfil#\smallhsk  & \strut\smallhsk \hfil#\smallhsk  &
\strut\smallhsk \hfil#\smallhsk  & \strut\smallhsk \hfil#\smallhsk  & \strut\smallhsk \hfil#\smallhsk  &
\strut\smallhsk \hfil#\smallhsk  & \strut\smallhsk \hfil#\smallhsk  & \strut\smallhsk \hfil#\smallhsk  &
\strut\smallhsk \hfil#\smallhsk  & \strut\smallhsk \hfil#\smallhsk  & \strut\smallhsk \hfil#\smallhsk  &
\vrule # \hskip 2pt & \strut\smallhsk \hfil#\smallhsk  & \vrule #
\cr
\HS
& $r$\phantom{al} & &
$1$ \hfil & $2$ \hfil & $3$ \hfil &
$4$ \hfil & $5$ \hfil & $6$ \hfil &
$7$ \hfil & $8$ \hfil & $9$ \hfil &
$10$ \hfil & $11$ \hfil & $12$ \hfil &
$13$ \hfil & $14$ \hfil & $15$ \hfil &
$16$ \hfil & $17$ \hfil & $18$ \hfil &
 & total \hfil  & 
\cr
\HL
\HS
& $|\PPP\sp{[1]}\sb r|$ & &
$1$  & $2$  & $3$  &
$6$  & $9$  & $16$ &
$24$ & $39$  & $57$  &
$88$  & $127$  & $189$  &
$262$  & $360$  & $448$  &
$500$  & $416$  & $199$  &
 & $2746$ \hfil & 
\cr
\HS
& $|\PPP\sp{[2]}\sb r|$ & &
$0$ &$0$ &$0$ &$0$ &$0$ &$0$ &$0$ &$1$ &$2$ &
$6$ &$13$ &$29$ &$53$ &$92$ &$133$ &$164$ &$155$ &$84$ & &
$732$ \hfil&
\cr
\HS
& $|\PPP\sp{[3]}\sb r|$ & &
$0$ &$0$ &$0$ &$0$ &$0$ &$0$ &$0$ &$0$ &$0$ &
$0$ &$0$ &$1$ &$2$ &$6$ &$12$ &$21$ &$24$ &$19$ & &
$85$ \hfil&
\cr
\HS
& $|\PPP\sp{[4]}\sb r|$ & &
$0$ &$0$ &$0$ &$0$ &$0$ &$0$ &$0$ &$0$ &$0$ &
$0$ &$0$ &$0$ &$0$ &$1$ &$4$ &$10$ &$15$ &$11$ & &
$41$ \hfil&
\cr
\HS
& $|\PPP\sp{[2,2]}\sb r|$ & &
$0$ &$0$ &$0$ &$0$ &$0$ &$0$ &$0$ &$0$ &$0$ &
$0$ &$0$ &$1$ &$2$ &$5$ &$10$ &$16$ &$16$ &$11$ & &
$61$ \hfil&
\cr
\HS
}
\hrule
}
}
\end{table}
\par
\medskip
{\bf Case 2.} $G\in \{[5], [6], [4,2]\}$.
We simply give the table of $\PPP\sp G$.
In each box, the $ADE$-types are listed according to the rank and the lexicographical order.
\par
\medskip
\noindent
$G=[5]$:
\listup{%
$2A_{{9}}$, 
$A_{{9}}+2A_{{4}}+A_{{1}}$, 
$4A_{{4}}+2A_{{1}}$, 
$A_{{9}}+2A_{{4}}$, 
$4A_{{4}}+A_{{1}}$, 
$4A_{{4}}$.
}%
$G=[6]$:
\listup{%
$A_{{11}}+A_{{5}}+2A_{{1}}$, 
$A_{{11}}+A_{{3}}+2A_{{2}}$,
$A_{{11}}+2A_{{2}}+3A_{{1}}$, 
$3A_{{5}}+A_{{3}}$,
$2A_{{5}}+A_{{3}}+2A_{{2}}+A_{ {1}}$, 
$A_{{11}}+2A_{{2}}+2A_{{1}}$,
$3A_{{5}}+2A_{{1}}$, 
$2A_{{5}}+A _{{3}}+2A_{{2}}$,
$2A_{{5}}+2A_{{2}}+3A_{{1}}$, 
$2A_{{5}}+2A_{{2 }}+2A_{{1}}$.
}%
$G=[4,2]$:
\listup{%
$2A_{{7}}+4A_{{1}}$, 
$A_{{7}}+3A_{{3}}+2A_{{1}}$, 
$A_{{7}}+2A_{{3}}+4A_{{1}}$, 
$5A_{{3}}+2A_{{1}}$, 
$4A_{{3}}+4A_{{1}}$.
}%
\par
\medskip
{\bf Case 3.} $G\in \{[7], [8], [6, 2], [3, 3], [4, 4]\}$.
In this case,
the $ADE$-type determines the torsion of the  Mordell-Weil group uniquely.
\begin{theorem}
Let $f : X\to \P\sp 1$ be an elliptic $K3$ surface.
Then the following hold.
\begin{itemize}
\item
 $G\sb f\cong \Z/(7) \Longleftrightarrow  \Sigma\sb f = 3 A\sb 6$.
\item
 $G\sb f\cong \Z/(8) \Longleftrightarrow  \Sigma\sb f =  2 A\sb 7 + A\sb 3 + A\sb 1$.
\item
 $G\sb f\cong \Z/(6)\times \Z/(2) \Longleftrightarrow  \Sigma\sb f =  3 A\sb 5 + 3 A\sb 1$.
\item
 $G\sb f\cong \Z/(4)\times \Z /(4) \Longleftrightarrow  \Sigma\sb f = 6 A\sb 3$.
\item
 $G\sb f\cong \Z/(3)\times \Z /(3) \Longleftrightarrow  \Sigma\sb f \in \{  
2 A\sb 5 + 4 A\sb 2 , A\sb 5 + 6 A\sb 2, 8 A\sb 2 \}$.
\end{itemize}
\end{theorem}
\begin{remark}
Elliptic $K3$ surfaces 
with
$G\sb f=[7], [8], [6,2], [4,4]$
are constructed as elliptic modular surfaces
(cf.\ \cite{ShiodaEllipticModular, Sebbar}).
The corresponding 
congruence groups $\Gamma \subset \SL \sb 2 (\Z)$  are as follows.
$$
\def\spaceheight{height 3pt}
\def\smallhsk{\hskip 1pt}
\def\HS{
\spaceheight & \omit  & \spaceheight  & \spaceheight&
\omit  & \spaceheight &\omit  &\spaceheight & \omit  &\spaceheight &\omit  &\spaceheight
\cr
}
\def\HL{\HS\noalign{\hrule}\HS}
\centerline{
\vbox{\offinterlineskip
\hrule
\halign{
\vrule # \hskip 5pt
&  \strut\smallhsk \hfil#\smallhsk \hfil & \vrule # \hskip 1pt\hfill  & \vrule # \hskip 5pt\hfill  &
\hfil #\hfil &
\hskip 5pt\vrule #\hskip 5pt  & &
\hfil #\hfil &
\hskip 5pt\vrule #\hskip 5pt  &
\hfil #\hfil &
\hskip 5pt\vrule #\hskip 5pt  &
\hfil#\hfil &
\hskip 5pt\vrule #
\cr
\HS
&  $G\sb f$  &  &&
$[7]$ & & $[8]$ & & $[6, 2]$ & & $[4, 4]$ & 
\cr
\HL
\HS
&  $\Gamma$  &  &&
$\Gamma\sb 1 (7) $ & & $\Gamma\sb 1 (8) $ & & $\Gamma\sb 0 (3) \cap \Gamma (2) $ & & $\Gamma (4)$ & 
\cr
\HS
}
\hrule
}
}
$$
\end{remark}
\subsection{From $ADE$-types to configurations of singular fibers}
\label{subsec:recover}
\begin{table}
\caption{Singular fibers of elliptic fibration.}
\label{table:kodaira}
\def\spaceheight{height 3pt} %Here is the papameter of the vspace above and below of the hrule.
\def\smallspaceheight{height 1pt}
\def\HS{\spaceheight &
\omit& \spaceheight& 
\omit& \spaceheight& 
\omit& \spaceheight&
\omit& \spaceheight
\cr
}
\def\smallHS{\smallspaceheight &
\omit& \smallspaceheight& 
\omit& \smallspaceheight& 
\omit& \smallspaceheight&
\omit& \smallspaceheight
\cr
}
\def\BB{\phantom{\quad (b\ge 2)}}
\def\HL{\HS\noalign{\hrule}\HS}
\centerline{
\vbox{\offinterlineskip
\hrule
\halign{
\vrule # & 
\strut\quad\hfil#\hfil\quad &\vrule# &
\strut\quad\hfil#\hfil\quad &\vrule# &
\strut\quad\hfil#\hfil\quad &\vrule# &
\strut\quad\hfil#\hfil\quad &\vrule# 
\cr
\HS
 &
Singular fiber & &
$ADE$-type  & &
Euler number  & & Possible tortion parts&
\cr
\HS\noalign{\hrule}\smallHS\noalign{\hrule}\HS
&
${\rm I}\sb 0\BB$ &&
regular &&
$0$ &&
all &
\cr
\HL
&
${\rm I}\sb 1\BB$ &&
irreducible &&
$1$ &&
 &
\cr
\HS
&
\multispan{5}\hrulefill&&
\omit\vbox  to 0pt { \hbox{\hskip 68pt \raise 15pt\hbox{ $\diamondsuit$}}\vss }&
\cr
\HS
\cr
&
${\rm I}\sb b\quad (b\ge 2)$ &&
$A\sb{b-1}$ &&
$b$ &&
 &
\cr
\HL
&
${\rm I}\sp *\sb b \quad (b\ge 0)$ &&
$D\sb{4+b}$ &&
$6+b$ &&
\hbox{
$
\begin{cases}
[1], [2], [2,2] & \hbox{\rm if $b$ is even} \cr
[1], [2], [4] & \hbox{\rm if $b$ is odd} \cr
\end{cases}
$
}
&
\cr
\HL
&
${\rm II}\BB$ &&
irreducible &&
$2$ &&
[1] &
\cr
\HL
&
\rlap{${\rm II}\sp *$}\phantom{${\rm II}\BB$} &&
$E\sb 8$ &&
$10$ &&
[1] &
\cr
\HL
&
${\rm III}\BB$ &&
$A\sb 1$ &&
$3$ &&
$[1], [2]$ &
\cr
\HL
&
\rlap{${\rm III}\sp *$}\phantom{${\rm III}\BB$} &&
$E\sb 7$ &&
$9$ &&
$[1], [2]$ &
\cr
\HL
&
${\rm IV}\BB$ &&
$A\sb 2$ &&
$4$ &&
$[1], [3]$ &
\cr
\HL
&
\rlap{${\rm IV}\sp *$}\phantom{${\rm IV}\BB$} &&
$E\sb 6$ &&
$8$ &&
$[1], [3]$ &
\cr
\HS
}
\hrule
}
}
$$
\diamondsuit\left\{
\vcenter{
\hbox{
\vbox{
\hbox{ $[a]$ is possible for $a=1, \dots, 8$,}\vskip 3pt
\hbox{ $[2a, 2]$ is possible for $a=1, \dots, 3$ if and only if $b=0 \bmod 2$,}\vskip 3pt
\hbox{ $[3,3]$ is possible if and only if $b=0 \bmod 3$,}\vskip 3pt
\hbox{ $[4,4]$ is possible  if and only if $b=0 \bmod 4$.}
}
}
}
\right.
$$
\vskip 12pt 
\hrule
\end{table}
The correspondence between the type (in the  notation of Kodaira) 
of a singular fiber of an elliptic fibration
and an  $ADE$-type is shown in Table~\ref{table:kodaira}.
There are following ambiguities in recovering the
configurations of singular fibers 
from its $ADE$-type.
\begin{itemize}
\item
An irreducible singular fiber is of type either ${\rm I}\sb 1$ or ${\rm II}$. 
\item
A singular fiber of $ADE$-type $A\sb 1$ is of type either ${\rm I}\sb 2$ or ${\rm III}$. 
\item
A singular fiber of $ADE$-type $A\sb 2$ is of type either ${\rm I}\sb 3$ or ${\rm IV}$. 
\end{itemize}
We present some restrictions
on the possibilities of
configuration of singular fibers of
an elliptic $K3$ surface $f: X\to \P\sp 1$
with a given $ADE$-type.
\par
\medskip
Let $i\sb b$ be the number of singular fibers of $f$ of type ${\rm I}\sb b$.
We define similarly $i\sb b \sp *$, $ii$, $ii\sp *$, $iii$,  $iii\sp *$,
$iv$, $iv\sp *$.
Miranda and Persson gave a formula
for the degree of the modulus function
$J\sb f : \P\sp 1 \to \P\sp 1 :=\UpperHalfPlane/\SL\sb 2 (\Z)$
associated with $f : X\to \P\sp 1$:
$$
\deg J\sb f:=\sum\sb{b\ge 1} b(i\sb b + i\sb b\sp *).
$$
By the Hurwitz formula,
they obtained the following necessary condition for configurarions;
if $\deg J\sb f >0$, then 
$$
\deg J\sb f \le 6\sum\sb{b\ge 1} (i\sb b + i\sp *\sb b) + 4(ii+iv\sp *)+3(iii+iii\sp *)+2(iv+ii\sp *)-12.
$$
See \cite[\S 3]{MirandaPerssonRational} for the proof.
\par
\medskip
The euler number $24$ of the $K3$ surface $X$ is equal to
the sum of euler numbers of singular fibers of $f$.
The third column of Table~\ref{table:kodaira} shows the euler number of
a singular fiber of each type.
We define the euler number $\euler (\Sigma)$ of
an $ADE$-type $\Sigma :=\sum a\sb l A\sb l + \sum d\sb m D\sb m + \sum e\sb n E\sb n$
by
$$
\euler (\Sigma) :=  \sum a\sb l \cdot (l+1)  + \sum d\sb m\cdot (m+2) + \sum e\sb n \cdot (n+2).
$$
Then $\euler (\Sigma\sb f)$ is less than or equal to 
the sum of euler numbers of reducible singular fibers.
Hence we always have
$$
\euler (\Sigma\sb f)\le 24.
$$
We can deduce from Table~\ref{table:kodaira},
for example,
that,
if $\euler (\Sigma\sb f) =24$,
then$f: X\to \P\sp 1$  has no
irreducible fibers
nor fibers of type ${\rm III}$ or ${\rm IV}$.
\par
\medskip
When $G\sb f$ is non-trivial,
certain types of singular fibers cannot appear.
Let $g : S\to \Delta$ be an elliptic fibration
over an open unit disk $\Delta$
such that $g$ is smooth over $\Delta\sp\times :=\Delta\setminus\{0\}$,
and let $E:=g\inv (p)$  be the fiber over a point $p\in \Delta\sp\times$.
Looking at the monodromy action of $\pi\sb 1 (\Delta\sp\times, p)$ on the set of
torsion points of $E$,
we can determine whether a finite abelian group
can be embedded into the Mordell-Weil group of $g$.
The fourth column
of Table~\ref{table:kodaira} shows
the groups among the list~\eqref{Gset} that can be 
isomorphic to the torsion part of the Mordell-Weil group
of an elliptic surface having 
the singular fiber.
We see, for example,
that,
if $G\sb f$ is non-trivial,
then every irreducible singular fiber must be
of type ${\rm I}\sb 1$.
\subsection{Miscellaneous facts}
For an integer $r$ with $1\le r\le 18$,
we put as follows.
\begin{align*}
\RRR\sb r &:= \{ \; \Sigma \; ;  \;  \text{\rm $\Sigma$ is an $ADE$-type with $\rank (\Sigma)=r$ } \; \}, \\
\EEE\sb r &:= \{ \; \Sigma \in \RRR\sb r  \; ;  \;  \euler (\Sigma)\le 24 \; \}, \quad\text{\rm and } \\
\PPP\sb r &:= \{ \; \Sigma \in \EEE\sb r  \; ;  \; 
\text{\rm there exists an elliptic $K3$ surface $f : X\to \P\sp 1$ with $\Sigma\sb f=\Sigma$} \; \}
=\cup\sb G \PPP\sb r \sp G.
\end{align*}
For $\Sigma \in \cup\sb {r=1}\sp{18}\;  \PPP\sb r$, we 
denote by $\GGG (\Sigma)$ the set of isomorphism classes
of finite abelian groups $G$ such that
$(\Sigma, G) \in \PPP$.
For each $r$, we denote by $\TTT\sb r$ 
 the set of $\Sigma \in \PPP\sb r$
such that $\GGG (\Sigma ) $ consists of only the trivial group $[1]$.
The cardinalities of these sets are given in Table~\ref{table:REK}.
\begin{table}
\caption{Cardinalities of $\RRR\sb r$, $\EEE\sb r$ and $\PPP\sb r$}
\label{table:REK}
\def\spaceheight{height 3pt}
\def\smallhsk{\hskip 1pt}
\def\HS{
\spaceheight & \omit  & \spaceheight &
\omit  & \omit  & \omit  &
\omit  & \omit  & \omit  &
\omit  & \omit  & \omit  &
\omit  & \omit  & \omit  &
\omit  & \omit  & \omit  &
\omit  & \omit  & \omit  &
 \spaceheight & \omit &\spaceheight
\cr
}
\def\HL{\HS\noalign{\hrule}\HS}
\centerline{
\vbox{\offinterlineskip
\hrule
\halign{
\vrule # \hskip 2pt& \strut\smallhsk \hfil#\smallhsk  & \vrule # \hskip 2pt &
\strut\smallhsk \hfil#\smallhsk  & \strut\smallhsk \hfil#\smallhsk  & \strut\smallhsk \hfil#\smallhsk  &
\strut\smallhsk \hfil#\smallhsk  & \strut\smallhsk \hfil#\smallhsk  & \strut\smallhsk \hfil#\smallhsk  &
\strut\smallhsk \hfil#\smallhsk  & \strut\smallhsk \hfil#\smallhsk  & \strut\smallhsk \hfil#\smallhsk  &
\strut\smallhsk \hfil#\smallhsk  & \strut\smallhsk \hfil#\smallhsk  & \strut\smallhsk \hfil#\smallhsk  &
\strut\smallhsk \hfil#\smallhsk  & \strut\smallhsk \hfil#\smallhsk  & \strut\smallhsk \hfil#\smallhsk  &
\strut\smallhsk \hfil#\smallhsk  & \strut\smallhsk \hfil#\smallhsk  & \strut\smallhsk \hfil#\smallhsk  &
\vrule # \hskip 2pt & \strut\smallhsk \hfil#\smallhsk  & \vrule #
\cr
\HS
& $r$\phantom{al} & &
$1$ \hfil & $2$ \hfil & $3$ \hfil &
$4$ \hfil & $5$ \hfil & $6$ \hfil &
$7$ \hfil & $8$ \hfil & $9$ \hfil &
$10$ \hfil & $11$ \hfil & $12$ \hfil &
$13$ \hfil & $14$ \hfil & $15$ \hfil &
$16$ \hfil & $17$ \hfil & $18$ \hfil &
 & total \hfil  & 
\cr
\HL
& $|\RRR\sb r|$ & &
$1$  & $2$  & $3$  &
$6$  & $9$  & $16$ &
$24$ & $39$  & $57$  &
$88$  & $128$  & $193$  &
$276$  & $403$  & $570$  &
$815$  & $1137$  & $1599$  &
 & $5366 $ &
\cr
\HS
& $|\hskip 2pt \EEE\sb r \hskip .8pt |$ & &
$1$  & $2$  & $3$  &
$6$  & $9$  & $16$ &
$24$ & $39$  & $57$  &
$88$  & $128$  & $193$  &
$274$  & $393$  & $531$  &
$688$  & $773$  & $712$  &
 & $3937$ & 
\cr
\HS
& $|\hskip .6pt \PPP\sb r|$ & &
$1$  & $2$  & $3$  &
$6$  & $9$  & $16$ &
$24$ & $39$  & $57$  &
$88$  & $128$  & $193$  &
$274$  & $392$  & $518$  &
$624$  & $580$  & $325$  &
 & $3279$ &  
\cr
\HS
& $|\hskip 1.4pt \TTT\sb r|$ & &
$1$  & $2$  & $3$  &
$6$  & $9$  & $16$ &
$24$ & $38$  & $55$  &
$82$  & $115$  & $162$  &
$217$  & $289$  & $362$  &
$419$  & $372$  & $188$  &
 & $2360$ & 
\cr
\HS
}
\hrule
}
}
\end{table}
Note that, if $\rank (\Sigma) \le 12$, then
$\euler (\Sigma) \le 24$ holds automatically.
\begin{theorem}
Let $\Sigma$ be an $ADE$-type with $\euler (\Sigma) \le 24$.
Suppose that 
$\rank (\Sigma) \le 13$.
Then there exists an elliptic $K3$ surface $f:X\to \P\sp 1$ with $\Sigma\sb f=\Sigma$.
\end{theorem}
\begin{remark}
The complement of $ \PPP\sb{14}$ in $\EEE\sb{14}$ consists of a single element $E\sb 6 + 8 A\sb 1$.
Hence, when $\euler (\Sigma) \le 24$ and $\rank (\Sigma) =14$,
there exists an elliptic $K3$ surface $f:X\to \P\sp 1$ with $\Sigma\sb f=\Sigma$
if and only if $\Sigma \ne  E\sb 6 + 8\; A\sb 1$.
\end{remark}
\begin{theorem}
Suppose that 
$\rank (\Sigma) \le 10$.
Then there exists an elliptic $K3$ surface $f:X\to \P\sp 1$ 
with $G\sb f = [1]$   and 
$\Sigma\sb f=\Sigma$.
\end{theorem}
\begin{remark}
The complement  $\PPP\sb{11}\setminus \PPP\sp{[1]} \sb{11} $ consists of a single element $11 A\sb 1$.
We have $\GGG (11A\sb 1) =\{ [2] \}$.
\end{remark}
\begin{theorem}
Let $f : X\to \P\sp 1$ be an elliptic $K3$ surface.
If $\rank (\Sigma\sb f) \le 7$, then $G\sb f$ must be trivial.
\end{theorem}
\begin{remark}
The complement  $\PPP\sb{8}\sp{[1]}\setminus \TTT \sb{8} $ consists of a single element $8 A\sb 1$, and 
the complement  $\PPP\sb{9}\sp{[1]} \setminus \TTT\sb{9} $ consists of two  elements $9 A\sb 1$ and $A\sb 3 + 6
A\sb 1$. We have 
$$
\GGG (8 A\sb 1) = \GGG (9 A\sb 1) = \GGG (A\sb 3 + 6 A\sb 1) =\{[1], [2]\}.
$$
\end{remark}
\begin{remark}
There are several  $ADE$-types  $\Sigma$ 
with $|\GGG (\Sigma)| \ge 3$.
For example,
$$
\GGG(2 A\sb 5 + 2 A\sb 2 + 2 A\sb 1 )=
\GGG(A\sb{11}+2 A\sb 2 + 2A\sb 1)
=\{ [1], [2], [3], [6]\}.
$$
\end{remark}
\section{Local invariants of lattices}\label{sec:local}
First we fix some terminologies about lattices.
\par
\medskip
Let $R$ be either $\Z$ or $\Zp$.
A {\it lattice} over $R$ is, by definition,
a free $R$-module $L$  of finite rank equipped with a non-degenerate 
symmetric bilinear form $(\phantom{a}, \phantom{a}) : L\times L \to R$.
For $\alpha \in R \setminus \{ 0 \}$, let $\alpha L$ denote the lattice obtained from $L$ by
multiplying the symmetric bilinear form by $\alpha$.
We will  denote $\opp{L}$ for $(-1)\thinspace L$.
We often express a lattice by the intersection matrix
with respect to a certain basis of $L$.
For example, $(a)$ is the lattice of rank $1$ generated by a vector $e$
such that $(e, e)=a$.
A sublattice $N$ of $L$ is said to be  {\it primitive} 
if $L/N$ is torsion free.
A lattice $L$ over $R$ is said to be  {\it even}
if $(v, v)\in 2\thinspace R$ holds for any $v\in L$.
Note that, when $R$ is $\Zp$ with $p$ an odd prime,
every lattice over $R$ is even.
The {\it discriminant} $\disc (L)$ of a lattice $L$
is considered as an element of $(R\setminus \{ 0 \})/\utsq{R}$.
A lattice $L$ is said to be  {\it unimodular} if $\disc (L) \in R\sp\times/\utsq{R}$.
\par
\medskip
Suppose that $R=\Zp$.
Then we have $\disc (L) = p\sp{\nu} u$ for some $\nu \ge 0$, where $u \in \Zp\sp\times /\utsq{\Zp}$.
We denote the element $u$ by $\reddisc (L)$ 
and call it the  {\it reduced discriminant} of $L$.
\par
\medskip
Let $k$ be the quotient field of $R$.
The $k$-vector space $L\otimes\sb R k$ 
has a natural symmetric bilinear form with values in $k$.
We denote by 
$\dual{L}$
the $R$-submodule of $L\otimes\sb R  k$ consisting of
all vectors $v$ such that $(v, w)\in R$ holds
for every $w\in L$,
and call it the  {\it dual lattice} of $L$.
An $R$-submodule $M$ of $\dual{L}$ is said to be 
an {\it overlattice} of $L$ if $M$ contains $L$
and the symmetric  bilinear form restricted to $M$ takes values in $R$.
Two lattices $L$ and $M$ over $R$ are said to be {\it  $k$-equivalent}
if $L\otimes\sb R k$ and $M\otimes\sb R k$ together
with their symmetric $k$-valued bilinear forms are isomorphic.
\par
\medskip
For a detailed  account of the following definitions and theorems,
see Conway and Sloane~\cite[Chapter 15]{ConwaySloane} and  Cassels~\cite[Chapters 8 and 9]{Cassels}.
\subsection{Local invariants}
Let $\Lambda$ be a lattice over $\Zp$.
Then $\Lambda$ is decomposed into the orthogonal direct 
sum
$\Lambda=\bigoplus\sb{\nu\ge 0} p\sp\nu \Lambda\sb\nu$
with each  $\Lambda\sb\nu$ being unimodular.
This  decomposition is called a {\it Jordan decomposition}
of $\Lambda$,
and each $p\sp\nu \Lambda\sb\nu$ is called a {\it Jordan component} of $\Lambda$.
Note that the reduced discriminant of $\Lambda$ is the product
of the discriminants of $\Lambda\sb\nu$.
\par
\medskip
Suppose that $p$ is odd.
Then a lattice $\Lambda$ over $\Zp$ is
isomorphic to  an orthogonal direct sum
$\oplus\sb{i} p\sp{\nu\sb i} (a\sb i) $,
where $a\sb i \in \Zp\sp\times$.
The {\it $p$-excess} of $\Lambda$ is defined to be
$$
-\rank (\Lambda) + 4 m + \sum\sb i  p\sp{\nu\sb i} \in \Z/ (8),
$$
where $m$ is the number of  orthogonal direct summands $p\sp{\nu\sb i} (a\sb i)$
such that $\nu\sb i $ is odd and that $a\sb i $ is not square in $\Zp\sp\times$.
It is known that the $p$-excess is a well-defined invariant of $\Qp$-equivalence classes of lattices over $\Zp$.
\par
\medskip
Suppose that $p=2$.
We put
$$
U:=\begin{pmatrix} 0 & 1 \\ 1 & 0 \end{pmatrix} \quad\text{\rm and }\quad
V:=\begin{pmatrix} 2 & 1 \\ 1 & 2 \end{pmatrix},
$$
both of which are even unimodular lattices of rank $2$ over $\Zt$.
Then a lattice over $\Zt$ is decomposed into the orthogonal direct sum of
 lattices such that each direct summand is isomorphic to 
$2\sp\nu (a)$ $(a\in \Zt\sp\times)$,
$2\sp\nu  U$ or $2\sp\nu  V$.
We define the $2$-excesses of these lattices
by
\begin{align*}
\texcess (2\sp\nu (a)) &=
\begin{cases}
1-a \bmod 8 & \text {\rm if $\nu$ is even or $a=\pm 1 \bmod 8$, } \\
5-a  \bmod 8 & \text {\rm if $\nu$ is odd and $a=\pm 3 \bmod 8$, }
\end{cases}
\\
\texcess (2\sp\nu U) &=\;\;2  \bmod 8 , \\
\texcess (2\sp\nu V) &= 
\begin{cases}
2  \bmod 8 & \text {\rm if $\nu$ is even, } \\
6  \bmod 8 & \text {\rm if $\nu$ is odd.}
\end{cases}
\end{align*}
Then we define the $2$-excess of
\begin{equation}\label{twodecomp}
\Lambda \cong \bigoplus\sb i 2\sp{\nu\sb i} (a\sb i) \oplus \bigoplus\sb j 2\sp{\nu\sb j} U 
\oplus \bigoplus\sb k 2\sp{\nu\sb k} V
\end{equation}
to be the sum of the $2$-excesses of direct summands in the decomposition~\eqref{twodecomp}.
Even though the decomposition~\eqref{twodecomp} is not unique in general,
it turns out that  the $2$-excess is a well-defined invariant of $\Qt$-equivalence classes of lattices over $\Zt$.
(Note that $U$ and $V$ are $\Qt$-equivalent to $2 (1)\oplus 2 (7)$ and
$2 (1)\oplus 2 (3)$, respectively.)
\subsection{Existence of lattices over $\Z$ with given local data}\label{sec:existencelattice}
Combining~\cite[Chapter 15, Theorem~5]{ConwaySloane} and~\cite[Chapter~9, Theorem~1.2]{Cassels},
we obtain the following:
\begin{theorem}\label{existence}
Let $d$ be a non-zero integer,
and $(r, s)$ a pair of non-negative integers such that  $n:=r+s$ is positive and that $d=(-1)\sp s |d|$ holds.
Suppose that, for each prime divisor $p$ of $2d$,
a lattice $\Lambda\sp{(p)}$ of rank $n$ over $\Zp$
is given.
Then there exists a lattice $L$ over $\Z$ with discriminant $d$ and signature $(r, s)$ such that
$L\otimes\sb{\Z} \Zp$ is isomorphic to $\Lambda\sp{(p)}$ for each $p$ if and only if
the following two conditions are satisfied{\rm :}
\begin{enumerate}
\renewcommand{\labelenumi}{(\roman{enumi})}
\item $\disc (\Lambda\sp{(p)})$ is equal to $d\cdot \utsq{\Zp}$ for each $p$, 
and
\item $r-s+\sum\sb{p|2d} \pexcess (\Lambda\sp{(p)}) = n \bmod 8$ holds. \hfill \qed
\end{enumerate}
\end{theorem}
\section{Theory of discriminant forms}\label{sec:theorydiscfrm}
\subsection{Definitions}
Let $R$ and $k$ be as above.
Let $D$ be a finite abelian group.
A finite symmetric bilinear form on $D$ with values in $k/R$ is,
by definition, a homomorphism $b : D\times D \to k/R$
such that $b(x, y)=b(y, x)$ holds for any $x, y \in D$.
A finite quadratic form on $D$ with values in $k/2R$ is a map $q : D \to k/2R$ with the following properties:
\begin{itemize}
\item[(i)] $q(nx)=n\sp 2 q(x)$ for $n\in \Z$ and $x\in D$,
and 
\item[(ii)] the map $b[q] : D\times D \to k/R$ defined by
$(x, y) \mapsto  (q(x+y)-q(x)-q(y))/2$ 
is a finite symmetric bilinear form.
\end{itemize}
Let $H$ be a subgroup of $D$. The orthogonal complement $H\sp\bot$ of $H$
with respect to $q$ is the subgroup of $D$ consisting of
elements $y$ such that $b[q] (x, y)=0$ holds for any $x\in H$.
We say that $q$ is non-degenerate if $D\sp\bot=(0)$.
Note that, if $D=H\oplus H\sp\bot$, then $q$ is written as $q|\sb H \oplus q|\sb{H\sp\bot}$,
because the homomorphism  $a\mapsto a/2$ from $k/2R$ to $k/R$ is injective.
\par
\medskip
The length of $D$ is,
by definition,
the minimal number of generators of $D$.
A subset $\{\gamma\sb 1, \dots, \gamma\sb l\}$
of $D$ is said to be a {\it reduced set of generators} of $D$
if $l$ is the length of $D$ and 
$D=\langle \gamma\sb 1 \rangle \times \cdots \times \langle \gamma\sb l \rangle$
holds.
Let  $\{\gamma\sb 1, \dots, \gamma\sb l\}$
 be a  reduced set of generators of $D$.
Then a finite quadratic form $q$ on $D$ is expressed 
by a symmetric $l\times l$ matrix
whose diagonal entries are $q(\gamma\sb i) \in k/2R$ and 
whose off-diagonal entries are $b[q](\gamma\sb i, \gamma\sb j) \in k/R$.
\par
\medskip
Let $L$ be a lattice over $R$.
The discriminant group $D\sb L$ of $L$ is,
by definition, the quotient group $\dual{L} /L$.
We denote by $\Psi\sb L : \dual{L}\to D\sb L$ the natural projection.
Suppose that $L$ is even.
Then we can define a finite quadratic form $q\sb L$ on $D\sb L$ with values in $k/2R$ by
$q\sb L (x):= (x\sp\prime, x\sp\prime )\bmod 2R$,
where $x\sp\prime$ is a vector of $\dual {L}$ such that $\Psi \sb L (x\sp\prime)=x$.
We call $q\sb L$ the 
{\it discriminant form} of $L$.
Because $L$ is non-degenerate, 
$q\sb L$ is also non-degenerate.
By definition,
we have
$(D\sb{L\oplus M}, q\sb{L\oplus M})=(D\sb L, q\sb L) \oplus (D\sb M, q\sb M)$.
\subsection{Discriminant forms and overlattices}
The following two propositions, due to Nikulin, play a central role in making the list $\PPP$.
\begin{proposition}[\cite{Nikulin}~Proposition 1.4.1]\label{nikulin1}
Let $L$ be an even lattice over $\Z$.
\par
{\rm (1)}
If $H\subset D\sb L$ is a  subgroup isotopic
with respect to $q\sb L$,
then $M:=\Psi\sb L \sp{-1} (H)$ is an even overlattice of $L$,
and
the discriminant form of $M$ is isomorphic to $(H\sp\bot /H, q\sb L |\sb {H\sp\bot /H})$.
\par
{\rm (2)} The map $H\mapsto \Psi\sb L \sp{-1} (H)$ establishes a bijection
between the set of isotopic subgroups of $(D\sb L, q\sb L)$ and the set of even
overlattices of $L$.
\qed
\end{proposition}
\begin{proposition}[\cite{Nikulin}~Proposition 1.6.1]\label{nikulin2}
Let $L$ and $M$ be even lattices over $\Z$.
Then the following are equivalent.
\begin{itemize}
\item[{\rm (i)}]
The two  finite quadratic forms $(D\sb L, q\sb L)$ and $(D\sb M, -q\sb M)$ are isomorphic.
\item[{\rm (ii)}]
 There exists an even unimodular overlattice of $L\oplus M$ into
which $L$ and $M$ are  embedded primitively. 
\end{itemize}
\qed
\end{proposition}
\subsection{Localization and  discriminant form}
Let $L$ be an even lattice over $\Z$.
We decompose $D\sb L$ into the direct sum of 
its $p$-Sylow subgroups $D\sb L \sp{(p)}$,
where $p$ runs through the set of prime divisors of
$|D\sb L | = | \disc (L) |$.
These $p$-parts are orthogonal to each other with 
respect to $q\sb L$,
and hence $q\sb L$ is also decomposed into the $p$-parts;
$q\sb L = \oplus \sb p q\sb L \sp{(p)}$,
where $q\sb L \sp{(p)}$ is the restriction of $q\sb L$ to $D\sb L \sp{(p)}$.
By the definition of the discriminant form,
we can easily prove the following:
\begin{lemma}\label{qdecomp}
The image of $q\sb L \sp{(p)}$ is contained in $2\thinspace\Z [1/p] / 2 \thinspace\Z \subset \Q/ 2 \thinspace \Z$.
The natural inclusion $ 2\thinspace \Z [1/p] \hookrightarrow \Qp$
 induces an isomorphism 
$2\thinspace\Z [1/p] / 2 \thinspace\Z \cong \Qp / 2\thinspace\Zp$
Under this identification,
$(D\sb L\sp{(p)}, q\sb L\sp{(p)})$ is isomorphic to
$ (D\sb{L\otimes \Zp}, q\sb{L\otimes\Zp})$.
\qed
\end{lemma}
The discriminant form  of an even  lattice $\Lambda$ over $\Zp$ is
calculated by  Table~\ref{table:discfrmZp}.
\begin{table}
\caption{Discriminant forms of even lattices over $\Zp$} %even
\label{table:discfrmZp}
\centerline{
\vbox{\offinterlineskip
\hrule
\halign{
\vrule # & 
\strut\quad\hfil#\hfil\quad &\vrule# \hskip 1pt &\vrule# &
\strut\quad\hfil#\hfil\quad &\vrule# &
\strut\quad\hfil#\hfil\quad &\vrule# &
\strut\quad\hfil#\hfil\quad &\vrule# 
\cr
height 2pt &
\omit& height 2pt& height 2pt &
\omit& height 2pt&
\omit& height 2pt&
\omit& height 2pt
\cr
 &
$\Lambda$ & &  &
$p\sp\nu (a) $ & &
$2\sp\nu U$ & &
$2\sp\nu V$ & 
\cr
height 2pt &
\omit& height 2pt& height 2pt &
\omit& height 2pt&
\omit& height 2pt&
\omit& height 2pt
\cr
\noalign{\hrule}
height 1pt &
\omit& height 1pt& height 1pt &
\omit& height 1pt&
\omit& height 1pt&
\omit& height 1pt
\cr
\noalign{\hrule}
height 3pt &
\omit& height 3pt& height 3pt &
\omit& height 3pt&
\omit& height 3pt&
\omit& height 3pt
\cr
 &
$D\sb \Lambda$ & &  &
$\Z/ (p\sp\nu) $
& &
$ (\Z/ (2\sp\nu))\sp{\oplus 2} $ &
&
$ (\Z/ (2\sp\nu))\sp{\oplus 2} $ & 
\cr
height 3pt &
\omit& height 3pt& height 3pt &
\omit& height 3pt&
\omit& height 3pt&
\omit& height 3pt
\cr
\noalign{\hrule}
height 3pt &
\omit& height 3pt& height 3pt &
\omit& height 3pt&
\omit& height 3pt&
\omit& height 3pt
\cr
 &
$q\sb \Lambda$ & &  &
$
\begin{bmatrix}  \bigfrac{a\sp{\phantom{}} }{ {p\sp \nu} }\end{bmatrix}$
& &
$
\bigfrac {1\sp{\phantom{}} }{ 2\sp\nu}
\begin{bmatrix}  0 & 1 \\ 1 & 0\end{bmatrix} $ &
&
$
\bigfrac{1\sp{\phantom{}} }{ 2\sp\nu}
\begin{bmatrix}  2 & 1 \\ 1 & 2\end{bmatrix}  $ & 
\cr
height 3pt &
\omit& height 3pt& height 3pt &
\omit& height 3pt&
\omit& height 3pt&
\omit& height 3pt
\cr
}
\hrule
}
}
\end{table}
In particular, $D\sb{\Lambda}$ is a $p$-group of length
 equal to $\rank (\Lambda )-\rank (\Lambda\sb 0)$,
where $\Lambda \sb 0$ is the first Jordan component of $\Lambda$.
We also have
$\disc (\Lambda) = |D\sb \Lambda | \cdot\reddisc (\Lambda)$.
\section{Existence of even lattices with a given discriminant form}\label{sec:givendiscfrm}%even
\subsection{Over $\Zp$}
Suppose that a finite abelian $p$-group $D$ and a non-degenerate finite quadratic form
$q : D\to \Qp /2\thinspace \Zp$ are given.
It is known that, if $n\ge \length (D)$, then there exists an even lattice $\Lambda$
of rank $n$ over $\Zp$ such that $(D\sb\Lambda, q\sb\Lambda)$ is isomorphic to $(D, q)$.
The purpose of this subsection
is to describe a method to determine the set $\LLL\sp{(p)} (n, D, q)$
of all $[\sigma, u] \in \Z /(8) \times \Zp\sp\times/ \utsq{\Zp}$
such that there exists an even lattice $\Lambda$ of rank $n$
over $\Zp$ with $(D\sb\Lambda, q\sb\Lambda)\cong (D, q)$,
$\pexcess (\Lambda)=\sigma$ and $\reddisc (\Lambda)=u$.
\par
\medskip
Note that
\begin{align*}
\pexcess (\Lambda\sb 1 \oplus \Lambda\sb 2 ) &= \pexcess (\Lambda\sb 1) + \pexcess (\Lambda\sb 2),\quad\text{\rm and}\\
\reddisc (\Lambda\sb 1 \oplus \Lambda\sb 2 ) &= \reddisc(\Lambda\sb 1) \cdot \reddisc (\Lambda\sb 2).
\end{align*}
Taking these into account, 
for sets $\LLL$ and $\LLL\sp\prime$ of elements of $\Z /(8) \times \Zp\sp\times/ \utsq{\Zp}$,  we define
$\LLL * \LLL\sp\prime$ to be the set
$$
\{\; [\sigma +\sigma\sp\prime, u u\sp\prime] \; ; \; 
[\sigma, u] \in \LLL, \;  [\sigma\sp\prime, u\sp\prime] \in \LLL\sp\prime
\; \}.
$$
We also put $\LLL\sp{(p)}\sb 0 := \{ [0, 1] \} $.
Then $\LLL* \LLL\sp{(p)}\sb 0 = \LLL$ holds for any $\LLL$.
\begin{lemma}\label{L0}
Let $l$ be the length of $D$.
Then we have
\begin{equation}\label{L0_1}
\LLL\sp{(p)} (n, D, q) =\LLL\sp{(p)} (n-l, (0), [0]) * \LLL\sp{(p)} (l, D, q).
\end{equation}
If $p$ is odd, then
\begin{equation}\label{L0_2}
  \LLL\sp{(p)} (n-l, (0), [0])=
\begin{cases}
\emptyset & \text{\rm if $n<l$,} \\
\LLL\sp{(p)}\sb 0 & \text{\rm if $n=l$,} \\
\{[0, 1], [0, v\sb p] \} & \text{\rm if $n>l$,} 
\end{cases}
\end{equation}
where $v\sb p$ is the unique non-trivial element of $\Zp\sp\times / \utsq{\Zp}$.
If $p=2$, then
\begin{equation*}\label{L0_3}
  \LLL\sp{(2)} (n-l, (0), [0])=
\begin{cases}
\emptyset & \text{\rm if $n<l$ or $n-l \bmod 2 =1$,} \\
\LLL\sp{(2)}\sb 0 & \text{\rm if $n=l$,} \\
\{[n-l, 1], [n-l, 5] \} & \text{\rm if $n>l$ and $n-l \bmod 4 =0$,} \\
\{[n-l, 3], [n-l, 7]\} & \text{\rm if $n>l$ and $n-l \bmod 4 =2$}. 
\end{cases}
\end{equation*}
\end{lemma}
\begin{proof}
Let $\Lambda = \Lambda\sb 0 \oplus \bigoplus\sb{\nu>0} p\sp\nu \Lambda\sb \nu$
be a Jordan decomposition of an even lattice $\Lambda$ over $\Zp$ %even
with $(D\sb{\Lambda}, q\sb{\Lambda})\cong (D, q) $.
We put $\Lambda\sb{>0} := \Lambda\sb 0\sp\bot= \bigoplus\sb{\nu>0} p\sp\nu \Lambda\sb \nu $.
Then we have $\rank (\Lambda\sb{>0}) =l$,
$(D\sb{\Lambda\sb 0}, q\sb{\Lambda\sb 0})=((0), [0]) $ and 
$ (D\sb{\Lambda\sb{> 0}}, q\sb{\Lambda\sb {> 0}})=(D\sb\Lambda, q\sb{\Lambda}) \cong (D, q) $.
Hence~\eqref{L0_1} holds.
The statement~\eqref{L0_2} is obvious.
A lattice $\Lambda$ over $\Zt$ is even if and only if
$\Lambda\sb 0$ is of even rank and is isomorphic
to an orthogonal direct sum of copies of $U$ and $V$.
Because of 
$$
[\texcess (U), \reddisc (U)]=[2, 7] \quad\text{\rm and} \quad
[\texcess (V), \reddisc (V)]=[2, 3], 
$$
we can easily prove the last statement.
\end{proof}
 \begin{lemma}\label{L1}
Let  $p$ be  a prime integer,  and $\nu$  a positive integer.

{\rm (1)}
Suppose that $p$ is an odd prime,
and let $v_p$ be the unique non-trivial element of $\Zp\sp\times / \utsq{\Zp}$.
Let $u$ be an integer prime to $p$.
We put
$$
\chi_p(u):=\left(\frac{u}{p}\right)\;\;\in\;\;\{\pm 1\}.
$$
Then
$$
\LLL\sp{(p)} (1, \Z/ (p\sp\nu), \begin{bmatrix}\bigfrac{u\sp{\phantom{}} }{ p\sp\nu} \end{bmatrix} ) =
\begin{cases}
\{ [ p\sp\nu -1, 1]\} & \text{\rm if $\chi_p(u)=1$,} \\
\{ [ p\sp\nu -1, v_p]\} & \text{\rm if $\nu$ is even and $\chi_p(u)=-1$,} \\
\{ [ p\sp\nu +3, v_p]\} & \text{\rm if $\nu$ is odd and $\chi_p(u)=-1$.} 
\end{cases}
$$

{\rm (2)}
Suppose that $p=2$ and $a$ is an odd integer.
Then we have
$$
\LLL\sp{(2)} (1, \Z/ (2\sp\nu), \begin{bmatrix}\bigfrac{a\sp{\phantom{}} }{ 2\sp\nu} \end{bmatrix} ) =
\begin{cases}
\{ [1-a,a]\} & \text{\rm if $\nu$ is even,} \\
\{ [1-a,a]\} & \text{\rm if $\nu$ is odd, $\nu\ge 2$, and $a\equiv\pm 1\bmod 8$,} \\
\{ [5-a,a]\} & \text{\rm if $\nu$ is odd, $\nu\ge 2$, and $a\equiv\pm 3\bmod 8$,} \\
\{ [1-a,a],[1-a,5a]\} & \text{\rm if $\nu=1$ and $a\equiv\pm 1\bmod 8$,} \\
\{  [5-a,a],[5-a,5a]\} & \text{\rm if $\nu=1$ and $a\equiv\pm 3\bmod 8$.} 
\end{cases}
$$
Let $u$, $v$ and $w$ be integers with $v$ being odd.
Then
$$
\LLL\sp{(2)} (2, (\Z/ (2\sp\nu))\sp{\oplus 2} , 
\bigfrac{1\sp{\phantom{}} }{ 2\sp\nu} \begin{bmatrix}2u & v \\ v &2w  \end{bmatrix} )
=
\begin{cases}
\{ [2, 7]\} & \text{\rm if $\nu$ is even and $uw$ is even,} \\
\{ [2, 3]\} & \text{\rm if $\nu$ is even and $uw$  is odd,} \\
\{ [2, 7]\} & \text{\rm if $\nu$ is odd and $uw$ is even,} \\
\{ [6, 3]\} & \text{\rm if $\nu$ is odd and $uw$  is odd.} 
\end{cases}
$$
\end{lemma}
\begin{proof}
Two non-degenerate quadratic forms $[u/p\sp{\nu}]$ and $[u\sp\prime/p\sp{\nu}]$
on $\Z/ (p\sp\nu )$ with values in $\Qp / 2\thinspace \Zp$ are isomorphic if and only if
$$
uu\sp\prime \in \utsq{\Zp}, \qquad\text{\rm or}\qquad
(\quad p=2, \quad \nu=1, \quad\text{\rm and}\quad u = u\sp\prime \bmod 4 \quad)
$$
is satisfied.
On the other hand, two lattices $p\sp\nu (u)$ and $p\sp\nu (u\sp\prime)$ 
with $u$, $u\sp\prime \in \Zp\sp\times$
of rank $1$ over $\Zp$ are isomorphic  if and only if $uu\sp\prime \in \utsq{\Zp}$ holds.
Therefore the first statement follows.
The finite quadratic form
$$
q=\bigfrac{1}{2\sp\nu} 
\begin{bmatrix}
2 u & v \\ v & 2w 
\end{bmatrix}\qquad (v : \text{\rm odd})
$$
on $(\Z /(2\sp\nu))\sp{\oplus 2} $ 
with values in $\Qt/2\thinspace \Zt$ 
 is isomorphic to
$q\sb {2\sp\nu U}$
(resp.
$q\sb {2\sp\nu V}$ )
if and only if $uw \bmod 2 =0$ (resp. $uw \bmod 2 =1$).
These two forms can never be isomorphic to
$[u\sp\prime /2\sp\nu]\oplus [w\sp\prime /2\sp\nu]$
with $u\sp\prime$ and $w\sp\prime$ being odd.
Thus the second statement follows.
\end{proof}
\par
\medskip
Now we state an algorithm to calculate  $\LLL\sp{(p)} (n, D, q)$.
By Lemma~\ref{L0},
it is enough to determine $\LLL\sp{(p)} (l, D, q)$.
Let $\{ \gamma\sb 1, \dots, \gamma\sb l\}$
be a reduced set of generators of $D$.
We denote  the order of $\gamma\sb i$ by  $p\sp{\nu\sb i}$, and 
 arrange the generators in such a way that
$\nu\sb 1 \ge \dots\ge \nu\sb l$ holds.
For an element $\alpha \in \Qp/\Zp$, 
we define $\phi\sb p (\alpha)$ to be the integer
such that the order of $\alpha$ is $p\sp{\phi\sb p (\alpha)}$.
Note that $\phi\sb p (b [q] (\gamma\sb i, \gamma\sb j)) \le \min (\nu\sb i, \nu\sb j)$
holds for any $\gamma\sb i$ and $\gamma\sb j$.
\par
\medskip
When $l=1$,
$\LLL\sp{(p)} (l, D, q)$ is given by Lemma~\ref{L1}.
Suppose that $l>1$.
\par
\medskip
{\bf Case 1.}
Suppose that there exists a generator $\gamma\sb i$
such that $\phi\sb p (b[q] (\gamma\sb i, \gamma\sb i) )=\nu\sb 1$.
Then we have $\nu\sb i =\nu\sb 1$.
Interchanging $\gamma\sb 1$ and $\gamma\sb i$,
we will assume that $\phi\sb p (b[q] (\gamma\sb 1, \gamma\sb 1))=\nu\sb 1$.
Let $u$ be  an integer
such that
$b[q](\gamma\sb 1, \gamma\sb 1)=u/p\sp{\nu\sb 1} \bmod \Zp$.
Then $u$ is prime to $p$,
and hence there is 
an integer $v$ such that $uv=1 \bmod p\sp{\nu\sb 1}$ holds.
Since $\phi\sb p (b[q] (\gamma\sb j, \gamma\sb 1) ) \le \min (\nu\sb j, \nu\sb 1)=\nu \sb j$,
we can write $b[q] (\gamma\sb j, \gamma\sb 1)$ in the form  $w\sb j /p\sp{\nu\sb 1} \bmod \Zp$
by some integer $w\sb j$ that is divisible by $p\sp{\nu\sb 1-\nu\sb j}$.
For $j\ge 2$,
we put $\gamma\sb j\sp\prime :=\gamma\sb j-vw\sb j \gamma\sb 1$.
Because $\gamma\sb 1$ is of order $p\sp{\nu\sb 1}$ in $D$,
$\gamma\sb j\sp\prime$ is 
independent of the choice of $u$,  $v$ and $w\sb j$.
Moreover, $\gamma\sb j\sp\prime$ is of order $p\sp{\nu\sb j}$,
and $\{\gamma\sb 1, \gamma\sb 2\sp\prime, \dots, \gamma\sb l\sp\prime \}$
is again a reduced set of generators.
By definition, we have $b[q](\gamma\sb j\sp\prime, \gamma\sb 1)=0$ for any $j\ge 2$.
We put 
$$
(D\sb 1, q\sb 1):=(\langle \gamma\sb 1 \rangle , q|\sb{\langle \gamma\sb 1 \rangle })
\cong (\Z/ (p\sp{\nu\sb 1}), [u/p\sp{\nu\sb 1}])
$$ 
and 
$(D\sb 2, q\sb 2):=(\langle \gamma\sb 2\sp\prime , \dots, \gamma\sb l\sp\prime \rangle , 
q|\sb{\langle \gamma\sb 2\sp\prime , \dots, \gamma\sb l\sp\prime \rangle })$.
Then $(D, q)$ is decomposed into the orthogonal direct sum
of  $(D\sb 1, q\sb 1)$ and $(D\sb 2, q\sb 2)$.
\par
\medskip
Let $\Lambda$ be an even  lattice of rank $l$ over $\Zp$ such that %even
there exists  an isomorphism $h : (D\sb\Lambda, q\sb \Lambda)\isom (D, q)$.
Let $e\sp * \in \dual{\Lambda}$ be a vector 
such that $h\circ \Psi\sb\Lambda (e\sp *) =\gamma\sb 1$,
and $\Lambda\sb 1\sp\prime\subset \dual{\Lambda}$
the $\Zp$-submodule generated by $e\sp *$.
Then $\Lambda\sb 1 :=\Lambda\sb 1 \sp\prime \cap \Lambda$ is 
a sublattice of rank $1$ generated by $e:=p\sp{\nu\sb 1} e\sp *$.
Let $x$ be an arbitrary vector of $\Lambda$.
Because of $\ord{p} ( (x, e) )\ge \nu\sb 1 = \ord{p} ( (e, e) )$,
the vector
$$
x\sp\prime := x-\bigfrac{(x, e)}{(e, e)} e
$$
is in $\Lambda$ and orthogonal to $\Lambda\sb 1$.
Hence we obtain an orthogonal decomposition $\Lambda=\Lambda\sb 1 \oplus \Lambda\sb 1 \sp\bot$.
The homomorphism
$h\circ \Psi\sb{\Lambda} : \dual{\Lambda} = \dual{\Lambda\sb 1}\oplus \dual{(\Lambda\sb 1\sp{\bot})} \to D$
induces isomorphisms $(D\sb{\Lambda\sb 1}, q\sb{\Lambda\sb 1})\cong(D\sb{1}, q\sb{1})$ and 
$(D\sb{\Lambda\sb 1\sp\bot}, q\sb{\Lambda\sb 1\sp\bot})\cong(D\sb{2}, q\sb{2})$.
It follows that 
$$
\LLL\sp{(p)} (l, D, q) =\LLL\sp{(p)} (1, D\sb 1, q\sb 1) * \LLL\sp{(p)} (l-1, D\sb 2, q\sb 2).
$$
Thus $\LLL\sp{(p)} (l, D, q)$ is calculated by Lemma~\ref{L1} and the induction hypothesis on $l$.
\par
\medskip
{\bf Case 2.}
Suppose that $\phi\sb p (b[q] (\gamma\sb i, \gamma\sb i))< \nu\sb 1$
holds for any  generator  $\gamma\sb i$.
Since $q$ is non-degenerate,
there exists at least one generator  $\gamma \sb k$ 
that satisfies $\phi\sb p (b[q] (\gamma\sb1 , \gamma\sb k))=\nu\sb 1$.
Because of  $\phi\sb p (b[q] (\gamma\sb 1, \gamma\sb k))\le \nu\sb k$,
we have $\nu\sb k=\nu\sb 1$.
\par
\smallskip
{\bf Case  2.1.}
Suppose that $p$ is odd.
We replace $\gamma\sb 1$ by  $\gamma\sb 1\sp\prime :=\gamma\sb 1 + \gamma\sb k$,
which is an element  of order $p\sp{\nu\sb 1}$.
It is obvious that $\{\gamma\sb 1\sp\prime, \gamma\sb 2, \dots, \gamma\sb l \}$
is again a reduced set of generators of $D$.
Moreover we have $\phi\sb p (b[q] (\gamma\sb 1\sp\prime, \gamma\sb 1\sp\prime))= \nu\sb 1$.
Therefore  we are led to Case~1.
\par
\smallskip
{\bf Case  2.2.}
Suppose that $p=2$.
We replace $\gamma\sb 2$ by $\gamma\sb k$.
There exist integers $u, v$ and $w$ with $v$ being odd
such that
$$
b[q] (\gamma\sb 1, \gamma\sb 1)=\bigfrac{2\thinspace u}{2\sp{\nu\sb 1}},
\quad
b[q] (\gamma\sb 1, \gamma\sb 2)=\bigfrac{v}{2\sp{\nu\sb 1}},
\quad\text{\rm and}\quad
b[q] (\gamma\sb 2, \gamma\sb 2)=\bigfrac{2\thinspace w}{2\sp{\nu\sb 1}}
$$
hold modulo $\Zt$.
Note that
$q (\gamma\sb 1)={2\tilde u}/{2\sp{\nu\sb 1}}$ and 
$q (\gamma\sb 2)={2\tilde w}/{2\sp{\nu\sb 1}}$
hold modulo $2\thinspace \Zt$  for some integers $\tilde u$ and $\tilde w$
with $u=\tilde u \bmod 2\sp{\nu\sb 1-1}$ and  $w=\tilde w \bmod 2\sp{\nu\sb 1-1}$.
\par
\medskip
If $l=2$, then $\LLL \sp{(2)} (l, D, q)$ is determined by Lemma~\ref{L1}.
Suppose that $l\ge 3$.
There exists an integer $t$ such that
$(4uw -v\sp 2) t=1 \bmod 2\sp{\nu\sb 1}$ holds. %correct on 11 Apr 2006
For each $j\ge 3$,
we choose integers $s\sb {j1}$ and $s\sb {j2}$
such that
$ b[q] (\gamma\sb j, \gamma\sb 1) =s\sb {j1}/ 2\sp{\nu\sb 1} \bmod \Zt$
and 
$ b[q] (\gamma\sb j, \gamma\sb 2) =s\sb {j2}/ 2\sp{\nu\sb 1} \bmod \Zt$
hold,
and calculate
$$
\begin{pmatrix}
\raise 1pt \hbox{$\beta\sb{j1}$} \\ \raise -1pt \hbox{$\beta\sb{j2}$}
\end{pmatrix} 
:=
t\cdot
\begin{pmatrix}
2w & -v \\ 
-v & 2u
\end{pmatrix}
\cdot
\begin{pmatrix}
\raise 1pt \hbox{$s\sb{j1}$} \\ 
\raise -1pt \hbox{$s\sb{j2}$}
\end{pmatrix}
.
$$
Then $\beta\sb{j1}$ and $\beta\sb{j2}$ are divisible by $2\sp{\nu\sb 1-\nu\sb j}$.
Hence $\gamma\sb j \sp\prime :=\gamma\sb j -\beta \sb{j1} \gamma\sb 1-\beta \sb{j2} \gamma\sb 2$
is an element  of  order $2\sp{\nu\sb j}$, which is independent of the choice of the integers.
The set $\{\gamma\sb 1, \gamma\sb 2, \gamma\sb 3\sp\prime, \dots, \gamma\sb l\sp\prime \}$
is again a reduced set of generators of $D$,
and the two subgroups
$\langle \gamma\sb1, \gamma\sb 2\rangle$
and 
$\langle \gamma\sb 3\sp\prime, \dots, \gamma\sb l\sp\prime \rangle$ of $D$ 
are orthogonal with respect to $q$.
Therefore, putting
$$
(D\sb 1, q\sb 1):=(\langle \gamma\sb1, \gamma\sb 2\rangle, q|\sb{\langle \gamma\sb1, \gamma\sb 2\rangle}) 
\cong \Bigl( (\Z/(p\sp{\nu\sb 1}))\sp{\oplus 2}, 
\bigfrac{1}{2\sp{\nu\sb 1}}
\begin{bmatrix}
2\tilde u & v \\ v & 2 \tilde w
\end{bmatrix}
  \Bigr)
$$
and $(D\sb 2, q\sb 2):=(\langle \gamma\sb 3\sp\prime, \dots, \gamma\sb l\sp\prime \rangle, 
q |\sb{\langle \gamma\sb 3\sp\prime, \dots,
\gamma\sb l\sp\prime \rangle})$, we obtain an orthogonal decomposition
$(D, q)=(D\sb 1, q\sb 1)\oplus (D\sb 2, q\sb 2) $.
\par
\medskip
Let $\Lambda$ be an even  lattice of rank $l$ over $\Zt$ %even
such that there exists an isomorphism $h : (D\sb\Lambda, q\sb\Lambda)\isom (D, q)$.
We pick up two vectors
$e\sb 1\sp *, e\sb 2\sp * \in \dual{\Lambda}$
such that $h\circ \Psi\sb \Lambda (e\sb 1 \sp *) =\gamma\sb 1$
and $h\circ \Psi\sb \Lambda (e\sb 2 \sp *) =\gamma\sb 2$.
Let $\Lambda\sb 1\sp\prime \subset \dual{\Lambda}$ be the $\Zt$-submodule of $\dual{\Lambda}$
generated by $e\sb 1\sp *$ and $e\sb 2\sp *$.
Then $\Lambda\sb 1:=\Lambda\sb 1\sp\prime \cap \Lambda$
is a sublattice of $\Lambda$ generated by $e\sb 1:=2\sp{\nu\sb 1}e\sb 1\sp *$
and $e\sb 2:=2\sp{\nu\sb 1}e\sb 2\sp *$.
The intersection matrix $M\sb1 $ of $\Lambda\sb 1$
with respect to $e\sb 1$ and $e\sb 2$ 
satisfies $\ord{2} (\det M\sb 1\sp{-1})=-\nu\sb 1$.
Because $\ord{2} ( (x, e\sb 1) ) \ge \nu\sb 1$ and $\ord{2} ( (x, e\sb 2) ) \ge \nu\sb 1$
hold for any vector $x\in \Lambda$,
we have $( (x, e\sb 1), (x, e\sb 2))\cdot M\sb 1 \sp{-1} \in \Zt\sp{\oplus 2}$.
Therefore $\Lambda$ is decomposed into
the orthogonal direct sum of $\Lambda\sb 1$ and $\Lambda\sb 1\sp\bot$.
The homomorphism
$h\circ \Psi\sb{\Lambda}$
induces isomorphisms $(D\sb{\Lambda\sb 1}, q\sb{\Lambda\sb 1})\cong(D\sb{1}, q\sb{1})$ and 
$(D\sb{\Lambda\sb 1\sp\bot}, q\sb{\Lambda\sb 1\sp\bot})\cong(D\sb{2}, q\sb{2})$.
It follows that 
$$
\LLL\sp{(2)} (l, D, q) =\LLL\sp{(2)} (2, D\sb 1, q\sb 1) * \LLL\sp{(2)} (l-2, D\sb 2, q\sb 2).
$$
Thus $\LLL\sp{(2)} (l, D, q)$ is calculated by Lemma~\ref{L1} and the induction hypothesis on $l$.
\subsection{Over $\Z$}
Let $D$ be a finite abelian group,
and $q : D\to \Q/2\thinspace \Z$ a non-degenerate finite quadratic form.
Let $(r, s)$ be a pair of non-negative integers 
such that $n:=r+s >0$.
We will describe a criterion to determine 
whether there exists an even  lattice $L$  over $\Z$ %even
with signature $(r, s)$
such that $(D\sb L, q\sb L)$ is isomorphic to the given
finite quadratic form  $(D, q)$.
\par
\medskip
We put $d:=(-1)\sp s |D|$.
Let $P$ be the set of prime divisors of $2\thinspace d$,
and let
$(D, q)=\oplus\sb{p\in P} (D\sp{(p)}, q\sp{(p)} )$ be the orthogonal decomposition of $(D, q)$ into
the $p$-parts. If $d$ is odd, then we put $(D\sp{(2)}, q\sp{(2)} )=((0), [0])$.
By Lemma~\ref{qdecomp} and Theorem~\ref{existence},
an even  lattice $L$ over $\Z$ with signature $(r, s)$ and $(D\sb L, q\sb L)\cong (D, q)$ %even
exists
if and only if the following claim is verified:
\par
\medskip
$(\sharp)$
For each $p\in P$,
there exists an even  lattice $\Lambda\sp{(p)}$ of rank $n$ over $\Zp$
such that
\begin{itemize}
\item[(i)] $\disc (\Lambda\sp{(p)}) = d\cdot \utsq{\Zp}$ and 
\item[(ii)] $(D\sb {\Lambda\sp{(p)}}, q\sb {\Lambda\sp{(p)}})\cong (D\sp{(p)}, q\sp{(p)})$ hold,
\end{itemize}
and they satisfy
$$
r-s +\sum\sb{p\in P} \pexcess (\Lambda\sp{(p)}) =n \quad\bmod 8.
$$
We put $\delta\sb p :=d/ p\sp{\ord{p} (d)} \in \Z$.
Under the condition (ii),
which implies $| D\sb{\Lambda\sp{(p)}} | = d/\delta\sb p$,
the condition (i) is equivalent to the condition 
 $$
\reddisc (\Lambda\sp{(p)} ) = \delta\sb p \cdot \utsq{\Zp}.
$$
Therefore we can check the claim $(\sharp)$ by the following method.
First we calculate $\LLL\sp{(p)} (n, D\sp{(p)}, q\sp{(p)})$
for each $p \in P$.
Then we search for an element
$(\;[\sigma\sb p, u\sb p]\; ; \; p\in P\;)$ of the Cartesian product of the sets $\LLL\sp{(p)} (n, D\sp{(p)}, q\sp{(p)})$
that satisfies
$u\sb p = \delta\sb p \cdot \utsq{\Zp}$ for each $p\in P$
and  $r-s +\sum \sigma\sb p =n \bmod 8$.
The claim $(\sharp)$ is true if and only if we find such an element.
\section{Roots}\label{sec:roots}
For the following, we refer to~\cite{Bourbaki},  \cite[Chapter 4]{ConwaySloane} or~\cite{Nishiyama}.
\subsection{Root system of a  positive-definite even lattice over $\Z$}
Let $L$ be a positive-definite  even lattice over $\Z$.
A vector of $ L$ is said to be a {\it root} if its norm is $2$.
We denote by $L\sb{\rm root}$ the sublattice of $L$ generated by roots.
A lattice $L$ is said to be a  {\it root lattice} if $L=L\sb{\rm root}$ holds.
Let $\Roots (L)$ be the set of roots of $L$.
We define $\sim$ to be the finest equivalence relation on $\Roots (L)$
that satisfies 
$(v, w)\ne 0 \Longrightarrow v\sim w$.
Let $I\sb 1$, \dots, $I\sb k$ be the equivalence classes of roots
under the relation $\sim$,
and let $L\sb i$ be the sublattice of $L\sb{\rm root}$ generated by $I\sb i$.
There exists a basis  $B\sb i \subset I\sb i$
such that the intersection matrix of $L\sb i$
with respect to  $B\sb i$ is the Cartan matrix corresponding to  a Dynkin diagram
of type $A\sb l$, $D\sb m $ or $E\sb n$.
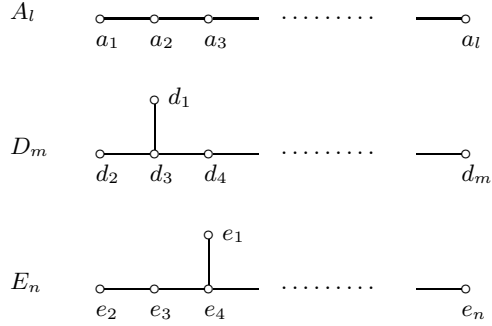
\begin{figure} \label{dynkin}
\caption{Dynkin diagram}
\def\ha{40}
\def\hav{37}
\def\hd{25}
\def\hdv{22}
\def\he{10}
\def\hev{7}
\setlength{\unitlength}{1.2mm}
\centerline{
{\small
\begin{picture}(100, 37)(-20, 7)
\put(0, \ha){$A\sb l$}
\put(10, \ha){\circle{1}}
\put(9.5, \hav){$a\sb 1$}
\put(10.5, \ha){\line(5, 0){5}}
\put(16, \ha){\circle{1}}
\put(15.5, \hav){$a\sb 2$}
\put(16.5, \ha){\line(5, 0){5}}
\put(22, \ha){\circle{1}}
\put(21.5, \hav){$a\sb 3$}
\put(22.5, \ha){\line(5, 0){5}}
\put(30, \ha){$\dots\dots\dots$}
\put(45, \ha){\line(5, 0){5}}
\put(50.5, \ha){\circle{1}}
\put(50, \hav){$a\sb {l}$}
\put(0, \hd){$D\sb m$}
\put(10, \hd){\circle{1}}
\put(9.5, \hdv){$d\sb 2$}
\put(10.5, \hd){\line(5, 0){5}}
\put(16, 31){\circle{1}}
\put(17.5, 30.5){$d\sb 1$}
\put(16, 25.5){\line(0,1){5}}
\put(16, \hd){\circle{1}}
\put(15.5, \hdv){$d\sb 3$}
\put(16.5, \hd){\line(5, 0){5}}
\put(22, \hd){\circle{1}}
\put(21.5, \hdv){$d\sb 4$}
\put(22.5, \hd){\line(5, 0){5}}
\put(30, \hd){$\dots\dots\dots$}
\put(45, \hd){\line(5, 0){5}}
\put(50.5, \hd){\circle{1}}
\put(50, \hdv){$d\sb {m}$}
\put(0, \he){$E\sb n$}
\put(10, \he){\circle{1}}
\put(9.5, \hev){$e\sb 2$}
\put(10.5, \he){\line(5, 0){5}}
\put(16, \he){\circle{1}}
\put(15.5, \hev){$e\sb 3$}
\put(22, 16){\circle{1}}
\put(23.5, 15.5){$e\sb 1$}
\put(22, 10.5){\line(0,1){5}}
\put(16.5, \he){\line(5, 0){5}}
\put(22, \he){\circle{1}}
\put(21.5, \hev){$e\sb 4$}
\put(22.5, \he){\line(5, 0){5}}
\put(30, \he){$\dots\dots\dots$}
\put(45, \he){\line(5, 0){5}}
\put(50.5, \he){\circle{1}}
\put(50, \hev){$e\sb {n}$}
\end{picture}
}
}
\vskip 10pt
\end{figure}
Let $\tau\sb i$ be the type of the Dynkin diagram
of the intersection matrix of $L\sb i$.
We define the root type of $L$ to be 
$\sum\sb{i=1}\sp k \tau\sb i$.
Conversely, for an $ADE$-type $\Sigma$,
there exists a root lattice $L(\Sigma)$, unique up to isomorphism,
whose root type is $\Sigma$.
\par
\medskip
The root type of a positive-definite  even lattice $L$  over $\Z$ is therefore determined by the 
following procedure.
\begin{itemize}
\item[(1)] Create the list $\Roots (L)$, and decompose it into $I\sb 1$, \dots, $I\sb k$.
\item[(2)] Calculate the rank of $L\sb i$ for $i=1, \dots, k$.
\item[(3)] Determine the type $\tau\sb i$ from $\rank (L\sb i)$ and $|I\sb i|$
by using Table~\ref{table:discfrmroot}. 
\end{itemize}
\subsection{Discriminant forms of root lattices}
The discriminant form  $(D\sb{L(\tau)}, q\sb{L(\tau)})$,
where $\tau$ is $A\sb l$, $D\sb m$ or $E\sb n$,
 is indicated in Table~\ref{table:discfrmroot}.
In this table, for example,  $\{ a\sb 1\sp *, \dots, a\sb l\sp *\}$ is the basis 
of $\dual{L(A\sb l)}$ dual to the basis
$\{ a\sb 1, \dots, a\sb l\}$ of $L (A\sb l)$ given in Figure~\ref{dynkin},
and $\bar a\sb i\sp * \in D\sb{L(A\sb l)}$ is the image of $a\sb i\sp *$ by the 
homomorphism $\Psi\sb{L(A\sb l)} :\dual{L(A\sb l)} \to  D\sb{L(A\sb l)} $.
\begin{table}
\caption{Number of roots and discriminant forms of root lattices}
\label{table:discfrmroot}
\def\hstable{\hskip 4pt}
\def\spaceheight{height 3pt}
\def\HS{
\spaceheight & 
\omit & \spaceheight & \spaceheight&
\omit & \spaceheight & \spaceheight&
\omit & \spaceheight &
\omit & \spaceheight 
\cr}
\def\HL{\HS
\noalign{\hrule}
\HS
\cr}
\def\HHL{\HS
\noalign{\hrule}
height 1pt &
\omit& height 1pt& height 1pt &
\omit& height 1pt& height 1pt &
\omit& height 1pt &
\omit& height 1pt
\cr
\noalign{\hrule}
\HS}
\centerline{
\vbox{\offinterlineskip
\hrule
\halign{
\vrule # & 
\strut\hstable\hfil#\hfil\hstable &\vrule# &\hskip 1pt \vrule #&
\strut\hstable\hfil#\hfil\hstable &\vrule# &\hskip 1pt \vrule #&
\strut\hstable\hfil#\hfil\hstable &\vrule# &
\strut\hstable\hfil#\hfil\hstable &\vrule# 
\cr
\spaceheight &
\omit& \spaceheight& \spaceheight &
\omit& \spaceheight& \spaceheight &
\omit& \spaceheight&
\omit& \spaceheight
\cr
 &
$\tau$ & & &
$|\textrm{Roots} (L (\tau)) | $ & & & 
$D\sb{L(\tau)}$ & &
$q\sb{L(\tau)}$ & 
\cr
\HHL
 &
 $A\sb l$ & & &
$ l (l+1)  $ & & & 
 $\langle \bar a \sb l\sp * \rangle \cong \Z/ (l+1)$ & &
 $\begin{bmatrix}  \bigfrac{l}{l+1}\end{bmatrix}$& 
\cr
\HL
 &
 $D\sb m$ $(m: \text{\rm even})$ & & &
$ 2 m (m-1)  $ & & & 
 $\langle \bar d\sb 1\sp * \rangle \oplus \langle  \bar d\sb m\sp *\rangle 
\cong (\Z/ (2))\sp{\oplus 2}$ & &
 $\begin{bmatrix}  m/4 & 1/2 \\ 1/2 & 1 \end{bmatrix}$& 
\cr
\HL
 &
 $D\sb m$ $(m: \text{\rm odd})$ & & &
$ 2 m (m-1)  $ & & &
 $\langle \bar d\sb 1\sp * \rangle \cong \Z/ (4)$ & &
 $\begin{bmatrix}  m/4 \end{bmatrix}$& 
\cr
\HL
 &
 $E\sb 6$ & & &
$ 72  $ & & &
 $\langle \bar e\sb 6\sp * \rangle \cong \Z/ (3)$ & &
 $\begin{bmatrix}  4/3 \end{bmatrix}$& 
\cr
\HL
 &
 $E\sb 7$ & & &
$ 126  $ & & &
 $\langle \bar e\sb 7\sp * \rangle \cong \Z/ (2)$ & &
 $\begin{bmatrix}  3/2 \end{bmatrix}$& 
\cr
\HL
 &
 $E\sb 8$ & & &
$ 240  $ & & &
 $(0)$ & &
 $\begin{bmatrix} 0 \end{bmatrix}$& 
\cr
\spaceheight &
\omit& \spaceheight& \spaceheight &
\omit& \spaceheight& \spaceheight &
\omit& \spaceheight&
\omit& \spaceheight
\cr
}
\hrule
}
}
\end{table}
\par
\medskip
Let $\Gamma (\tau)$ denote the image 
of the natural homomorphism from the orthogonal group $ O (L (\tau))$ 
of the lattice $L(\tau)$ to
$\Aut (D\sb{L(\tau)}, q\sb{L (\tau)})$.
The structure of $\Gamma (\tau)$ is given as follows.
\begin{itemize}
\item
If $\tau=A\sb 1$ or $\tau=E\sb 7$, then $\Gamma (\tau)$ is trivial.
\item
If $\tau=A\sb l $ $(l>1)$ or $\tau=D\sb m$ $(m: \text{\rm odd})$ or
 $\tau=E\sb 6$, then $\Gamma (\tau)$ is isomorphic to $\Z / (2)$ generated by
the multiplication by $-1$.
\item
If $\tau=D\sb m$ with $m$ being even and $>4$,
then $\Gamma (\tau)$ is isomorphic to $\Z / (2)$ generated by
$$
\bar d\sb 1 \sp * \mapsto \bar d\sb 1 \sp * + \bar d\sb m \sp * , \quad
\bar d\sb m \sp * \mapsto \bar d\sb m \sp *.
$$
\item
If $\tau=D\sb 4$,
then $\Gamma (\tau)$ is isomorphic to 
the full symmetric group acting on 
the set $\{ \bar d\sb 1\sp *, \bar d\sb4 \sp *, \bar d\sb 1\sp *+ \bar d\sb4 \sp * \}$
of non-trivial elements of $D\sb{L (\tau)}$.
\end{itemize}
\section{Existence of an elliptic $K3$ surface with given data}\label{sec:existenceK3}
\begin{theorem}\label{thm:existenceK3}
Let $\Sigma$ be an $ADE$-type with $\rank (\Sigma)\le 18$,
and $G$ a finite abelian group.
There exists an elliptic $K3$ surface 
$f : X\to \P\sp 1$
with $\Sigma\sb f = \Sigma$ and $G\sb f\cong G$
if and only if the root lattice $L (\Sigma)$
has an even overlattice $M$
with the following properties.
\begin{enumerate}
\renewcommand{\labelenumi}{(\roman{enumi})}
\item
$M/ L(\Sigma) \cong G$,
\par
\item
there exists an even lattice $N$ of 
signature $(2, 18-\rank (\Sigma))$ such that $(D\sb N, q\sb N)$ is isomorphic to 
$(D\sb M, q\sb M)$, and
\item
the sublattice  $M\sb{\rm root}$ of $M$ coincides with $L (\Sigma)$.
\end{enumerate}
\end{theorem}
\begin{proof}
Suppose that a  pair  $(\Sigma, G)$
satisfies the condition of Theorem.
By Proposition~\ref{nikulin2},
the property (ii) implies that there exists
an even unimodular overlattice $K\sp\prime$
of $\opp{M}\oplus N$
into which $\opp{M}$ and $N$ are primitively embedded.
Let $H$ denote the hyperbolic lattice;
$$
H:=\begin{pmatrix} 0 & 1 \\ 1 & 0 \end{pmatrix}.
$$
Then $K:=K\sp\prime\oplus H$
is an even unimodular lattice with signature $(3, 19)$.
Hence $K$ is isomorphic to the $K3$ lattice $\opp{L (2 E\sb 8)} \oplus H\sp{\oplus 2}$
by Milnor's structure theorem
(cf.\ \cite{Serre}).
There exists a $2$-dimensional linear subspace $V$ of $N\otimes\sb{\Z} \R$
such that the bilinear form is positive-definite on $V$
and that, if $N\sp\prime \subset N$ is a sublattice such that 
$N\sp\prime\otimes\sb{\Z} \R$ contains $V$, then $N\sp\prime$ 
coincides with  $N$.
By the surjectivity of the period map on 
the moduli of $K3$ surfaces,
there exists a complex $K3$ surface $X$ 
and an isomorphism $\alpha : H\sp 2 (X; \Z) \isom K$
of lattices
such that
$$
\alpha\sb{\R}\sp{-1} (V) = (H\sp{0,2} (X) \oplus H\sp{2, 0} (X) ) \cap H\sp 2 (X; \R)
$$
holds, where $\alpha\sb{\R}:=\alpha\otimes\sb{\Z} {\R}$. Then we have 
\begin{equation}\label{eq:alpha}
\alpha\sp{-1} (\opp{M}\oplus H) = \NS\sb X. 
\end{equation}
By Kondo's lemma~\cite[Lemma 2.1]{Kondo},
there exists a structure of the  elliptic fibration $f : X\to \P\sp 1$
with a section $O : \P\sp 1 \to X$ such that,
if $F$ denotes the cohomology class of a general fiber of $f$,
then 
\begin{equation}\label{eq:kondo}
\Z [F]\sp\bot /\Z [F] \cong \opp{M}
\end{equation}
holds, where $\Z [F]\sp\bot$ is the orthogonal complement of $\Z [F]$ in 
the N\'eron-Severi lattice $\NS\sb X$ of $X$.
Let $H\sb f$ be  the sublattice of $\NS\sb X$
spanned by the cohomology classes
of the zero section and a general fiber of $f$, 
$S\sb f$ the sublattice of  $\NS\sb X$
defined  in  \S\ \ref{sec:intro},
and $W\sb f$ the orthogonal complement 
of $H\sb f$ in $\NS\sb X$.
The lattice $H\sb f$ is isomorphic to
the hyperbolic lattice $H$, 
and is orthogonal to $S\sb f$.
By abuse of notation,
we denote  $(W\sb f)\sb{\rm root}$ for 
the sublattice of $W\sb f$ generated by the vectors of norm $-2$.
From~\eqref{eq:kondo},
we have
\begin{equation}\label{eq:kondo2}
W\sb f \cong \opp{M}.
\end{equation}
On the other hand, 
by   Nishiyama's lemma~\cite[Lemma 6.1]{Nishiyama},
we have
\begin{align}
W\sb f / (W\sb f)\sb{\rm root} &\cong  MW\sb f, \quad\text{\rm and}\label{eq:nishiyama1} \\
S\sb f & =   (W\sb f)\sb{\rm root} \label{eq:nishiyama2}.
\end{align}
Combining these with the properties (i) and (iii) of $M$ and 
the isomorphism~\eqref{eq:kondo2},
we have $S\sb f \cong  \opp{L (\Sigma)}$ and $MW\sb f \cong G$.
Hence $\Sigma = \Sigma\sb f$ and $G\cong G\sb f$ hold.
\par
\medskip
Conversely,
suppose that there exists an elliptic $K3$ surface $f :X\to \P\sp 1$
with $\Sigma\sb f=\Sigma$ and $G\sb f\cong G$.
Using  Nishiyama's lemma again, we see that
the primitive closure $\closure{S\sb f}$ of $S\sb f$ in $\NS\sb X$  satisfies
$\closure{S\sb f}/S\sb f \cong G$ and $(\closure{S\sb f})\sb{\rm root}=S\sb f$.
We have an isomorphism $S\sb f \cong \opp{L(\Sigma)}$.
Let $\opp{M}$ be the overlattice of $\opp{L(\Sigma)}$
corresponding to $\closure{S\sb f}$ via this isomorphism.
Then $M:=\opp{(\opp{M})}$ is an overlattice of $L (\Sigma)$ that possess
 the properties (i) and (iii).
Moreover, $\closure{S\sb f} \oplus H\sb f$ is primitive in the even unimodular lattice $H\sp 2 (X; \Z)$,
and hence Proposition~\ref{nikulin2} implies that the orthogonal complement $N\sb f$ of
$\closure{S\sb f} \oplus H\sb f$ in $H\sp 2 (X; \Z)$
satisfies $(D\sb{N\sb f}, q\sb{N\sb f})\cong 
(D\sb{\closure{S\sb f}}, -q\sb{\closure{S\sb f}})\cong
(D\sb M, q\sb M)$.
Because the signature of $N\sb f$ is $(2, 18-\rank (\Sigma))$,
the overlattice $M$ has  the property (ii).
\end{proof}
\section{Making the list}\label{sec:making}
Recall that,
in order for an $ADE$-type $\Sigma$ to be an $ADE$-type of an elliptic $K3$ surface,
it is necessary that
$\rank (\Sigma) \le 18$ and $\euler (\Sigma ) \le 24$.
It is obvious that
the torsion part of the Mordell-Weil group
of an elliptic surface is
of length $\le 2$.
\par
\medskip
First we list up all $ADE$-types $\Sigma$ with $\rank (\Sigma)\le 18$
and $\euler (\Sigma ) \le 24$.
There are $3937$ such $ADE$-types.
For each 
$$
\Sigma :=\sum a\sb l A\sb l + \sum d\sb m D\sb m + \sum e\sb n E\sb n
$$
in this list, we carry out the following calculation.
\par
\medskip
{\it Step 1.}
We calculate the discriminant form $(D\sb{L(\Sigma)}, q\sb{L(\Sigma)})$
using Table~\ref{table:discfrmroot}.
Note that the product of the wreath products
$$
\prod\sb{a\sb l >0} (\Gamma (A\sb l) \wr \Sym\sb{a\sb l}) \times 
\prod\sb{d\sb m >0} (\Gamma (D\sb m) \wr \Sym\sb{d\sb m})
\times \prod\sb{e\sb n >0} (\Gamma (E\sb n) \wr \Sym\sb{e\sb n}) 
$$
acts on $(D\sb{L(\Sigma)}, q\sb{L(\Sigma)})$.
Here, for example, the full symmetric group $\Sym\sb{a\sb l}$ acts on $D\sb{L (\Sigma) }$ 
as the permutation group on the $a\sb l$ components of $D\sb{L (\Sigma)}$ 
isomorphic to $D\sb{L (A\sb l )}$. 
We denote this group by $\Gamma(\Sigma)$.
\par
\medskip
{\it Step 2.}
We make a complete list of representatives
of the quotient set $D\sb{L(\Sigma)}/\Gamma (\Sigma)$
and pick up from this list 
elements isotopic with respect to $q\sb{L(\Sigma)}$.
Let $\VVV\sb{\Sigma}=\{\bar v\sb 1, \dots, \bar v\sb N \}$
be the list of isotopic  elements of $D\sb{L(\Sigma)}$ modulo $\Gamma(\Sigma)$.
For each $\bar v\sb i \in \VVV\sb{\Sigma}$,
we calculate the stabilizer subgroup $St (\Gamma(\Sigma), \bar v\sb i)$
of $\bar v\sb i$ in $\Gamma(\Sigma)$.
Then we make a complete list of representatives of
$D\sb {L (\Sigma)} / St (\Gamma(\Sigma), \bar v\sb i)$,
and pick up from this list elements isotopic with respect to $q\sb{L(\Sigma)}$
and orthogonal to $\bar v\sb i$ with respect to $b[q\sb{L(\Sigma)}]$.
Let $\WWW\sb{\Sigma, i}$ be the list of isotopic elements 
orthogonal to $\bar v\sb i$ modulo $St (\Gamma(\Sigma), \bar v\sb i)$.
\par
Next we make the list
$\GGG\sp\prime\sb{\Sigma}$
of all pairs $[\bar v\sb i, \bar w\sb j]$
of  $\bar v\sb i \in \VVV\sb{\Sigma}$ and $\bar w\sb j \in \WWW\sb{\Sigma, i}$.
Then every isotopic subgroup 
of $(D\sb{L (\Sigma)}, q\sb{L(\Sigma)})$ with length $\le 2$ is conjugate under the action of $\Gamma (\Sigma)$ 
to a subgroup $\langle \bar v\sb i, \bar w\sb j \rangle$
generated by 
$\bar v\sb i$ and $ \bar w\sb j$ for some $[ \bar v\sb i, \bar w\sb j ] \in \GGG\sp\prime\sb{\Sigma}$.
Of course, there are several different pairs that generate a same subgroup.
We remove this redundancy from $\GGG\sp\prime\sb{\Sigma}$,
and make a list $\GGG\sb{\Sigma}$.
\par
\medskip
{\it Step 3.}
For each $[\bar v, \bar w] \in \GGG\sb{\Sigma}$,
we calculate the subgroup $G:=\langle \bar v, \bar w \rangle$ of $D\sb{L (\Sigma)}$,
its orthogonal complement $G\sp\bot$ in $(D\sb{ L (\Sigma)}, q\sb{ L (\Sigma)})$,
and the finite quadratic  form $(D\sb G, q\sb G) := (G\sp\bot / G, q\sb{L (\Sigma) } |\sb{G\sp\bot / G})$.
\par
\smallskip
{\it Step 3.1.}
By the algorithm described in \S\ \ref{sec:givendiscfrm},
we determine whether there exists an even lattice $N$ over $\Z$ of signature $(2, 18-\rank (\Sigma))$
such that $(D\sb N, q\sb N)\cong (D\sb G, q\sb G)$.
If the answer is affirmative, we go to the next step.
\par
\smallskip
{\it Step 3.2.}
We calculate the intersection matrix of the even overlattice $M\sb G$ of $L (\Sigma)$
generated by $L (\Sigma)$ and $v$, $w$ in $\dual{ L(\Sigma)}$,
where $v$ and $w$ are vectors of $\dual{L (\Sigma)}$ such that $\Psi \sb{L (\Sigma)} (v)=\bar v$
and $\Psi \sb{L (\Sigma)} (w)=\bar w$.
Then we calculate the root type of $M\sb G$
by the algorithm described in \S\ \ref{sec:roots}.
If this root type coincides with the initial $ADE$-type  $\Sigma$,
then we let the pair  $(\Sigma, G)$ be a member of the list $\PPP$.
\par
\medskip
By Theorem~\ref{thm:existenceK3},
the list $\PPP$ thus made is the complete list of 
the data of elliptic $K3$ surfaces.
\par
\medskip
The following  remarks are useful in checking the program.
\begin{remark}
Note that neither  $\euler (\Sigma ) \le 24$ 
nor $\length (G) \le 2$ is contained in the conditions of Theorem~\ref{thm:existenceK3}.
Therefore,
if we input $\Sigma$ with $\euler (\Sigma)> 24$
into the program,
then it should return no subgroups $G$ of $D\sb{L (\Sigma )}$
such that $(\Sigma, G)$ can be a member of  the list $\PPP$.
If we change Step 2 of the program so that it lists up
all isotopic subgroups of length $\ge 3$,
then the result should also be an empty set.
\end{remark}
\begin{remark}
Suppose that the root type $\Sigma\sp\prime$ of $M\sb G$
is not  equal to $\Sigma$ in Step~3.2  of the program.
Let $G\sp\prime$ be the finite abelian group $M\sb G / (M\sb G)\sb{\rm root}$.
Then $(\Sigma\sp\prime, G\sp\prime)$  appears in $\PPP$.
\end{remark}
\begin{remark}
For each $(\Sigma, G)\in \PPP$,
there should be
at least one configuration
that satisfies the conditions given in \S\ \ref{subsec:recover}.
\end{remark}

\vfill
\eject

This is the table of all $ADE$-types of singular fibers of 
complex elliptic 
$K3$ surfaces with a zero section 
and the torsion parts of their
Mordell-Weil groups.
\par
In Table 1,
the  $ADE$-types of the singular fibers are listed 
according to the rank and the lexicographic order.
For each $ADE$-type $\Sigma$,
the third column shows the list of all abelian groups that can be realized
as the torsion part of the Mordell-Weil group
of an elliptic $K3$ surface $f : X\to \P\sp 1$ 
with $\Sigma\sb f=\Sigma$.
Here, $[1]$ is the trivial group, 
$[a]$ is the cyclic group $\Z / (a)$, and 
$[a, b]$ is the group $\Z / (a)\times \Z / (b)$.
\par
Table 2 shows, for each abelian group $G$ with order $\ge 3$
that appears as $G\sb f$ of some elliptic $K3$ surface,
the list of all $ADE$-types $\Sigma$ such that 
the pair $(\Sigma, G)$ appears in Table 1.

\begin{center}
Table 1.
\end{center}

{\small
%%%%%%%%%%%%%%%%%%%%%%%%%%%%%%%
%
% Table 1
%
%%%%%%%%%%%%%%%%%%%%%%%%%%%%%%%%

\vsr \elldata{No.}{rank}{$ADE$-type}{$G$}

\vsrs \elldata{0}{$1$}{$A_{1}$}{$[1]$}

\vsr \elldata{No.}{rank}{$ADE$-type}{$G$}

\vsrs \elldata{1}{$2$}{$A_{2}$}{$[1]$}
\elldata{2}{$2$}{$2\,A_{1}$}{$[1]$}

\vsr \elldata{No.}{rank}{$ADE$-type}{$G$}

\vsrs \elldata{3}{$3$}{$A_{3}$}{$[1]$}
\elldata{4}{$3$}{$A_{2}\,+\,A_{1}$}{$[1]$}
\elldata{5}{$3$}{$3\,A_{1}$}{$[1]$}

\vsr \elldata{No.}{rank}{$ADE$-type}{$G$}

\vsrs \elldata{6}{$4$}{$D_{4}$}{$[1]$}
\elldata{7}{$4$}{$A_{4}$}{$[1]$}
\elldata{8}{$4$}{$A_{3}\,+\,A_{1}$}{$[1]$}
\elldata{9}{$4$}{$2\,A_{2}$}{$[1]$}
\elldata{10}{$4$}{$A_{2}\,+\,2\,A_{1}$}{$[1]$}
\elldata{11}{$4$}{$4\,A_{1}$}{$[1]$}

\vsr \elldata{No.}{rank}{$ADE$-type}{$G$}

\vsrs \elldata{12}{$5$}{$D_{5}$}{$[1]$}
\elldata{13}{$5$}{$D_{4}\,+\,A_{1}$}{$[1]$}
\elldata{14}{$5$}{$A_{5}$}{$[1]$}
\elldata{15}{$5$}{$A_{4}\,+\,A_{1}$}{$[1]$}
\elldata{16}{$5$}{$A_{3}\,+\,A_{2}$}{$[1]$}
\elldata{17}{$5$}{$A_{3}\,+\,2\,A_{1}$}{$[1]$}
\elldata{18}{$5$}{$2\,A_{2}\,+\,A_{1}$}{$[1]$}
\elldata{19}{$5$}{$A_{2}\,+\,3\,A_{1}$}{$[1]$}
\elldata{20}{$5$}{$5\,A_{1}$}{$[1]$}

\vsr \elldata{No.}{rank}{$ADE$-type}{$G$}

\vsrs \elldata{21}{$6$}{$E_{6}$}{$[1]$}
\elldata{22}{$6$}{$D_{6}$}{$[1]$}
\elldata{23}{$6$}{$D_{5}\,+\,A_{1}$}{$[1]$}
\elldata{24}{$6$}{$D_{4}\,+\,A_{2}$}{$[1]$}
\elldata{25}{$6$}{$D_{4}\,+\,2\,A_{1}$}{$[1]$}
\elldata{26}{$6$}{$A_{6}$}{$[1]$}
\elldata{27}{$6$}{$A_{5}\,+\,A_{1}$}{$[1]$}
\elldata{28}{$6$}{$A_{4}\,+\,A_{2}$}{$[1]$}
\elldata{29}{$6$}{$A_{4}\,+\,2\,A_{1}$}{$[1]$}
\elldata{30}{$6$}{$2\,A_{3}$}{$[1]$}
\elldata{31}{$6$}{$A_{3}\,+\,A_{2}\,+\,A_{1}$}{$[1]$}
\elldata{32}{$6$}{$A_{3}\,+\,3\,A_{1}$}{$[1]$}
\elldata{33}{$6$}{$3\,A_{2}$}{$[1]$}
\elldata{34}{$6$}{$2\,A_{2}\,+\,2\,A_{1}$}{$[1]$}
\elldata{35}{$6$}{$A_{2}\,+\,4\,A_{1}$}{$[1]$}
\elldata{36}{$6$}{$6\,A_{1}$}{$[1]$}

\vsr \elldata{No.}{rank}{$ADE$-type}{$G$}

\vsrs \elldata{37}{$7$}{$E_{7}$}{$[1]$}
\elldata{38}{$7$}{$E_{6}\,+\,A_{1}$}{$[1]$}
\elldata{39}{$7$}{$D_{7}$}{$[1]$}
\elldata{40}{$7$}{$D_{6}\,+\,A_{1}$}{$[1]$}
\elldata{41}{$7$}{$D_{5}\,+\,A_{2}$}{$[1]$}
\elldata{42}{$7$}{$D_{5}\,+\,2\,A_{1}$}{$[1]$}
\elldata{43}{$7$}{$D_{4}\,+\,A_{3}$}{$[1]$}
\elldata{44}{$7$}{$D_{4}\,+\,A_{2}\,+\,A_{1}$}{$[1]$}
\elldata{45}{$7$}{$D_{4}\,+\,3\,A_{1}$}{$[1]$}
\elldata{46}{$7$}{$A_{7}$}{$[1]$}
\elldata{47}{$7$}{$A_{6}\,+\,A_{1}$}{$[1]$}
\elldata{48}{$7$}{$A_{5}\,+\,A_{2}$}{$[1]$}
\elldata{49}{$7$}{$A_{5}\,+\,2\,A_{1}$}{$[1]$}
\elldata{50}{$7$}{$A_{4}\,+\,A_{3}$}{$[1]$}
\elldata{51}{$7$}{$A_{4}\,+\,A_{2}\,+\,A_{1}$}{$[1]$}
\elldata{52}{$7$}{$A_{4}\,+\,3\,A_{1}$}{$[1]$}
\elldata{53}{$7$}{$2\,A_{3}\,+\,A_{1}$}{$[1]$}
\elldata{54}{$7$}{$A_{3}\,+\,2\,A_{2}$}{$[1]$}
\elldata{55}{$7$}{$A_{3}\,+\,A_{2}\,+\,2\,A_{1}$}{$[1]$}
\elldata{56}{$7$}{$A_{3}\,+\,4\,A_{1}$}{$[1]$}
\elldata{57}{$7$}{$3\,A_{2}\,+\,A_{1}$}{$[1]$}
\elldata{58}{$7$}{$2\,A_{2}\,+\,3\,A_{1}$}{$[1]$}
\elldata{59}{$7$}{$A_{2}\,+\,5\,A_{1}$}{$[1]$}
\elldata{60}{$7$}{$7\,A_{1}$}{$[1]$}

\vsr \elldata{No.}{rank}{$ADE$-type}{$G$}

\vsrs \elldata{61}{$8$}{$E_{8}$}{$[1]$}
\elldata{62}{$8$}{$E_{7}\,+\,A_{1}$}{$[1]$}
\elldata{63}{$8$}{$E_{6}\,+\,A_{2}$}{$[1]$}
\elldata{64}{$8$}{$E_{6}\,+\,2\,A_{1}$}{$[1]$}
\elldata{65}{$8$}{$D_{8}$}{$[1]$}
\elldata{66}{$8$}{$D_{7}\,+\,A_{1}$}{$[1]$}
\elldata{67}{$8$}{$D_{6}\,+\,A_{2}$}{$[1]$}
\elldata{68}{$8$}{$D_{6}\,+\,2\,A_{1}$}{$[1]$}
\elldata{69}{$8$}{$D_{5}\,+\,A_{3}$}{$[1]$}
\elldata{70}{$8$}{$D_{5}\,+\,A_{2}\,+\,A_{1}$}{$[1]$}
\elldata{71}{$8$}{$D_{5}\,+\,3\,A_{1}$}{$[1]$}
\elldata{72}{$8$}{$2\,D_{4}$}{$[1]$}
\elldata{73}{$8$}{$D_{4}\,+\,A_{4}$}{$[1]$}
\elldata{74}{$8$}{$D_{4}\,+\,A_{3}\,+\,A_{1}$}{$[1]$}
\elldata{75}{$8$}{$D_{4}\,+\,2\,A_{2}$}{$[1]$}
\elldata{76}{$8$}{$D_{4}\,+\,A_{2}\,+\,2\,A_{1}$}{$[1]$}
\elldata{77}{$8$}{$D_{4}\,+\,4\,A_{1}$}{$[1]$}
\elldata{78}{$8$}{$A_{8}$}{$[1]$}
\elldata{79}{$8$}{$A_{7}\,+\,A_{1}$}{$[1]$}
\elldata{80}{$8$}{$A_{6}\,+\,A_{2}$}{$[1]$}
\elldata{81}{$8$}{$A_{6}\,+\,2\,A_{1}$}{$[1]$}
\elldata{82}{$8$}{$A_{5}\,+\,A_{3}$}{$[1]$}
\elldata{83}{$8$}{$A_{5}\,+\,A_{2}\,+\,A_{1}$}{$[1]$}
\elldata{84}{$8$}{$A_{5}\,+\,3\,A_{1}$}{$[1]$}
\elldata{85}{$8$}{$2\,A_{4}$}{$[1]$}
\elldata{86}{$8$}{$A_{4}\,+\,A_{3}\,+\,A_{1}$}{$[1]$}
\elldata{87}{$8$}{$A_{4}\,+\,2\,A_{2}$}{$[1]$}
\elldata{88}{$8$}{$A_{4}\,+\,A_{2}\,+\,2\,A_{1}$}{$[1]$}
\elldata{89}{$8$}{$A_{4}\,+\,4\,A_{1}$}{$[1]$}
\elldata{90}{$8$}{$2\,A_{3}\,+\,A_{2}$}{$[1]$}
\elldata{91}{$8$}{$2\,A_{3}\,+\,2\,A_{1}$}{$[1]$}
\elldata{92}{$8$}{$A_{3}\,+\,2\,A_{2}\,+\,A_{1}$}{$[1]$}
\elldata{93}{$8$}{$A_{3}\,+\,A_{2}\,+\,3\,A_{1}$}{$[1]$}
\elldata{94}{$8$}{$A_{3}\,+\,5\,A_{1}$}{$[1]$}
\elldata{95}{$8$}{$4\,A_{2}$}{$[1]$}
\elldata{96}{$8$}{$3\,A_{2}\,+\,2\,A_{1}$}{$[1]$}
\elldata{97}{$8$}{$2\,A_{2}\,+\,4\,A_{1}$}{$[1]$}
\elldata{98}{$8$}{$A_{2}\,+\,6\,A_{1}$}{$[1]$}
\elldata{99}{$8$}{$8\,A_{1}$}{$[2], \,[1]$}

\vsr \elldata{No.}{rank}{$ADE$-type}{$G$}

\vsrs \elldata{100}{$9$}{$E_{8}\,+\,A_{1}$}{$[1]$}
\elldata{101}{$9$}{$E_{7}\,+\,A_{2}$}{$[1]$}
\elldata{102}{$9$}{$E_{7}\,+\,2\,A_{1}$}{$[1]$}
\elldata{103}{$9$}{$E_{6}\,+\,A_{3}$}{$[1]$}
\elldata{104}{$9$}{$E_{6}\,+\,A_{2}\,+\,A_{1}$}{$[1]$}
\elldata{105}{$9$}{$E_{6}\,+\,3\,A_{1}$}{$[1]$}
\elldata{106}{$9$}{$D_{9}$}{$[1]$}
\elldata{107}{$9$}{$D_{8}\,+\,A_{1}$}{$[1]$}
\elldata{108}{$9$}{$D_{7}\,+\,A_{2}$}{$[1]$}
\elldata{109}{$9$}{$D_{7}\,+\,2\,A_{1}$}{$[1]$}
\elldata{110}{$9$}{$D_{6}\,+\,A_{3}$}{$[1]$}
\elldata{111}{$9$}{$D_{6}\,+\,A_{2}\,+\,A_{1}$}{$[1]$}
\elldata{112}{$9$}{$D_{6}\,+\,3\,A_{1}$}{$[1]$}
\elldata{113}{$9$}{$D_{5}\,+\,D_{4}$}{$[1]$}
\elldata{114}{$9$}{$D_{5}\,+\,A_{4}$}{$[1]$}
\elldata{115}{$9$}{$D_{5}\,+\,A_{3}\,+\,A_{1}$}{$[1]$}
\elldata{116}{$9$}{$D_{5}\,+\,2\,A_{2}$}{$[1]$}
\elldata{117}{$9$}{$D_{5}\,+\,A_{2}\,+\,2\,A_{1}$}{$[1]$}
\elldata{118}{$9$}{$D_{5}\,+\,4\,A_{1}$}{$[1]$}
\elldata{119}{$9$}{$2\,D_{4}\,+\,A_{1}$}{$[1]$}
\elldata{120}{$9$}{$D_{4}\,+\,A_{5}$}{$[1]$}
\elldata{121}{$9$}{$D_{4}\,+\,A_{4}\,+\,A_{1}$}{$[1]$}
\elldata{122}{$9$}{$D_{4}\,+\,A_{3}\,+\,A_{2}$}{$[1]$}
\elldata{123}{$9$}{$D_{4}\,+\,A_{3}\,+\,2\,A_{1}$}{$[1]$}
\elldata{124}{$9$}{$D_{4}\,+\,2\,A_{2}\,+\,A_{1}$}{$[1]$}
\elldata{125}{$9$}{$D_{4}\,+\,A_{2}\,+\,3\,A_{1}$}{$[1]$}
\elldata{126}{$9$}{$D_{4}\,+\,5\,A_{1}$}{$[1]$}
\elldata{127}{$9$}{$A_{9}$}{$[1]$}
\elldata{128}{$9$}{$A_{8}\,+\,A_{1}$}{$[1]$}
\elldata{129}{$9$}{$A_{7}\,+\,A_{2}$}{$[1]$}
\elldata{130}{$9$}{$A_{7}\,+\,2\,A_{1}$}{$[1]$}
\elldata{131}{$9$}{$A_{6}\,+\,A_{3}$}{$[1]$}
\elldata{132}{$9$}{$A_{6}\,+\,A_{2}\,+\,A_{1}$}{$[1]$}
\elldata{133}{$9$}{$A_{6}\,+\,3\,A_{1}$}{$[1]$}
\elldata{134}{$9$}{$A_{5}\,+\,A_{4}$}{$[1]$}
\elldata{135}{$9$}{$A_{5}\,+\,A_{3}\,+\,A_{1}$}{$[1]$}
\elldata{136}{$9$}{$A_{5}\,+\,2\,A_{2}$}{$[1]$}
\elldata{137}{$9$}{$A_{5}\,+\,A_{2}\,+\,2\,A_{1}$}{$[1]$}
\elldata{138}{$9$}{$A_{5}\,+\,4\,A_{1}$}{$[1]$}
\elldata{139}{$9$}{$2\,A_{4}\,+\,A_{1}$}{$[1]$}
\elldata{140}{$9$}{$A_{4}\,+\,A_{3}\,+\,A_{2}$}{$[1]$}
\elldata{141}{$9$}{$A_{4}\,+\,A_{3}\,+\,2\,A_{1}$}{$[1]$}
\elldata{142}{$9$}{$A_{4}\,+\,2\,A_{2}\,+\,A_{1}$}{$[1]$}
\elldata{143}{$9$}{$A_{4}\,+\,A_{2}\,+\,3\,A_{1}$}{$[1]$}
\elldata{144}{$9$}{$A_{4}\,+\,5\,A_{1}$}{$[1]$}
\elldata{145}{$9$}{$3\,A_{3}$}{$[1]$}
\elldata{146}{$9$}{$2\,A_{3}\,+\,A_{2}\,+\,A_{1}$}{$[1]$}
\elldata{147}{$9$}{$2\,A_{3}\,+\,3\,A_{1}$}{$[1]$}
\elldata{148}{$9$}{$A_{3}\,+\,3\,A_{2}$}{$[1]$}
\elldata{149}{$9$}{$A_{3}\,+\,2\,A_{2}\,+\,2\,A_{1}$}{$[1]$}
\elldata{150}{$9$}{$A_{3}\,+\,A_{2}\,+\,4\,A_{1}$}{$[1]$}
\elldata{151}{$9$}{$A_{3}\,+\,6\,A_{1}$}{$[2], \,[1]$}
\elldata{152}{$9$}{$4\,A_{2}\,+\,A_{1}$}{$[1]$}
\elldata{153}{$9$}{$3\,A_{2}\,+\,3\,A_{1}$}{$[1]$}
\elldata{154}{$9$}{$2\,A_{2}\,+\,5\,A_{1}$}{$[1]$}
\elldata{155}{$9$}{$A_{2}\,+\,7\,A_{1}$}{$[1]$}
\elldata{156}{$9$}{$9\,A_{1}$}{$[2], \,[1]$}

\vsr \elldata{No.}{rank}{$ADE$-type}{$G$}

\vsrs \elldata{157}{$10$}{$E_{8}\,+\,A_{2}$}{$[1]$}
\elldata{158}{$10$}{$E_{8}\,+\,2\,A_{1}$}{$[1]$}
\elldata{159}{$10$}{$E_{7}\,+\,A_{3}$}{$[1]$}
\elldata{160}{$10$}{$E_{7}\,+\,A_{2}\,+\,A_{1}$}{$[1]$}
\elldata{161}{$10$}{$E_{7}\,+\,3\,A_{1}$}{$[1]$}
\elldata{162}{$10$}{$E_{6}\,+\,D_{4}$}{$[1]$}
\elldata{163}{$10$}{$E_{6}\,+\,A_{4}$}{$[1]$}
\elldata{164}{$10$}{$E_{6}\,+\,A_{3}\,+\,A_{1}$}{$[1]$}
\elldata{165}{$10$}{$E_{6}\,+\,2\,A_{2}$}{$[1]$}
\elldata{166}{$10$}{$E_{6}\,+\,A_{2}\,+\,2\,A_{1}$}{$[1]$}
\elldata{167}{$10$}{$E_{6}\,+\,4\,A_{1}$}{$[1]$}
\elldata{168}{$10$}{$D_{10}$}{$[1]$}
\elldata{169}{$10$}{$D_{9}\,+\,A_{1}$}{$[1]$}
\elldata{170}{$10$}{$D_{8}\,+\,A_{2}$}{$[1]$}
\elldata{171}{$10$}{$D_{8}\,+\,2\,A_{1}$}{$[1]$}
\elldata{172}{$10$}{$D_{7}\,+\,A_{3}$}{$[1]$}
\elldata{173}{$10$}{$D_{7}\,+\,A_{2}\,+\,A_{1}$}{$[1]$}
\elldata{174}{$10$}{$D_{7}\,+\,3\,A_{1}$}{$[1]$}
\elldata{175}{$10$}{$D_{6}\,+\,D_{4}$}{$[1]$}
\elldata{176}{$10$}{$D_{6}\,+\,A_{4}$}{$[1]$}
\elldata{177}{$10$}{$D_{6}\,+\,A_{3}\,+\,A_{1}$}{$[1]$}
\elldata{178}{$10$}{$D_{6}\,+\,2\,A_{2}$}{$[1]$}
\elldata{179}{$10$}{$D_{6}\,+\,A_{2}\,+\,2\,A_{1}$}{$[1]$}
\elldata{180}{$10$}{$D_{6}\,+\,4\,A_{1}$}{$[1]$}
\elldata{181}{$10$}{$2\,D_{5}$}{$[1]$}
\elldata{182}{$10$}{$D_{5}\,+\,D_{4}\,+\,A_{1}$}{$[1]$}
\elldata{183}{$10$}{$D_{5}\,+\,A_{5}$}{$[1]$}
\elldata{184}{$10$}{$D_{5}\,+\,A_{4}\,+\,A_{1}$}{$[1]$}
\elldata{185}{$10$}{$D_{5}\,+\,A_{3}\,+\,A_{2}$}{$[1]$}
\elldata{186}{$10$}{$D_{5}\,+\,A_{3}\,+\,2\,A_{1}$}{$[1]$}
\elldata{187}{$10$}{$D_{5}\,+\,2\,A_{2}\,+\,A_{1}$}{$[1]$}
\elldata{188}{$10$}{$D_{5}\,+\,A_{2}\,+\,3\,A_{1}$}{$[1]$}
\elldata{189}{$10$}{$D_{5}\,+\,5\,A_{1}$}{$[1]$}
\elldata{190}{$10$}{$2\,D_{4}\,+\,A_{2}$}{$[1]$}
\elldata{191}{$10$}{$2\,D_{4}\,+\,2\,A_{1}$}{$[1]$}
\elldata{192}{$10$}{$D_{4}\,+\,A_{6}$}{$[1]$}
\elldata{193}{$10$}{$D_{4}\,+\,A_{5}\,+\,A_{1}$}{$[1]$}
\elldata{194}{$10$}{$D_{4}\,+\,A_{4}\,+\,A_{2}$}{$[1]$}
\elldata{195}{$10$}{$D_{4}\,+\,A_{4}\,+\,2\,A_{1}$}{$[1]$}
\elldata{196}{$10$}{$D_{4}\,+\,2\,A_{3}$}{$[1]$}
\elldata{197}{$10$}{$D_{4}\,+\,A_{3}\,+\,A_{2}\,+\,A_{1}$}{$[1]$}
\elldata{198}{$10$}{$D_{4}\,+\,A_{3}\,+\,3\,A_{1}$}{$[1]$}
\elldata{199}{$10$}{$D_{4}\,+\,3\,A_{2}$}{$[1]$}
\elldata{200}{$10$}{$D_{4}\,+\,2\,A_{2}\,+\,2\,A_{1}$}{$[1]$}
\elldata{201}{$10$}{$D_{4}\,+\,A_{2}\,+\,4\,A_{1}$}{$[1]$}
\elldata{202}{$10$}{$D_{4}\,+\,6\,A_{1}$}{$[2], \,[1]$}
\elldata{203}{$10$}{$A_{10}$}{$[1]$}
\elldata{204}{$10$}{$A_{9}\,+\,A_{1}$}{$[1]$}
\elldata{205}{$10$}{$A_{8}\,+\,A_{2}$}{$[1]$}
\elldata{206}{$10$}{$A_{8}\,+\,2\,A_{1}$}{$[1]$}
\elldata{207}{$10$}{$A_{7}\,+\,A_{3}$}{$[1]$}
\elldata{208}{$10$}{$A_{7}\,+\,A_{2}\,+\,A_{1}$}{$[1]$}
\elldata{209}{$10$}{$A_{7}\,+\,3\,A_{1}$}{$[1]$}
\elldata{210}{$10$}{$A_{6}\,+\,A_{4}$}{$[1]$}
\elldata{211}{$10$}{$A_{6}\,+\,A_{3}\,+\,A_{1}$}{$[1]$}
\elldata{212}{$10$}{$A_{6}\,+\,2\,A_{2}$}{$[1]$}
\elldata{213}{$10$}{$A_{6}\,+\,A_{2}\,+\,2\,A_{1}$}{$[1]$}
\elldata{214}{$10$}{$A_{6}\,+\,4\,A_{1}$}{$[1]$}
\elldata{215}{$10$}{$2\,A_{5}$}{$[1]$}
\elldata{216}{$10$}{$A_{5}\,+\,A_{4}\,+\,A_{1}$}{$[1]$}
\elldata{217}{$10$}{$A_{5}\,+\,A_{3}\,+\,A_{2}$}{$[1]$}
\elldata{218}{$10$}{$A_{5}\,+\,A_{3}\,+\,2\,A_{1}$}{$[1]$}
\elldata{219}{$10$}{$A_{5}\,+\,2\,A_{2}\,+\,A_{1}$}{$[1]$}
\elldata{220}{$10$}{$A_{5}\,+\,A_{2}\,+\,3\,A_{1}$}{$[1]$}
\elldata{221}{$10$}{$A_{5}\,+\,5\,A_{1}$}{$[2], \,[1]$}
\elldata{222}{$10$}{$2\,A_{4}\,+\,A_{2}$}{$[1]$}
\elldata{223}{$10$}{$2\,A_{4}\,+\,2\,A_{1}$}{$[1]$}
\elldata{224}{$10$}{$A_{4}\,+\,2\,A_{3}$}{$[1]$}
\elldata{225}{$10$}{$A_{4}\,+\,A_{3}\,+\,A_{2}\,+\,A_{1}$}{$[1]$}
\elldata{226}{$10$}{$A_{4}\,+\,A_{3}\,+\,3\,A_{1}$}{$[1]$}
\elldata{227}{$10$}{$A_{4}\,+\,3\,A_{2}$}{$[1]$}
\elldata{228}{$10$}{$A_{4}\,+\,2\,A_{2}\,+\,2\,A_{1}$}{$[1]$}
\elldata{229}{$10$}{$A_{4}\,+\,A_{2}\,+\,4\,A_{1}$}{$[1]$}
\elldata{230}{$10$}{$A_{4}\,+\,6\,A_{1}$}{$[1]$}
\elldata{231}{$10$}{$3\,A_{3}\,+\,A_{1}$}{$[1]$}
\elldata{232}{$10$}{$2\,A_{3}\,+\,2\,A_{2}$}{$[1]$}
\elldata{233}{$10$}{$2\,A_{3}\,+\,A_{2}\,+\,2\,A_{1}$}{$[1]$}
\elldata{234}{$10$}{$2\,A_{3}\,+\,4\,A_{1}$}{$[2], \,[1]$}
\elldata{235}{$10$}{$A_{3}\,+\,3\,A_{2}\,+\,A_{1}$}{$[1]$}
\elldata{236}{$10$}{$A_{3}\,+\,2\,A_{2}\,+\,3\,A_{1}$}{$[1]$}
\elldata{237}{$10$}{$A_{3}\,+\,A_{2}\,+\,5\,A_{1}$}{$[1]$}
\elldata{238}{$10$}{$A_{3}\,+\,7\,A_{1}$}{$[2], \,[1]$}
\elldata{239}{$10$}{$5\,A_{2}$}{$[1]$}
\elldata{240}{$10$}{$4\,A_{2}\,+\,2\,A_{1}$}{$[1]$}
\elldata{241}{$10$}{$3\,A_{2}\,+\,4\,A_{1}$}{$[1]$}
\elldata{242}{$10$}{$2\,A_{2}\,+\,6\,A_{1}$}{$[1]$}
\elldata{243}{$10$}{$A_{2}\,+\,8\,A_{1}$}{$[2], \,[1]$}
\elldata{244}{$10$}{$10\,A_{1}$}{$[2], \,[1]$}

\vsr \elldata{No.}{rank}{$ADE$-type}{$G$}

\vsrs \elldata{245}{$11$}{$E_{8}\,+\,A_{3}$}{$[1]$}
\elldata{246}{$11$}{$E_{8}\,+\,A_{2}\,+\,A_{1}$}{$[1]$}
\elldata{247}{$11$}{$E_{8}\,+\,3\,A_{1}$}{$[1]$}
\elldata{248}{$11$}{$E_{7}\,+\,D_{4}$}{$[1]$}
\elldata{249}{$11$}{$E_{7}\,+\,A_{4}$}{$[1]$}
\elldata{250}{$11$}{$E_{7}\,+\,A_{3}\,+\,A_{1}$}{$[1]$}
\elldata{251}{$11$}{$E_{7}\,+\,2\,A_{2}$}{$[1]$}
\elldata{252}{$11$}{$E_{7}\,+\,A_{2}\,+\,2\,A_{1}$}{$[1]$}
\elldata{253}{$11$}{$E_{7}\,+\,4\,A_{1}$}{$[1]$}
\elldata{254}{$11$}{$E_{6}\,+\,D_{5}$}{$[1]$}
\elldata{255}{$11$}{$E_{6}\,+\,D_{4}\,+\,A_{1}$}{$[1]$}
\elldata{256}{$11$}{$E_{6}\,+\,A_{5}$}{$[1]$}
\elldata{257}{$11$}{$E_{6}\,+\,A_{4}\,+\,A_{1}$}{$[1]$}
\elldata{258}{$11$}{$E_{6}\,+\,A_{3}\,+\,A_{2}$}{$[1]$}
\elldata{259}{$11$}{$E_{6}\,+\,A_{3}\,+\,2\,A_{1}$}{$[1]$}
\elldata{260}{$11$}{$E_{6}\,+\,2\,A_{2}\,+\,A_{1}$}{$[1]$}
\elldata{261}{$11$}{$E_{6}\,+\,A_{2}\,+\,3\,A_{1}$}{$[1]$}
\elldata{262}{$11$}{$E_{6}\,+\,5\,A_{1}$}{$[1]$}
\elldata{263}{$11$}{$D_{11}$}{$[1]$}
\elldata{264}{$11$}{$D_{10}\,+\,A_{1}$}{$[1]$}
\elldata{265}{$11$}{$D_{9}\,+\,A_{2}$}{$[1]$}
\elldata{266}{$11$}{$D_{9}\,+\,2\,A_{1}$}{$[1]$}
\elldata{267}{$11$}{$D_{8}\,+\,A_{3}$}{$[1]$}
\elldata{268}{$11$}{$D_{8}\,+\,A_{2}\,+\,A_{1}$}{$[1]$}
\elldata{269}{$11$}{$D_{8}\,+\,3\,A_{1}$}{$[1]$}
\elldata{270}{$11$}{$D_{7}\,+\,D_{4}$}{$[1]$}
\elldata{271}{$11$}{$D_{7}\,+\,A_{4}$}{$[1]$}
\elldata{272}{$11$}{$D_{7}\,+\,A_{3}\,+\,A_{1}$}{$[1]$}
\elldata{273}{$11$}{$D_{7}\,+\,2\,A_{2}$}{$[1]$}
\elldata{274}{$11$}{$D_{7}\,+\,A_{2}\,+\,2\,A_{1}$}{$[1]$}
\elldata{275}{$11$}{$D_{7}\,+\,4\,A_{1}$}{$[1]$}
\elldata{276}{$11$}{$D_{6}\,+\,D_{5}$}{$[1]$}
\elldata{277}{$11$}{$D_{6}\,+\,D_{4}\,+\,A_{1}$}{$[1]$}
\elldata{278}{$11$}{$D_{6}\,+\,A_{5}$}{$[1]$}
\elldata{279}{$11$}{$D_{6}\,+\,A_{4}\,+\,A_{1}$}{$[1]$}
\elldata{280}{$11$}{$D_{6}\,+\,A_{3}\,+\,A_{2}$}{$[1]$}
\elldata{281}{$11$}{$D_{6}\,+\,A_{3}\,+\,2\,A_{1}$}{$[1]$}
\elldata{282}{$11$}{$D_{6}\,+\,2\,A_{2}\,+\,A_{1}$}{$[1]$}
\elldata{283}{$11$}{$D_{6}\,+\,A_{2}\,+\,3\,A_{1}$}{$[1]$}
\elldata{284}{$11$}{$D_{6}\,+\,5\,A_{1}$}{$[2], \,[1]$}
\elldata{285}{$11$}{$2\,D_{5}\,+\,A_{1}$}{$[1]$}
\elldata{286}{$11$}{$D_{5}\,+\,D_{4}\,+\,A_{2}$}{$[1]$}
\elldata{287}{$11$}{$D_{5}\,+\,D_{4}\,+\,2\,A_{1}$}{$[1]$}
\elldata{288}{$11$}{$D_{5}\,+\,A_{6}$}{$[1]$}
\elldata{289}{$11$}{$D_{5}\,+\,A_{5}\,+\,A_{1}$}{$[1]$}
\elldata{290}{$11$}{$D_{5}\,+\,A_{4}\,+\,A_{2}$}{$[1]$}
\elldata{291}{$11$}{$D_{5}\,+\,A_{4}\,+\,2\,A_{1}$}{$[1]$}
\elldata{292}{$11$}{$D_{5}\,+\,2\,A_{3}$}{$[1]$}
\elldata{293}{$11$}{$D_{5}\,+\,A_{3}\,+\,A_{2}\,+\,A_{1}$}{$[1]$}
\elldata{294}{$11$}{$D_{5}\,+\,A_{3}\,+\,3\,A_{1}$}{$[1]$}
\elldata{295}{$11$}{$D_{5}\,+\,3\,A_{2}$}{$[1]$}
\elldata{296}{$11$}{$D_{5}\,+\,2\,A_{2}\,+\,2\,A_{1}$}{$[1]$}
\elldata{297}{$11$}{$D_{5}\,+\,A_{2}\,+\,4\,A_{1}$}{$[1]$}
\elldata{298}{$11$}{$D_{5}\,+\,6\,A_{1}$}{$[2], \,[1]$}
\elldata{299}{$11$}{$2\,D_{4}\,+\,A_{3}$}{$[1]$}
\elldata{300}{$11$}{$2\,D_{4}\,+\,A_{2}\,+\,A_{1}$}{$[1]$}
\elldata{301}{$11$}{$2\,D_{4}\,+\,3\,A_{1}$}{$[1]$}
\elldata{302}{$11$}{$D_{4}\,+\,A_{7}$}{$[1]$}
\elldata{303}{$11$}{$D_{4}\,+\,A_{6}\,+\,A_{1}$}{$[1]$}
\elldata{304}{$11$}{$D_{4}\,+\,A_{5}\,+\,A_{2}$}{$[1]$}
\elldata{305}{$11$}{$D_{4}\,+\,A_{5}\,+\,2\,A_{1}$}{$[1]$}
\elldata{306}{$11$}{$D_{4}\,+\,A_{4}\,+\,A_{3}$}{$[1]$}
\elldata{307}{$11$}{$D_{4}\,+\,A_{4}\,+\,A_{2}\,+\,A_{1}$}{$[1]$}
\elldata{308}{$11$}{$D_{4}\,+\,A_{4}\,+\,3\,A_{1}$}{$[1]$}
\elldata{309}{$11$}{$D_{4}\,+\,2\,A_{3}\,+\,A_{1}$}{$[1]$}
\elldata{310}{$11$}{$D_{4}\,+\,A_{3}\,+\,2\,A_{2}$}{$[1]$}
\elldata{311}{$11$}{$D_{4}\,+\,A_{3}\,+\,A_{2}\,+\,2\,A_{1}$}{$[1]$}
\elldata{312}{$11$}{$D_{4}\,+\,A_{3}\,+\,4\,A_{1}$}{$[2], \,[1]$}
\elldata{313}{$11$}{$D_{4}\,+\,3\,A_{2}\,+\,A_{1}$}{$[1]$}
\elldata{314}{$11$}{$D_{4}\,+\,2\,A_{2}\,+\,3\,A_{1}$}{$[1]$}
\elldata{315}{$11$}{$D_{4}\,+\,A_{2}\,+\,5\,A_{1}$}{$[1]$}
\elldata{316}{$11$}{$D_{4}\,+\,7\,A_{1}$}{$[2], \,[1]$}
\elldata{317}{$11$}{$A_{11}$}{$[1]$}
\elldata{318}{$11$}{$A_{10}\,+\,A_{1}$}{$[1]$}
\elldata{319}{$11$}{$A_{9}\,+\,A_{2}$}{$[1]$}
\elldata{320}{$11$}{$A_{9}\,+\,2\,A_{1}$}{$[1]$}
\elldata{321}{$11$}{$A_{8}\,+\,A_{3}$}{$[1]$}
\elldata{322}{$11$}{$A_{8}\,+\,A_{2}\,+\,A_{1}$}{$[1]$}
\elldata{323}{$11$}{$A_{8}\,+\,3\,A_{1}$}{$[1]$}
\elldata{324}{$11$}{$A_{7}\,+\,A_{4}$}{$[1]$}
\elldata{325}{$11$}{$A_{7}\,+\,A_{3}\,+\,A_{1}$}{$[1]$}
\elldata{326}{$11$}{$A_{7}\,+\,2\,A_{2}$}{$[1]$}
\elldata{327}{$11$}{$A_{7}\,+\,A_{2}\,+\,2\,A_{1}$}{$[1]$}
\elldata{328}{$11$}{$A_{7}\,+\,4\,A_{1}$}{$[2], \,[1]$}
\elldata{329}{$11$}{$A_{6}\,+\,A_{5}$}{$[1]$}
\elldata{330}{$11$}{$A_{6}\,+\,A_{4}\,+\,A_{1}$}{$[1]$}
\elldata{331}{$11$}{$A_{6}\,+\,A_{3}\,+\,A_{2}$}{$[1]$}
\elldata{332}{$11$}{$A_{6}\,+\,A_{3}\,+\,2\,A_{1}$}{$[1]$}
\elldata{333}{$11$}{$A_{6}\,+\,2\,A_{2}\,+\,A_{1}$}{$[1]$}
\elldata{334}{$11$}{$A_{6}\,+\,A_{2}\,+\,3\,A_{1}$}{$[1]$}
\elldata{335}{$11$}{$A_{6}\,+\,5\,A_{1}$}{$[1]$}
\elldata{336}{$11$}{$2\,A_{5}\,+\,A_{1}$}{$[1]$}
\elldata{337}{$11$}{$A_{5}\,+\,A_{4}\,+\,A_{2}$}{$[1]$}
\elldata{338}{$11$}{$A_{5}\,+\,A_{4}\,+\,2\,A_{1}$}{$[1]$}
\elldata{339}{$11$}{$A_{5}\,+\,2\,A_{3}$}{$[1]$}
\elldata{340}{$11$}{$A_{5}\,+\,A_{3}\,+\,A_{2}\,+\,A_{1}$}{$[1]$}
\elldata{341}{$11$}{$A_{5}\,+\,A_{3}\,+\,3\,A_{1}$}{$[2], \,[1]$}
\elldata{342}{$11$}{$A_{5}\,+\,3\,A_{2}$}{$[1]$}
\elldata{343}{$11$}{$A_{5}\,+\,2\,A_{2}\,+\,2\,A_{1}$}{$[1]$}
\elldata{344}{$11$}{$A_{5}\,+\,A_{2}\,+\,4\,A_{1}$}{$[1]$}
\elldata{345}{$11$}{$A_{5}\,+\,6\,A_{1}$}{$[2], \,[1]$}
\elldata{346}{$11$}{$2\,A_{4}\,+\,A_{3}$}{$[1]$}
\elldata{347}{$11$}{$2\,A_{4}\,+\,A_{2}\,+\,A_{1}$}{$[1]$}
\elldata{348}{$11$}{$2\,A_{4}\,+\,3\,A_{1}$}{$[1]$}
\elldata{349}{$11$}{$A_{4}\,+\,2\,A_{3}\,+\,A_{1}$}{$[1]$}
\elldata{350}{$11$}{$A_{4}\,+\,A_{3}\,+\,2\,A_{2}$}{$[1]$}
\elldata{351}{$11$}{$A_{4}\,+\,A_{3}\,+\,A_{2}\,+\,2\,A_{1}$}{$[1]$}
\elldata{352}{$11$}{$A_{4}\,+\,A_{3}\,+\,4\,A_{1}$}{$[1]$}
\elldata{353}{$11$}{$A_{4}\,+\,3\,A_{2}\,+\,A_{1}$}{$[1]$}
\elldata{354}{$11$}{$A_{4}\,+\,2\,A_{2}\,+\,3\,A_{1}$}{$[1]$}
\elldata{355}{$11$}{$A_{4}\,+\,A_{2}\,+\,5\,A_{1}$}{$[1]$}
\elldata{356}{$11$}{$A_{4}\,+\,7\,A_{1}$}{$[1]$}
\elldata{357}{$11$}{$3\,A_{3}\,+\,A_{2}$}{$[1]$}
\elldata{358}{$11$}{$3\,A_{3}\,+\,2\,A_{1}$}{$[2], \,[1]$}
\elldata{359}{$11$}{$2\,A_{3}\,+\,2\,A_{2}\,+\,A_{1}$}{$[1]$}
\elldata{360}{$11$}{$2\,A_{3}\,+\,A_{2}\,+\,3\,A_{1}$}{$[1]$}
\elldata{361}{$11$}{$2\,A_{3}\,+\,5\,A_{1}$}{$[2], \,[1]$}
\elldata{362}{$11$}{$A_{3}\,+\,4\,A_{2}$}{$[1]$}
\elldata{363}{$11$}{$A_{3}\,+\,3\,A_{2}\,+\,2\,A_{1}$}{$[1]$}
\elldata{364}{$11$}{$A_{3}\,+\,2\,A_{2}\,+\,4\,A_{1}$}{$[1]$}
\elldata{365}{$11$}{$A_{3}\,+\,A_{2}\,+\,6\,A_{1}$}{$[2], \,[1]$}
\elldata{366}{$11$}{$A_{3}\,+\,8\,A_{1}$}{$[2], \,[1]$}
\elldata{367}{$11$}{$5\,A_{2}\,+\,A_{1}$}{$[1]$}
\elldata{368}{$11$}{$4\,A_{2}\,+\,3\,A_{1}$}{$[1]$}
\elldata{369}{$11$}{$3\,A_{2}\,+\,5\,A_{1}$}{$[1]$}
\elldata{370}{$11$}{$2\,A_{2}\,+\,7\,A_{1}$}{$[1]$}
\elldata{371}{$11$}{$A_{2}\,+\,9\,A_{1}$}{$[2], \,[1]$}
\elldata{372}{$11$}{$11\,A_{1}$}{$[2]$}

\vsr \elldata{No.}{rank}{$ADE$-type}{$G$}

\vsrs \elldata{373}{$12$}{$E_{8}\,+\,D_{4}$}{$[1]$}
\elldata{374}{$12$}{$E_{8}\,+\,A_{4}$}{$[1]$}
\elldata{375}{$12$}{$E_{8}\,+\,A_{3}\,+\,A_{1}$}{$[1]$}
\elldata{376}{$12$}{$E_{8}\,+\,2\,A_{2}$}{$[1]$}
\elldata{377}{$12$}{$E_{8}\,+\,A_{2}\,+\,2\,A_{1}$}{$[1]$}
\elldata{378}{$12$}{$E_{8}\,+\,4\,A_{1}$}{$[1]$}
\elldata{379}{$12$}{$E_{7}\,+\,D_{5}$}{$[1]$}
\elldata{380}{$12$}{$E_{7}\,+\,D_{4}\,+\,A_{1}$}{$[1]$}
\elldata{381}{$12$}{$E_{7}\,+\,A_{5}$}{$[1]$}
\elldata{382}{$12$}{$E_{7}\,+\,A_{4}\,+\,A_{1}$}{$[1]$}
\elldata{383}{$12$}{$E_{7}\,+\,A_{3}\,+\,A_{2}$}{$[1]$}
\elldata{384}{$12$}{$E_{7}\,+\,A_{3}\,+\,2\,A_{1}$}{$[1]$}
\elldata{385}{$12$}{$E_{7}\,+\,2\,A_{2}\,+\,A_{1}$}{$[1]$}
\elldata{386}{$12$}{$E_{7}\,+\,A_{2}\,+\,3\,A_{1}$}{$[1]$}
\elldata{387}{$12$}{$E_{7}\,+\,5\,A_{1}$}{$[2], \,[1]$}
\elldata{388}{$12$}{$2\,E_{6}$}{$[1]$}
\elldata{389}{$12$}{$E_{6}\,+\,D_{6}$}{$[1]$}
\elldata{390}{$12$}{$E_{6}\,+\,D_{5}\,+\,A_{1}$}{$[1]$}
\elldata{391}{$12$}{$E_{6}\,+\,D_{4}\,+\,A_{2}$}{$[1]$}
\elldata{392}{$12$}{$E_{6}\,+\,D_{4}\,+\,2\,A_{1}$}{$[1]$}
\elldata{393}{$12$}{$E_{6}\,+\,A_{6}$}{$[1]$}
\elldata{394}{$12$}{$E_{6}\,+\,A_{5}\,+\,A_{1}$}{$[1]$}
\elldata{395}{$12$}{$E_{6}\,+\,A_{4}\,+\,A_{2}$}{$[1]$}
\elldata{396}{$12$}{$E_{6}\,+\,A_{4}\,+\,2\,A_{1}$}{$[1]$}
\elldata{397}{$12$}{$E_{6}\,+\,2\,A_{3}$}{$[1]$}
\elldata{398}{$12$}{$E_{6}\,+\,A_{3}\,+\,A_{2}\,+\,A_{1}$}{$[1]$}
\elldata{399}{$12$}{$E_{6}\,+\,A_{3}\,+\,3\,A_{1}$}{$[1]$}
\elldata{400}{$12$}{$E_{6}\,+\,3\,A_{2}$}{$[1]$}
\elldata{401}{$12$}{$E_{6}\,+\,2\,A_{2}\,+\,2\,A_{1}$}{$[1]$}
\elldata{402}{$12$}{$E_{6}\,+\,A_{2}\,+\,4\,A_{1}$}{$[1]$}
\elldata{403}{$12$}{$E_{6}\,+\,6\,A_{1}$}{$[1]$}
\elldata{404}{$12$}{$D_{12}$}{$[1]$}
\elldata{405}{$12$}{$D_{11}\,+\,A_{1}$}{$[1]$}
\elldata{406}{$12$}{$D_{10}\,+\,A_{2}$}{$[1]$}
\elldata{407}{$12$}{$D_{10}\,+\,2\,A_{1}$}{$[1]$}
\elldata{408}{$12$}{$D_{9}\,+\,A_{3}$}{$[1]$}
\elldata{409}{$12$}{$D_{9}\,+\,A_{2}\,+\,A_{1}$}{$[1]$}
\elldata{410}{$12$}{$D_{9}\,+\,3\,A_{1}$}{$[1]$}
\elldata{411}{$12$}{$D_{8}\,+\,D_{4}$}{$[1]$}
\elldata{412}{$12$}{$D_{8}\,+\,A_{4}$}{$[1]$}
\elldata{413}{$12$}{$D_{8}\,+\,A_{3}\,+\,A_{1}$}{$[1]$}
\elldata{414}{$12$}{$D_{8}\,+\,2\,A_{2}$}{$[1]$}
\elldata{415}{$12$}{$D_{8}\,+\,A_{2}\,+\,2\,A_{1}$}{$[1]$}
\elldata{416}{$12$}{$D_{8}\,+\,4\,A_{1}$}{$[2], \,[1]$}
\elldata{417}{$12$}{$D_{7}\,+\,D_{5}$}{$[1]$}
\elldata{418}{$12$}{$D_{7}\,+\,D_{4}\,+\,A_{1}$}{$[1]$}
\elldata{419}{$12$}{$D_{7}\,+\,A_{5}$}{$[1]$}
\elldata{420}{$12$}{$D_{7}\,+\,A_{4}\,+\,A_{1}$}{$[1]$}
\elldata{421}{$12$}{$D_{7}\,+\,A_{3}\,+\,A_{2}$}{$[1]$}
\elldata{422}{$12$}{$D_{7}\,+\,A_{3}\,+\,2\,A_{1}$}{$[1]$}
\elldata{423}{$12$}{$D_{7}\,+\,2\,A_{2}\,+\,A_{1}$}{$[1]$}
\elldata{424}{$12$}{$D_{7}\,+\,A_{2}\,+\,3\,A_{1}$}{$[1]$}
\elldata{425}{$12$}{$D_{7}\,+\,5\,A_{1}$}{$[1]$}
\elldata{426}{$12$}{$2\,D_{6}$}{$[1]$}
\elldata{427}{$12$}{$D_{6}\,+\,D_{5}\,+\,A_{1}$}{$[1]$}
\elldata{428}{$12$}{$D_{6}\,+\,D_{4}\,+\,A_{2}$}{$[1]$}
\elldata{429}{$12$}{$D_{6}\,+\,D_{4}\,+\,2\,A_{1}$}{$[1]$}
\elldata{430}{$12$}{$D_{6}\,+\,A_{6}$}{$[1]$}
\elldata{431}{$12$}{$D_{6}\,+\,A_{5}\,+\,A_{1}$}{$[1]$}
\elldata{432}{$12$}{$D_{6}\,+\,A_{4}\,+\,A_{2}$}{$[1]$}
\elldata{433}{$12$}{$D_{6}\,+\,A_{4}\,+\,2\,A_{1}$}{$[1]$}
\elldata{434}{$12$}{$D_{6}\,+\,2\,A_{3}$}{$[1]$}
\elldata{435}{$12$}{$D_{6}\,+\,A_{3}\,+\,A_{2}\,+\,A_{1}$}{$[1]$}
\elldata{436}{$12$}{$D_{6}\,+\,A_{3}\,+\,3\,A_{1}$}{$[2], \,[1]$}
\elldata{437}{$12$}{$D_{6}\,+\,3\,A_{2}$}{$[1]$}
\elldata{438}{$12$}{$D_{6}\,+\,2\,A_{2}\,+\,2\,A_{1}$}{$[1]$}
\elldata{439}{$12$}{$D_{6}\,+\,A_{2}\,+\,4\,A_{1}$}{$[1]$}
\elldata{440}{$12$}{$D_{6}\,+\,6\,A_{1}$}{$[2], \,[1]$}
\elldata{441}{$12$}{$2\,D_{5}\,+\,A_{2}$}{$[1]$}
\elldata{442}{$12$}{$2\,D_{5}\,+\,2\,A_{1}$}{$[1]$}
\elldata{443}{$12$}{$D_{5}\,+\,D_{4}\,+\,A_{3}$}{$[1]$}
\elldata{444}{$12$}{$D_{5}\,+\,D_{4}\,+\,A_{2}\,+\,A_{1}$}{$[1]$}
\elldata{445}{$12$}{$D_{5}\,+\,D_{4}\,+\,3\,A_{1}$}{$[1]$}
\elldata{446}{$12$}{$D_{5}\,+\,A_{7}$}{$[1]$}
\elldata{447}{$12$}{$D_{5}\,+\,A_{6}\,+\,A_{1}$}{$[1]$}
\elldata{448}{$12$}{$D_{5}\,+\,A_{5}\,+\,A_{2}$}{$[1]$}
\elldata{449}{$12$}{$D_{5}\,+\,A_{5}\,+\,2\,A_{1}$}{$[1]$}
\elldata{450}{$12$}{$D_{5}\,+\,A_{4}\,+\,A_{3}$}{$[1]$}
\elldata{451}{$12$}{$D_{5}\,+\,A_{4}\,+\,A_{2}\,+\,A_{1}$}{$[1]$}
\elldata{452}{$12$}{$D_{5}\,+\,A_{4}\,+\,3\,A_{1}$}{$[1]$}
\elldata{453}{$12$}{$D_{5}\,+\,2\,A_{3}\,+\,A_{1}$}{$[1]$}
\elldata{454}{$12$}{$D_{5}\,+\,A_{3}\,+\,2\,A_{2}$}{$[1]$}
\elldata{455}{$12$}{$D_{5}\,+\,A_{3}\,+\,A_{2}\,+\,2\,A_{1}$}{$[1]$}
\elldata{456}{$12$}{$D_{5}\,+\,A_{3}\,+\,4\,A_{1}$}{$[2], \,[1]$}
\elldata{457}{$12$}{$D_{5}\,+\,3\,A_{2}\,+\,A_{1}$}{$[1]$}
\elldata{458}{$12$}{$D_{5}\,+\,2\,A_{2}\,+\,3\,A_{1}$}{$[1]$}
\elldata{459}{$12$}{$D_{5}\,+\,A_{2}\,+\,5\,A_{1}$}{$[1]$}
\elldata{460}{$12$}{$D_{5}\,+\,7\,A_{1}$}{$[2], \,[1]$}
\elldata{461}{$12$}{$3\,D_{4}$}{$[1]$}
\elldata{462}{$12$}{$2\,D_{4}\,+\,A_{4}$}{$[1]$}
\elldata{463}{$12$}{$2\,D_{4}\,+\,A_{3}\,+\,A_{1}$}{$[1]$}
\elldata{464}{$12$}{$2\,D_{4}\,+\,2\,A_{2}$}{$[1]$}
\elldata{465}{$12$}{$2\,D_{4}\,+\,A_{2}\,+\,2\,A_{1}$}{$[1]$}
\elldata{466}{$12$}{$2\,D_{4}\,+\,4\,A_{1}$}{$[2], \,[1]$}
\elldata{467}{$12$}{$D_{4}\,+\,A_{8}$}{$[1]$}
\elldata{468}{$12$}{$D_{4}\,+\,A_{7}\,+\,A_{1}$}{$[1]$}
\elldata{469}{$12$}{$D_{4}\,+\,A_{6}\,+\,A_{2}$}{$[1]$}
\elldata{470}{$12$}{$D_{4}\,+\,A_{6}\,+\,2\,A_{1}$}{$[1]$}
\elldata{471}{$12$}{$D_{4}\,+\,A_{5}\,+\,A_{3}$}{$[1]$}
\elldata{472}{$12$}{$D_{4}\,+\,A_{5}\,+\,A_{2}\,+\,A_{1}$}{$[1]$}
\elldata{473}{$12$}{$D_{4}\,+\,A_{5}\,+\,3\,A_{1}$}{$[2], \,[1]$}
\elldata{474}{$12$}{$D_{4}\,+\,2\,A_{4}$}{$[1]$}
\elldata{475}{$12$}{$D_{4}\,+\,A_{4}\,+\,A_{3}\,+\,A_{1}$}{$[1]$}
\elldata{476}{$12$}{$D_{4}\,+\,A_{4}\,+\,2\,A_{2}$}{$[1]$}
\elldata{477}{$12$}{$D_{4}\,+\,A_{4}\,+\,A_{2}\,+\,2\,A_{1}$}{$[1]$}
\elldata{478}{$12$}{$D_{4}\,+\,A_{4}\,+\,4\,A_{1}$}{$[1]$}
\elldata{479}{$12$}{$D_{4}\,+\,2\,A_{3}\,+\,A_{2}$}{$[1]$}
\elldata{480}{$12$}{$D_{4}\,+\,2\,A_{3}\,+\,2\,A_{1}$}{$[2], \,[1]$}
\elldata{481}{$12$}{$D_{4}\,+\,A_{3}\,+\,2\,A_{2}\,+\,A_{1}$}{$[1]$}
\elldata{482}{$12$}{$D_{4}\,+\,A_{3}\,+\,A_{2}\,+\,3\,A_{1}$}{$[1]$}
\elldata{483}{$12$}{$D_{4}\,+\,A_{3}\,+\,5\,A_{1}$}{$[2], \,[1]$}
\elldata{484}{$12$}{$D_{4}\,+\,4\,A_{2}$}{$[1]$}
\elldata{485}{$12$}{$D_{4}\,+\,3\,A_{2}\,+\,2\,A_{1}$}{$[1]$}
\elldata{486}{$12$}{$D_{4}\,+\,2\,A_{2}\,+\,4\,A_{1}$}{$[1]$}
\elldata{487}{$12$}{$D_{4}\,+\,A_{2}\,+\,6\,A_{1}$}{$[2], \,[1]$}
\elldata{488}{$12$}{$D_{4}\,+\,8\,A_{1}$}{$[2]$}
\elldata{489}{$12$}{$A_{12}$}{$[1]$}
\elldata{490}{$12$}{$A_{11}\,+\,A_{1}$}{$[1]$}
\elldata{491}{$12$}{$A_{10}\,+\,A_{2}$}{$[1]$}
\elldata{492}{$12$}{$A_{10}\,+\,2\,A_{1}$}{$[1]$}
\elldata{493}{$12$}{$A_{9}\,+\,A_{3}$}{$[1]$}
\elldata{494}{$12$}{$A_{9}\,+\,A_{2}\,+\,A_{1}$}{$[1]$}
\elldata{495}{$12$}{$A_{9}\,+\,3\,A_{1}$}{$[2], \,[1]$}
\elldata{496}{$12$}{$A_{8}\,+\,A_{4}$}{$[1]$}
\elldata{497}{$12$}{$A_{8}\,+\,A_{3}\,+\,A_{1}$}{$[1]$}
\elldata{498}{$12$}{$A_{8}\,+\,2\,A_{2}$}{$[1]$}
\elldata{499}{$12$}{$A_{8}\,+\,A_{2}\,+\,2\,A_{1}$}{$[1]$}
\elldata{500}{$12$}{$A_{8}\,+\,4\,A_{1}$}{$[1]$}
\elldata{501}{$12$}{$A_{7}\,+\,A_{5}$}{$[1]$}
\elldata{502}{$12$}{$A_{7}\,+\,A_{4}\,+\,A_{1}$}{$[1]$}
\elldata{503}{$12$}{$A_{7}\,+\,A_{3}\,+\,A_{2}$}{$[1]$}
\elldata{504}{$12$}{$A_{7}\,+\,A_{3}\,+\,2\,A_{1}$}{$[2], \,[1]$}
\elldata{505}{$12$}{$A_{7}\,+\,2\,A_{2}\,+\,A_{1}$}{$[1]$}
\elldata{506}{$12$}{$A_{7}\,+\,A_{2}\,+\,3\,A_{1}$}{$[1]$}
\elldata{507}{$12$}{$A_{7}\,+\,5\,A_{1}$}{$[2], \,[1]$}
\elldata{508}{$12$}{$2\,A_{6}$}{$[1]$}
\elldata{509}{$12$}{$A_{6}\,+\,A_{5}\,+\,A_{1}$}{$[1]$}
\elldata{510}{$12$}{$A_{6}\,+\,A_{4}\,+\,A_{2}$}{$[1]$}
\elldata{511}{$12$}{$A_{6}\,+\,A_{4}\,+\,2\,A_{1}$}{$[1]$}
\elldata{512}{$12$}{$A_{6}\,+\,2\,A_{3}$}{$[1]$}
\elldata{513}{$12$}{$A_{6}\,+\,A_{3}\,+\,A_{2}\,+\,A_{1}$}{$[1]$}
\elldata{514}{$12$}{$A_{6}\,+\,A_{3}\,+\,3\,A_{1}$}{$[1]$}
\elldata{515}{$12$}{$A_{6}\,+\,3\,A_{2}$}{$[1]$}
\elldata{516}{$12$}{$A_{6}\,+\,2\,A_{2}\,+\,2\,A_{1}$}{$[1]$}
\elldata{517}{$12$}{$A_{6}\,+\,A_{2}\,+\,4\,A_{1}$}{$[1]$}
\elldata{518}{$12$}{$A_{6}\,+\,6\,A_{1}$}{$[1]$}
\elldata{519}{$12$}{$2\,A_{5}\,+\,A_{2}$}{$[1]$}
\elldata{520}{$12$}{$2\,A_{5}\,+\,2\,A_{1}$}{$[2], \,[1]$}
\elldata{521}{$12$}{$A_{5}\,+\,A_{4}\,+\,A_{3}$}{$[1]$}
\elldata{522}{$12$}{$A_{5}\,+\,A_{4}\,+\,A_{2}\,+\,A_{1}$}{$[1]$}
\elldata{523}{$12$}{$A_{5}\,+\,A_{4}\,+\,3\,A_{1}$}{$[1]$}
\elldata{524}{$12$}{$A_{5}\,+\,2\,A_{3}\,+\,A_{1}$}{$[2], \,[1]$}
\elldata{525}{$12$}{$A_{5}\,+\,A_{3}\,+\,2\,A_{2}$}{$[1]$}
\elldata{526}{$12$}{$A_{5}\,+\,A_{3}\,+\,A_{2}\,+\,2\,A_{1}$}{$[1]$}
\elldata{527}{$12$}{$A_{5}\,+\,A_{3}\,+\,4\,A_{1}$}{$[2], \,[1]$}
\elldata{528}{$12$}{$A_{5}\,+\,3\,A_{2}\,+\,A_{1}$}{$[1]$}
\elldata{529}{$12$}{$A_{5}\,+\,2\,A_{2}\,+\,3\,A_{1}$}{$[1]$}
\elldata{530}{$12$}{$A_{5}\,+\,A_{2}\,+\,5\,A_{1}$}{$[2], \,[1]$}
\elldata{531}{$12$}{$A_{5}\,+\,7\,A_{1}$}{$[2], \,[1]$}
\elldata{532}{$12$}{$3\,A_{4}$}{$[1]$}
\elldata{533}{$12$}{$2\,A_{4}\,+\,A_{3}\,+\,A_{1}$}{$[1]$}
\elldata{534}{$12$}{$2\,A_{4}\,+\,2\,A_{2}$}{$[1]$}
\elldata{535}{$12$}{$2\,A_{4}\,+\,A_{2}\,+\,2\,A_{1}$}{$[1]$}
\elldata{536}{$12$}{$2\,A_{4}\,+\,4\,A_{1}$}{$[1]$}
\elldata{537}{$12$}{$A_{4}\,+\,2\,A_{3}\,+\,A_{2}$}{$[1]$}
\elldata{538}{$12$}{$A_{4}\,+\,2\,A_{3}\,+\,2\,A_{1}$}{$[1]$}
\elldata{539}{$12$}{$A_{4}\,+\,A_{3}\,+\,2\,A_{2}\,+\,A_{1}$}{$[1]$}
\elldata{540}{$12$}{$A_{4}\,+\,A_{3}\,+\,A_{2}\,+\,3\,A_{1}$}{$[1]$}
\elldata{541}{$12$}{$A_{4}\,+\,A_{3}\,+\,5\,A_{1}$}{$[1]$}
\elldata{542}{$12$}{$A_{4}\,+\,4\,A_{2}$}{$[1]$}
\elldata{543}{$12$}{$A_{4}\,+\,3\,A_{2}\,+\,2\,A_{1}$}{$[1]$}
\elldata{544}{$12$}{$A_{4}\,+\,2\,A_{2}\,+\,4\,A_{1}$}{$[1]$}
\elldata{545}{$12$}{$A_{4}\,+\,A_{2}\,+\,6\,A_{1}$}{$[1]$}
\elldata{546}{$12$}{$A_{4}\,+\,8\,A_{1}$}{$[2], \,[1]$}
\elldata{547}{$12$}{$4\,A_{3}$}{$[2], \,[1]$}
\elldata{548}{$12$}{$3\,A_{3}\,+\,A_{2}\,+\,A_{1}$}{$[1]$}
\elldata{549}{$12$}{$3\,A_{3}\,+\,3\,A_{1}$}{$[2], \,[1]$}
\elldata{550}{$12$}{$2\,A_{3}\,+\,3\,A_{2}$}{$[1]$}
\elldata{551}{$12$}{$2\,A_{3}\,+\,2\,A_{2}\,+\,2\,A_{1}$}{$[1]$}
\elldata{552}{$12$}{$2\,A_{3}\,+\,A_{2}\,+\,4\,A_{1}$}{$[2], \,[1]$}
\elldata{553}{$12$}{$2\,A_{3}\,+\,6\,A_{1}$}{$[2], \,[1]$}
\elldata{554}{$12$}{$A_{3}\,+\,4\,A_{2}\,+\,A_{1}$}{$[1]$}
\elldata{555}{$12$}{$A_{3}\,+\,3\,A_{2}\,+\,3\,A_{1}$}{$[1]$}
\elldata{556}{$12$}{$A_{3}\,+\,2\,A_{2}\,+\,5\,A_{1}$}{$[1]$}
\elldata{557}{$12$}{$A_{3}\,+\,A_{2}\,+\,7\,A_{1}$}{$[2], \,[1]$}
\elldata{558}{$12$}{$A_{3}\,+\,9\,A_{1}$}{$[2]$}
\elldata{559}{$12$}{$6\,A_{2}$}{$[3], \,[1]$}
\elldata{560}{$12$}{$5\,A_{2}\,+\,2\,A_{1}$}{$[1]$}
\elldata{561}{$12$}{$4\,A_{2}\,+\,4\,A_{1}$}{$[1]$}
\elldata{562}{$12$}{$3\,A_{2}\,+\,6\,A_{1}$}{$[1]$}
\elldata{563}{$12$}{$2\,A_{2}\,+\,8\,A_{1}$}{$[2], \,[1]$}
\elldata{564}{$12$}{$A_{2}\,+\,10\,A_{1}$}{$[2]$}
\elldata{565}{$12$}{$12\,A_{1}$}{$[2, 2]$}

\vsr \elldata{No.}{rank}{$ADE$-type}{$G$}

\vsrs \elldata{566}{$13$}{$E_{8}\,+\,D_{5}$}{$[1]$}
\elldata{567}{$13$}{$E_{8}\,+\,D_{4}\,+\,A_{1}$}{$[1]$}
\elldata{568}{$13$}{$E_{8}\,+\,A_{5}$}{$[1]$}
\elldata{569}{$13$}{$E_{8}\,+\,A_{4}\,+\,A_{1}$}{$[1]$}
\elldata{570}{$13$}{$E_{8}\,+\,A_{3}\,+\,A_{2}$}{$[1]$}
\elldata{571}{$13$}{$E_{8}\,+\,A_{3}\,+\,2\,A_{1}$}{$[1]$}
\elldata{572}{$13$}{$E_{8}\,+\,2\,A_{2}\,+\,A_{1}$}{$[1]$}
\elldata{573}{$13$}{$E_{8}\,+\,A_{2}\,+\,3\,A_{1}$}{$[1]$}
\elldata{574}{$13$}{$E_{8}\,+\,5\,A_{1}$}{$[1]$}
\elldata{575}{$13$}{$E_{7}\,+\,E_{6}$}{$[1]$}
\elldata{576}{$13$}{$E_{7}\,+\,D_{6}$}{$[1]$}
\elldata{577}{$13$}{$E_{7}\,+\,D_{5}\,+\,A_{1}$}{$[1]$}
\elldata{578}{$13$}{$E_{7}\,+\,D_{4}\,+\,A_{2}$}{$[1]$}
\elldata{579}{$13$}{$E_{7}\,+\,D_{4}\,+\,2\,A_{1}$}{$[1]$}
\elldata{580}{$13$}{$E_{7}\,+\,A_{6}$}{$[1]$}
\elldata{581}{$13$}{$E_{7}\,+\,A_{5}\,+\,A_{1}$}{$[1]$}
\elldata{582}{$13$}{$E_{7}\,+\,A_{4}\,+\,A_{2}$}{$[1]$}
\elldata{583}{$13$}{$E_{7}\,+\,A_{4}\,+\,2\,A_{1}$}{$[1]$}
\elldata{584}{$13$}{$E_{7}\,+\,2\,A_{3}$}{$[1]$}
\elldata{585}{$13$}{$E_{7}\,+\,A_{3}\,+\,A_{2}\,+\,A_{1}$}{$[1]$}
\elldata{586}{$13$}{$E_{7}\,+\,A_{3}\,+\,3\,A_{1}$}{$[2], \,[1]$}
\elldata{587}{$13$}{$E_{7}\,+\,3\,A_{2}$}{$[1]$}
\elldata{588}{$13$}{$E_{7}\,+\,2\,A_{2}\,+\,2\,A_{1}$}{$[1]$}
\elldata{589}{$13$}{$E_{7}\,+\,A_{2}\,+\,4\,A_{1}$}{$[1]$}
\elldata{590}{$13$}{$E_{7}\,+\,6\,A_{1}$}{$[2], \,[1]$}
\elldata{591}{$13$}{$2\,E_{6}\,+\,A_{1}$}{$[1]$}
\elldata{592}{$13$}{$E_{6}\,+\,D_{7}$}{$[1]$}
\elldata{593}{$13$}{$E_{6}\,+\,D_{6}\,+\,A_{1}$}{$[1]$}
\elldata{594}{$13$}{$E_{6}\,+\,D_{5}\,+\,A_{2}$}{$[1]$}
\elldata{595}{$13$}{$E_{6}\,+\,D_{5}\,+\,2\,A_{1}$}{$[1]$}
\elldata{596}{$13$}{$E_{6}\,+\,D_{4}\,+\,A_{3}$}{$[1]$}
\elldata{597}{$13$}{$E_{6}\,+\,D_{4}\,+\,A_{2}\,+\,A_{1}$}{$[1]$}
\elldata{598}{$13$}{$E_{6}\,+\,D_{4}\,+\,3\,A_{1}$}{$[1]$}
\elldata{599}{$13$}{$E_{6}\,+\,A_{7}$}{$[1]$}
\elldata{600}{$13$}{$E_{6}\,+\,A_{6}\,+\,A_{1}$}{$[1]$}
\elldata{601}{$13$}{$E_{6}\,+\,A_{5}\,+\,A_{2}$}{$[1]$}
\elldata{602}{$13$}{$E_{6}\,+\,A_{5}\,+\,2\,A_{1}$}{$[1]$}
\elldata{603}{$13$}{$E_{6}\,+\,A_{4}\,+\,A_{3}$}{$[1]$}
\elldata{604}{$13$}{$E_{6}\,+\,A_{4}\,+\,A_{2}\,+\,A_{1}$}{$[1]$}
\elldata{605}{$13$}{$E_{6}\,+\,A_{4}\,+\,3\,A_{1}$}{$[1]$}
\elldata{606}{$13$}{$E_{6}\,+\,2\,A_{3}\,+\,A_{1}$}{$[1]$}
\elldata{607}{$13$}{$E_{6}\,+\,A_{3}\,+\,2\,A_{2}$}{$[1]$}
\elldata{608}{$13$}{$E_{6}\,+\,A_{3}\,+\,A_{2}\,+\,2\,A_{1}$}{$[1]$}
\elldata{609}{$13$}{$E_{6}\,+\,A_{3}\,+\,4\,A_{1}$}{$[1]$}
\elldata{610}{$13$}{$E_{6}\,+\,3\,A_{2}\,+\,A_{1}$}{$[1]$}
\elldata{611}{$13$}{$E_{6}\,+\,2\,A_{2}\,+\,3\,A_{1}$}{$[1]$}
\elldata{612}{$13$}{$E_{6}\,+\,A_{2}\,+\,5\,A_{1}$}{$[1]$}
\elldata{613}{$13$}{$E_{6}\,+\,7\,A_{1}$}{$[1]$}
\elldata{614}{$13$}{$D_{13}$}{$[1]$}
\elldata{615}{$13$}{$D_{12}\,+\,A_{1}$}{$[1]$}
\elldata{616}{$13$}{$D_{11}\,+\,A_{2}$}{$[1]$}
\elldata{617}{$13$}{$D_{11}\,+\,2\,A_{1}$}{$[1]$}
\elldata{618}{$13$}{$D_{10}\,+\,A_{3}$}{$[1]$}
\elldata{619}{$13$}{$D_{10}\,+\,A_{2}\,+\,A_{1}$}{$[1]$}
\elldata{620}{$13$}{$D_{10}\,+\,3\,A_{1}$}{$[2], \,[1]$}
\elldata{621}{$13$}{$D_{9}\,+\,D_{4}$}{$[1]$}
\elldata{622}{$13$}{$D_{9}\,+\,A_{4}$}{$[1]$}
\elldata{623}{$13$}{$D_{9}\,+\,A_{3}\,+\,A_{1}$}{$[1]$}
\elldata{624}{$13$}{$D_{9}\,+\,2\,A_{2}$}{$[1]$}
\elldata{625}{$13$}{$D_{9}\,+\,A_{2}\,+\,2\,A_{1}$}{$[1]$}
\elldata{626}{$13$}{$D_{9}\,+\,4\,A_{1}$}{$[1]$}
\elldata{627}{$13$}{$D_{8}\,+\,D_{5}$}{$[1]$}
\elldata{628}{$13$}{$D_{8}\,+\,D_{4}\,+\,A_{1}$}{$[1]$}
\elldata{629}{$13$}{$D_{8}\,+\,A_{5}$}{$[1]$}
\elldata{630}{$13$}{$D_{8}\,+\,A_{4}\,+\,A_{1}$}{$[1]$}
\elldata{631}{$13$}{$D_{8}\,+\,A_{3}\,+\,A_{2}$}{$[1]$}
\elldata{632}{$13$}{$D_{8}\,+\,A_{3}\,+\,2\,A_{1}$}{$[2], \,[1]$}
\elldata{633}{$13$}{$D_{8}\,+\,2\,A_{2}\,+\,A_{1}$}{$[1]$}
\elldata{634}{$13$}{$D_{8}\,+\,A_{2}\,+\,3\,A_{1}$}{$[1]$}
\elldata{635}{$13$}{$D_{8}\,+\,5\,A_{1}$}{$[2], \,[1]$}
\elldata{636}{$13$}{$D_{7}\,+\,D_{6}$}{$[1]$}
\elldata{637}{$13$}{$D_{7}\,+\,D_{5}\,+\,A_{1}$}{$[1]$}
\elldata{638}{$13$}{$D_{7}\,+\,D_{4}\,+\,A_{2}$}{$[1]$}
\elldata{639}{$13$}{$D_{7}\,+\,D_{4}\,+\,2\,A_{1}$}{$[1]$}
\elldata{640}{$13$}{$D_{7}\,+\,A_{6}$}{$[1]$}
\elldata{641}{$13$}{$D_{7}\,+\,A_{5}\,+\,A_{1}$}{$[1]$}
\elldata{642}{$13$}{$D_{7}\,+\,A_{4}\,+\,A_{2}$}{$[1]$}
\elldata{643}{$13$}{$D_{7}\,+\,A_{4}\,+\,2\,A_{1}$}{$[1]$}
\elldata{644}{$13$}{$D_{7}\,+\,2\,A_{3}$}{$[1]$}
\elldata{645}{$13$}{$D_{7}\,+\,A_{3}\,+\,A_{2}\,+\,A_{1}$}{$[1]$}
\elldata{646}{$13$}{$D_{7}\,+\,A_{3}\,+\,3\,A_{1}$}{$[1]$}
\elldata{647}{$13$}{$D_{7}\,+\,3\,A_{2}$}{$[1]$}
\elldata{648}{$13$}{$D_{7}\,+\,2\,A_{2}\,+\,2\,A_{1}$}{$[1]$}
\elldata{649}{$13$}{$D_{7}\,+\,A_{2}\,+\,4\,A_{1}$}{$[1]$}
\elldata{650}{$13$}{$D_{7}\,+\,6\,A_{1}$}{$[2], \,[1]$}
\elldata{651}{$13$}{$2\,D_{6}\,+\,A_{1}$}{$[1]$}
\elldata{652}{$13$}{$D_{6}\,+\,D_{5}\,+\,A_{2}$}{$[1]$}
\elldata{653}{$13$}{$D_{6}\,+\,D_{5}\,+\,2\,A_{1}$}{$[1]$}
\elldata{654}{$13$}{$D_{6}\,+\,D_{4}\,+\,A_{3}$}{$[1]$}
\elldata{655}{$13$}{$D_{6}\,+\,D_{4}\,+\,A_{2}\,+\,A_{1}$}{$[1]$}
\elldata{656}{$13$}{$D_{6}\,+\,D_{4}\,+\,3\,A_{1}$}{$[2], \,[1]$}
\elldata{657}{$13$}{$D_{6}\,+\,A_{7}$}{$[1]$}
\elldata{658}{$13$}{$D_{6}\,+\,A_{6}\,+\,A_{1}$}{$[1]$}
\elldata{659}{$13$}{$D_{6}\,+\,A_{5}\,+\,A_{2}$}{$[1]$}
\elldata{660}{$13$}{$D_{6}\,+\,A_{5}\,+\,2\,A_{1}$}{$[2], \,[1]$}
\elldata{661}{$13$}{$D_{6}\,+\,A_{4}\,+\,A_{3}$}{$[1]$}
\elldata{662}{$13$}{$D_{6}\,+\,A_{4}\,+\,A_{2}\,+\,A_{1}$}{$[1]$}
\elldata{663}{$13$}{$D_{6}\,+\,A_{4}\,+\,3\,A_{1}$}{$[1]$}
\elldata{664}{$13$}{$D_{6}\,+\,2\,A_{3}\,+\,A_{1}$}{$[2], \,[1]$}
\elldata{665}{$13$}{$D_{6}\,+\,A_{3}\,+\,2\,A_{2}$}{$[1]$}
\elldata{666}{$13$}{$D_{6}\,+\,A_{3}\,+\,A_{2}\,+\,2\,A_{1}$}{$[1]$}
\elldata{667}{$13$}{$D_{6}\,+\,A_{3}\,+\,4\,A_{1}$}{$[2], \,[1]$}
\elldata{668}{$13$}{$D_{6}\,+\,3\,A_{2}\,+\,A_{1}$}{$[1]$}
\elldata{669}{$13$}{$D_{6}\,+\,2\,A_{2}\,+\,3\,A_{1}$}{$[1]$}
\elldata{670}{$13$}{$D_{6}\,+\,A_{2}\,+\,5\,A_{1}$}{$[2], \,[1]$}
\elldata{671}{$13$}{$D_{6}\,+\,7\,A_{1}$}{$[2]$}
\elldata{672}{$13$}{$2\,D_{5}\,+\,A_{3}$}{$[1]$}
\elldata{673}{$13$}{$2\,D_{5}\,+\,A_{2}\,+\,A_{1}$}{$[1]$}
\elldata{674}{$13$}{$2\,D_{5}\,+\,3\,A_{1}$}{$[1]$}
\elldata{675}{$13$}{$D_{5}\,+\,2\,D_{4}$}{$[1]$}
\elldata{676}{$13$}{$D_{5}\,+\,D_{4}\,+\,A_{4}$}{$[1]$}
\elldata{677}{$13$}{$D_{5}\,+\,D_{4}\,+\,A_{3}\,+\,A_{1}$}{$[1]$}
\elldata{678}{$13$}{$D_{5}\,+\,D_{4}\,+\,2\,A_{2}$}{$[1]$}
\elldata{679}{$13$}{$D_{5}\,+\,D_{4}\,+\,A_{2}\,+\,2\,A_{1}$}{$[1]$}
\elldata{680}{$13$}{$D_{5}\,+\,D_{4}\,+\,4\,A_{1}$}{$[2], \,[1]$}
\elldata{681}{$13$}{$D_{5}\,+\,A_{8}$}{$[1]$}
\elldata{682}{$13$}{$D_{5}\,+\,A_{7}\,+\,A_{1}$}{$[1]$}
\elldata{683}{$13$}{$D_{5}\,+\,A_{6}\,+\,A_{2}$}{$[1]$}
\elldata{684}{$13$}{$D_{5}\,+\,A_{6}\,+\,2\,A_{1}$}{$[1]$}
\elldata{685}{$13$}{$D_{5}\,+\,A_{5}\,+\,A_{3}$}{$[1]$}
\elldata{686}{$13$}{$D_{5}\,+\,A_{5}\,+\,A_{2}\,+\,A_{1}$}{$[1]$}
\elldata{687}{$13$}{$D_{5}\,+\,A_{5}\,+\,3\,A_{1}$}{$[2], \,[1]$}
\elldata{688}{$13$}{$D_{5}\,+\,2\,A_{4}$}{$[1]$}
\elldata{689}{$13$}{$D_{5}\,+\,A_{4}\,+\,A_{3}\,+\,A_{1}$}{$[1]$}
\elldata{690}{$13$}{$D_{5}\,+\,A_{4}\,+\,2\,A_{2}$}{$[1]$}
\elldata{691}{$13$}{$D_{5}\,+\,A_{4}\,+\,A_{2}\,+\,2\,A_{1}$}{$[1]$}
\elldata{692}{$13$}{$D_{5}\,+\,A_{4}\,+\,4\,A_{1}$}{$[1]$}
\elldata{693}{$13$}{$D_{5}\,+\,2\,A_{3}\,+\,A_{2}$}{$[1]$}
\elldata{694}{$13$}{$D_{5}\,+\,2\,A_{3}\,+\,2\,A_{1}$}{$[2], \,[1]$}
\elldata{695}{$13$}{$D_{5}\,+\,A_{3}\,+\,2\,A_{2}\,+\,A_{1}$}{$[1]$}
\elldata{696}{$13$}{$D_{5}\,+\,A_{3}\,+\,A_{2}\,+\,3\,A_{1}$}{$[1]$}
\elldata{697}{$13$}{$D_{5}\,+\,A_{3}\,+\,5\,A_{1}$}{$[2], \,[1]$}
\elldata{698}{$13$}{$D_{5}\,+\,4\,A_{2}$}{$[1]$}
\elldata{699}{$13$}{$D_{5}\,+\,3\,A_{2}\,+\,2\,A_{1}$}{$[1]$}
\elldata{700}{$13$}{$D_{5}\,+\,2\,A_{2}\,+\,4\,A_{1}$}{$[1]$}
\elldata{701}{$13$}{$D_{5}\,+\,A_{2}\,+\,6\,A_{1}$}{$[2], \,[1]$}
\elldata{702}{$13$}{$D_{5}\,+\,8\,A_{1}$}{$[2]$}
\elldata{703}{$13$}{$3\,D_{4}\,+\,A_{1}$}{$[1]$}
\elldata{704}{$13$}{$2\,D_{4}\,+\,A_{5}$}{$[1]$}
\elldata{705}{$13$}{$2\,D_{4}\,+\,A_{4}\,+\,A_{1}$}{$[1]$}
\elldata{706}{$13$}{$2\,D_{4}\,+\,A_{3}\,+\,A_{2}$}{$[1]$}
\elldata{707}{$13$}{$2\,D_{4}\,+\,A_{3}\,+\,2\,A_{1}$}{$[2], \,[1]$}
\elldata{708}{$13$}{$2\,D_{4}\,+\,2\,A_{2}\,+\,A_{1}$}{$[1]$}
\elldata{709}{$13$}{$2\,D_{4}\,+\,A_{2}\,+\,3\,A_{1}$}{$[1]$}
\elldata{710}{$13$}{$2\,D_{4}\,+\,5\,A_{1}$}{$[2]$}
\elldata{711}{$13$}{$D_{4}\,+\,A_{9}$}{$[1]$}
\elldata{712}{$13$}{$D_{4}\,+\,A_{8}\,+\,A_{1}$}{$[1]$}
\elldata{713}{$13$}{$D_{4}\,+\,A_{7}\,+\,A_{2}$}{$[1]$}
\elldata{714}{$13$}{$D_{4}\,+\,A_{7}\,+\,2\,A_{1}$}{$[2], \,[1]$}
\elldata{715}{$13$}{$D_{4}\,+\,A_{6}\,+\,A_{3}$}{$[1]$}
\elldata{716}{$13$}{$D_{4}\,+\,A_{6}\,+\,A_{2}\,+\,A_{1}$}{$[1]$}
\elldata{717}{$13$}{$D_{4}\,+\,A_{6}\,+\,3\,A_{1}$}{$[1]$}
\elldata{718}{$13$}{$D_{4}\,+\,A_{5}\,+\,A_{4}$}{$[1]$}
\elldata{719}{$13$}{$D_{4}\,+\,A_{5}\,+\,A_{3}\,+\,A_{1}$}{$[2], \,[1]$}
\elldata{720}{$13$}{$D_{4}\,+\,A_{5}\,+\,2\,A_{2}$}{$[1]$}
\elldata{721}{$13$}{$D_{4}\,+\,A_{5}\,+\,A_{2}\,+\,2\,A_{1}$}{$[1]$}
\elldata{722}{$13$}{$D_{4}\,+\,A_{5}\,+\,4\,A_{1}$}{$[2], \,[1]$}
\elldata{723}{$13$}{$D_{4}\,+\,2\,A_{4}\,+\,A_{1}$}{$[1]$}
\elldata{724}{$13$}{$D_{4}\,+\,A_{4}\,+\,A_{3}\,+\,A_{2}$}{$[1]$}
\elldata{725}{$13$}{$D_{4}\,+\,A_{4}\,+\,A_{3}\,+\,2\,A_{1}$}{$[1]$}
\elldata{726}{$13$}{$D_{4}\,+\,A_{4}\,+\,2\,A_{2}\,+\,A_{1}$}{$[1]$}
\elldata{727}{$13$}{$D_{4}\,+\,A_{4}\,+\,A_{2}\,+\,3\,A_{1}$}{$[1]$}
\elldata{728}{$13$}{$D_{4}\,+\,A_{4}\,+\,5\,A_{1}$}{$[1]$}
\elldata{729}{$13$}{$D_{4}\,+\,3\,A_{3}$}{$[2], \,[1]$}
\elldata{730}{$13$}{$D_{4}\,+\,2\,A_{3}\,+\,A_{2}\,+\,A_{1}$}{$[1]$}
\elldata{731}{$13$}{$D_{4}\,+\,2\,A_{3}\,+\,3\,A_{1}$}{$[2], \,[1]$}
\elldata{732}{$13$}{$D_{4}\,+\,A_{3}\,+\,3\,A_{2}$}{$[1]$}
\elldata{733}{$13$}{$D_{4}\,+\,A_{3}\,+\,2\,A_{2}\,+\,2\,A_{1}$}{$[1]$}
\elldata{734}{$13$}{$D_{4}\,+\,A_{3}\,+\,A_{2}\,+\,4\,A_{1}$}{$[2], \,[1]$}
\elldata{735}{$13$}{$D_{4}\,+\,A_{3}\,+\,6\,A_{1}$}{$[2]$}
\elldata{736}{$13$}{$D_{4}\,+\,4\,A_{2}\,+\,A_{1}$}{$[1]$}
\elldata{737}{$13$}{$D_{4}\,+\,3\,A_{2}\,+\,3\,A_{1}$}{$[1]$}
\elldata{738}{$13$}{$D_{4}\,+\,2\,A_{2}\,+\,5\,A_{1}$}{$[1]$}
\elldata{739}{$13$}{$D_{4}\,+\,A_{2}\,+\,7\,A_{1}$}{$[2]$}
\elldata{740}{$13$}{$D_{4}\,+\,9\,A_{1}$}{$[2, 2]$}
\elldata{741}{$13$}{$A_{13}$}{$[1]$}
\elldata{742}{$13$}{$A_{12}\,+\,A_{1}$}{$[1]$}
\elldata{743}{$13$}{$A_{11}\,+\,A_{2}$}{$[1]$}
\elldata{744}{$13$}{$A_{11}\,+\,2\,A_{1}$}{$[2], \,[1]$}
\elldata{745}{$13$}{$A_{10}\,+\,A_{3}$}{$[1]$}
\elldata{746}{$13$}{$A_{10}\,+\,A_{2}\,+\,A_{1}$}{$[1]$}
\elldata{747}{$13$}{$A_{10}\,+\,3\,A_{1}$}{$[1]$}
\elldata{748}{$13$}{$A_{9}\,+\,A_{4}$}{$[1]$}
\elldata{749}{$13$}{$A_{9}\,+\,A_{3}\,+\,A_{1}$}{$[2], \,[1]$}
\elldata{750}{$13$}{$A_{9}\,+\,2\,A_{2}$}{$[1]$}
\elldata{751}{$13$}{$A_{9}\,+\,A_{2}\,+\,2\,A_{1}$}{$[1]$}
\elldata{752}{$13$}{$A_{9}\,+\,4\,A_{1}$}{$[2], \,[1]$}
\elldata{753}{$13$}{$A_{8}\,+\,A_{5}$}{$[1]$}
\elldata{754}{$13$}{$A_{8}\,+\,A_{4}\,+\,A_{1}$}{$[1]$}
\elldata{755}{$13$}{$A_{8}\,+\,A_{3}\,+\,A_{2}$}{$[1]$}
\elldata{756}{$13$}{$A_{8}\,+\,A_{3}\,+\,2\,A_{1}$}{$[1]$}
\elldata{757}{$13$}{$A_{8}\,+\,2\,A_{2}\,+\,A_{1}$}{$[1]$}
\elldata{758}{$13$}{$A_{8}\,+\,A_{2}\,+\,3\,A_{1}$}{$[1]$}
\elldata{759}{$13$}{$A_{8}\,+\,5\,A_{1}$}{$[1]$}
\elldata{760}{$13$}{$A_{7}\,+\,A_{6}$}{$[1]$}
\elldata{761}{$13$}{$A_{7}\,+\,A_{5}\,+\,A_{1}$}{$[2], \,[1]$}
\elldata{762}{$13$}{$A_{7}\,+\,A_{4}\,+\,A_{2}$}{$[1]$}
\elldata{763}{$13$}{$A_{7}\,+\,A_{4}\,+\,2\,A_{1}$}{$[1]$}
\elldata{764}{$13$}{$A_{7}\,+\,2\,A_{3}$}{$[2], \,[1]$}
\elldata{765}{$13$}{$A_{7}\,+\,A_{3}\,+\,A_{2}\,+\,A_{1}$}{$[1]$}
\elldata{766}{$13$}{$A_{7}\,+\,A_{3}\,+\,3\,A_{1}$}{$[2], \,[1]$}
\elldata{767}{$13$}{$A_{7}\,+\,3\,A_{2}$}{$[1]$}
\elldata{768}{$13$}{$A_{7}\,+\,2\,A_{2}\,+\,2\,A_{1}$}{$[1]$}
\elldata{769}{$13$}{$A_{7}\,+\,A_{2}\,+\,4\,A_{1}$}{$[2], \,[1]$}
\elldata{770}{$13$}{$A_{7}\,+\,6\,A_{1}$}{$[2], \,[1]$}
\elldata{771}{$13$}{$2\,A_{6}\,+\,A_{1}$}{$[1]$}
\elldata{772}{$13$}{$A_{6}\,+\,A_{5}\,+\,A_{2}$}{$[1]$}
\elldata{773}{$13$}{$A_{6}\,+\,A_{5}\,+\,2\,A_{1}$}{$[1]$}
\elldata{774}{$13$}{$A_{6}\,+\,A_{4}\,+\,A_{3}$}{$[1]$}
\elldata{775}{$13$}{$A_{6}\,+\,A_{4}\,+\,A_{2}\,+\,A_{1}$}{$[1]$}
\elldata{776}{$13$}{$A_{6}\,+\,A_{4}\,+\,3\,A_{1}$}{$[1]$}
\elldata{777}{$13$}{$A_{6}\,+\,2\,A_{3}\,+\,A_{1}$}{$[1]$}
\elldata{778}{$13$}{$A_{6}\,+\,A_{3}\,+\,2\,A_{2}$}{$[1]$}
\elldata{779}{$13$}{$A_{6}\,+\,A_{3}\,+\,A_{2}\,+\,2\,A_{1}$}{$[1]$}
\elldata{780}{$13$}{$A_{6}\,+\,A_{3}\,+\,4\,A_{1}$}{$[1]$}
\elldata{781}{$13$}{$A_{6}\,+\,3\,A_{2}\,+\,A_{1}$}{$[1]$}
\elldata{782}{$13$}{$A_{6}\,+\,2\,A_{2}\,+\,3\,A_{1}$}{$[1]$}
\elldata{783}{$13$}{$A_{6}\,+\,A_{2}\,+\,5\,A_{1}$}{$[1]$}
\elldata{784}{$13$}{$A_{6}\,+\,7\,A_{1}$}{$[1]$}
\elldata{785}{$13$}{$2\,A_{5}\,+\,A_{3}$}{$[2], \,[1]$}
\elldata{786}{$13$}{$2\,A_{5}\,+\,A_{2}\,+\,A_{1}$}{$[1]$}
\elldata{787}{$13$}{$2\,A_{5}\,+\,3\,A_{1}$}{$[2], \,[1]$}
\elldata{788}{$13$}{$A_{5}\,+\,2\,A_{4}$}{$[1]$}
\elldata{789}{$13$}{$A_{5}\,+\,A_{4}\,+\,A_{3}\,+\,A_{1}$}{$[1]$}
\elldata{790}{$13$}{$A_{5}\,+\,A_{4}\,+\,2\,A_{2}$}{$[1]$}
\elldata{791}{$13$}{$A_{5}\,+\,A_{4}\,+\,A_{2}\,+\,2\,A_{1}$}{$[1]$}
\elldata{792}{$13$}{$A_{5}\,+\,A_{4}\,+\,4\,A_{1}$}{$[1]$}
\elldata{793}{$13$}{$A_{5}\,+\,2\,A_{3}\,+\,A_{2}$}{$[1]$}
\elldata{794}{$13$}{$A_{5}\,+\,2\,A_{3}\,+\,2\,A_{1}$}{$[2], \,[1]$}
\elldata{795}{$13$}{$A_{5}\,+\,A_{3}\,+\,2\,A_{2}\,+\,A_{1}$}{$[1]$}
\elldata{796}{$13$}{$A_{5}\,+\,A_{3}\,+\,A_{2}\,+\,3\,A_{1}$}{$[2], \,[1]$}
\elldata{797}{$13$}{$A_{5}\,+\,A_{3}\,+\,5\,A_{1}$}{$[2], \,[1]$}
\elldata{798}{$13$}{$A_{5}\,+\,4\,A_{2}$}{$[3], \,[1]$}
\elldata{799}{$13$}{$A_{5}\,+\,3\,A_{2}\,+\,2\,A_{1}$}{$[1]$}
\elldata{800}{$13$}{$A_{5}\,+\,2\,A_{2}\,+\,4\,A_{1}$}{$[1]$}
\elldata{801}{$13$}{$A_{5}\,+\,A_{2}\,+\,6\,A_{1}$}{$[2], \,[1]$}
\elldata{802}{$13$}{$A_{5}\,+\,8\,A_{1}$}{$[2]$}
\elldata{803}{$13$}{$3\,A_{4}\,+\,A_{1}$}{$[1]$}
\elldata{804}{$13$}{$2\,A_{4}\,+\,A_{3}\,+\,A_{2}$}{$[1]$}
\elldata{805}{$13$}{$2\,A_{4}\,+\,A_{3}\,+\,2\,A_{1}$}{$[1]$}
\elldata{806}{$13$}{$2\,A_{4}\,+\,2\,A_{2}\,+\,A_{1}$}{$[1]$}
\elldata{807}{$13$}{$2\,A_{4}\,+\,A_{2}\,+\,3\,A_{1}$}{$[1]$}
\elldata{808}{$13$}{$2\,A_{4}\,+\,5\,A_{1}$}{$[1]$}
\elldata{809}{$13$}{$A_{4}\,+\,3\,A_{3}$}{$[1]$}
\elldata{810}{$13$}{$A_{4}\,+\,2\,A_{3}\,+\,A_{2}\,+\,A_{1}$}{$[1]$}
\elldata{811}{$13$}{$A_{4}\,+\,2\,A_{3}\,+\,3\,A_{1}$}{$[1]$}
\elldata{812}{$13$}{$A_{4}\,+\,A_{3}\,+\,3\,A_{2}$}{$[1]$}
\elldata{813}{$13$}{$A_{4}\,+\,A_{3}\,+\,2\,A_{2}\,+\,2\,A_{1}$}{$[1]$}
\elldata{814}{$13$}{$A_{4}\,+\,A_{3}\,+\,A_{2}\,+\,4\,A_{1}$}{$[1]$}
\elldata{815}{$13$}{$A_{4}\,+\,A_{3}\,+\,6\,A_{1}$}{$[2], \,[1]$}
\elldata{816}{$13$}{$A_{4}\,+\,4\,A_{2}\,+\,A_{1}$}{$[1]$}
\elldata{817}{$13$}{$A_{4}\,+\,3\,A_{2}\,+\,3\,A_{1}$}{$[1]$}
\elldata{818}{$13$}{$A_{4}\,+\,2\,A_{2}\,+\,5\,A_{1}$}{$[1]$}
\elldata{819}{$13$}{$A_{4}\,+\,A_{2}\,+\,7\,A_{1}$}{$[1]$}
\elldata{820}{$13$}{$A_{4}\,+\,9\,A_{1}$}{$[2]$}
\elldata{821}{$13$}{$4\,A_{3}\,+\,A_{1}$}{$[2], \,[1]$}
\elldata{822}{$13$}{$3\,A_{3}\,+\,2\,A_{2}$}{$[1]$}
\elldata{823}{$13$}{$3\,A_{3}\,+\,A_{2}\,+\,2\,A_{1}$}{$[2], \,[1]$}
\elldata{824}{$13$}{$3\,A_{3}\,+\,4\,A_{1}$}{$[2], \,[1]$}
\elldata{825}{$13$}{$2\,A_{3}\,+\,3\,A_{2}\,+\,A_{1}$}{$[1]$}
\elldata{826}{$13$}{$2\,A_{3}\,+\,2\,A_{2}\,+\,3\,A_{1}$}{$[1]$}
\elldata{827}{$13$}{$2\,A_{3}\,+\,A_{2}\,+\,5\,A_{1}$}{$[2], \,[1]$}
\elldata{828}{$13$}{$2\,A_{3}\,+\,7\,A_{1}$}{$[2]$}
\elldata{829}{$13$}{$A_{3}\,+\,5\,A_{2}$}{$[1]$}
\elldata{830}{$13$}{$A_{3}\,+\,4\,A_{2}\,+\,2\,A_{1}$}{$[1]$}
\elldata{831}{$13$}{$A_{3}\,+\,3\,A_{2}\,+\,4\,A_{1}$}{$[1]$}
\elldata{832}{$13$}{$A_{3}\,+\,2\,A_{2}\,+\,6\,A_{1}$}{$[2], \,[1]$}
\elldata{833}{$13$}{$A_{3}\,+\,A_{2}\,+\,8\,A_{1}$}{$[2]$}
\elldata{834}{$13$}{$A_{3}\,+\,10\,A_{1}$}{$[2, 2]$}
\elldata{835}{$13$}{$6\,A_{2}\,+\,A_{1}$}{$[3], \,[1]$}
\elldata{836}{$13$}{$5\,A_{2}\,+\,3\,A_{1}$}{$[1]$}
\elldata{837}{$13$}{$4\,A_{2}\,+\,5\,A_{1}$}{$[1]$}
\elldata{838}{$13$}{$3\,A_{2}\,+\,7\,A_{1}$}{$[1]$}
\elldata{839}{$13$}{$2\,A_{2}\,+\,9\,A_{1}$}{$[2]$}

\vsr \elldata{No.}{rank}{$ADE$-type}{$G$}

\vsrs \elldata{840}{$14$}{$E_{8}\,+\,E_{6}$}{$[1]$}
\elldata{841}{$14$}{$E_{8}\,+\,D_{6}$}{$[1]$}
\elldata{842}{$14$}{$E_{8}\,+\,D_{5}\,+\,A_{1}$}{$[1]$}
\elldata{843}{$14$}{$E_{8}\,+\,D_{4}\,+\,A_{2}$}{$[1]$}
\elldata{844}{$14$}{$E_{8}\,+\,D_{4}\,+\,2\,A_{1}$}{$[1]$}
\elldata{845}{$14$}{$E_{8}\,+\,A_{6}$}{$[1]$}
\elldata{846}{$14$}{$E_{8}\,+\,A_{5}\,+\,A_{1}$}{$[1]$}
\elldata{847}{$14$}{$E_{8}\,+\,A_{4}\,+\,A_{2}$}{$[1]$}
\elldata{848}{$14$}{$E_{8}\,+\,A_{4}\,+\,2\,A_{1}$}{$[1]$}
\elldata{849}{$14$}{$E_{8}\,+\,2\,A_{3}$}{$[1]$}
\elldata{850}{$14$}{$E_{8}\,+\,A_{3}\,+\,A_{2}\,+\,A_{1}$}{$[1]$}
\elldata{851}{$14$}{$E_{8}\,+\,A_{3}\,+\,3\,A_{1}$}{$[1]$}
\elldata{852}{$14$}{$E_{8}\,+\,3\,A_{2}$}{$[1]$}
\elldata{853}{$14$}{$E_{8}\,+\,2\,A_{2}\,+\,2\,A_{1}$}{$[1]$}
\elldata{854}{$14$}{$E_{8}\,+\,A_{2}\,+\,4\,A_{1}$}{$[1]$}
\elldata{855}{$14$}{$E_{8}\,+\,6\,A_{1}$}{$[1]$}
\elldata{856}{$14$}{$2\,E_{7}$}{$[1]$}
\elldata{857}{$14$}{$E_{7}\,+\,E_{6}\,+\,A_{1}$}{$[1]$}
\elldata{858}{$14$}{$E_{7}\,+\,D_{7}$}{$[1]$}
\elldata{859}{$14$}{$E_{7}\,+\,D_{6}\,+\,A_{1}$}{$[1]$}
\elldata{860}{$14$}{$E_{7}\,+\,D_{5}\,+\,A_{2}$}{$[1]$}
\elldata{861}{$14$}{$E_{7}\,+\,D_{5}\,+\,2\,A_{1}$}{$[1]$}
\elldata{862}{$14$}{$E_{7}\,+\,D_{4}\,+\,A_{3}$}{$[1]$}
\elldata{863}{$14$}{$E_{7}\,+\,D_{4}\,+\,A_{2}\,+\,A_{1}$}{$[1]$}
\elldata{864}{$14$}{$E_{7}\,+\,D_{4}\,+\,3\,A_{1}$}{$[2], \,[1]$}
\elldata{865}{$14$}{$E_{7}\,+\,A_{7}$}{$[1]$}
\elldata{866}{$14$}{$E_{7}\,+\,A_{6}\,+\,A_{1}$}{$[1]$}
\elldata{867}{$14$}{$E_{7}\,+\,A_{5}\,+\,A_{2}$}{$[1]$}
\elldata{868}{$14$}{$E_{7}\,+\,A_{5}\,+\,2\,A_{1}$}{$[2], \,[1]$}
\elldata{869}{$14$}{$E_{7}\,+\,A_{4}\,+\,A_{3}$}{$[1]$}
\elldata{870}{$14$}{$E_{7}\,+\,A_{4}\,+\,A_{2}\,+\,A_{1}$}{$[1]$}
\elldata{871}{$14$}{$E_{7}\,+\,A_{4}\,+\,3\,A_{1}$}{$[1]$}
\elldata{872}{$14$}{$E_{7}\,+\,2\,A_{3}\,+\,A_{1}$}{$[2], \,[1]$}
\elldata{873}{$14$}{$E_{7}\,+\,A_{3}\,+\,2\,A_{2}$}{$[1]$}
\elldata{874}{$14$}{$E_{7}\,+\,A_{3}\,+\,A_{2}\,+\,2\,A_{1}$}{$[1]$}
\elldata{875}{$14$}{$E_{7}\,+\,A_{3}\,+\,4\,A_{1}$}{$[2], \,[1]$}
\elldata{876}{$14$}{$E_{7}\,+\,3\,A_{2}\,+\,A_{1}$}{$[1]$}
\elldata{877}{$14$}{$E_{7}\,+\,2\,A_{2}\,+\,3\,A_{1}$}{$[1]$}
\elldata{878}{$14$}{$E_{7}\,+\,A_{2}\,+\,5\,A_{1}$}{$[2], \,[1]$}
\elldata{879}{$14$}{$E_{7}\,+\,7\,A_{1}$}{$[2]$}
\elldata{880}{$14$}{$2\,E_{6}\,+\,A_{2}$}{$[1]$}
\elldata{881}{$14$}{$2\,E_{6}\,+\,2\,A_{1}$}{$[1]$}
\elldata{882}{$14$}{$E_{6}\,+\,D_{8}$}{$[1]$}
\elldata{883}{$14$}{$E_{6}\,+\,D_{7}\,+\,A_{1}$}{$[1]$}
\elldata{884}{$14$}{$E_{6}\,+\,D_{6}\,+\,A_{2}$}{$[1]$}
\elldata{885}{$14$}{$E_{6}\,+\,D_{6}\,+\,2\,A_{1}$}{$[1]$}
\elldata{886}{$14$}{$E_{6}\,+\,D_{5}\,+\,A_{3}$}{$[1]$}
\elldata{887}{$14$}{$E_{6}\,+\,D_{5}\,+\,A_{2}\,+\,A_{1}$}{$[1]$}
\elldata{888}{$14$}{$E_{6}\,+\,D_{5}\,+\,3\,A_{1}$}{$[1]$}
\elldata{889}{$14$}{$E_{6}\,+\,2\,D_{4}$}{$[1]$}
\elldata{890}{$14$}{$E_{6}\,+\,D_{4}\,+\,A_{4}$}{$[1]$}
\elldata{891}{$14$}{$E_{6}\,+\,D_{4}\,+\,A_{3}\,+\,A_{1}$}{$[1]$}
\elldata{892}{$14$}{$E_{6}\,+\,D_{4}\,+\,2\,A_{2}$}{$[1]$}
\elldata{893}{$14$}{$E_{6}\,+\,D_{4}\,+\,A_{2}\,+\,2\,A_{1}$}{$[1]$}
\elldata{894}{$14$}{$E_{6}\,+\,D_{4}\,+\,4\,A_{1}$}{$[1]$}
\elldata{895}{$14$}{$E_{6}\,+\,A_{8}$}{$[1]$}
\elldata{896}{$14$}{$E_{6}\,+\,A_{7}\,+\,A_{1}$}{$[1]$}
\elldata{897}{$14$}{$E_{6}\,+\,A_{6}\,+\,A_{2}$}{$[1]$}
\elldata{898}{$14$}{$E_{6}\,+\,A_{6}\,+\,2\,A_{1}$}{$[1]$}
\elldata{899}{$14$}{$E_{6}\,+\,A_{5}\,+\,A_{3}$}{$[1]$}
\elldata{900}{$14$}{$E_{6}\,+\,A_{5}\,+\,A_{2}\,+\,A_{1}$}{$[1]$}
\elldata{901}{$14$}{$E_{6}\,+\,A_{5}\,+\,3\,A_{1}$}{$[1]$}
\elldata{902}{$14$}{$E_{6}\,+\,2\,A_{4}$}{$[1]$}
\elldata{903}{$14$}{$E_{6}\,+\,A_{4}\,+\,A_{3}\,+\,A_{1}$}{$[1]$}
\elldata{904}{$14$}{$E_{6}\,+\,A_{4}\,+\,2\,A_{2}$}{$[1]$}
\elldata{905}{$14$}{$E_{6}\,+\,A_{4}\,+\,A_{2}\,+\,2\,A_{1}$}{$[1]$}
\elldata{906}{$14$}{$E_{6}\,+\,A_{4}\,+\,4\,A_{1}$}{$[1]$}
\elldata{907}{$14$}{$E_{6}\,+\,2\,A_{3}\,+\,A_{2}$}{$[1]$}
\elldata{908}{$14$}{$E_{6}\,+\,2\,A_{3}\,+\,2\,A_{1}$}{$[1]$}
\elldata{909}{$14$}{$E_{6}\,+\,A_{3}\,+\,2\,A_{2}\,+\,A_{1}$}{$[1]$}
\elldata{910}{$14$}{$E_{6}\,+\,A_{3}\,+\,A_{2}\,+\,3\,A_{1}$}{$[1]$}
\elldata{911}{$14$}{$E_{6}\,+\,A_{3}\,+\,5\,A_{1}$}{$[1]$}
\elldata{912}{$14$}{$E_{6}\,+\,4\,A_{2}$}{$[3], \,[1]$}
\elldata{913}{$14$}{$E_{6}\,+\,3\,A_{2}\,+\,2\,A_{1}$}{$[1]$}
\elldata{914}{$14$}{$E_{6}\,+\,2\,A_{2}\,+\,4\,A_{1}$}{$[1]$}
\elldata{915}{$14$}{$E_{6}\,+\,A_{2}\,+\,6\,A_{1}$}{$[1]$}
\elldata{916}{$14$}{$D_{14}$}{$[1]$}
\elldata{917}{$14$}{$D_{13}\,+\,A_{1}$}{$[1]$}
\elldata{918}{$14$}{$D_{12}\,+\,A_{2}$}{$[1]$}
\elldata{919}{$14$}{$D_{12}\,+\,2\,A_{1}$}{$[2], \,[1]$}
\elldata{920}{$14$}{$D_{11}\,+\,A_{3}$}{$[1]$}
\elldata{921}{$14$}{$D_{11}\,+\,A_{2}\,+\,A_{1}$}{$[1]$}
\elldata{922}{$14$}{$D_{11}\,+\,3\,A_{1}$}{$[1]$}
\elldata{923}{$14$}{$D_{10}\,+\,D_{4}$}{$[1]$}
\elldata{924}{$14$}{$D_{10}\,+\,A_{4}$}{$[1]$}
\elldata{925}{$14$}{$D_{10}\,+\,A_{3}\,+\,A_{1}$}{$[2], \,[1]$}
\elldata{926}{$14$}{$D_{10}\,+\,2\,A_{2}$}{$[1]$}
\elldata{927}{$14$}{$D_{10}\,+\,A_{2}\,+\,2\,A_{1}$}{$[1]$}
\elldata{928}{$14$}{$D_{10}\,+\,4\,A_{1}$}{$[2], \,[1]$}
\elldata{929}{$14$}{$D_{9}\,+\,D_{5}$}{$[1]$}
\elldata{930}{$14$}{$D_{9}\,+\,D_{4}\,+\,A_{1}$}{$[1]$}
\elldata{931}{$14$}{$D_{9}\,+\,A_{5}$}{$[1]$}
\elldata{932}{$14$}{$D_{9}\,+\,A_{4}\,+\,A_{1}$}{$[1]$}
\elldata{933}{$14$}{$D_{9}\,+\,A_{3}\,+\,A_{2}$}{$[1]$}
\elldata{934}{$14$}{$D_{9}\,+\,A_{3}\,+\,2\,A_{1}$}{$[1]$}
\elldata{935}{$14$}{$D_{9}\,+\,2\,A_{2}\,+\,A_{1}$}{$[1]$}
\elldata{936}{$14$}{$D_{9}\,+\,A_{2}\,+\,3\,A_{1}$}{$[1]$}
\elldata{937}{$14$}{$D_{9}\,+\,5\,A_{1}$}{$[1]$}
\elldata{938}{$14$}{$D_{8}\,+\,D_{6}$}{$[1]$}
\elldata{939}{$14$}{$D_{8}\,+\,D_{5}\,+\,A_{1}$}{$[1]$}
\elldata{940}{$14$}{$D_{8}\,+\,D_{4}\,+\,A_{2}$}{$[1]$}
\elldata{941}{$14$}{$D_{8}\,+\,D_{4}\,+\,2\,A_{1}$}{$[2], \,[1]$}
\elldata{942}{$14$}{$D_{8}\,+\,A_{6}$}{$[1]$}
\elldata{943}{$14$}{$D_{8}\,+\,A_{5}\,+\,A_{1}$}{$[2], \,[1]$}
\elldata{944}{$14$}{$D_{8}\,+\,A_{4}\,+\,A_{2}$}{$[1]$}
\elldata{945}{$14$}{$D_{8}\,+\,A_{4}\,+\,2\,A_{1}$}{$[1]$}
\elldata{946}{$14$}{$D_{8}\,+\,2\,A_{3}$}{$[2], \,[1]$}
\elldata{947}{$14$}{$D_{8}\,+\,A_{3}\,+\,A_{2}\,+\,A_{1}$}{$[1]$}
\elldata{948}{$14$}{$D_{8}\,+\,A_{3}\,+\,3\,A_{1}$}{$[2], \,[1]$}
\elldata{949}{$14$}{$D_{8}\,+\,3\,A_{2}$}{$[1]$}
\elldata{950}{$14$}{$D_{8}\,+\,2\,A_{2}\,+\,2\,A_{1}$}{$[1]$}
\elldata{951}{$14$}{$D_{8}\,+\,A_{2}\,+\,4\,A_{1}$}{$[2], \,[1]$}
\elldata{952}{$14$}{$D_{8}\,+\,6\,A_{1}$}{$[2]$}
\elldata{953}{$14$}{$2\,D_{7}$}{$[1]$}
\elldata{954}{$14$}{$D_{7}\,+\,D_{6}\,+\,A_{1}$}{$[1]$}
\elldata{955}{$14$}{$D_{7}\,+\,D_{5}\,+\,A_{2}$}{$[1]$}
\elldata{956}{$14$}{$D_{7}\,+\,D_{5}\,+\,2\,A_{1}$}{$[1]$}
\elldata{957}{$14$}{$D_{7}\,+\,D_{4}\,+\,A_{3}$}{$[1]$}
\elldata{958}{$14$}{$D_{7}\,+\,D_{4}\,+\,A_{2}\,+\,A_{1}$}{$[1]$}
\elldata{959}{$14$}{$D_{7}\,+\,D_{4}\,+\,3\,A_{1}$}{$[1]$}
\elldata{960}{$14$}{$D_{7}\,+\,A_{7}$}{$[1]$}
\elldata{961}{$14$}{$D_{7}\,+\,A_{6}\,+\,A_{1}$}{$[1]$}
\elldata{962}{$14$}{$D_{7}\,+\,A_{5}\,+\,A_{2}$}{$[1]$}
\elldata{963}{$14$}{$D_{7}\,+\,A_{5}\,+\,2\,A_{1}$}{$[1]$}
\elldata{964}{$14$}{$D_{7}\,+\,A_{4}\,+\,A_{3}$}{$[1]$}
\elldata{965}{$14$}{$D_{7}\,+\,A_{4}\,+\,A_{2}\,+\,A_{1}$}{$[1]$}
\elldata{966}{$14$}{$D_{7}\,+\,A_{4}\,+\,3\,A_{1}$}{$[1]$}
\elldata{967}{$14$}{$D_{7}\,+\,2\,A_{3}\,+\,A_{1}$}{$[1]$}
\elldata{968}{$14$}{$D_{7}\,+\,A_{3}\,+\,2\,A_{2}$}{$[1]$}
\elldata{969}{$14$}{$D_{7}\,+\,A_{3}\,+\,A_{2}\,+\,2\,A_{1}$}{$[1]$}
\elldata{970}{$14$}{$D_{7}\,+\,A_{3}\,+\,4\,A_{1}$}{$[2], \,[1]$}
\elldata{971}{$14$}{$D_{7}\,+\,3\,A_{2}\,+\,A_{1}$}{$[1]$}
\elldata{972}{$14$}{$D_{7}\,+\,2\,A_{2}\,+\,3\,A_{1}$}{$[1]$}
\elldata{973}{$14$}{$D_{7}\,+\,A_{2}\,+\,5\,A_{1}$}{$[1]$}
\elldata{974}{$14$}{$D_{7}\,+\,7\,A_{1}$}{$[2]$}
\elldata{975}{$14$}{$2\,D_{6}\,+\,A_{2}$}{$[1]$}
\elldata{976}{$14$}{$2\,D_{6}\,+\,2\,A_{1}$}{$[2], \,[1]$}
\elldata{977}{$14$}{$D_{6}\,+\,D_{5}\,+\,A_{3}$}{$[1]$}
\elldata{978}{$14$}{$D_{6}\,+\,D_{5}\,+\,A_{2}\,+\,A_{1}$}{$[1]$}
\elldata{979}{$14$}{$D_{6}\,+\,D_{5}\,+\,3\,A_{1}$}{$[2], \,[1]$}
\elldata{980}{$14$}{$D_{6}\,+\,2\,D_{4}$}{$[1]$}
\elldata{981}{$14$}{$D_{6}\,+\,D_{4}\,+\,A_{4}$}{$[1]$}
\elldata{982}{$14$}{$D_{6}\,+\,D_{4}\,+\,A_{3}\,+\,A_{1}$}{$[2], \,[1]$}
\elldata{983}{$14$}{$D_{6}\,+\,D_{4}\,+\,2\,A_{2}$}{$[1]$}
\elldata{984}{$14$}{$D_{6}\,+\,D_{4}\,+\,A_{2}\,+\,2\,A_{1}$}{$[1]$}
\elldata{985}{$14$}{$D_{6}\,+\,D_{4}\,+\,4\,A_{1}$}{$[2]$}
\elldata{986}{$14$}{$D_{6}\,+\,A_{8}$}{$[1]$}
\elldata{987}{$14$}{$D_{6}\,+\,A_{7}\,+\,A_{1}$}{$[2], \,[1]$}
\elldata{988}{$14$}{$D_{6}\,+\,A_{6}\,+\,A_{2}$}{$[1]$}
\elldata{989}{$14$}{$D_{6}\,+\,A_{6}\,+\,2\,A_{1}$}{$[1]$}
\elldata{990}{$14$}{$D_{6}\,+\,A_{5}\,+\,A_{3}$}{$[2], \,[1]$}
\elldata{991}{$14$}{$D_{6}\,+\,A_{5}\,+\,A_{2}\,+\,A_{1}$}{$[1]$}
\elldata{992}{$14$}{$D_{6}\,+\,A_{5}\,+\,3\,A_{1}$}{$[2], \,[1]$}
\elldata{993}{$14$}{$D_{6}\,+\,2\,A_{4}$}{$[1]$}
\elldata{994}{$14$}{$D_{6}\,+\,A_{4}\,+\,A_{3}\,+\,A_{1}$}{$[1]$}
\elldata{995}{$14$}{$D_{6}\,+\,A_{4}\,+\,2\,A_{2}$}{$[1]$}
\elldata{996}{$14$}{$D_{6}\,+\,A_{4}\,+\,A_{2}\,+\,2\,A_{1}$}{$[1]$}
\elldata{997}{$14$}{$D_{6}\,+\,A_{4}\,+\,4\,A_{1}$}{$[1]$}
\elldata{998}{$14$}{$D_{6}\,+\,2\,A_{3}\,+\,A_{2}$}{$[1]$}
\elldata{999}{$14$}{$D_{6}\,+\,2\,A_{3}\,+\,2\,A_{1}$}{$[2], \,[1]$}
\elldata{1000}{$14$}{$D_{6}\,+\,A_{3}\,+\,2\,A_{2}\,+\,A_{1}$}{$[1]$}
\elldata{1001}{$14$}{$D_{6}\,+\,A_{3}\,+\,A_{2}\,+\,3\,A_{1}$}{$[2], \,[1]$}
\elldata{1002}{$14$}{$D_{6}\,+\,A_{3}\,+\,5\,A_{1}$}{$[2]$}
\elldata{1003}{$14$}{$D_{6}\,+\,4\,A_{2}$}{$[1]$}
\elldata{1004}{$14$}{$D_{6}\,+\,3\,A_{2}\,+\,2\,A_{1}$}{$[1]$}
\elldata{1005}{$14$}{$D_{6}\,+\,2\,A_{2}\,+\,4\,A_{1}$}{$[1]$}
\elldata{1006}{$14$}{$D_{6}\,+\,A_{2}\,+\,6\,A_{1}$}{$[2]$}
\elldata{1007}{$14$}{$D_{6}\,+\,8\,A_{1}$}{$[2, 2]$}
\elldata{1008}{$14$}{$2\,D_{5}\,+\,D_{4}$}{$[1]$}
\elldata{1009}{$14$}{$2\,D_{5}\,+\,A_{4}$}{$[1]$}
\elldata{1010}{$14$}{$2\,D_{5}\,+\,A_{3}\,+\,A_{1}$}{$[1]$}
\elldata{1011}{$14$}{$2\,D_{5}\,+\,2\,A_{2}$}{$[1]$}
\elldata{1012}{$14$}{$2\,D_{5}\,+\,A_{2}\,+\,2\,A_{1}$}{$[1]$}
\elldata{1013}{$14$}{$2\,D_{5}\,+\,4\,A_{1}$}{$[2], \,[1]$}
\elldata{1014}{$14$}{$D_{5}\,+\,2\,D_{4}\,+\,A_{1}$}{$[1]$}
\elldata{1015}{$14$}{$D_{5}\,+\,D_{4}\,+\,A_{5}$}{$[1]$}
\elldata{1016}{$14$}{$D_{5}\,+\,D_{4}\,+\,A_{4}\,+\,A_{1}$}{$[1]$}
\elldata{1017}{$14$}{$D_{5}\,+\,D_{4}\,+\,A_{3}\,+\,A_{2}$}{$[1]$}
\elldata{1018}{$14$}{$D_{5}\,+\,D_{4}\,+\,A_{3}\,+\,2\,A_{1}$}{$[2], \,[1]$}
\elldata{1019}{$14$}{$D_{5}\,+\,D_{4}\,+\,2\,A_{2}\,+\,A_{1}$}{$[1]$}
\elldata{1020}{$14$}{$D_{5}\,+\,D_{4}\,+\,A_{2}\,+\,3\,A_{1}$}{$[1]$}
\elldata{1021}{$14$}{$D_{5}\,+\,D_{4}\,+\,5\,A_{1}$}{$[2]$}
\elldata{1022}{$14$}{$D_{5}\,+\,A_{9}$}{$[1]$}
\elldata{1023}{$14$}{$D_{5}\,+\,A_{8}\,+\,A_{1}$}{$[1]$}
\elldata{1024}{$14$}{$D_{5}\,+\,A_{7}\,+\,A_{2}$}{$[1]$}
\elldata{1025}{$14$}{$D_{5}\,+\,A_{7}\,+\,2\,A_{1}$}{$[2], \,[1]$}
\elldata{1026}{$14$}{$D_{5}\,+\,A_{6}\,+\,A_{3}$}{$[1]$}
\elldata{1027}{$14$}{$D_{5}\,+\,A_{6}\,+\,A_{2}\,+\,A_{1}$}{$[1]$}
\elldata{1028}{$14$}{$D_{5}\,+\,A_{6}\,+\,3\,A_{1}$}{$[1]$}
\elldata{1029}{$14$}{$D_{5}\,+\,A_{5}\,+\,A_{4}$}{$[1]$}
\elldata{1030}{$14$}{$D_{5}\,+\,A_{5}\,+\,A_{3}\,+\,A_{1}$}{$[2], \,[1]$}
\elldata{1031}{$14$}{$D_{5}\,+\,A_{5}\,+\,2\,A_{2}$}{$[1]$}
\elldata{1032}{$14$}{$D_{5}\,+\,A_{5}\,+\,A_{2}\,+\,2\,A_{1}$}{$[1]$}
\elldata{1033}{$14$}{$D_{5}\,+\,A_{5}\,+\,4\,A_{1}$}{$[2], \,[1]$}
\elldata{1034}{$14$}{$D_{5}\,+\,2\,A_{4}\,+\,A_{1}$}{$[1]$}
\elldata{1035}{$14$}{$D_{5}\,+\,A_{4}\,+\,A_{3}\,+\,A_{2}$}{$[1]$}
\elldata{1036}{$14$}{$D_{5}\,+\,A_{4}\,+\,A_{3}\,+\,2\,A_{1}$}{$[1]$}
\elldata{1037}{$14$}{$D_{5}\,+\,A_{4}\,+\,2\,A_{2}\,+\,A_{1}$}{$[1]$}
\elldata{1038}{$14$}{$D_{5}\,+\,A_{4}\,+\,A_{2}\,+\,3\,A_{1}$}{$[1]$}
\elldata{1039}{$14$}{$D_{5}\,+\,A_{4}\,+\,5\,A_{1}$}{$[1]$}
\elldata{1040}{$14$}{$D_{5}\,+\,3\,A_{3}$}{$[2], \,[1]$}
\elldata{1041}{$14$}{$D_{5}\,+\,2\,A_{3}\,+\,A_{2}\,+\,A_{1}$}{$[1]$}
\elldata{1042}{$14$}{$D_{5}\,+\,2\,A_{3}\,+\,3\,A_{1}$}{$[2], \,[1]$}
\elldata{1043}{$14$}{$D_{5}\,+\,A_{3}\,+\,3\,A_{2}$}{$[1]$}
\elldata{1044}{$14$}{$D_{5}\,+\,A_{3}\,+\,2\,A_{2}\,+\,2\,A_{1}$}{$[1]$}
\elldata{1045}{$14$}{$D_{5}\,+\,A_{3}\,+\,A_{2}\,+\,4\,A_{1}$}{$[2], \,[1]$}
\elldata{1046}{$14$}{$D_{5}\,+\,A_{3}\,+\,6\,A_{1}$}{$[2]$}
\elldata{1047}{$14$}{$D_{5}\,+\,4\,A_{2}\,+\,A_{1}$}{$[1]$}
\elldata{1048}{$14$}{$D_{5}\,+\,3\,A_{2}\,+\,3\,A_{1}$}{$[1]$}
\elldata{1049}{$14$}{$D_{5}\,+\,2\,A_{2}\,+\,5\,A_{1}$}{$[1]$}
\elldata{1050}{$14$}{$D_{5}\,+\,A_{2}\,+\,7\,A_{1}$}{$[2]$}
\elldata{1051}{$14$}{$3\,D_{4}\,+\,A_{2}$}{$[1]$}
\elldata{1052}{$14$}{$3\,D_{4}\,+\,2\,A_{1}$}{$[2]$}
\elldata{1053}{$14$}{$2\,D_{4}\,+\,A_{6}$}{$[1]$}
\elldata{1054}{$14$}{$2\,D_{4}\,+\,A_{5}\,+\,A_{1}$}{$[2], \,[1]$}
\elldata{1055}{$14$}{$2\,D_{4}\,+\,A_{4}\,+\,A_{2}$}{$[1]$}
\elldata{1056}{$14$}{$2\,D_{4}\,+\,A_{4}\,+\,2\,A_{1}$}{$[1]$}
\elldata{1057}{$14$}{$2\,D_{4}\,+\,2\,A_{3}$}{$[2], \,[1]$}
\elldata{1058}{$14$}{$2\,D_{4}\,+\,A_{3}\,+\,A_{2}\,+\,A_{1}$}{$[1]$}
\elldata{1059}{$14$}{$2\,D_{4}\,+\,A_{3}\,+\,3\,A_{1}$}{$[2]$}
\elldata{1060}{$14$}{$2\,D_{4}\,+\,3\,A_{2}$}{$[1]$}
\elldata{1061}{$14$}{$2\,D_{4}\,+\,2\,A_{2}\,+\,2\,A_{1}$}{$[1]$}
\elldata{1062}{$14$}{$2\,D_{4}\,+\,A_{2}\,+\,4\,A_{1}$}{$[2]$}
\elldata{1063}{$14$}{$2\,D_{4}\,+\,6\,A_{1}$}{$[2, 2]$}
\elldata{1064}{$14$}{$D_{4}\,+\,A_{10}$}{$[1]$}
\elldata{1065}{$14$}{$D_{4}\,+\,A_{9}\,+\,A_{1}$}{$[2], \,[1]$}
\elldata{1066}{$14$}{$D_{4}\,+\,A_{8}\,+\,A_{2}$}{$[1]$}
\elldata{1067}{$14$}{$D_{4}\,+\,A_{8}\,+\,2\,A_{1}$}{$[1]$}
\elldata{1068}{$14$}{$D_{4}\,+\,A_{7}\,+\,A_{3}$}{$[2], \,[1]$}
\elldata{1069}{$14$}{$D_{4}\,+\,A_{7}\,+\,A_{2}\,+\,A_{1}$}{$[1]$}
\elldata{1070}{$14$}{$D_{4}\,+\,A_{7}\,+\,3\,A_{1}$}{$[2], \,[1]$}
\elldata{1071}{$14$}{$D_{4}\,+\,A_{6}\,+\,A_{4}$}{$[1]$}
\elldata{1072}{$14$}{$D_{4}\,+\,A_{6}\,+\,A_{3}\,+\,A_{1}$}{$[1]$}
\elldata{1073}{$14$}{$D_{4}\,+\,A_{6}\,+\,2\,A_{2}$}{$[1]$}
\elldata{1074}{$14$}{$D_{4}\,+\,A_{6}\,+\,A_{2}\,+\,2\,A_{1}$}{$[1]$}
\elldata{1075}{$14$}{$D_{4}\,+\,A_{6}\,+\,4\,A_{1}$}{$[1]$}
\elldata{1076}{$14$}{$D_{4}\,+\,2\,A_{5}$}{$[2], \,[1]$}
\elldata{1077}{$14$}{$D_{4}\,+\,A_{5}\,+\,A_{4}\,+\,A_{1}$}{$[1]$}
\elldata{1078}{$14$}{$D_{4}\,+\,A_{5}\,+\,A_{3}\,+\,A_{2}$}{$[1]$}
\elldata{1079}{$14$}{$D_{4}\,+\,A_{5}\,+\,A_{3}\,+\,2\,A_{1}$}{$[2], \,[1]$}
\elldata{1080}{$14$}{$D_{4}\,+\,A_{5}\,+\,2\,A_{2}\,+\,A_{1}$}{$[1]$}
\elldata{1081}{$14$}{$D_{4}\,+\,A_{5}\,+\,A_{2}\,+\,3\,A_{1}$}{$[2], \,[1]$}
\elldata{1082}{$14$}{$D_{4}\,+\,A_{5}\,+\,5\,A_{1}$}{$[2]$}
\elldata{1083}{$14$}{$D_{4}\,+\,2\,A_{4}\,+\,A_{2}$}{$[1]$}
\elldata{1084}{$14$}{$D_{4}\,+\,2\,A_{4}\,+\,2\,A_{1}$}{$[1]$}
\elldata{1085}{$14$}{$D_{4}\,+\,A_{4}\,+\,2\,A_{3}$}{$[1]$}
\elldata{1086}{$14$}{$D_{4}\,+\,A_{4}\,+\,A_{3}\,+\,A_{2}\,+\,A_{1}$}{$[1]$}
\elldata{1087}{$14$}{$D_{4}\,+\,A_{4}\,+\,A_{3}\,+\,3\,A_{1}$}{$[1]$}
\elldata{1088}{$14$}{$D_{4}\,+\,A_{4}\,+\,3\,A_{2}$}{$[1]$}
\elldata{1089}{$14$}{$D_{4}\,+\,A_{4}\,+\,2\,A_{2}\,+\,2\,A_{1}$}{$[1]$}
\elldata{1090}{$14$}{$D_{4}\,+\,A_{4}\,+\,A_{2}\,+\,4\,A_{1}$}{$[1]$}
\elldata{1091}{$14$}{$D_{4}\,+\,A_{4}\,+\,6\,A_{1}$}{$[2]$}
\elldata{1092}{$14$}{$D_{4}\,+\,3\,A_{3}\,+\,A_{1}$}{$[2], \,[1]$}
\elldata{1093}{$14$}{$D_{4}\,+\,2\,A_{3}\,+\,2\,A_{2}$}{$[1]$}
\elldata{1094}{$14$}{$D_{4}\,+\,2\,A_{3}\,+\,A_{2}\,+\,2\,A_{1}$}{$[2], \,[1]$}
\elldata{1095}{$14$}{$D_{4}\,+\,2\,A_{3}\,+\,4\,A_{1}$}{$[2]$}
\elldata{1096}{$14$}{$D_{4}\,+\,A_{3}\,+\,3\,A_{2}\,+\,A_{1}$}{$[1]$}
\elldata{1097}{$14$}{$D_{4}\,+\,A_{3}\,+\,2\,A_{2}\,+\,3\,A_{1}$}{$[1]$}
\elldata{1098}{$14$}{$D_{4}\,+\,A_{3}\,+\,A_{2}\,+\,5\,A_{1}$}{$[2]$}
\elldata{1099}{$14$}{$D_{4}\,+\,A_{3}\,+\,7\,A_{1}$}{$[2, 2]$}
\elldata{1100}{$14$}{$D_{4}\,+\,5\,A_{2}$}{$[1]$}
\elldata{1101}{$14$}{$D_{4}\,+\,4\,A_{2}\,+\,2\,A_{1}$}{$[1]$}
\elldata{1102}{$14$}{$D_{4}\,+\,3\,A_{2}\,+\,4\,A_{1}$}{$[1]$}
\elldata{1103}{$14$}{$D_{4}\,+\,2\,A_{2}\,+\,6\,A_{1}$}{$[2]$}
\elldata{1104}{$14$}{$A_{14}$}{$[1]$}
\elldata{1105}{$14$}{$A_{13}\,+\,A_{1}$}{$[2], \,[1]$}
\elldata{1106}{$14$}{$A_{12}\,+\,A_{2}$}{$[1]$}
\elldata{1107}{$14$}{$A_{12}\,+\,2\,A_{1}$}{$[1]$}
\elldata{1108}{$14$}{$A_{11}\,+\,A_{3}$}{$[2], \,[1]$}
\elldata{1109}{$14$}{$A_{11}\,+\,A_{2}\,+\,A_{1}$}{$[1]$}
\elldata{1110}{$14$}{$A_{11}\,+\,3\,A_{1}$}{$[2], \,[1]$}
\elldata{1111}{$14$}{$A_{10}\,+\,A_{4}$}{$[1]$}
\elldata{1112}{$14$}{$A_{10}\,+\,A_{3}\,+\,A_{1}$}{$[1]$}
\elldata{1113}{$14$}{$A_{10}\,+\,2\,A_{2}$}{$[1]$}
\elldata{1114}{$14$}{$A_{10}\,+\,A_{2}\,+\,2\,A_{1}$}{$[1]$}
\elldata{1115}{$14$}{$A_{10}\,+\,4\,A_{1}$}{$[1]$}
\elldata{1116}{$14$}{$A_{9}\,+\,A_{5}$}{$[2], \,[1]$}
\elldata{1117}{$14$}{$A_{9}\,+\,A_{4}\,+\,A_{1}$}{$[1]$}
\elldata{1118}{$14$}{$A_{9}\,+\,A_{3}\,+\,A_{2}$}{$[1]$}
\elldata{1119}{$14$}{$A_{9}\,+\,A_{3}\,+\,2\,A_{1}$}{$[2], \,[1]$}
\elldata{1120}{$14$}{$A_{9}\,+\,2\,A_{2}\,+\,A_{1}$}{$[1]$}
\elldata{1121}{$14$}{$A_{9}\,+\,A_{2}\,+\,3\,A_{1}$}{$[2], \,[1]$}
\elldata{1122}{$14$}{$A_{9}\,+\,5\,A_{1}$}{$[2], \,[1]$}
\elldata{1123}{$14$}{$A_{8}\,+\,A_{6}$}{$[1]$}
\elldata{1124}{$14$}{$A_{8}\,+\,A_{5}\,+\,A_{1}$}{$[1]$}
\elldata{1125}{$14$}{$A_{8}\,+\,A_{4}\,+\,A_{2}$}{$[1]$}
\elldata{1126}{$14$}{$A_{8}\,+\,A_{4}\,+\,2\,A_{1}$}{$[1]$}
\elldata{1127}{$14$}{$A_{8}\,+\,2\,A_{3}$}{$[1]$}
\elldata{1128}{$14$}{$A_{8}\,+\,A_{3}\,+\,A_{2}\,+\,A_{1}$}{$[1]$}
\elldata{1129}{$14$}{$A_{8}\,+\,A_{3}\,+\,3\,A_{1}$}{$[1]$}
\elldata{1130}{$14$}{$A_{8}\,+\,3\,A_{2}$}{$[3], \,[1]$}
\elldata{1131}{$14$}{$A_{8}\,+\,2\,A_{2}\,+\,2\,A_{1}$}{$[1]$}
\elldata{1132}{$14$}{$A_{8}\,+\,A_{2}\,+\,4\,A_{1}$}{$[1]$}
\elldata{1133}{$14$}{$A_{8}\,+\,6\,A_{1}$}{$[1]$}
\elldata{1134}{$14$}{$2\,A_{7}$}{$[2], \,[1]$}
\elldata{1135}{$14$}{$A_{7}\,+\,A_{6}\,+\,A_{1}$}{$[1]$}
\elldata{1136}{$14$}{$A_{7}\,+\,A_{5}\,+\,A_{2}$}{$[1]$}
\elldata{1137}{$14$}{$A_{7}\,+\,A_{5}\,+\,2\,A_{1}$}{$[2], \,[1]$}
\elldata{1138}{$14$}{$A_{7}\,+\,A_{4}\,+\,A_{3}$}{$[1]$}
\elldata{1139}{$14$}{$A_{7}\,+\,A_{4}\,+\,A_{2}\,+\,A_{1}$}{$[1]$}
\elldata{1140}{$14$}{$A_{7}\,+\,A_{4}\,+\,3\,A_{1}$}{$[1]$}
\elldata{1141}{$14$}{$A_{7}\,+\,2\,A_{3}\,+\,A_{1}$}{$[2], \,[1]$}
\elldata{1142}{$14$}{$A_{7}\,+\,A_{3}\,+\,2\,A_{2}$}{$[1]$}
\elldata{1143}{$14$}{$A_{7}\,+\,A_{3}\,+\,A_{2}\,+\,2\,A_{1}$}{$[2], \,[1]$}
\elldata{1144}{$14$}{$A_{7}\,+\,A_{3}\,+\,4\,A_{1}$}{$[2], \,[1]$}
\elldata{1145}{$14$}{$A_{7}\,+\,3\,A_{2}\,+\,A_{1}$}{$[1]$}
\elldata{1146}{$14$}{$A_{7}\,+\,2\,A_{2}\,+\,3\,A_{1}$}{$[1]$}
\elldata{1147}{$14$}{$A_{7}\,+\,A_{2}\,+\,5\,A_{1}$}{$[2], \,[1]$}
\elldata{1148}{$14$}{$A_{7}\,+\,7\,A_{1}$}{$[2]$}
\elldata{1149}{$14$}{$2\,A_{6}\,+\,A_{2}$}{$[1]$}
\elldata{1150}{$14$}{$2\,A_{6}\,+\,2\,A_{1}$}{$[1]$}
\elldata{1151}{$14$}{$A_{6}\,+\,A_{5}\,+\,A_{3}$}{$[1]$}
\elldata{1152}{$14$}{$A_{6}\,+\,A_{5}\,+\,A_{2}\,+\,A_{1}$}{$[1]$}
\elldata{1153}{$14$}{$A_{6}\,+\,A_{5}\,+\,3\,A_{1}$}{$[1]$}
\elldata{1154}{$14$}{$A_{6}\,+\,2\,A_{4}$}{$[1]$}
\elldata{1155}{$14$}{$A_{6}\,+\,A_{4}\,+\,A_{3}\,+\,A_{1}$}{$[1]$}
\elldata{1156}{$14$}{$A_{6}\,+\,A_{4}\,+\,2\,A_{2}$}{$[1]$}
\elldata{1157}{$14$}{$A_{6}\,+\,A_{4}\,+\,A_{2}\,+\,2\,A_{1}$}{$[1]$}
\elldata{1158}{$14$}{$A_{6}\,+\,A_{4}\,+\,4\,A_{1}$}{$[1]$}
\elldata{1159}{$14$}{$A_{6}\,+\,2\,A_{3}\,+\,A_{2}$}{$[1]$}
\elldata{1160}{$14$}{$A_{6}\,+\,2\,A_{3}\,+\,2\,A_{1}$}{$[1]$}
\elldata{1161}{$14$}{$A_{6}\,+\,A_{3}\,+\,2\,A_{2}\,+\,A_{1}$}{$[1]$}
\elldata{1162}{$14$}{$A_{6}\,+\,A_{3}\,+\,A_{2}\,+\,3\,A_{1}$}{$[1]$}
\elldata{1163}{$14$}{$A_{6}\,+\,A_{3}\,+\,5\,A_{1}$}{$[1]$}
\elldata{1164}{$14$}{$A_{6}\,+\,4\,A_{2}$}{$[1]$}
\elldata{1165}{$14$}{$A_{6}\,+\,3\,A_{2}\,+\,2\,A_{1}$}{$[1]$}
\elldata{1166}{$14$}{$A_{6}\,+\,2\,A_{2}\,+\,4\,A_{1}$}{$[1]$}
\elldata{1167}{$14$}{$A_{6}\,+\,A_{2}\,+\,6\,A_{1}$}{$[1]$}
\elldata{1168}{$14$}{$A_{6}\,+\,8\,A_{1}$}{$[2]$}
\elldata{1169}{$14$}{$2\,A_{5}\,+\,A_{4}$}{$[1]$}
\elldata{1170}{$14$}{$2\,A_{5}\,+\,A_{3}\,+\,A_{1}$}{$[2], \,[1]$}
\elldata{1171}{$14$}{$2\,A_{5}\,+\,2\,A_{2}$}{$[3], \,[1]$}
\elldata{1172}{$14$}{$2\,A_{5}\,+\,A_{2}\,+\,2\,A_{1}$}{$[2], \,[1]$}
\elldata{1173}{$14$}{$2\,A_{5}\,+\,4\,A_{1}$}{$[2], \,[1]$}
\elldata{1174}{$14$}{$A_{5}\,+\,2\,A_{4}\,+\,A_{1}$}{$[1]$}
\elldata{1175}{$14$}{$A_{5}\,+\,A_{4}\,+\,A_{3}\,+\,A_{2}$}{$[1]$}
\elldata{1176}{$14$}{$A_{5}\,+\,A_{4}\,+\,A_{3}\,+\,2\,A_{1}$}{$[1]$}
\elldata{1177}{$14$}{$A_{5}\,+\,A_{4}\,+\,2\,A_{2}\,+\,A_{1}$}{$[1]$}
\elldata{1178}{$14$}{$A_{5}\,+\,A_{4}\,+\,A_{2}\,+\,3\,A_{1}$}{$[1]$}
\elldata{1179}{$14$}{$A_{5}\,+\,A_{4}\,+\,5\,A_{1}$}{$[2], \,[1]$}
\elldata{1180}{$14$}{$A_{5}\,+\,3\,A_{3}$}{$[1]$}
\elldata{1181}{$14$}{$A_{5}\,+\,2\,A_{3}\,+\,A_{2}\,+\,A_{1}$}{$[2], \,[1]$}
\elldata{1182}{$14$}{$A_{5}\,+\,2\,A_{3}\,+\,3\,A_{1}$}{$[2], \,[1]$}
\elldata{1183}{$14$}{$A_{5}\,+\,A_{3}\,+\,3\,A_{2}$}{$[1]$}
\elldata{1184}{$14$}{$A_{5}\,+\,A_{3}\,+\,2\,A_{2}\,+\,2\,A_{1}$}{$[1]$}
\elldata{1185}{$14$}{$A_{5}\,+\,A_{3}\,+\,A_{2}\,+\,4\,A_{1}$}{$[2], \,[1]$}
\elldata{1186}{$14$}{$A_{5}\,+\,A_{3}\,+\,6\,A_{1}$}{$[2]$}
\elldata{1187}{$14$}{$A_{5}\,+\,4\,A_{2}\,+\,A_{1}$}{$[3], \,[1]$}
\elldata{1188}{$14$}{$A_{5}\,+\,3\,A_{2}\,+\,3\,A_{1}$}{$[1]$}
\elldata{1189}{$14$}{$A_{5}\,+\,2\,A_{2}\,+\,5\,A_{1}$}{$[2], \,[1]$}
\elldata{1190}{$14$}{$A_{5}\,+\,A_{2}\,+\,7\,A_{1}$}{$[2]$}
\elldata{1191}{$14$}{$A_{5}\,+\,9\,A_{1}$}{$[2, 2]$}
\elldata{1192}{$14$}{$3\,A_{4}\,+\,A_{2}$}{$[1]$}
\elldata{1193}{$14$}{$3\,A_{4}\,+\,2\,A_{1}$}{$[1]$}
\elldata{1194}{$14$}{$2\,A_{4}\,+\,2\,A_{3}$}{$[1]$}
\elldata{1195}{$14$}{$2\,A_{4}\,+\,A_{3}\,+\,A_{2}\,+\,A_{1}$}{$[1]$}
\elldata{1196}{$14$}{$2\,A_{4}\,+\,A_{3}\,+\,3\,A_{1}$}{$[1]$}
\elldata{1197}{$14$}{$2\,A_{4}\,+\,3\,A_{2}$}{$[1]$}
\elldata{1198}{$14$}{$2\,A_{4}\,+\,2\,A_{2}\,+\,2\,A_{1}$}{$[1]$}
\elldata{1199}{$14$}{$2\,A_{4}\,+\,A_{2}\,+\,4\,A_{1}$}{$[1]$}
\elldata{1200}{$14$}{$2\,A_{4}\,+\,6\,A_{1}$}{$[1]$}
\elldata{1201}{$14$}{$A_{4}\,+\,3\,A_{3}\,+\,A_{1}$}{$[1]$}
\elldata{1202}{$14$}{$A_{4}\,+\,2\,A_{3}\,+\,2\,A_{2}$}{$[1]$}
\elldata{1203}{$14$}{$A_{4}\,+\,2\,A_{3}\,+\,A_{2}\,+\,2\,A_{1}$}{$[1]$}
\elldata{1204}{$14$}{$A_{4}\,+\,2\,A_{3}\,+\,4\,A_{1}$}{$[2], \,[1]$}
\elldata{1205}{$14$}{$A_{4}\,+\,A_{3}\,+\,3\,A_{2}\,+\,A_{1}$}{$[1]$}
\elldata{1206}{$14$}{$A_{4}\,+\,A_{3}\,+\,2\,A_{2}\,+\,3\,A_{1}$}{$[1]$}
\elldata{1207}{$14$}{$A_{4}\,+\,A_{3}\,+\,A_{2}\,+\,5\,A_{1}$}{$[1]$}
\elldata{1208}{$14$}{$A_{4}\,+\,A_{3}\,+\,7\,A_{1}$}{$[2]$}
\elldata{1209}{$14$}{$A_{4}\,+\,5\,A_{2}$}{$[1]$}
\elldata{1210}{$14$}{$A_{4}\,+\,4\,A_{2}\,+\,2\,A_{1}$}{$[1]$}
\elldata{1211}{$14$}{$A_{4}\,+\,3\,A_{2}\,+\,4\,A_{1}$}{$[1]$}
\elldata{1212}{$14$}{$A_{4}\,+\,2\,A_{2}\,+\,6\,A_{1}$}{$[1]$}
\elldata{1213}{$14$}{$A_{4}\,+\,A_{2}\,+\,8\,A_{1}$}{$[2]$}
\elldata{1214}{$14$}{$4\,A_{3}\,+\,A_{2}$}{$[2], \,[1]$}
\elldata{1215}{$14$}{$4\,A_{3}\,+\,2\,A_{1}$}{$[4], \,[2], \,[1]$}
\elldata{1216}{$14$}{$3\,A_{3}\,+\,2\,A_{2}\,+\,A_{1}$}{$[1]$}
\elldata{1217}{$14$}{$3\,A_{3}\,+\,A_{2}\,+\,3\,A_{1}$}{$[2], \,[1]$}
\elldata{1218}{$14$}{$3\,A_{3}\,+\,5\,A_{1}$}{$[2]$}
\elldata{1219}{$14$}{$2\,A_{3}\,+\,4\,A_{2}$}{$[1]$}
\elldata{1220}{$14$}{$2\,A_{3}\,+\,3\,A_{2}\,+\,2\,A_{1}$}{$[1]$}
\elldata{1221}{$14$}{$2\,A_{3}\,+\,2\,A_{2}\,+\,4\,A_{1}$}{$[2], \,[1]$}
\elldata{1222}{$14$}{$2\,A_{3}\,+\,A_{2}\,+\,6\,A_{1}$}{$[2]$}
\elldata{1223}{$14$}{$2\,A_{3}\,+\,8\,A_{1}$}{$[2, 2]$}
\elldata{1224}{$14$}{$A_{3}\,+\,5\,A_{2}\,+\,A_{1}$}{$[1]$}
\elldata{1225}{$14$}{$A_{3}\,+\,4\,A_{2}\,+\,3\,A_{1}$}{$[1]$}
\elldata{1226}{$14$}{$A_{3}\,+\,3\,A_{2}\,+\,5\,A_{1}$}{$[1]$}
\elldata{1227}{$14$}{$A_{3}\,+\,2\,A_{2}\,+\,7\,A_{1}$}{$[2]$}
\elldata{1228}{$14$}{$7\,A_{2}$}{$[3]$}
\elldata{1229}{$14$}{$6\,A_{2}\,+\,2\,A_{1}$}{$[3], \,[1]$}
\elldata{1230}{$14$}{$5\,A_{2}\,+\,4\,A_{1}$}{$[1]$}
\elldata{1231}{$14$}{$4\,A_{2}\,+\,6\,A_{1}$}{$[1]$}

\vsr \elldata{No.}{rank}{$ADE$-type}{$G$}

\vsrs \elldata{1232}{$15$}{$E_{8}\,+\,E_{7}$}{$[1]$}
\elldata{1233}{$15$}{$E_{8}\,+\,E_{6}\,+\,A_{1}$}{$[1]$}
\elldata{1234}{$15$}{$E_{8}\,+\,D_{7}$}{$[1]$}
\elldata{1235}{$15$}{$E_{8}\,+\,D_{6}\,+\,A_{1}$}{$[1]$}
\elldata{1236}{$15$}{$E_{8}\,+\,D_{5}\,+\,A_{2}$}{$[1]$}
\elldata{1237}{$15$}{$E_{8}\,+\,D_{5}\,+\,2\,A_{1}$}{$[1]$}
\elldata{1238}{$15$}{$E_{8}\,+\,D_{4}\,+\,A_{3}$}{$[1]$}
\elldata{1239}{$15$}{$E_{8}\,+\,D_{4}\,+\,A_{2}\,+\,A_{1}$}{$[1]$}
\elldata{1240}{$15$}{$E_{8}\,+\,D_{4}\,+\,3\,A_{1}$}{$[1]$}
\elldata{1241}{$15$}{$E_{8}\,+\,A_{7}$}{$[1]$}
\elldata{1242}{$15$}{$E_{8}\,+\,A_{6}\,+\,A_{1}$}{$[1]$}
\elldata{1243}{$15$}{$E_{8}\,+\,A_{5}\,+\,A_{2}$}{$[1]$}
\elldata{1244}{$15$}{$E_{8}\,+\,A_{5}\,+\,2\,A_{1}$}{$[1]$}
\elldata{1245}{$15$}{$E_{8}\,+\,A_{4}\,+\,A_{3}$}{$[1]$}
\elldata{1246}{$15$}{$E_{8}\,+\,A_{4}\,+\,A_{2}\,+\,A_{1}$}{$[1]$}
\elldata{1247}{$15$}{$E_{8}\,+\,A_{4}\,+\,3\,A_{1}$}{$[1]$}
\elldata{1248}{$15$}{$E_{8}\,+\,2\,A_{3}\,+\,A_{1}$}{$[1]$}
\elldata{1249}{$15$}{$E_{8}\,+\,A_{3}\,+\,2\,A_{2}$}{$[1]$}
\elldata{1250}{$15$}{$E_{8}\,+\,A_{3}\,+\,A_{2}\,+\,2\,A_{1}$}{$[1]$}
\elldata{1251}{$15$}{$E_{8}\,+\,A_{3}\,+\,4\,A_{1}$}{$[1]$}
\elldata{1252}{$15$}{$E_{8}\,+\,3\,A_{2}\,+\,A_{1}$}{$[1]$}
\elldata{1253}{$15$}{$E_{8}\,+\,2\,A_{2}\,+\,3\,A_{1}$}{$[1]$}
\elldata{1254}{$15$}{$E_{8}\,+\,A_{2}\,+\,5\,A_{1}$}{$[1]$}
\elldata{1255}{$15$}{$2\,E_{7}\,+\,A_{1}$}{$[1]$}
\elldata{1256}{$15$}{$E_{7}\,+\,E_{6}\,+\,A_{2}$}{$[1]$}
\elldata{1257}{$15$}{$E_{7}\,+\,E_{6}\,+\,2\,A_{1}$}{$[1]$}
\elldata{1258}{$15$}{$E_{7}\,+\,D_{8}$}{$[1]$}
\elldata{1259}{$15$}{$E_{7}\,+\,D_{7}\,+\,A_{1}$}{$[1]$}
\elldata{1260}{$15$}{$E_{7}\,+\,D_{6}\,+\,A_{2}$}{$[1]$}
\elldata{1261}{$15$}{$E_{7}\,+\,D_{6}\,+\,2\,A_{1}$}{$[2], \,[1]$}
\elldata{1262}{$15$}{$E_{7}\,+\,D_{5}\,+\,A_{3}$}{$[1]$}
\elldata{1263}{$15$}{$E_{7}\,+\,D_{5}\,+\,A_{2}\,+\,A_{1}$}{$[1]$}
\elldata{1264}{$15$}{$E_{7}\,+\,D_{5}\,+\,3\,A_{1}$}{$[2], \,[1]$}
\elldata{1265}{$15$}{$E_{7}\,+\,2\,D_{4}$}{$[1]$}
\elldata{1266}{$15$}{$E_{7}\,+\,D_{4}\,+\,A_{4}$}{$[1]$}
\elldata{1267}{$15$}{$E_{7}\,+\,D_{4}\,+\,A_{3}\,+\,A_{1}$}{$[2], \,[1]$}
\elldata{1268}{$15$}{$E_{7}\,+\,D_{4}\,+\,2\,A_{2}$}{$[1]$}
\elldata{1269}{$15$}{$E_{7}\,+\,D_{4}\,+\,A_{2}\,+\,2\,A_{1}$}{$[1]$}
\elldata{1270}{$15$}{$E_{7}\,+\,D_{4}\,+\,4\,A_{1}$}{$[2]$}
\elldata{1271}{$15$}{$E_{7}\,+\,A_{8}$}{$[1]$}
\elldata{1272}{$15$}{$E_{7}\,+\,A_{7}\,+\,A_{1}$}{$[2], \,[1]$}
\elldata{1273}{$15$}{$E_{7}\,+\,A_{6}\,+\,A_{2}$}{$[1]$}
\elldata{1274}{$15$}{$E_{7}\,+\,A_{6}\,+\,2\,A_{1}$}{$[1]$}
\elldata{1275}{$15$}{$E_{7}\,+\,A_{5}\,+\,A_{3}$}{$[2], \,[1]$}
\elldata{1276}{$15$}{$E_{7}\,+\,A_{5}\,+\,A_{2}\,+\,A_{1}$}{$[1]$}
\elldata{1277}{$15$}{$E_{7}\,+\,A_{5}\,+\,3\,A_{1}$}{$[2], \,[1]$}
\elldata{1278}{$15$}{$E_{7}\,+\,2\,A_{4}$}{$[1]$}
\elldata{1279}{$15$}{$E_{7}\,+\,A_{4}\,+\,A_{3}\,+\,A_{1}$}{$[1]$}
\elldata{1280}{$15$}{$E_{7}\,+\,A_{4}\,+\,2\,A_{2}$}{$[1]$}
\elldata{1281}{$15$}{$E_{7}\,+\,A_{4}\,+\,A_{2}\,+\,2\,A_{1}$}{$[1]$}
\elldata{1282}{$15$}{$E_{7}\,+\,A_{4}\,+\,4\,A_{1}$}{$[1]$}
\elldata{1283}{$15$}{$E_{7}\,+\,2\,A_{3}\,+\,A_{2}$}{$[1]$}
\elldata{1284}{$15$}{$E_{7}\,+\,2\,A_{3}\,+\,2\,A_{1}$}{$[2], \,[1]$}
\elldata{1285}{$15$}{$E_{7}\,+\,A_{3}\,+\,2\,A_{2}\,+\,A_{1}$}{$[1]$}
\elldata{1286}{$15$}{$E_{7}\,+\,A_{3}\,+\,A_{2}\,+\,3\,A_{1}$}{$[2], \,[1]$}
\elldata{1287}{$15$}{$E_{7}\,+\,A_{3}\,+\,5\,A_{1}$}{$[2]$}
\elldata{1288}{$15$}{$E_{7}\,+\,4\,A_{2}$}{$[1]$}
\elldata{1289}{$15$}{$E_{7}\,+\,3\,A_{2}\,+\,2\,A_{1}$}{$[1]$}
\elldata{1290}{$15$}{$E_{7}\,+\,2\,A_{2}\,+\,4\,A_{1}$}{$[1]$}
\elldata{1291}{$15$}{$E_{7}\,+\,A_{2}\,+\,6\,A_{1}$}{$[2]$}
\elldata{1292}{$15$}{$2\,E_{6}\,+\,A_{3}$}{$[1]$}
\elldata{1293}{$15$}{$2\,E_{6}\,+\,A_{2}\,+\,A_{1}$}{$[1]$}
\elldata{1294}{$15$}{$2\,E_{6}\,+\,3\,A_{1}$}{$[1]$}
\elldata{1295}{$15$}{$E_{6}\,+\,D_{9}$}{$[1]$}
\elldata{1296}{$15$}{$E_{6}\,+\,D_{8}\,+\,A_{1}$}{$[1]$}
\elldata{1297}{$15$}{$E_{6}\,+\,D_{7}\,+\,A_{2}$}{$[1]$}
\elldata{1298}{$15$}{$E_{6}\,+\,D_{7}\,+\,2\,A_{1}$}{$[1]$}
\elldata{1299}{$15$}{$E_{6}\,+\,D_{6}\,+\,A_{3}$}{$[1]$}
\elldata{1300}{$15$}{$E_{6}\,+\,D_{6}\,+\,A_{2}\,+\,A_{1}$}{$[1]$}
\elldata{1301}{$15$}{$E_{6}\,+\,D_{6}\,+\,3\,A_{1}$}{$[1]$}
\elldata{1302}{$15$}{$E_{6}\,+\,D_{5}\,+\,D_{4}$}{$[1]$}
\elldata{1303}{$15$}{$E_{6}\,+\,D_{5}\,+\,A_{4}$}{$[1]$}
\elldata{1304}{$15$}{$E_{6}\,+\,D_{5}\,+\,A_{3}\,+\,A_{1}$}{$[1]$}
\elldata{1305}{$15$}{$E_{6}\,+\,D_{5}\,+\,2\,A_{2}$}{$[1]$}
\elldata{1306}{$15$}{$E_{6}\,+\,D_{5}\,+\,A_{2}\,+\,2\,A_{1}$}{$[1]$}
\elldata{1307}{$15$}{$E_{6}\,+\,D_{5}\,+\,4\,A_{1}$}{$[1]$}
\elldata{1308}{$15$}{$E_{6}\,+\,2\,D_{4}\,+\,A_{1}$}{$[1]$}
\elldata{1309}{$15$}{$E_{6}\,+\,D_{4}\,+\,A_{5}$}{$[1]$}
\elldata{1310}{$15$}{$E_{6}\,+\,D_{4}\,+\,A_{4}\,+\,A_{1}$}{$[1]$}
\elldata{1311}{$15$}{$E_{6}\,+\,D_{4}\,+\,A_{3}\,+\,A_{2}$}{$[1]$}
\elldata{1312}{$15$}{$E_{6}\,+\,D_{4}\,+\,A_{3}\,+\,2\,A_{1}$}{$[1]$}
\elldata{1313}{$15$}{$E_{6}\,+\,D_{4}\,+\,2\,A_{2}\,+\,A_{1}$}{$[1]$}
\elldata{1314}{$15$}{$E_{6}\,+\,D_{4}\,+\,A_{2}\,+\,3\,A_{1}$}{$[1]$}
\elldata{1315}{$15$}{$E_{6}\,+\,A_{9}$}{$[1]$}
\elldata{1316}{$15$}{$E_{6}\,+\,A_{8}\,+\,A_{1}$}{$[1]$}
\elldata{1317}{$15$}{$E_{6}\,+\,A_{7}\,+\,A_{2}$}{$[1]$}
\elldata{1318}{$15$}{$E_{6}\,+\,A_{7}\,+\,2\,A_{1}$}{$[1]$}
\elldata{1319}{$15$}{$E_{6}\,+\,A_{6}\,+\,A_{3}$}{$[1]$}
\elldata{1320}{$15$}{$E_{6}\,+\,A_{6}\,+\,A_{2}\,+\,A_{1}$}{$[1]$}
\elldata{1321}{$15$}{$E_{6}\,+\,A_{6}\,+\,3\,A_{1}$}{$[1]$}
\elldata{1322}{$15$}{$E_{6}\,+\,A_{5}\,+\,A_{4}$}{$[1]$}
\elldata{1323}{$15$}{$E_{6}\,+\,A_{5}\,+\,A_{3}\,+\,A_{1}$}{$[1]$}
\elldata{1324}{$15$}{$E_{6}\,+\,A_{5}\,+\,2\,A_{2}$}{$[3], \,[1]$}
\elldata{1325}{$15$}{$E_{6}\,+\,A_{5}\,+\,A_{2}\,+\,2\,A_{1}$}{$[1]$}
\elldata{1326}{$15$}{$E_{6}\,+\,A_{5}\,+\,4\,A_{1}$}{$[1]$}
\elldata{1327}{$15$}{$E_{6}\,+\,2\,A_{4}\,+\,A_{1}$}{$[1]$}
\elldata{1328}{$15$}{$E_{6}\,+\,A_{4}\,+\,A_{3}\,+\,A_{2}$}{$[1]$}
\elldata{1329}{$15$}{$E_{6}\,+\,A_{4}\,+\,A_{3}\,+\,2\,A_{1}$}{$[1]$}
\elldata{1330}{$15$}{$E_{6}\,+\,A_{4}\,+\,2\,A_{2}\,+\,A_{1}$}{$[1]$}
\elldata{1331}{$15$}{$E_{6}\,+\,A_{4}\,+\,A_{2}\,+\,3\,A_{1}$}{$[1]$}
\elldata{1332}{$15$}{$E_{6}\,+\,A_{4}\,+\,5\,A_{1}$}{$[1]$}
\elldata{1333}{$15$}{$E_{6}\,+\,3\,A_{3}$}{$[1]$}
\elldata{1334}{$15$}{$E_{6}\,+\,2\,A_{3}\,+\,A_{2}\,+\,A_{1}$}{$[1]$}
\elldata{1335}{$15$}{$E_{6}\,+\,2\,A_{3}\,+\,3\,A_{1}$}{$[1]$}
\elldata{1336}{$15$}{$E_{6}\,+\,A_{3}\,+\,3\,A_{2}$}{$[1]$}
\elldata{1337}{$15$}{$E_{6}\,+\,A_{3}\,+\,2\,A_{2}\,+\,2\,A_{1}$}{$[1]$}
\elldata{1338}{$15$}{$E_{6}\,+\,A_{3}\,+\,A_{2}\,+\,4\,A_{1}$}{$[1]$}
\elldata{1339}{$15$}{$E_{6}\,+\,4\,A_{2}\,+\,A_{1}$}{$[3], \,[1]$}
\elldata{1340}{$15$}{$E_{6}\,+\,3\,A_{2}\,+\,3\,A_{1}$}{$[1]$}
\elldata{1341}{$15$}{$E_{6}\,+\,2\,A_{2}\,+\,5\,A_{1}$}{$[1]$}
\elldata{1342}{$15$}{$D_{15}$}{$[1]$}
\elldata{1343}{$15$}{$D_{14}\,+\,A_{1}$}{$[2], \,[1]$}
\elldata{1344}{$15$}{$D_{13}\,+\,A_{2}$}{$[1]$}
\elldata{1345}{$15$}{$D_{13}\,+\,2\,A_{1}$}{$[1]$}
\elldata{1346}{$15$}{$D_{12}\,+\,A_{3}$}{$[2], \,[1]$}
\elldata{1347}{$15$}{$D_{12}\,+\,A_{2}\,+\,A_{1}$}{$[1]$}
\elldata{1348}{$15$}{$D_{12}\,+\,3\,A_{1}$}{$[2], \,[1]$}
\elldata{1349}{$15$}{$D_{11}\,+\,D_{4}$}{$[1]$}
\elldata{1350}{$15$}{$D_{11}\,+\,A_{4}$}{$[1]$}
\elldata{1351}{$15$}{$D_{11}\,+\,A_{3}\,+\,A_{1}$}{$[1]$}
\elldata{1352}{$15$}{$D_{11}\,+\,2\,A_{2}$}{$[1]$}
\elldata{1353}{$15$}{$D_{11}\,+\,A_{2}\,+\,2\,A_{1}$}{$[1]$}
\elldata{1354}{$15$}{$D_{11}\,+\,4\,A_{1}$}{$[1]$}
\elldata{1355}{$15$}{$D_{10}\,+\,D_{5}$}{$[1]$}
\elldata{1356}{$15$}{$D_{10}\,+\,D_{4}\,+\,A_{1}$}{$[2], \,[1]$}
\elldata{1357}{$15$}{$D_{10}\,+\,A_{5}$}{$[2], \,[1]$}
\elldata{1358}{$15$}{$D_{10}\,+\,A_{4}\,+\,A_{1}$}{$[1]$}
\elldata{1359}{$15$}{$D_{10}\,+\,A_{3}\,+\,A_{2}$}{$[1]$}
\elldata{1360}{$15$}{$D_{10}\,+\,A_{3}\,+\,2\,A_{1}$}{$[2], \,[1]$}
\elldata{1361}{$15$}{$D_{10}\,+\,2\,A_{2}\,+\,A_{1}$}{$[1]$}
\elldata{1362}{$15$}{$D_{10}\,+\,A_{2}\,+\,3\,A_{1}$}{$[2], \,[1]$}
\elldata{1363}{$15$}{$D_{10}\,+\,5\,A_{1}$}{$[2]$}
\elldata{1364}{$15$}{$D_{9}\,+\,D_{6}$}{$[1]$}
\elldata{1365}{$15$}{$D_{9}\,+\,D_{5}\,+\,A_{1}$}{$[1]$}
\elldata{1366}{$15$}{$D_{9}\,+\,D_{4}\,+\,A_{2}$}{$[1]$}
\elldata{1367}{$15$}{$D_{9}\,+\,D_{4}\,+\,2\,A_{1}$}{$[1]$}
\elldata{1368}{$15$}{$D_{9}\,+\,A_{6}$}{$[1]$}
\elldata{1369}{$15$}{$D_{9}\,+\,A_{5}\,+\,A_{1}$}{$[1]$}
\elldata{1370}{$15$}{$D_{9}\,+\,A_{4}\,+\,A_{2}$}{$[1]$}
\elldata{1371}{$15$}{$D_{9}\,+\,A_{4}\,+\,2\,A_{1}$}{$[1]$}
\elldata{1372}{$15$}{$D_{9}\,+\,2\,A_{3}$}{$[1]$}
\elldata{1373}{$15$}{$D_{9}\,+\,A_{3}\,+\,A_{2}\,+\,A_{1}$}{$[1]$}
\elldata{1374}{$15$}{$D_{9}\,+\,A_{3}\,+\,3\,A_{1}$}{$[1]$}
\elldata{1375}{$15$}{$D_{9}\,+\,3\,A_{2}$}{$[1]$}
\elldata{1376}{$15$}{$D_{9}\,+\,2\,A_{2}\,+\,2\,A_{1}$}{$[1]$}
\elldata{1377}{$15$}{$D_{9}\,+\,A_{2}\,+\,4\,A_{1}$}{$[1]$}
\elldata{1378}{$15$}{$D_{9}\,+\,6\,A_{1}$}{$[2]$}
\elldata{1379}{$15$}{$D_{8}\,+\,D_{7}$}{$[1]$}
\elldata{1380}{$15$}{$D_{8}\,+\,D_{6}\,+\,A_{1}$}{$[2], \,[1]$}
\elldata{1381}{$15$}{$D_{8}\,+\,D_{5}\,+\,A_{2}$}{$[1]$}
\elldata{1382}{$15$}{$D_{8}\,+\,D_{5}\,+\,2\,A_{1}$}{$[2], \,[1]$}
\elldata{1383}{$15$}{$D_{8}\,+\,D_{4}\,+\,A_{3}$}{$[2], \,[1]$}
\elldata{1384}{$15$}{$D_{8}\,+\,D_{4}\,+\,A_{2}\,+\,A_{1}$}{$[1]$}
\elldata{1385}{$15$}{$D_{8}\,+\,D_{4}\,+\,3\,A_{1}$}{$[2]$}
\elldata{1386}{$15$}{$D_{8}\,+\,A_{7}$}{$[2], \,[1]$}
\elldata{1387}{$15$}{$D_{8}\,+\,A_{6}\,+\,A_{1}$}{$[1]$}
\elldata{1388}{$15$}{$D_{8}\,+\,A_{5}\,+\,A_{2}$}{$[1]$}
\elldata{1389}{$15$}{$D_{8}\,+\,A_{5}\,+\,2\,A_{1}$}{$[2], \,[1]$}
\elldata{1390}{$15$}{$D_{8}\,+\,A_{4}\,+\,A_{3}$}{$[1]$}
\elldata{1391}{$15$}{$D_{8}\,+\,A_{4}\,+\,A_{2}\,+\,A_{1}$}{$[1]$}
\elldata{1392}{$15$}{$D_{8}\,+\,A_{4}\,+\,3\,A_{1}$}{$[1]$}
\elldata{1393}{$15$}{$D_{8}\,+\,2\,A_{3}\,+\,A_{1}$}{$[2], \,[1]$}
\elldata{1394}{$15$}{$D_{8}\,+\,A_{3}\,+\,2\,A_{2}$}{$[1]$}
\elldata{1395}{$15$}{$D_{8}\,+\,A_{3}\,+\,A_{2}\,+\,2\,A_{1}$}{$[2], \,[1]$}
\elldata{1396}{$15$}{$D_{8}\,+\,A_{3}\,+\,4\,A_{1}$}{$[2]$}
\elldata{1397}{$15$}{$D_{8}\,+\,3\,A_{2}\,+\,A_{1}$}{$[1]$}
\elldata{1398}{$15$}{$D_{8}\,+\,2\,A_{2}\,+\,3\,A_{1}$}{$[1]$}
\elldata{1399}{$15$}{$D_{8}\,+\,A_{2}\,+\,5\,A_{1}$}{$[2]$}
\elldata{1400}{$15$}{$D_{8}\,+\,7\,A_{1}$}{$[2, 2]$}
\elldata{1401}{$15$}{$2\,D_{7}\,+\,A_{1}$}{$[1]$}
\elldata{1402}{$15$}{$D_{7}\,+\,D_{6}\,+\,A_{2}$}{$[1]$}
\elldata{1403}{$15$}{$D_{7}\,+\,D_{6}\,+\,2\,A_{1}$}{$[1]$}
\elldata{1404}{$15$}{$D_{7}\,+\,D_{5}\,+\,A_{3}$}{$[1]$}
\elldata{1405}{$15$}{$D_{7}\,+\,D_{5}\,+\,A_{2}\,+\,A_{1}$}{$[1]$}
\elldata{1406}{$15$}{$D_{7}\,+\,D_{5}\,+\,3\,A_{1}$}{$[1]$}
\elldata{1407}{$15$}{$D_{7}\,+\,2\,D_{4}$}{$[1]$}
\elldata{1408}{$15$}{$D_{7}\,+\,D_{4}\,+\,A_{4}$}{$[1]$}
\elldata{1409}{$15$}{$D_{7}\,+\,D_{4}\,+\,A_{3}\,+\,A_{1}$}{$[1]$}
\elldata{1410}{$15$}{$D_{7}\,+\,D_{4}\,+\,2\,A_{2}$}{$[1]$}
\elldata{1411}{$15$}{$D_{7}\,+\,D_{4}\,+\,A_{2}\,+\,2\,A_{1}$}{$[1]$}
\elldata{1412}{$15$}{$D_{7}\,+\,D_{4}\,+\,4\,A_{1}$}{$[2]$}
\elldata{1413}{$15$}{$D_{7}\,+\,A_{8}$}{$[1]$}
\elldata{1414}{$15$}{$D_{7}\,+\,A_{7}\,+\,A_{1}$}{$[1]$}
\elldata{1415}{$15$}{$D_{7}\,+\,A_{6}\,+\,A_{2}$}{$[1]$}
\elldata{1416}{$15$}{$D_{7}\,+\,A_{6}\,+\,2\,A_{1}$}{$[1]$}
\elldata{1417}{$15$}{$D_{7}\,+\,A_{5}\,+\,A_{3}$}{$[1]$}
\elldata{1418}{$15$}{$D_{7}\,+\,A_{5}\,+\,A_{2}\,+\,A_{1}$}{$[1]$}
\elldata{1419}{$15$}{$D_{7}\,+\,A_{5}\,+\,3\,A_{1}$}{$[2], \,[1]$}
\elldata{1420}{$15$}{$D_{7}\,+\,2\,A_{4}$}{$[1]$}
\elldata{1421}{$15$}{$D_{7}\,+\,A_{4}\,+\,A_{3}\,+\,A_{1}$}{$[1]$}
\elldata{1422}{$15$}{$D_{7}\,+\,A_{4}\,+\,2\,A_{2}$}{$[1]$}
\elldata{1423}{$15$}{$D_{7}\,+\,A_{4}\,+\,A_{2}\,+\,2\,A_{1}$}{$[1]$}
\elldata{1424}{$15$}{$D_{7}\,+\,A_{4}\,+\,4\,A_{1}$}{$[1]$}
\elldata{1425}{$15$}{$D_{7}\,+\,2\,A_{3}\,+\,A_{2}$}{$[1]$}
\elldata{1426}{$15$}{$D_{7}\,+\,2\,A_{3}\,+\,2\,A_{1}$}{$[2], \,[1]$}
\elldata{1427}{$15$}{$D_{7}\,+\,A_{3}\,+\,2\,A_{2}\,+\,A_{1}$}{$[1]$}
\elldata{1428}{$15$}{$D_{7}\,+\,A_{3}\,+\,A_{2}\,+\,3\,A_{1}$}{$[1]$}
\elldata{1429}{$15$}{$D_{7}\,+\,A_{3}\,+\,5\,A_{1}$}{$[2]$}
\elldata{1430}{$15$}{$D_{7}\,+\,4\,A_{2}$}{$[1]$}
\elldata{1431}{$15$}{$D_{7}\,+\,3\,A_{2}\,+\,2\,A_{1}$}{$[1]$}
\elldata{1432}{$15$}{$D_{7}\,+\,2\,A_{2}\,+\,4\,A_{1}$}{$[1]$}
\elldata{1433}{$15$}{$D_{7}\,+\,A_{2}\,+\,6\,A_{1}$}{$[2]$}
\elldata{1434}{$15$}{$2\,D_{6}\,+\,A_{3}$}{$[2], \,[1]$}
\elldata{1435}{$15$}{$2\,D_{6}\,+\,A_{2}\,+\,A_{1}$}{$[1]$}
\elldata{1436}{$15$}{$2\,D_{6}\,+\,3\,A_{1}$}{$[2]$}
\elldata{1437}{$15$}{$D_{6}\,+\,D_{5}\,+\,D_{4}$}{$[1]$}
\elldata{1438}{$15$}{$D_{6}\,+\,D_{5}\,+\,A_{4}$}{$[1]$}
\elldata{1439}{$15$}{$D_{6}\,+\,D_{5}\,+\,A_{3}\,+\,A_{1}$}{$[2], \,[1]$}
\elldata{1440}{$15$}{$D_{6}\,+\,D_{5}\,+\,2\,A_{2}$}{$[1]$}
\elldata{1441}{$15$}{$D_{6}\,+\,D_{5}\,+\,A_{2}\,+\,2\,A_{1}$}{$[1]$}
\elldata{1442}{$15$}{$D_{6}\,+\,D_{5}\,+\,4\,A_{1}$}{$[2]$}
\elldata{1443}{$15$}{$D_{6}\,+\,2\,D_{4}\,+\,A_{1}$}{$[2]$}
\elldata{1444}{$15$}{$D_{6}\,+\,D_{4}\,+\,A_{5}$}{$[2], \,[1]$}
\elldata{1445}{$15$}{$D_{6}\,+\,D_{4}\,+\,A_{4}\,+\,A_{1}$}{$[1]$}
\elldata{1446}{$15$}{$D_{6}\,+\,D_{4}\,+\,A_{3}\,+\,A_{2}$}{$[1]$}
\elldata{1447}{$15$}{$D_{6}\,+\,D_{4}\,+\,A_{3}\,+\,2\,A_{1}$}{$[2]$}
\elldata{1448}{$15$}{$D_{6}\,+\,D_{4}\,+\,2\,A_{2}\,+\,A_{1}$}{$[1]$}
\elldata{1449}{$15$}{$D_{6}\,+\,D_{4}\,+\,A_{2}\,+\,3\,A_{1}$}{$[2]$}
\elldata{1450}{$15$}{$D_{6}\,+\,D_{4}\,+\,5\,A_{1}$}{$[2, 2]$}
\elldata{1451}{$15$}{$D_{6}\,+\,A_{9}$}{$[2], \,[1]$}
\elldata{1452}{$15$}{$D_{6}\,+\,A_{8}\,+\,A_{1}$}{$[1]$}
\elldata{1453}{$15$}{$D_{6}\,+\,A_{7}\,+\,A_{2}$}{$[1]$}
\elldata{1454}{$15$}{$D_{6}\,+\,A_{7}\,+\,2\,A_{1}$}{$[2], \,[1]$}
\elldata{1455}{$15$}{$D_{6}\,+\,A_{6}\,+\,A_{3}$}{$[1]$}
\elldata{1456}{$15$}{$D_{6}\,+\,A_{6}\,+\,A_{2}\,+\,A_{1}$}{$[1]$}
\elldata{1457}{$15$}{$D_{6}\,+\,A_{6}\,+\,3\,A_{1}$}{$[1]$}
\elldata{1458}{$15$}{$D_{6}\,+\,A_{5}\,+\,A_{4}$}{$[1]$}
\elldata{1459}{$15$}{$D_{6}\,+\,A_{5}\,+\,A_{3}\,+\,A_{1}$}{$[2], \,[1]$}
\elldata{1460}{$15$}{$D_{6}\,+\,A_{5}\,+\,2\,A_{2}$}{$[1]$}
\elldata{1461}{$15$}{$D_{6}\,+\,A_{5}\,+\,A_{2}\,+\,2\,A_{1}$}{$[2], \,[1]$}
\elldata{1462}{$15$}{$D_{6}\,+\,A_{5}\,+\,4\,A_{1}$}{$[2]$}
\elldata{1463}{$15$}{$D_{6}\,+\,2\,A_{4}\,+\,A_{1}$}{$[1]$}
\elldata{1464}{$15$}{$D_{6}\,+\,A_{4}\,+\,A_{3}\,+\,A_{2}$}{$[1]$}
\elldata{1465}{$15$}{$D_{6}\,+\,A_{4}\,+\,A_{3}\,+\,2\,A_{1}$}{$[1]$}
\elldata{1466}{$15$}{$D_{6}\,+\,A_{4}\,+\,2\,A_{2}\,+\,A_{1}$}{$[1]$}
\elldata{1467}{$15$}{$D_{6}\,+\,A_{4}\,+\,A_{2}\,+\,3\,A_{1}$}{$[1]$}
\elldata{1468}{$15$}{$D_{6}\,+\,A_{4}\,+\,5\,A_{1}$}{$[2]$}
\elldata{1469}{$15$}{$D_{6}\,+\,3\,A_{3}$}{$[2], \,[1]$}
\elldata{1470}{$15$}{$D_{6}\,+\,2\,A_{3}\,+\,A_{2}\,+\,A_{1}$}{$[2], \,[1]$}
\elldata{1471}{$15$}{$D_{6}\,+\,2\,A_{3}\,+\,3\,A_{1}$}{$[2]$}
\elldata{1472}{$15$}{$D_{6}\,+\,A_{3}\,+\,3\,A_{2}$}{$[1]$}
\elldata{1473}{$15$}{$D_{6}\,+\,A_{3}\,+\,2\,A_{2}\,+\,2\,A_{1}$}{$[1]$}
\elldata{1474}{$15$}{$D_{6}\,+\,A_{3}\,+\,A_{2}\,+\,4\,A_{1}$}{$[2]$}
\elldata{1475}{$15$}{$D_{6}\,+\,A_{3}\,+\,6\,A_{1}$}{$[2, 2]$}
\elldata{1476}{$15$}{$D_{6}\,+\,4\,A_{2}\,+\,A_{1}$}{$[1]$}
\elldata{1477}{$15$}{$D_{6}\,+\,3\,A_{2}\,+\,3\,A_{1}$}{$[1]$}
\elldata{1478}{$15$}{$D_{6}\,+\,2\,A_{2}\,+\,5\,A_{1}$}{$[2]$}
\elldata{1479}{$15$}{$3\,D_{5}$}{$[1]$}
\elldata{1480}{$15$}{$2\,D_{5}\,+\,D_{4}\,+\,A_{1}$}{$[1]$}
\elldata{1481}{$15$}{$2\,D_{5}\,+\,A_{5}$}{$[1]$}
\elldata{1482}{$15$}{$2\,D_{5}\,+\,A_{4}\,+\,A_{1}$}{$[1]$}
\elldata{1483}{$15$}{$2\,D_{5}\,+\,A_{3}\,+\,A_{2}$}{$[1]$}
\elldata{1484}{$15$}{$2\,D_{5}\,+\,A_{3}\,+\,2\,A_{1}$}{$[2], \,[1]$}
\elldata{1485}{$15$}{$2\,D_{5}\,+\,2\,A_{2}\,+\,A_{1}$}{$[1]$}
\elldata{1486}{$15$}{$2\,D_{5}\,+\,A_{2}\,+\,3\,A_{1}$}{$[1]$}
\elldata{1487}{$15$}{$2\,D_{5}\,+\,5\,A_{1}$}{$[2]$}
\elldata{1488}{$15$}{$D_{5}\,+\,2\,D_{4}\,+\,A_{2}$}{$[1]$}
\elldata{1489}{$15$}{$D_{5}\,+\,2\,D_{4}\,+\,2\,A_{1}$}{$[2]$}
\elldata{1490}{$15$}{$D_{5}\,+\,D_{4}\,+\,A_{6}$}{$[1]$}
\elldata{1491}{$15$}{$D_{5}\,+\,D_{4}\,+\,A_{5}\,+\,A_{1}$}{$[2], \,[1]$}
\elldata{1492}{$15$}{$D_{5}\,+\,D_{4}\,+\,A_{4}\,+\,A_{2}$}{$[1]$}
\elldata{1493}{$15$}{$D_{5}\,+\,D_{4}\,+\,A_{4}\,+\,2\,A_{1}$}{$[1]$}
\elldata{1494}{$15$}{$D_{5}\,+\,D_{4}\,+\,2\,A_{3}$}{$[2], \,[1]$}
\elldata{1495}{$15$}{$D_{5}\,+\,D_{4}\,+\,A_{3}\,+\,A_{2}\,+\,A_{1}$}{$[1]$}
\elldata{1496}{$15$}{$D_{5}\,+\,D_{4}\,+\,A_{3}\,+\,3\,A_{1}$}{$[2]$}
\elldata{1497}{$15$}{$D_{5}\,+\,D_{4}\,+\,3\,A_{2}$}{$[1]$}
\elldata{1498}{$15$}{$D_{5}\,+\,D_{4}\,+\,2\,A_{2}\,+\,2\,A_{1}$}{$[1]$}
\elldata{1499}{$15$}{$D_{5}\,+\,D_{4}\,+\,A_{2}\,+\,4\,A_{1}$}{$[2]$}
\elldata{1500}{$15$}{$D_{5}\,+\,A_{10}$}{$[1]$}
\elldata{1501}{$15$}{$D_{5}\,+\,A_{9}\,+\,A_{1}$}{$[2], \,[1]$}
\elldata{1502}{$15$}{$D_{5}\,+\,A_{8}\,+\,A_{2}$}{$[1]$}
\elldata{1503}{$15$}{$D_{5}\,+\,A_{8}\,+\,2\,A_{1}$}{$[1]$}
\elldata{1504}{$15$}{$D_{5}\,+\,A_{7}\,+\,A_{3}$}{$[2], \,[1]$}
\elldata{1505}{$15$}{$D_{5}\,+\,A_{7}\,+\,A_{2}\,+\,A_{1}$}{$[1]$}
\elldata{1506}{$15$}{$D_{5}\,+\,A_{7}\,+\,3\,A_{1}$}{$[2], \,[1]$}
\elldata{1507}{$15$}{$D_{5}\,+\,A_{6}\,+\,A_{4}$}{$[1]$}
\elldata{1508}{$15$}{$D_{5}\,+\,A_{6}\,+\,A_{3}\,+\,A_{1}$}{$[1]$}
\elldata{1509}{$15$}{$D_{5}\,+\,A_{6}\,+\,2\,A_{2}$}{$[1]$}
\elldata{1510}{$15$}{$D_{5}\,+\,A_{6}\,+\,A_{2}\,+\,2\,A_{1}$}{$[1]$}
\elldata{1511}{$15$}{$D_{5}\,+\,A_{6}\,+\,4\,A_{1}$}{$[1]$}
\elldata{1512}{$15$}{$D_{5}\,+\,2\,A_{5}$}{$[2], \,[1]$}
\elldata{1513}{$15$}{$D_{5}\,+\,A_{5}\,+\,A_{4}\,+\,A_{1}$}{$[1]$}
\elldata{1514}{$15$}{$D_{5}\,+\,A_{5}\,+\,A_{3}\,+\,A_{2}$}{$[1]$}
\elldata{1515}{$15$}{$D_{5}\,+\,A_{5}\,+\,A_{3}\,+\,2\,A_{1}$}{$[2], \,[1]$}
\elldata{1516}{$15$}{$D_{5}\,+\,A_{5}\,+\,2\,A_{2}\,+\,A_{1}$}{$[1]$}
\elldata{1517}{$15$}{$D_{5}\,+\,A_{5}\,+\,A_{2}\,+\,3\,A_{1}$}{$[2], \,[1]$}
\elldata{1518}{$15$}{$D_{5}\,+\,A_{5}\,+\,5\,A_{1}$}{$[2]$}
\elldata{1519}{$15$}{$D_{5}\,+\,2\,A_{4}\,+\,A_{2}$}{$[1]$}
\elldata{1520}{$15$}{$D_{5}\,+\,2\,A_{4}\,+\,2\,A_{1}$}{$[1]$}
\elldata{1521}{$15$}{$D_{5}\,+\,A_{4}\,+\,2\,A_{3}$}{$[1]$}
\elldata{1522}{$15$}{$D_{5}\,+\,A_{4}\,+\,A_{3}\,+\,A_{2}\,+\,A_{1}$}{$[1]$}
\elldata{1523}{$15$}{$D_{5}\,+\,A_{4}\,+\,A_{3}\,+\,3\,A_{1}$}{$[1]$}
\elldata{1524}{$15$}{$D_{5}\,+\,A_{4}\,+\,3\,A_{2}$}{$[1]$}
\elldata{1525}{$15$}{$D_{5}\,+\,A_{4}\,+\,2\,A_{2}\,+\,2\,A_{1}$}{$[1]$}
\elldata{1526}{$15$}{$D_{5}\,+\,A_{4}\,+\,A_{2}\,+\,4\,A_{1}$}{$[1]$}
\elldata{1527}{$15$}{$D_{5}\,+\,A_{4}\,+\,6\,A_{1}$}{$[2]$}
\elldata{1528}{$15$}{$D_{5}\,+\,3\,A_{3}\,+\,A_{1}$}{$[4], \,[2], \,[1]$}
\elldata{1529}{$15$}{$D_{5}\,+\,2\,A_{3}\,+\,2\,A_{2}$}{$[1]$}
\elldata{1530}{$15$}{$D_{5}\,+\,2\,A_{3}\,+\,A_{2}\,+\,2\,A_{1}$}{$[2], \,[1]$}
\elldata{1531}{$15$}{$D_{5}\,+\,2\,A_{3}\,+\,4\,A_{1}$}{$[2]$}
\elldata{1532}{$15$}{$D_{5}\,+\,A_{3}\,+\,3\,A_{2}\,+\,A_{1}$}{$[1]$}
\elldata{1533}{$15$}{$D_{5}\,+\,A_{3}\,+\,2\,A_{2}\,+\,3\,A_{1}$}{$[1]$}
\elldata{1534}{$15$}{$D_{5}\,+\,A_{3}\,+\,A_{2}\,+\,5\,A_{1}$}{$[2]$}
\elldata{1535}{$15$}{$D_{5}\,+\,5\,A_{2}$}{$[1]$}
\elldata{1536}{$15$}{$D_{5}\,+\,4\,A_{2}\,+\,2\,A_{1}$}{$[1]$}
\elldata{1537}{$15$}{$D_{5}\,+\,3\,A_{2}\,+\,4\,A_{1}$}{$[1]$}
\elldata{1538}{$15$}{$3\,D_{4}\,+\,A_{3}$}{$[2]$}
\elldata{1539}{$15$}{$3\,D_{4}\,+\,3\,A_{1}$}{$[2, 2]$}
\elldata{1540}{$15$}{$2\,D_{4}\,+\,A_{7}$}{$[2], \,[1]$}
\elldata{1541}{$15$}{$2\,D_{4}\,+\,A_{6}\,+\,A_{1}$}{$[1]$}
\elldata{1542}{$15$}{$2\,D_{4}\,+\,A_{5}\,+\,A_{2}$}{$[1]$}
\elldata{1543}{$15$}{$2\,D_{4}\,+\,A_{5}\,+\,2\,A_{1}$}{$[2]$}
\elldata{1544}{$15$}{$2\,D_{4}\,+\,A_{4}\,+\,A_{3}$}{$[1]$}
\elldata{1545}{$15$}{$2\,D_{4}\,+\,A_{4}\,+\,A_{2}\,+\,A_{1}$}{$[1]$}
\elldata{1546}{$15$}{$2\,D_{4}\,+\,2\,A_{3}\,+\,A_{1}$}{$[2]$}
\elldata{1547}{$15$}{$2\,D_{4}\,+\,A_{3}\,+\,A_{2}\,+\,2\,A_{1}$}{$[2]$}
\elldata{1548}{$15$}{$2\,D_{4}\,+\,A_{3}\,+\,4\,A_{1}$}{$[2, 2]$}
\elldata{1549}{$15$}{$2\,D_{4}\,+\,3\,A_{2}\,+\,A_{1}$}{$[1]$}
\elldata{1550}{$15$}{$D_{4}\,+\,A_{11}$}{$[2], \,[1]$}
\elldata{1551}{$15$}{$D_{4}\,+\,A_{10}\,+\,A_{1}$}{$[1]$}
\elldata{1552}{$15$}{$D_{4}\,+\,A_{9}\,+\,A_{2}$}{$[1]$}
\elldata{1553}{$15$}{$D_{4}\,+\,A_{9}\,+\,2\,A_{1}$}{$[2], \,[1]$}
\elldata{1554}{$15$}{$D_{4}\,+\,A_{8}\,+\,A_{3}$}{$[1]$}
\elldata{1555}{$15$}{$D_{4}\,+\,A_{8}\,+\,A_{2}\,+\,A_{1}$}{$[1]$}
\elldata{1556}{$15$}{$D_{4}\,+\,A_{8}\,+\,3\,A_{1}$}{$[1]$}
\elldata{1557}{$15$}{$D_{4}\,+\,A_{7}\,+\,A_{4}$}{$[1]$}
\elldata{1558}{$15$}{$D_{4}\,+\,A_{7}\,+\,A_{3}\,+\,A_{1}$}{$[2], \,[1]$}
\elldata{1559}{$15$}{$D_{4}\,+\,A_{7}\,+\,2\,A_{2}$}{$[1]$}
\elldata{1560}{$15$}{$D_{4}\,+\,A_{7}\,+\,A_{2}\,+\,2\,A_{1}$}{$[2], \,[1]$}
\elldata{1561}{$15$}{$D_{4}\,+\,A_{7}\,+\,4\,A_{1}$}{$[2]$}
\elldata{1562}{$15$}{$D_{4}\,+\,A_{6}\,+\,A_{5}$}{$[1]$}
\elldata{1563}{$15$}{$D_{4}\,+\,A_{6}\,+\,A_{4}\,+\,A_{1}$}{$[1]$}
\elldata{1564}{$15$}{$D_{4}\,+\,A_{6}\,+\,A_{3}\,+\,A_{2}$}{$[1]$}
\elldata{1565}{$15$}{$D_{4}\,+\,A_{6}\,+\,A_{3}\,+\,2\,A_{1}$}{$[1]$}
\elldata{1566}{$15$}{$D_{4}\,+\,A_{6}\,+\,2\,A_{2}\,+\,A_{1}$}{$[1]$}
\elldata{1567}{$15$}{$D_{4}\,+\,A_{6}\,+\,A_{2}\,+\,3\,A_{1}$}{$[1]$}
\elldata{1568}{$15$}{$D_{4}\,+\,2\,A_{5}\,+\,A_{1}$}{$[2], \,[1]$}
\elldata{1569}{$15$}{$D_{4}\,+\,A_{5}\,+\,A_{4}\,+\,A_{2}$}{$[1]$}
\elldata{1570}{$15$}{$D_{4}\,+\,A_{5}\,+\,A_{4}\,+\,2\,A_{1}$}{$[1]$}
\elldata{1571}{$15$}{$D_{4}\,+\,A_{5}\,+\,2\,A_{3}$}{$[1]$}
\elldata{1572}{$15$}{$D_{4}\,+\,A_{5}\,+\,A_{3}\,+\,A_{2}\,+\,A_{1}$}{$[2], \,[1]$}
\elldata{1573}{$15$}{$D_{4}\,+\,A_{5}\,+\,A_{3}\,+\,3\,A_{1}$}{$[2]$}
\elldata{1574}{$15$}{$D_{4}\,+\,A_{5}\,+\,3\,A_{2}$}{$[1]$}
\elldata{1575}{$15$}{$D_{4}\,+\,A_{5}\,+\,2\,A_{2}\,+\,2\,A_{1}$}{$[1]$}
\elldata{1576}{$15$}{$D_{4}\,+\,A_{5}\,+\,A_{2}\,+\,4\,A_{1}$}{$[2]$}
\elldata{1577}{$15$}{$D_{4}\,+\,A_{5}\,+\,6\,A_{1}$}{$[2, 2]$}
\elldata{1578}{$15$}{$D_{4}\,+\,2\,A_{4}\,+\,A_{3}$}{$[1]$}
\elldata{1579}{$15$}{$D_{4}\,+\,2\,A_{4}\,+\,A_{2}\,+\,A_{1}$}{$[1]$}
\elldata{1580}{$15$}{$D_{4}\,+\,2\,A_{4}\,+\,3\,A_{1}$}{$[1]$}
\elldata{1581}{$15$}{$D_{4}\,+\,A_{4}\,+\,2\,A_{3}\,+\,A_{1}$}{$[1]$}
\elldata{1582}{$15$}{$D_{4}\,+\,A_{4}\,+\,A_{3}\,+\,2\,A_{2}$}{$[1]$}
\elldata{1583}{$15$}{$D_{4}\,+\,A_{4}\,+\,A_{3}\,+\,A_{2}\,+\,2\,A_{1}$}{$[1]$}
\elldata{1584}{$15$}{$D_{4}\,+\,A_{4}\,+\,A_{3}\,+\,4\,A_{1}$}{$[2]$}
\elldata{1585}{$15$}{$D_{4}\,+\,A_{4}\,+\,3\,A_{2}\,+\,A_{1}$}{$[1]$}
\elldata{1586}{$15$}{$D_{4}\,+\,A_{4}\,+\,2\,A_{2}\,+\,3\,A_{1}$}{$[1]$}
\elldata{1587}{$15$}{$D_{4}\,+\,3\,A_{3}\,+\,A_{2}$}{$[2]$}
\elldata{1588}{$15$}{$D_{4}\,+\,3\,A_{3}\,+\,2\,A_{1}$}{$[2]$}
\elldata{1589}{$15$}{$D_{4}\,+\,2\,A_{3}\,+\,2\,A_{2}\,+\,A_{1}$}{$[1]$}
\elldata{1590}{$15$}{$D_{4}\,+\,2\,A_{3}\,+\,A_{2}\,+\,3\,A_{1}$}{$[2]$}
\elldata{1591}{$15$}{$D_{4}\,+\,2\,A_{3}\,+\,5\,A_{1}$}{$[2, 2]$}
\elldata{1592}{$15$}{$D_{4}\,+\,A_{3}\,+\,4\,A_{2}$}{$[1]$}
\elldata{1593}{$15$}{$D_{4}\,+\,A_{3}\,+\,3\,A_{2}\,+\,2\,A_{1}$}{$[1]$}
\elldata{1594}{$15$}{$D_{4}\,+\,A_{3}\,+\,2\,A_{2}\,+\,4\,A_{1}$}{$[2]$}
\elldata{1595}{$15$}{$D_{4}\,+\,4\,A_{2}\,+\,3\,A_{1}$}{$[1]$}
\elldata{1596}{$15$}{$A_{15}$}{$[2], \,[1]$}
\elldata{1597}{$15$}{$A_{14}\,+\,A_{1}$}{$[1]$}
\elldata{1598}{$15$}{$A_{13}\,+\,A_{2}$}{$[1]$}
\elldata{1599}{$15$}{$A_{13}\,+\,2\,A_{1}$}{$[2], \,[1]$}
\elldata{1600}{$15$}{$A_{12}\,+\,A_{3}$}{$[1]$}
\elldata{1601}{$15$}{$A_{12}\,+\,A_{2}\,+\,A_{1}$}{$[1]$}
\elldata{1602}{$15$}{$A_{12}\,+\,3\,A_{1}$}{$[1]$}
\elldata{1603}{$15$}{$A_{11}\,+\,A_{4}$}{$[1]$}
\elldata{1604}{$15$}{$A_{11}\,+\,A_{3}\,+\,A_{1}$}{$[2], \,[1]$}
\elldata{1605}{$15$}{$A_{11}\,+\,2\,A_{2}$}{$[3], \,[1]$}
\elldata{1606}{$15$}{$A_{11}\,+\,A_{2}\,+\,2\,A_{1}$}{$[2], \,[1]$}
\elldata{1607}{$15$}{$A_{11}\,+\,4\,A_{1}$}{$[2], \,[1]$}
\elldata{1608}{$15$}{$A_{10}\,+\,A_{5}$}{$[1]$}
\elldata{1609}{$15$}{$A_{10}\,+\,A_{4}\,+\,A_{1}$}{$[1]$}
\elldata{1610}{$15$}{$A_{10}\,+\,A_{3}\,+\,A_{2}$}{$[1]$}
\elldata{1611}{$15$}{$A_{10}\,+\,A_{3}\,+\,2\,A_{1}$}{$[1]$}
\elldata{1612}{$15$}{$A_{10}\,+\,2\,A_{2}\,+\,A_{1}$}{$[1]$}
\elldata{1613}{$15$}{$A_{10}\,+\,A_{2}\,+\,3\,A_{1}$}{$[1]$}
\elldata{1614}{$15$}{$A_{10}\,+\,5\,A_{1}$}{$[1]$}
\elldata{1615}{$15$}{$A_{9}\,+\,A_{6}$}{$[1]$}
\elldata{1616}{$15$}{$A_{9}\,+\,A_{5}\,+\,A_{1}$}{$[2], \,[1]$}
\elldata{1617}{$15$}{$A_{9}\,+\,A_{4}\,+\,A_{2}$}{$[1]$}
\elldata{1618}{$15$}{$A_{9}\,+\,A_{4}\,+\,2\,A_{1}$}{$[1]$}
\elldata{1619}{$15$}{$A_{9}\,+\,2\,A_{3}$}{$[1]$}
\elldata{1620}{$15$}{$A_{9}\,+\,A_{3}\,+\,A_{2}\,+\,A_{1}$}{$[2], \,[1]$}
\elldata{1621}{$15$}{$A_{9}\,+\,A_{3}\,+\,3\,A_{1}$}{$[2], \,[1]$}
\elldata{1622}{$15$}{$A_{9}\,+\,3\,A_{2}$}{$[1]$}
\elldata{1623}{$15$}{$A_{9}\,+\,2\,A_{2}\,+\,2\,A_{1}$}{$[1]$}
\elldata{1624}{$15$}{$A_{9}\,+\,A_{2}\,+\,4\,A_{1}$}{$[2], \,[1]$}
\elldata{1625}{$15$}{$A_{9}\,+\,6\,A_{1}$}{$[2]$}
\elldata{1626}{$15$}{$A_{8}\,+\,A_{7}$}{$[1]$}
\elldata{1627}{$15$}{$A_{8}\,+\,A_{6}\,+\,A_{1}$}{$[1]$}
\elldata{1628}{$15$}{$A_{8}\,+\,A_{5}\,+\,A_{2}$}{$[3], \,[1]$}
\elldata{1629}{$15$}{$A_{8}\,+\,A_{5}\,+\,2\,A_{1}$}{$[1]$}
\elldata{1630}{$15$}{$A_{8}\,+\,A_{4}\,+\,A_{3}$}{$[1]$}
\elldata{1631}{$15$}{$A_{8}\,+\,A_{4}\,+\,A_{2}\,+\,A_{1}$}{$[1]$}
\elldata{1632}{$15$}{$A_{8}\,+\,A_{4}\,+\,3\,A_{1}$}{$[1]$}
\elldata{1633}{$15$}{$A_{8}\,+\,2\,A_{3}\,+\,A_{1}$}{$[1]$}
\elldata{1634}{$15$}{$A_{8}\,+\,A_{3}\,+\,2\,A_{2}$}{$[1]$}
\elldata{1635}{$15$}{$A_{8}\,+\,A_{3}\,+\,A_{2}\,+\,2\,A_{1}$}{$[1]$}
\elldata{1636}{$15$}{$A_{8}\,+\,A_{3}\,+\,4\,A_{1}$}{$[1]$}
\elldata{1637}{$15$}{$A_{8}\,+\,3\,A_{2}\,+\,A_{1}$}{$[3], \,[1]$}
\elldata{1638}{$15$}{$A_{8}\,+\,2\,A_{2}\,+\,3\,A_{1}$}{$[1]$}
\elldata{1639}{$15$}{$A_{8}\,+\,A_{2}\,+\,5\,A_{1}$}{$[1]$}
\elldata{1640}{$15$}{$2\,A_{7}\,+\,A_{1}$}{$[2], \,[1]$}
\elldata{1641}{$15$}{$A_{7}\,+\,A_{6}\,+\,A_{2}$}{$[1]$}
\elldata{1642}{$15$}{$A_{7}\,+\,A_{6}\,+\,2\,A_{1}$}{$[1]$}
\elldata{1643}{$15$}{$A_{7}\,+\,A_{5}\,+\,A_{3}$}{$[1]$}
\elldata{1644}{$15$}{$A_{7}\,+\,A_{5}\,+\,A_{2}\,+\,A_{1}$}{$[2], \,[1]$}
\elldata{1645}{$15$}{$A_{7}\,+\,A_{5}\,+\,3\,A_{1}$}{$[2], \,[1]$}
\elldata{1646}{$15$}{$A_{7}\,+\,2\,A_{4}$}{$[1]$}
\elldata{1647}{$15$}{$A_{7}\,+\,A_{4}\,+\,A_{3}\,+\,A_{1}$}{$[1]$}
\elldata{1648}{$15$}{$A_{7}\,+\,A_{4}\,+\,2\,A_{2}$}{$[1]$}
\elldata{1649}{$15$}{$A_{7}\,+\,A_{4}\,+\,A_{2}\,+\,2\,A_{1}$}{$[1]$}
\elldata{1650}{$15$}{$A_{7}\,+\,A_{4}\,+\,4\,A_{1}$}{$[2], \,[1]$}
\elldata{1651}{$15$}{$A_{7}\,+\,2\,A_{3}\,+\,A_{2}$}{$[2], \,[1]$}
\elldata{1652}{$15$}{$A_{7}\,+\,2\,A_{3}\,+\,2\,A_{1}$}{$[4], \,[2], \,[1]$}
\elldata{1653}{$15$}{$A_{7}\,+\,A_{3}\,+\,2\,A_{2}\,+\,A_{1}$}{$[1]$}
\elldata{1654}{$15$}{$A_{7}\,+\,A_{3}\,+\,A_{2}\,+\,3\,A_{1}$}{$[2], \,[1]$}
\elldata{1655}{$15$}{$A_{7}\,+\,A_{3}\,+\,5\,A_{1}$}{$[2]$}
\elldata{1656}{$15$}{$A_{7}\,+\,4\,A_{2}$}{$[1]$}
\elldata{1657}{$15$}{$A_{7}\,+\,3\,A_{2}\,+\,2\,A_{1}$}{$[1]$}
\elldata{1658}{$15$}{$A_{7}\,+\,2\,A_{2}\,+\,4\,A_{1}$}{$[2], \,[1]$}
\elldata{1659}{$15$}{$A_{7}\,+\,A_{2}\,+\,6\,A_{1}$}{$[2]$}
\elldata{1660}{$15$}{$A_{7}\,+\,8\,A_{1}$}{$[2, 2]$}
\elldata{1661}{$15$}{$2\,A_{6}\,+\,A_{3}$}{$[1]$}
\elldata{1662}{$15$}{$2\,A_{6}\,+\,A_{2}\,+\,A_{1}$}{$[1]$}
\elldata{1663}{$15$}{$2\,A_{6}\,+\,3\,A_{1}$}{$[1]$}
\elldata{1664}{$15$}{$A_{6}\,+\,A_{5}\,+\,A_{4}$}{$[1]$}
\elldata{1665}{$15$}{$A_{6}\,+\,A_{5}\,+\,A_{3}\,+\,A_{1}$}{$[1]$}
\elldata{1666}{$15$}{$A_{6}\,+\,A_{5}\,+\,2\,A_{2}$}{$[1]$}
\elldata{1667}{$15$}{$A_{6}\,+\,A_{5}\,+\,A_{2}\,+\,2\,A_{1}$}{$[1]$}
\elldata{1668}{$15$}{$A_{6}\,+\,A_{5}\,+\,4\,A_{1}$}{$[1]$}
\elldata{1669}{$15$}{$A_{6}\,+\,2\,A_{4}\,+\,A_{1}$}{$[1]$}
\elldata{1670}{$15$}{$A_{6}\,+\,A_{4}\,+\,A_{3}\,+\,A_{2}$}{$[1]$}
\elldata{1671}{$15$}{$A_{6}\,+\,A_{4}\,+\,A_{3}\,+\,2\,A_{1}$}{$[1]$}
\elldata{1672}{$15$}{$A_{6}\,+\,A_{4}\,+\,2\,A_{2}\,+\,A_{1}$}{$[1]$}
\elldata{1673}{$15$}{$A_{6}\,+\,A_{4}\,+\,A_{2}\,+\,3\,A_{1}$}{$[1]$}
\elldata{1674}{$15$}{$A_{6}\,+\,A_{4}\,+\,5\,A_{1}$}{$[1]$}
\elldata{1675}{$15$}{$A_{6}\,+\,3\,A_{3}$}{$[1]$}
\elldata{1676}{$15$}{$A_{6}\,+\,2\,A_{3}\,+\,A_{2}\,+\,A_{1}$}{$[1]$}
\elldata{1677}{$15$}{$A_{6}\,+\,2\,A_{3}\,+\,3\,A_{1}$}{$[1]$}
\elldata{1678}{$15$}{$A_{6}\,+\,A_{3}\,+\,3\,A_{2}$}{$[1]$}
\elldata{1679}{$15$}{$A_{6}\,+\,A_{3}\,+\,2\,A_{2}\,+\,2\,A_{1}$}{$[1]$}
\elldata{1680}{$15$}{$A_{6}\,+\,A_{3}\,+\,A_{2}\,+\,4\,A_{1}$}{$[1]$}
\elldata{1681}{$15$}{$A_{6}\,+\,A_{3}\,+\,6\,A_{1}$}{$[2]$}
\elldata{1682}{$15$}{$A_{6}\,+\,4\,A_{2}\,+\,A_{1}$}{$[1]$}
\elldata{1683}{$15$}{$A_{6}\,+\,3\,A_{2}\,+\,3\,A_{1}$}{$[1]$}
\elldata{1684}{$15$}{$A_{6}\,+\,2\,A_{2}\,+\,5\,A_{1}$}{$[1]$}
\elldata{1685}{$15$}{$3\,A_{5}$}{$[3], \,[1]$}
\elldata{1686}{$15$}{$2\,A_{5}\,+\,A_{4}\,+\,A_{1}$}{$[1]$}
\elldata{1687}{$15$}{$2\,A_{5}\,+\,A_{3}\,+\,A_{2}$}{$[2], \,[1]$}
\elldata{1688}{$15$}{$2\,A_{5}\,+\,A_{3}\,+\,2\,A_{1}$}{$[2], \,[1]$}
\elldata{1689}{$15$}{$2\,A_{5}\,+\,2\,A_{2}\,+\,A_{1}$}{$[3], \,[1]$}
\elldata{1690}{$15$}{$2\,A_{5}\,+\,A_{2}\,+\,3\,A_{1}$}{$[2], \,[1]$}
\elldata{1691}{$15$}{$2\,A_{5}\,+\,5\,A_{1}$}{$[2]$}
\elldata{1692}{$15$}{$A_{5}\,+\,2\,A_{4}\,+\,A_{2}$}{$[1]$}
\elldata{1693}{$15$}{$A_{5}\,+\,2\,A_{4}\,+\,2\,A_{1}$}{$[1]$}
\elldata{1694}{$15$}{$A_{5}\,+\,A_{4}\,+\,2\,A_{3}$}{$[1]$}
\elldata{1695}{$15$}{$A_{5}\,+\,A_{4}\,+\,A_{3}\,+\,A_{2}\,+\,A_{1}$}{$[1]$}
\elldata{1696}{$15$}{$A_{5}\,+\,A_{4}\,+\,A_{3}\,+\,3\,A_{1}$}{$[2], \,[1]$}
\elldata{1697}{$15$}{$A_{5}\,+\,A_{4}\,+\,3\,A_{2}$}{$[1]$}
\elldata{1698}{$15$}{$A_{5}\,+\,A_{4}\,+\,2\,A_{2}\,+\,2\,A_{1}$}{$[1]$}
\elldata{1699}{$15$}{$A_{5}\,+\,A_{4}\,+\,A_{2}\,+\,4\,A_{1}$}{$[1]$}
\elldata{1700}{$15$}{$A_{5}\,+\,A_{4}\,+\,6\,A_{1}$}{$[2]$}
\elldata{1701}{$15$}{$A_{5}\,+\,3\,A_{3}\,+\,A_{1}$}{$[2], \,[1]$}
\elldata{1702}{$15$}{$A_{5}\,+\,2\,A_{3}\,+\,2\,A_{2}$}{$[1]$}
\elldata{1703}{$15$}{$A_{5}\,+\,2\,A_{3}\,+\,A_{2}\,+\,2\,A_{1}$}{$[2], \,[1]$}
\elldata{1704}{$15$}{$A_{5}\,+\,2\,A_{3}\,+\,4\,A_{1}$}{$[2]$}
\elldata{1705}{$15$}{$A_{5}\,+\,A_{3}\,+\,3\,A_{2}\,+\,A_{1}$}{$[1]$}
\elldata{1706}{$15$}{$A_{5}\,+\,A_{3}\,+\,2\,A_{2}\,+\,3\,A_{1}$}{$[2], \,[1]$}
\elldata{1707}{$15$}{$A_{5}\,+\,A_{3}\,+\,A_{2}\,+\,5\,A_{1}$}{$[2]$}
\elldata{1708}{$15$}{$A_{5}\,+\,A_{3}\,+\,7\,A_{1}$}{$[2, 2]$}
\elldata{1709}{$15$}{$A_{5}\,+\,5\,A_{2}$}{$[3]$}
\elldata{1710}{$15$}{$A_{5}\,+\,4\,A_{2}\,+\,2\,A_{1}$}{$[3], \,[1]$}
\elldata{1711}{$15$}{$A_{5}\,+\,3\,A_{2}\,+\,4\,A_{1}$}{$[1]$}
\elldata{1712}{$15$}{$A_{5}\,+\,2\,A_{2}\,+\,6\,A_{1}$}{$[2]$}
\elldata{1713}{$15$}{$3\,A_{4}\,+\,A_{3}$}{$[1]$}
\elldata{1714}{$15$}{$3\,A_{4}\,+\,A_{2}\,+\,A_{1}$}{$[1]$}
\elldata{1715}{$15$}{$3\,A_{4}\,+\,3\,A_{1}$}{$[1]$}
\elldata{1716}{$15$}{$2\,A_{4}\,+\,2\,A_{3}\,+\,A_{1}$}{$[1]$}
\elldata{1717}{$15$}{$2\,A_{4}\,+\,A_{3}\,+\,2\,A_{2}$}{$[1]$}
\elldata{1718}{$15$}{$2\,A_{4}\,+\,A_{3}\,+\,A_{2}\,+\,2\,A_{1}$}{$[1]$}
\elldata{1719}{$15$}{$2\,A_{4}\,+\,A_{3}\,+\,4\,A_{1}$}{$[1]$}
\elldata{1720}{$15$}{$2\,A_{4}\,+\,3\,A_{2}\,+\,A_{1}$}{$[1]$}
\elldata{1721}{$15$}{$2\,A_{4}\,+\,2\,A_{2}\,+\,3\,A_{1}$}{$[1]$}
\elldata{1722}{$15$}{$2\,A_{4}\,+\,A_{2}\,+\,5\,A_{1}$}{$[1]$}
\elldata{1723}{$15$}{$A_{4}\,+\,3\,A_{3}\,+\,A_{2}$}{$[1]$}
\elldata{1724}{$15$}{$A_{4}\,+\,3\,A_{3}\,+\,2\,A_{1}$}{$[2], \,[1]$}
\elldata{1725}{$15$}{$A_{4}\,+\,2\,A_{3}\,+\,2\,A_{2}\,+\,A_{1}$}{$[1]$}
\elldata{1726}{$15$}{$A_{4}\,+\,2\,A_{3}\,+\,A_{2}\,+\,3\,A_{1}$}{$[1]$}
\elldata{1727}{$15$}{$A_{4}\,+\,2\,A_{3}\,+\,5\,A_{1}$}{$[2]$}
\elldata{1728}{$15$}{$A_{4}\,+\,A_{3}\,+\,4\,A_{2}$}{$[1]$}
\elldata{1729}{$15$}{$A_{4}\,+\,A_{3}\,+\,3\,A_{2}\,+\,2\,A_{1}$}{$[1]$}
\elldata{1730}{$15$}{$A_{4}\,+\,A_{3}\,+\,2\,A_{2}\,+\,4\,A_{1}$}{$[1]$}
\elldata{1731}{$15$}{$A_{4}\,+\,A_{3}\,+\,A_{2}\,+\,6\,A_{1}$}{$[2]$}
\elldata{1732}{$15$}{$A_{4}\,+\,5\,A_{2}\,+\,A_{1}$}{$[1]$}
\elldata{1733}{$15$}{$A_{4}\,+\,4\,A_{2}\,+\,3\,A_{1}$}{$[1]$}
\elldata{1734}{$15$}{$A_{4}\,+\,3\,A_{2}\,+\,5\,A_{1}$}{$[1]$}
\elldata{1735}{$15$}{$5\,A_{3}$}{$[4]$}
\elldata{1736}{$15$}{$4\,A_{3}\,+\,A_{2}\,+\,A_{1}$}{$[2], \,[1]$}
\elldata{1737}{$15$}{$4\,A_{3}\,+\,3\,A_{1}$}{$[4], \,[2]$}
\elldata{1738}{$15$}{$3\,A_{3}\,+\,3\,A_{2}$}{$[1]$}
\elldata{1739}{$15$}{$3\,A_{3}\,+\,2\,A_{2}\,+\,2\,A_{1}$}{$[2], \,[1]$}
\elldata{1740}{$15$}{$3\,A_{3}\,+\,A_{2}\,+\,4\,A_{1}$}{$[2]$}
\elldata{1741}{$15$}{$3\,A_{3}\,+\,6\,A_{1}$}{$[2, 2]$}
\elldata{1742}{$15$}{$2\,A_{3}\,+\,4\,A_{2}\,+\,A_{1}$}{$[1]$}
\elldata{1743}{$15$}{$2\,A_{3}\,+\,3\,A_{2}\,+\,3\,A_{1}$}{$[1]$}
\elldata{1744}{$15$}{$2\,A_{3}\,+\,2\,A_{2}\,+\,5\,A_{1}$}{$[2]$}
\elldata{1745}{$15$}{$A_{3}\,+\,6\,A_{2}$}{$[3]$}
\elldata{1746}{$15$}{$A_{3}\,+\,5\,A_{2}\,+\,2\,A_{1}$}{$[1]$}
\elldata{1747}{$15$}{$A_{3}\,+\,4\,A_{2}\,+\,4\,A_{1}$}{$[1]$}
\elldata{1748}{$15$}{$7\,A_{2}\,+\,A_{1}$}{$[3]$}
\elldata{1749}{$15$}{$6\,A_{2}\,+\,3\,A_{1}$}{$[3]$}

\vsr \elldata{No.}{rank}{$ADE$-type}{$G$}

\vsrs \elldata{1750}{$16$}{$2\,E_{8}$}{$[1]$}
\elldata{1751}{$16$}{$E_{8}\,+\,E_{7}\,+\,A_{1}$}{$[1]$}
\elldata{1752}{$16$}{$E_{8}\,+\,E_{6}\,+\,A_{2}$}{$[1]$}
\elldata{1753}{$16$}{$E_{8}\,+\,E_{6}\,+\,2\,A_{1}$}{$[1]$}
\elldata{1754}{$16$}{$E_{8}\,+\,D_{8}$}{$[1]$}
\elldata{1755}{$16$}{$E_{8}\,+\,D_{7}\,+\,A_{1}$}{$[1]$}
\elldata{1756}{$16$}{$E_{8}\,+\,D_{6}\,+\,A_{2}$}{$[1]$}
\elldata{1757}{$16$}{$E_{8}\,+\,D_{6}\,+\,2\,A_{1}$}{$[1]$}
\elldata{1758}{$16$}{$E_{8}\,+\,D_{5}\,+\,A_{3}$}{$[1]$}
\elldata{1759}{$16$}{$E_{8}\,+\,D_{5}\,+\,A_{2}\,+\,A_{1}$}{$[1]$}
\elldata{1760}{$16$}{$E_{8}\,+\,D_{5}\,+\,3\,A_{1}$}{$[1]$}
\elldata{1761}{$16$}{$E_{8}\,+\,2\,D_{4}$}{$[1]$}
\elldata{1762}{$16$}{$E_{8}\,+\,D_{4}\,+\,A_{4}$}{$[1]$}
\elldata{1763}{$16$}{$E_{8}\,+\,D_{4}\,+\,A_{3}\,+\,A_{1}$}{$[1]$}
\elldata{1764}{$16$}{$E_{8}\,+\,D_{4}\,+\,2\,A_{2}$}{$[1]$}
\elldata{1765}{$16$}{$E_{8}\,+\,D_{4}\,+\,A_{2}\,+\,2\,A_{1}$}{$[1]$}
\elldata{1766}{$16$}{$E_{8}\,+\,A_{8}$}{$[1]$}
\elldata{1767}{$16$}{$E_{8}\,+\,A_{7}\,+\,A_{1}$}{$[1]$}
\elldata{1768}{$16$}{$E_{8}\,+\,A_{6}\,+\,A_{2}$}{$[1]$}
\elldata{1769}{$16$}{$E_{8}\,+\,A_{6}\,+\,2\,A_{1}$}{$[1]$}
\elldata{1770}{$16$}{$E_{8}\,+\,A_{5}\,+\,A_{3}$}{$[1]$}
\elldata{1771}{$16$}{$E_{8}\,+\,A_{5}\,+\,A_{2}\,+\,A_{1}$}{$[1]$}
\elldata{1772}{$16$}{$E_{8}\,+\,A_{5}\,+\,3\,A_{1}$}{$[1]$}
\elldata{1773}{$16$}{$E_{8}\,+\,2\,A_{4}$}{$[1]$}
\elldata{1774}{$16$}{$E_{8}\,+\,A_{4}\,+\,A_{3}\,+\,A_{1}$}{$[1]$}
\elldata{1775}{$16$}{$E_{8}\,+\,A_{4}\,+\,2\,A_{2}$}{$[1]$}
\elldata{1776}{$16$}{$E_{8}\,+\,A_{4}\,+\,A_{2}\,+\,2\,A_{1}$}{$[1]$}
\elldata{1777}{$16$}{$E_{8}\,+\,A_{4}\,+\,4\,A_{1}$}{$[1]$}
\elldata{1778}{$16$}{$E_{8}\,+\,2\,A_{3}\,+\,A_{2}$}{$[1]$}
\elldata{1779}{$16$}{$E_{8}\,+\,2\,A_{3}\,+\,2\,A_{1}$}{$[1]$}
\elldata{1780}{$16$}{$E_{8}\,+\,A_{3}\,+\,2\,A_{2}\,+\,A_{1}$}{$[1]$}
\elldata{1781}{$16$}{$E_{8}\,+\,A_{3}\,+\,A_{2}\,+\,3\,A_{1}$}{$[1]$}
\elldata{1782}{$16$}{$E_{8}\,+\,4\,A_{2}$}{$[1]$}
\elldata{1783}{$16$}{$E_{8}\,+\,3\,A_{2}\,+\,2\,A_{1}$}{$[1]$}
\elldata{1784}{$16$}{$E_{8}\,+\,2\,A_{2}\,+\,4\,A_{1}$}{$[1]$}
\elldata{1785}{$16$}{$2\,E_{7}\,+\,A_{2}$}{$[1]$}
\elldata{1786}{$16$}{$2\,E_{7}\,+\,2\,A_{1}$}{$[2], \,[1]$}
\elldata{1787}{$16$}{$E_{7}\,+\,E_{6}\,+\,A_{3}$}{$[1]$}
\elldata{1788}{$16$}{$E_{7}\,+\,E_{6}\,+\,A_{2}\,+\,A_{1}$}{$[1]$}
\elldata{1789}{$16$}{$E_{7}\,+\,E_{6}\,+\,3\,A_{1}$}{$[1]$}
\elldata{1790}{$16$}{$E_{7}\,+\,D_{9}$}{$[1]$}
\elldata{1791}{$16$}{$E_{7}\,+\,D_{8}\,+\,A_{1}$}{$[2], \,[1]$}
\elldata{1792}{$16$}{$E_{7}\,+\,D_{7}\,+\,A_{2}$}{$[1]$}
\elldata{1793}{$16$}{$E_{7}\,+\,D_{7}\,+\,2\,A_{1}$}{$[1]$}
\elldata{1794}{$16$}{$E_{7}\,+\,D_{6}\,+\,A_{3}$}{$[2], \,[1]$}
\elldata{1795}{$16$}{$E_{7}\,+\,D_{6}\,+\,A_{2}\,+\,A_{1}$}{$[1]$}
\elldata{1796}{$16$}{$E_{7}\,+\,D_{6}\,+\,3\,A_{1}$}{$[2]$}
\elldata{1797}{$16$}{$E_{7}\,+\,D_{5}\,+\,D_{4}$}{$[1]$}
\elldata{1798}{$16$}{$E_{7}\,+\,D_{5}\,+\,A_{4}$}{$[1]$}
\elldata{1799}{$16$}{$E_{7}\,+\,D_{5}\,+\,A_{3}\,+\,A_{1}$}{$[2], \,[1]$}
\elldata{1800}{$16$}{$E_{7}\,+\,D_{5}\,+\,2\,A_{2}$}{$[1]$}
\elldata{1801}{$16$}{$E_{7}\,+\,D_{5}\,+\,A_{2}\,+\,2\,A_{1}$}{$[1]$}
\elldata{1802}{$16$}{$E_{7}\,+\,D_{5}\,+\,4\,A_{1}$}{$[2]$}
\elldata{1803}{$16$}{$E_{7}\,+\,2\,D_{4}\,+\,A_{1}$}{$[2]$}
\elldata{1804}{$16$}{$E_{7}\,+\,D_{4}\,+\,A_{5}$}{$[2], \,[1]$}
\elldata{1805}{$16$}{$E_{7}\,+\,D_{4}\,+\,A_{4}\,+\,A_{1}$}{$[1]$}
\elldata{1806}{$16$}{$E_{7}\,+\,D_{4}\,+\,A_{3}\,+\,A_{2}$}{$[1]$}
\elldata{1807}{$16$}{$E_{7}\,+\,D_{4}\,+\,A_{3}\,+\,2\,A_{1}$}{$[2]$}
\elldata{1808}{$16$}{$E_{7}\,+\,D_{4}\,+\,2\,A_{2}\,+\,A_{1}$}{$[1]$}
\elldata{1809}{$16$}{$E_{7}\,+\,D_{4}\,+\,A_{2}\,+\,3\,A_{1}$}{$[2]$}
\elldata{1810}{$16$}{$E_{7}\,+\,A_{9}$}{$[2], \,[1]$}
\elldata{1811}{$16$}{$E_{7}\,+\,A_{8}\,+\,A_{1}$}{$[1]$}
\elldata{1812}{$16$}{$E_{7}\,+\,A_{7}\,+\,A_{2}$}{$[1]$}
\elldata{1813}{$16$}{$E_{7}\,+\,A_{7}\,+\,2\,A_{1}$}{$[2], \,[1]$}
\elldata{1814}{$16$}{$E_{7}\,+\,A_{6}\,+\,A_{3}$}{$[1]$}
\elldata{1815}{$16$}{$E_{7}\,+\,A_{6}\,+\,A_{2}\,+\,A_{1}$}{$[1]$}
\elldata{1816}{$16$}{$E_{7}\,+\,A_{6}\,+\,3\,A_{1}$}{$[1]$}
\elldata{1817}{$16$}{$E_{7}\,+\,A_{5}\,+\,A_{4}$}{$[1]$}
\elldata{1818}{$16$}{$E_{7}\,+\,A_{5}\,+\,A_{3}\,+\,A_{1}$}{$[2], \,[1]$}
\elldata{1819}{$16$}{$E_{7}\,+\,A_{5}\,+\,2\,A_{2}$}{$[1]$}
\elldata{1820}{$16$}{$E_{7}\,+\,A_{5}\,+\,A_{2}\,+\,2\,A_{1}$}{$[2], \,[1]$}
\elldata{1821}{$16$}{$E_{7}\,+\,A_{5}\,+\,4\,A_{1}$}{$[2]$}
\elldata{1822}{$16$}{$E_{7}\,+\,2\,A_{4}\,+\,A_{1}$}{$[1]$}
\elldata{1823}{$16$}{$E_{7}\,+\,A_{4}\,+\,A_{3}\,+\,A_{2}$}{$[1]$}
\elldata{1824}{$16$}{$E_{7}\,+\,A_{4}\,+\,A_{3}\,+\,2\,A_{1}$}{$[1]$}
\elldata{1825}{$16$}{$E_{7}\,+\,A_{4}\,+\,2\,A_{2}\,+\,A_{1}$}{$[1]$}
\elldata{1826}{$16$}{$E_{7}\,+\,A_{4}\,+\,A_{2}\,+\,3\,A_{1}$}{$[1]$}
\elldata{1827}{$16$}{$E_{7}\,+\,A_{4}\,+\,5\,A_{1}$}{$[2]$}
\elldata{1828}{$16$}{$E_{7}\,+\,3\,A_{3}$}{$[1]$}
\elldata{1829}{$16$}{$E_{7}\,+\,2\,A_{3}\,+\,A_{2}\,+\,A_{1}$}{$[2], \,[1]$}
\elldata{1830}{$16$}{$E_{7}\,+\,2\,A_{3}\,+\,3\,A_{1}$}{$[2]$}
\elldata{1831}{$16$}{$E_{7}\,+\,A_{3}\,+\,3\,A_{2}$}{$[1]$}
\elldata{1832}{$16$}{$E_{7}\,+\,A_{3}\,+\,2\,A_{2}\,+\,2\,A_{1}$}{$[1]$}
\elldata{1833}{$16$}{$E_{7}\,+\,A_{3}\,+\,A_{2}\,+\,4\,A_{1}$}{$[2]$}
\elldata{1834}{$16$}{$E_{7}\,+\,4\,A_{2}\,+\,A_{1}$}{$[1]$}
\elldata{1835}{$16$}{$E_{7}\,+\,3\,A_{2}\,+\,3\,A_{1}$}{$[1]$}
\elldata{1836}{$16$}{$2\,E_{6}\,+\,D_{4}$}{$[1]$}
\elldata{1837}{$16$}{$2\,E_{6}\,+\,A_{4}$}{$[1]$}
\elldata{1838}{$16$}{$2\,E_{6}\,+\,A_{3}\,+\,A_{1}$}{$[1]$}
\elldata{1839}{$16$}{$2\,E_{6}\,+\,2\,A_{2}$}{$[3], \,[1]$}
\elldata{1840}{$16$}{$2\,E_{6}\,+\,A_{2}\,+\,2\,A_{1}$}{$[1]$}
\elldata{1841}{$16$}{$2\,E_{6}\,+\,4\,A_{1}$}{$[1]$}
\elldata{1842}{$16$}{$E_{6}\,+\,D_{10}$}{$[1]$}
\elldata{1843}{$16$}{$E_{6}\,+\,D_{9}\,+\,A_{1}$}{$[1]$}
\elldata{1844}{$16$}{$E_{6}\,+\,D_{8}\,+\,A_{2}$}{$[1]$}
\elldata{1845}{$16$}{$E_{6}\,+\,D_{8}\,+\,2\,A_{1}$}{$[1]$}
\elldata{1846}{$16$}{$E_{6}\,+\,D_{7}\,+\,A_{3}$}{$[1]$}
\elldata{1847}{$16$}{$E_{6}\,+\,D_{7}\,+\,A_{2}\,+\,A_{1}$}{$[1]$}
\elldata{1848}{$16$}{$E_{6}\,+\,D_{7}\,+\,3\,A_{1}$}{$[1]$}
\elldata{1849}{$16$}{$E_{6}\,+\,D_{6}\,+\,D_{4}$}{$[1]$}
\elldata{1850}{$16$}{$E_{6}\,+\,D_{6}\,+\,A_{4}$}{$[1]$}
\elldata{1851}{$16$}{$E_{6}\,+\,D_{6}\,+\,A_{3}\,+\,A_{1}$}{$[1]$}
\elldata{1852}{$16$}{$E_{6}\,+\,D_{6}\,+\,2\,A_{2}$}{$[1]$}
\elldata{1853}{$16$}{$E_{6}\,+\,D_{6}\,+\,A_{2}\,+\,2\,A_{1}$}{$[1]$}
\elldata{1854}{$16$}{$E_{6}\,+\,2\,D_{5}$}{$[1]$}
\elldata{1855}{$16$}{$E_{6}\,+\,D_{5}\,+\,D_{4}\,+\,A_{1}$}{$[1]$}
\elldata{1856}{$16$}{$E_{6}\,+\,D_{5}\,+\,A_{5}$}{$[1]$}
\elldata{1857}{$16$}{$E_{6}\,+\,D_{5}\,+\,A_{4}\,+\,A_{1}$}{$[1]$}
\elldata{1858}{$16$}{$E_{6}\,+\,D_{5}\,+\,A_{3}\,+\,A_{2}$}{$[1]$}
\elldata{1859}{$16$}{$E_{6}\,+\,D_{5}\,+\,A_{3}\,+\,2\,A_{1}$}{$[1]$}
\elldata{1860}{$16$}{$E_{6}\,+\,D_{5}\,+\,2\,A_{2}\,+\,A_{1}$}{$[1]$}
\elldata{1861}{$16$}{$E_{6}\,+\,D_{5}\,+\,A_{2}\,+\,3\,A_{1}$}{$[1]$}
\elldata{1862}{$16$}{$E_{6}\,+\,2\,D_{4}\,+\,A_{2}$}{$[1]$}
\elldata{1863}{$16$}{$E_{6}\,+\,D_{4}\,+\,A_{6}$}{$[1]$}
\elldata{1864}{$16$}{$E_{6}\,+\,D_{4}\,+\,A_{5}\,+\,A_{1}$}{$[1]$}
\elldata{1865}{$16$}{$E_{6}\,+\,D_{4}\,+\,A_{4}\,+\,A_{2}$}{$[1]$}
\elldata{1866}{$16$}{$E_{6}\,+\,D_{4}\,+\,A_{4}\,+\,2\,A_{1}$}{$[1]$}
\elldata{1867}{$16$}{$E_{6}\,+\,D_{4}\,+\,2\,A_{3}$}{$[1]$}
\elldata{1868}{$16$}{$E_{6}\,+\,D_{4}\,+\,A_{3}\,+\,A_{2}\,+\,A_{1}$}{$[1]$}
\elldata{1869}{$16$}{$E_{6}\,+\,D_{4}\,+\,2\,A_{2}\,+\,2\,A_{1}$}{$[1]$}
\elldata{1870}{$16$}{$E_{6}\,+\,A_{10}$}{$[1]$}
\elldata{1871}{$16$}{$E_{6}\,+\,A_{9}\,+\,A_{1}$}{$[1]$}
\elldata{1872}{$16$}{$E_{6}\,+\,A_{8}\,+\,A_{2}$}{$[3], \,[1]$}
\elldata{1873}{$16$}{$E_{6}\,+\,A_{8}\,+\,2\,A_{1}$}{$[1]$}
\elldata{1874}{$16$}{$E_{6}\,+\,A_{7}\,+\,A_{3}$}{$[1]$}
\elldata{1875}{$16$}{$E_{6}\,+\,A_{7}\,+\,A_{2}\,+\,A_{1}$}{$[1]$}
\elldata{1876}{$16$}{$E_{6}\,+\,A_{7}\,+\,3\,A_{1}$}{$[1]$}
\elldata{1877}{$16$}{$E_{6}\,+\,A_{6}\,+\,A_{4}$}{$[1]$}
\elldata{1878}{$16$}{$E_{6}\,+\,A_{6}\,+\,A_{3}\,+\,A_{1}$}{$[1]$}
\elldata{1879}{$16$}{$E_{6}\,+\,A_{6}\,+\,2\,A_{2}$}{$[1]$}
\elldata{1880}{$16$}{$E_{6}\,+\,A_{6}\,+\,A_{2}\,+\,2\,A_{1}$}{$[1]$}
\elldata{1881}{$16$}{$E_{6}\,+\,A_{6}\,+\,4\,A_{1}$}{$[1]$}
\elldata{1882}{$16$}{$E_{6}\,+\,2\,A_{5}$}{$[3], \,[1]$}
\elldata{1883}{$16$}{$E_{6}\,+\,A_{5}\,+\,A_{4}\,+\,A_{1}$}{$[1]$}
\elldata{1884}{$16$}{$E_{6}\,+\,A_{5}\,+\,A_{3}\,+\,A_{2}$}{$[1]$}
\elldata{1885}{$16$}{$E_{6}\,+\,A_{5}\,+\,A_{3}\,+\,2\,A_{1}$}{$[1]$}
\elldata{1886}{$16$}{$E_{6}\,+\,A_{5}\,+\,2\,A_{2}\,+\,A_{1}$}{$[3], \,[1]$}
\elldata{1887}{$16$}{$E_{6}\,+\,A_{5}\,+\,A_{2}\,+\,3\,A_{1}$}{$[1]$}
\elldata{1888}{$16$}{$E_{6}\,+\,2\,A_{4}\,+\,A_{2}$}{$[1]$}
\elldata{1889}{$16$}{$E_{6}\,+\,2\,A_{4}\,+\,2\,A_{1}$}{$[1]$}
\elldata{1890}{$16$}{$E_{6}\,+\,A_{4}\,+\,2\,A_{3}$}{$[1]$}
\elldata{1891}{$16$}{$E_{6}\,+\,A_{4}\,+\,A_{3}\,+\,A_{2}\,+\,A_{1}$}{$[1]$}
\elldata{1892}{$16$}{$E_{6}\,+\,A_{4}\,+\,A_{3}\,+\,3\,A_{1}$}{$[1]$}
\elldata{1893}{$16$}{$E_{6}\,+\,A_{4}\,+\,3\,A_{2}$}{$[1]$}
\elldata{1894}{$16$}{$E_{6}\,+\,A_{4}\,+\,2\,A_{2}\,+\,2\,A_{1}$}{$[1]$}
\elldata{1895}{$16$}{$E_{6}\,+\,A_{4}\,+\,A_{2}\,+\,4\,A_{1}$}{$[1]$}
\elldata{1896}{$16$}{$E_{6}\,+\,3\,A_{3}\,+\,A_{1}$}{$[1]$}
\elldata{1897}{$16$}{$E_{6}\,+\,2\,A_{3}\,+\,2\,A_{2}$}{$[1]$}
\elldata{1898}{$16$}{$E_{6}\,+\,2\,A_{3}\,+\,A_{2}\,+\,2\,A_{1}$}{$[1]$}
\elldata{1899}{$16$}{$E_{6}\,+\,A_{3}\,+\,3\,A_{2}\,+\,A_{1}$}{$[1]$}
\elldata{1900}{$16$}{$E_{6}\,+\,A_{3}\,+\,2\,A_{2}\,+\,3\,A_{1}$}{$[1]$}
\elldata{1901}{$16$}{$E_{6}\,+\,5\,A_{2}$}{$[3]$}
\elldata{1902}{$16$}{$E_{6}\,+\,4\,A_{2}\,+\,2\,A_{1}$}{$[3]$}
\elldata{1903}{$16$}{$D_{16}$}{$[2], \,[1]$}
\elldata{1904}{$16$}{$D_{15}\,+\,A_{1}$}{$[1]$}
\elldata{1905}{$16$}{$D_{14}\,+\,A_{2}$}{$[1]$}
\elldata{1906}{$16$}{$D_{14}\,+\,2\,A_{1}$}{$[2], \,[1]$}
\elldata{1907}{$16$}{$D_{13}\,+\,A_{3}$}{$[1]$}
\elldata{1908}{$16$}{$D_{13}\,+\,A_{2}\,+\,A_{1}$}{$[1]$}
\elldata{1909}{$16$}{$D_{13}\,+\,3\,A_{1}$}{$[1]$}
\elldata{1910}{$16$}{$D_{12}\,+\,D_{4}$}{$[2], \,[1]$}
\elldata{1911}{$16$}{$D_{12}\,+\,A_{4}$}{$[1]$}
\elldata{1912}{$16$}{$D_{12}\,+\,A_{3}\,+\,A_{1}$}{$[2], \,[1]$}
\elldata{1913}{$16$}{$D_{12}\,+\,2\,A_{2}$}{$[1]$}
\elldata{1914}{$16$}{$D_{12}\,+\,A_{2}\,+\,2\,A_{1}$}{$[2], \,[1]$}
\elldata{1915}{$16$}{$D_{12}\,+\,4\,A_{1}$}{$[2]$}
\elldata{1916}{$16$}{$D_{11}\,+\,D_{5}$}{$[1]$}
\elldata{1917}{$16$}{$D_{11}\,+\,D_{4}\,+\,A_{1}$}{$[1]$}
\elldata{1918}{$16$}{$D_{11}\,+\,A_{5}$}{$[1]$}
\elldata{1919}{$16$}{$D_{11}\,+\,A_{4}\,+\,A_{1}$}{$[1]$}
\elldata{1920}{$16$}{$D_{11}\,+\,A_{3}\,+\,A_{2}$}{$[1]$}
\elldata{1921}{$16$}{$D_{11}\,+\,A_{3}\,+\,2\,A_{1}$}{$[1]$}
\elldata{1922}{$16$}{$D_{11}\,+\,2\,A_{2}\,+\,A_{1}$}{$[1]$}
\elldata{1923}{$16$}{$D_{11}\,+\,A_{2}\,+\,3\,A_{1}$}{$[1]$}
\elldata{1924}{$16$}{$D_{10}\,+\,D_{6}$}{$[2], \,[1]$}
\elldata{1925}{$16$}{$D_{10}\,+\,D_{5}\,+\,A_{1}$}{$[2], \,[1]$}
\elldata{1926}{$16$}{$D_{10}\,+\,D_{4}\,+\,A_{2}$}{$[1]$}
\elldata{1927}{$16$}{$D_{10}\,+\,D_{4}\,+\,2\,A_{1}$}{$[2]$}
\elldata{1928}{$16$}{$D_{10}\,+\,A_{6}$}{$[1]$}
\elldata{1929}{$16$}{$D_{10}\,+\,A_{5}\,+\,A_{1}$}{$[2], \,[1]$}
\elldata{1930}{$16$}{$D_{10}\,+\,A_{4}\,+\,A_{2}$}{$[1]$}
\elldata{1931}{$16$}{$D_{10}\,+\,A_{4}\,+\,2\,A_{1}$}{$[1]$}
\elldata{1932}{$16$}{$D_{10}\,+\,2\,A_{3}$}{$[1]$}
\elldata{1933}{$16$}{$D_{10}\,+\,A_{3}\,+\,A_{2}\,+\,A_{1}$}{$[2], \,[1]$}
\elldata{1934}{$16$}{$D_{10}\,+\,A_{3}\,+\,3\,A_{1}$}{$[2]$}
\elldata{1935}{$16$}{$D_{10}\,+\,3\,A_{2}$}{$[1]$}
\elldata{1936}{$16$}{$D_{10}\,+\,2\,A_{2}\,+\,2\,A_{1}$}{$[1]$}
\elldata{1937}{$16$}{$D_{10}\,+\,A_{2}\,+\,4\,A_{1}$}{$[2]$}
\elldata{1938}{$16$}{$D_{10}\,+\,6\,A_{1}$}{$[2, 2]$}
\elldata{1939}{$16$}{$D_{9}\,+\,D_{7}$}{$[1]$}
\elldata{1940}{$16$}{$D_{9}\,+\,D_{6}\,+\,A_{1}$}{$[1]$}
\elldata{1941}{$16$}{$D_{9}\,+\,D_{5}\,+\,A_{2}$}{$[1]$}
\elldata{1942}{$16$}{$D_{9}\,+\,D_{5}\,+\,2\,A_{1}$}{$[1]$}
\elldata{1943}{$16$}{$D_{9}\,+\,D_{4}\,+\,A_{3}$}{$[1]$}
\elldata{1944}{$16$}{$D_{9}\,+\,D_{4}\,+\,A_{2}\,+\,A_{1}$}{$[1]$}
\elldata{1945}{$16$}{$D_{9}\,+\,A_{7}$}{$[1]$}
\elldata{1946}{$16$}{$D_{9}\,+\,A_{6}\,+\,A_{1}$}{$[1]$}
\elldata{1947}{$16$}{$D_{9}\,+\,A_{5}\,+\,A_{2}$}{$[1]$}
\elldata{1948}{$16$}{$D_{9}\,+\,A_{5}\,+\,2\,A_{1}$}{$[1]$}
\elldata{1949}{$16$}{$D_{9}\,+\,A_{4}\,+\,A_{3}$}{$[1]$}
\elldata{1950}{$16$}{$D_{9}\,+\,A_{4}\,+\,A_{2}\,+\,A_{1}$}{$[1]$}
\elldata{1951}{$16$}{$D_{9}\,+\,A_{4}\,+\,3\,A_{1}$}{$[1]$}
\elldata{1952}{$16$}{$D_{9}\,+\,2\,A_{3}\,+\,A_{1}$}{$[1]$}
\elldata{1953}{$16$}{$D_{9}\,+\,A_{3}\,+\,2\,A_{2}$}{$[1]$}
\elldata{1954}{$16$}{$D_{9}\,+\,A_{3}\,+\,A_{2}\,+\,2\,A_{1}$}{$[1]$}
\elldata{1955}{$16$}{$D_{9}\,+\,A_{3}\,+\,4\,A_{1}$}{$[2]$}
\elldata{1956}{$16$}{$D_{9}\,+\,3\,A_{2}\,+\,A_{1}$}{$[1]$}
\elldata{1957}{$16$}{$D_{9}\,+\,2\,A_{2}\,+\,3\,A_{1}$}{$[1]$}
\elldata{1958}{$16$}{$2\,D_{8}$}{$[2], \,[1]$}
\elldata{1959}{$16$}{$D_{8}\,+\,D_{7}\,+\,A_{1}$}{$[1]$}
\elldata{1960}{$16$}{$D_{8}\,+\,D_{6}\,+\,A_{2}$}{$[1]$}
\elldata{1961}{$16$}{$D_{8}\,+\,D_{6}\,+\,2\,A_{1}$}{$[2]$}
\elldata{1962}{$16$}{$D_{8}\,+\,D_{5}\,+\,A_{3}$}{$[2], \,[1]$}
\elldata{1963}{$16$}{$D_{8}\,+\,D_{5}\,+\,A_{2}\,+\,A_{1}$}{$[1]$}
\elldata{1964}{$16$}{$D_{8}\,+\,D_{5}\,+\,3\,A_{1}$}{$[2]$}
\elldata{1965}{$16$}{$D_{8}\,+\,2\,D_{4}$}{$[2]$}
\elldata{1966}{$16$}{$D_{8}\,+\,D_{4}\,+\,A_{4}$}{$[1]$}
\elldata{1967}{$16$}{$D_{8}\,+\,D_{4}\,+\,A_{3}\,+\,A_{1}$}{$[2]$}
\elldata{1968}{$16$}{$D_{8}\,+\,D_{4}\,+\,A_{2}\,+\,2\,A_{1}$}{$[2]$}
\elldata{1969}{$16$}{$D_{8}\,+\,D_{4}\,+\,4\,A_{1}$}{$[2, 2]$}
\elldata{1970}{$16$}{$D_{8}\,+\,A_{8}$}{$[1]$}
\elldata{1971}{$16$}{$D_{8}\,+\,A_{7}\,+\,A_{1}$}{$[2], \,[1]$}
\elldata{1972}{$16$}{$D_{8}\,+\,A_{6}\,+\,A_{2}$}{$[1]$}
\elldata{1973}{$16$}{$D_{8}\,+\,A_{6}\,+\,2\,A_{1}$}{$[1]$}
\elldata{1974}{$16$}{$D_{8}\,+\,A_{5}\,+\,A_{3}$}{$[1]$}
\elldata{1975}{$16$}{$D_{8}\,+\,A_{5}\,+\,A_{2}\,+\,A_{1}$}{$[2], \,[1]$}
\elldata{1976}{$16$}{$D_{8}\,+\,A_{5}\,+\,3\,A_{1}$}{$[2]$}
\elldata{1977}{$16$}{$D_{8}\,+\,2\,A_{4}$}{$[1]$}
\elldata{1978}{$16$}{$D_{8}\,+\,A_{4}\,+\,A_{3}\,+\,A_{1}$}{$[1]$}
\elldata{1979}{$16$}{$D_{8}\,+\,A_{4}\,+\,2\,A_{2}$}{$[1]$}
\elldata{1980}{$16$}{$D_{8}\,+\,A_{4}\,+\,A_{2}\,+\,2\,A_{1}$}{$[1]$}
\elldata{1981}{$16$}{$D_{8}\,+\,A_{4}\,+\,4\,A_{1}$}{$[2]$}
\elldata{1982}{$16$}{$D_{8}\,+\,2\,A_{3}\,+\,A_{2}$}{$[2]$}
\elldata{1983}{$16$}{$D_{8}\,+\,2\,A_{3}\,+\,2\,A_{1}$}{$[2]$}
\elldata{1984}{$16$}{$D_{8}\,+\,A_{3}\,+\,2\,A_{2}\,+\,A_{1}$}{$[1]$}
\elldata{1985}{$16$}{$D_{8}\,+\,A_{3}\,+\,A_{2}\,+\,3\,A_{1}$}{$[2]$}
\elldata{1986}{$16$}{$D_{8}\,+\,A_{3}\,+\,5\,A_{1}$}{$[2, 2]$}
\elldata{1987}{$16$}{$D_{8}\,+\,4\,A_{2}$}{$[1]$}
\elldata{1988}{$16$}{$D_{8}\,+\,3\,A_{2}\,+\,2\,A_{1}$}{$[1]$}
\elldata{1989}{$16$}{$D_{8}\,+\,2\,A_{2}\,+\,4\,A_{1}$}{$[2]$}
\elldata{1990}{$16$}{$2\,D_{7}\,+\,A_{2}$}{$[1]$}
\elldata{1991}{$16$}{$2\,D_{7}\,+\,2\,A_{1}$}{$[1]$}
\elldata{1992}{$16$}{$D_{7}\,+\,D_{6}\,+\,A_{3}$}{$[1]$}
\elldata{1993}{$16$}{$D_{7}\,+\,D_{6}\,+\,A_{2}\,+\,A_{1}$}{$[1]$}
\elldata{1994}{$16$}{$D_{7}\,+\,D_{6}\,+\,3\,A_{1}$}{$[2]$}
\elldata{1995}{$16$}{$D_{7}\,+\,D_{5}\,+\,D_{4}$}{$[1]$}
\elldata{1996}{$16$}{$D_{7}\,+\,D_{5}\,+\,A_{4}$}{$[1]$}
\elldata{1997}{$16$}{$D_{7}\,+\,D_{5}\,+\,A_{3}\,+\,A_{1}$}{$[1]$}
\elldata{1998}{$16$}{$D_{7}\,+\,D_{5}\,+\,2\,A_{2}$}{$[1]$}
\elldata{1999}{$16$}{$D_{7}\,+\,D_{5}\,+\,A_{2}\,+\,2\,A_{1}$}{$[1]$}
\elldata{2000}{$16$}{$D_{7}\,+\,D_{5}\,+\,4\,A_{1}$}{$[2]$}
\elldata{2001}{$16$}{$D_{7}\,+\,D_{4}\,+\,A_{5}$}{$[1]$}
\elldata{2002}{$16$}{$D_{7}\,+\,D_{4}\,+\,A_{4}\,+\,A_{1}$}{$[1]$}
\elldata{2003}{$16$}{$D_{7}\,+\,D_{4}\,+\,A_{3}\,+\,2\,A_{1}$}{$[2]$}
\elldata{2004}{$16$}{$D_{7}\,+\,D_{4}\,+\,2\,A_{2}\,+\,A_{1}$}{$[1]$}
\elldata{2005}{$16$}{$D_{7}\,+\,A_{9}$}{$[1]$}
\elldata{2006}{$16$}{$D_{7}\,+\,A_{8}\,+\,A_{1}$}{$[1]$}
\elldata{2007}{$16$}{$D_{7}\,+\,A_{7}\,+\,A_{2}$}{$[1]$}
\elldata{2008}{$16$}{$D_{7}\,+\,A_{7}\,+\,2\,A_{1}$}{$[2], \,[1]$}
\elldata{2009}{$16$}{$D_{7}\,+\,A_{6}\,+\,A_{3}$}{$[1]$}
\elldata{2010}{$16$}{$D_{7}\,+\,A_{6}\,+\,A_{2}\,+\,A_{1}$}{$[1]$}
\elldata{2011}{$16$}{$D_{7}\,+\,A_{6}\,+\,3\,A_{1}$}{$[1]$}
\elldata{2012}{$16$}{$D_{7}\,+\,A_{5}\,+\,A_{4}$}{$[1]$}
\elldata{2013}{$16$}{$D_{7}\,+\,A_{5}\,+\,A_{3}\,+\,A_{1}$}{$[2], \,[1]$}
\elldata{2014}{$16$}{$D_{7}\,+\,A_{5}\,+\,2\,A_{2}$}{$[1]$}
\elldata{2015}{$16$}{$D_{7}\,+\,A_{5}\,+\,A_{2}\,+\,2\,A_{1}$}{$[1]$}
\elldata{2016}{$16$}{$D_{7}\,+\,A_{5}\,+\,4\,A_{1}$}{$[2]$}
\elldata{2017}{$16$}{$D_{7}\,+\,2\,A_{4}\,+\,A_{1}$}{$[1]$}
\elldata{2018}{$16$}{$D_{7}\,+\,A_{4}\,+\,A_{3}\,+\,A_{2}$}{$[1]$}
\elldata{2019}{$16$}{$D_{7}\,+\,A_{4}\,+\,A_{3}\,+\,2\,A_{1}$}{$[1]$}
\elldata{2020}{$16$}{$D_{7}\,+\,A_{4}\,+\,2\,A_{2}\,+\,A_{1}$}{$[1]$}
\elldata{2021}{$16$}{$D_{7}\,+\,A_{4}\,+\,A_{2}\,+\,3\,A_{1}$}{$[1]$}
\elldata{2022}{$16$}{$D_{7}\,+\,3\,A_{3}$}{$[4]$}
\elldata{2023}{$16$}{$D_{7}\,+\,2\,A_{3}\,+\,A_{2}\,+\,A_{1}$}{$[1]$}
\elldata{2024}{$16$}{$D_{7}\,+\,2\,A_{3}\,+\,3\,A_{1}$}{$[2]$}
\elldata{2025}{$16$}{$D_{7}\,+\,A_{3}\,+\,3\,A_{2}$}{$[1]$}
\elldata{2026}{$16$}{$D_{7}\,+\,A_{3}\,+\,2\,A_{2}\,+\,2\,A_{1}$}{$[1]$}
\elldata{2027}{$16$}{$D_{7}\,+\,A_{3}\,+\,A_{2}\,+\,4\,A_{1}$}{$[2]$}
\elldata{2028}{$16$}{$D_{7}\,+\,3\,A_{2}\,+\,3\,A_{1}$}{$[1]$}
\elldata{2029}{$16$}{$2\,D_{6}\,+\,D_{4}$}{$[2]$}
\elldata{2030}{$16$}{$2\,D_{6}\,+\,A_{4}$}{$[1]$}
\elldata{2031}{$16$}{$2\,D_{6}\,+\,A_{3}\,+\,A_{1}$}{$[2]$}
\elldata{2032}{$16$}{$2\,D_{6}\,+\,2\,A_{2}$}{$[1]$}
\elldata{2033}{$16$}{$2\,D_{6}\,+\,A_{2}\,+\,2\,A_{1}$}{$[2]$}
\elldata{2034}{$16$}{$2\,D_{6}\,+\,4\,A_{1}$}{$[2, 2]$}
\elldata{2035}{$16$}{$D_{6}\,+\,2\,D_{5}$}{$[1]$}
\elldata{2036}{$16$}{$D_{6}\,+\,D_{5}\,+\,D_{4}\,+\,A_{1}$}{$[2]$}
\elldata{2037}{$16$}{$D_{6}\,+\,D_{5}\,+\,A_{5}$}{$[2], \,[1]$}
\elldata{2038}{$16$}{$D_{6}\,+\,D_{5}\,+\,A_{4}\,+\,A_{1}$}{$[1]$}
\elldata{2039}{$16$}{$D_{6}\,+\,D_{5}\,+\,A_{3}\,+\,A_{2}$}{$[1]$}
\elldata{2040}{$16$}{$D_{6}\,+\,D_{5}\,+\,A_{3}\,+\,2\,A_{1}$}{$[2]$}
\elldata{2041}{$16$}{$D_{6}\,+\,D_{5}\,+\,2\,A_{2}\,+\,A_{1}$}{$[1]$}
\elldata{2042}{$16$}{$D_{6}\,+\,D_{5}\,+\,A_{2}\,+\,3\,A_{1}$}{$[2]$}
\elldata{2043}{$16$}{$D_{6}\,+\,2\,D_{4}\,+\,2\,A_{1}$}{$[2, 2]$}
\elldata{2044}{$16$}{$D_{6}\,+\,D_{4}\,+\,A_{6}$}{$[1]$}
\elldata{2045}{$16$}{$D_{6}\,+\,D_{4}\,+\,A_{5}\,+\,A_{1}$}{$[2]$}
\elldata{2046}{$16$}{$D_{6}\,+\,D_{4}\,+\,A_{4}\,+\,A_{2}$}{$[1]$}
\elldata{2047}{$16$}{$D_{6}\,+\,D_{4}\,+\,2\,A_{3}$}{$[2]$}
\elldata{2048}{$16$}{$D_{6}\,+\,D_{4}\,+\,A_{3}\,+\,A_{2}\,+\,A_{1}$}{$[2]$}
\elldata{2049}{$16$}{$D_{6}\,+\,D_{4}\,+\,A_{3}\,+\,3\,A_{1}$}{$[2, 2]$}
\elldata{2050}{$16$}{$D_{6}\,+\,D_{4}\,+\,3\,A_{2}$}{$[1]$}
\elldata{2051}{$16$}{$D_{6}\,+\,A_{10}$}{$[1]$}
\elldata{2052}{$16$}{$D_{6}\,+\,A_{9}\,+\,A_{1}$}{$[2], \,[1]$}
\elldata{2053}{$16$}{$D_{6}\,+\,A_{8}\,+\,A_{2}$}{$[1]$}
\elldata{2054}{$16$}{$D_{6}\,+\,A_{8}\,+\,2\,A_{1}$}{$[1]$}
\elldata{2055}{$16$}{$D_{6}\,+\,A_{7}\,+\,A_{3}$}{$[2], \,[1]$}
\elldata{2056}{$16$}{$D_{6}\,+\,A_{7}\,+\,A_{2}\,+\,A_{1}$}{$[2], \,[1]$}
\elldata{2057}{$16$}{$D_{6}\,+\,A_{7}\,+\,3\,A_{1}$}{$[2]$}
\elldata{2058}{$16$}{$D_{6}\,+\,A_{6}\,+\,A_{4}$}{$[1]$}
\elldata{2059}{$16$}{$D_{6}\,+\,A_{6}\,+\,A_{3}\,+\,A_{1}$}{$[1]$}
\elldata{2060}{$16$}{$D_{6}\,+\,A_{6}\,+\,2\,A_{2}$}{$[1]$}
\elldata{2061}{$16$}{$D_{6}\,+\,A_{6}\,+\,A_{2}\,+\,2\,A_{1}$}{$[1]$}
\elldata{2062}{$16$}{$D_{6}\,+\,2\,A_{5}$}{$[2], \,[1]$}
\elldata{2063}{$16$}{$D_{6}\,+\,A_{5}\,+\,A_{4}\,+\,A_{1}$}{$[1]$}
\elldata{2064}{$16$}{$D_{6}\,+\,A_{5}\,+\,A_{3}\,+\,A_{2}$}{$[2], \,[1]$}
\elldata{2065}{$16$}{$D_{6}\,+\,A_{5}\,+\,A_{3}\,+\,2\,A_{1}$}{$[2]$}
\elldata{2066}{$16$}{$D_{6}\,+\,A_{5}\,+\,2\,A_{2}\,+\,A_{1}$}{$[1]$}
\elldata{2067}{$16$}{$D_{6}\,+\,A_{5}\,+\,A_{2}\,+\,3\,A_{1}$}{$[2]$}
\elldata{2068}{$16$}{$D_{6}\,+\,A_{5}\,+\,5\,A_{1}$}{$[2, 2]$}
\elldata{2069}{$16$}{$D_{6}\,+\,2\,A_{4}\,+\,A_{2}$}{$[1]$}
\elldata{2070}{$16$}{$D_{6}\,+\,2\,A_{4}\,+\,2\,A_{1}$}{$[1]$}
\elldata{2071}{$16$}{$D_{6}\,+\,A_{4}\,+\,2\,A_{3}$}{$[1]$}
\elldata{2072}{$16$}{$D_{6}\,+\,A_{4}\,+\,A_{3}\,+\,A_{2}\,+\,A_{1}$}{$[1]$}
\elldata{2073}{$16$}{$D_{6}\,+\,A_{4}\,+\,A_{3}\,+\,3\,A_{1}$}{$[2]$}
\elldata{2074}{$16$}{$D_{6}\,+\,A_{4}\,+\,3\,A_{2}$}{$[1]$}
\elldata{2075}{$16$}{$D_{6}\,+\,A_{4}\,+\,2\,A_{2}\,+\,2\,A_{1}$}{$[1]$}
\elldata{2076}{$16$}{$D_{6}\,+\,3\,A_{3}\,+\,A_{1}$}{$[2]$}
\elldata{2077}{$16$}{$D_{6}\,+\,2\,A_{3}\,+\,2\,A_{2}$}{$[1]$}
\elldata{2078}{$16$}{$D_{6}\,+\,2\,A_{3}\,+\,A_{2}\,+\,2\,A_{1}$}{$[2]$}
\elldata{2079}{$16$}{$D_{6}\,+\,2\,A_{3}\,+\,4\,A_{1}$}{$[2, 2]$}
\elldata{2080}{$16$}{$D_{6}\,+\,A_{3}\,+\,3\,A_{2}\,+\,A_{1}$}{$[1]$}
\elldata{2081}{$16$}{$D_{6}\,+\,A_{3}\,+\,2\,A_{2}\,+\,3\,A_{1}$}{$[2]$}
\elldata{2082}{$16$}{$D_{6}\,+\,4\,A_{2}\,+\,2\,A_{1}$}{$[1]$}
\elldata{2083}{$16$}{$3\,D_{5}\,+\,A_{1}$}{$[1]$}
\elldata{2084}{$16$}{$2\,D_{5}\,+\,D_{4}\,+\,A_{2}$}{$[1]$}
\elldata{2085}{$16$}{$2\,D_{5}\,+\,D_{4}\,+\,2\,A_{1}$}{$[2]$}
\elldata{2086}{$16$}{$2\,D_{5}\,+\,A_{6}$}{$[1]$}
\elldata{2087}{$16$}{$2\,D_{5}\,+\,A_{5}\,+\,A_{1}$}{$[2], \,[1]$}
\elldata{2088}{$16$}{$2\,D_{5}\,+\,A_{4}\,+\,A_{2}$}{$[1]$}
\elldata{2089}{$16$}{$2\,D_{5}\,+\,A_{4}\,+\,2\,A_{1}$}{$[1]$}
\elldata{2090}{$16$}{$2\,D_{5}\,+\,2\,A_{3}$}{$[4], \,[2], \,[1]$}
\elldata{2091}{$16$}{$2\,D_{5}\,+\,A_{3}\,+\,A_{2}\,+\,A_{1}$}{$[1]$}
\elldata{2092}{$16$}{$2\,D_{5}\,+\,A_{3}\,+\,3\,A_{1}$}{$[2]$}
\elldata{2093}{$16$}{$2\,D_{5}\,+\,3\,A_{2}$}{$[1]$}
\elldata{2094}{$16$}{$2\,D_{5}\,+\,2\,A_{2}\,+\,2\,A_{1}$}{$[1]$}
\elldata{2095}{$16$}{$D_{5}\,+\,2\,D_{4}\,+\,A_{3}$}{$[2]$}
\elldata{2096}{$16$}{$D_{5}\,+\,D_{4}\,+\,A_{7}$}{$[2], \,[1]$}
\elldata{2097}{$16$}{$D_{5}\,+\,D_{4}\,+\,A_{6}\,+\,A_{1}$}{$[1]$}
\elldata{2098}{$16$}{$D_{5}\,+\,D_{4}\,+\,A_{5}\,+\,A_{2}$}{$[1]$}
\elldata{2099}{$16$}{$D_{5}\,+\,D_{4}\,+\,A_{5}\,+\,2\,A_{1}$}{$[2]$}
\elldata{2100}{$16$}{$D_{5}\,+\,D_{4}\,+\,A_{4}\,+\,A_{3}$}{$[1]$}
\elldata{2101}{$16$}{$D_{5}\,+\,D_{4}\,+\,A_{4}\,+\,A_{2}\,+\,A_{1}$}{$[1]$}
\elldata{2102}{$16$}{$D_{5}\,+\,D_{4}\,+\,2\,A_{3}\,+\,A_{1}$}{$[2]$}
\elldata{2103}{$16$}{$D_{5}\,+\,D_{4}\,+\,A_{3}\,+\,A_{2}\,+\,2\,A_{1}$}{$[2]$}
\elldata{2104}{$16$}{$D_{5}\,+\,D_{4}\,+\,3\,A_{2}\,+\,A_{1}$}{$[1]$}
\elldata{2105}{$16$}{$D_{5}\,+\,A_{11}$}{$[2], \,[1]$}
\elldata{2106}{$16$}{$D_{5}\,+\,A_{10}\,+\,A_{1}$}{$[1]$}
\elldata{2107}{$16$}{$D_{5}\,+\,A_{9}\,+\,A_{2}$}{$[1]$}
\elldata{2108}{$16$}{$D_{5}\,+\,A_{9}\,+\,2\,A_{1}$}{$[2], \,[1]$}
\elldata{2109}{$16$}{$D_{5}\,+\,A_{8}\,+\,A_{3}$}{$[1]$}
\elldata{2110}{$16$}{$D_{5}\,+\,A_{8}\,+\,A_{2}\,+\,A_{1}$}{$[1]$}
\elldata{2111}{$16$}{$D_{5}\,+\,A_{8}\,+\,3\,A_{1}$}{$[1]$}
\elldata{2112}{$16$}{$D_{5}\,+\,A_{7}\,+\,A_{4}$}{$[1]$}
\elldata{2113}{$16$}{$D_{5}\,+\,A_{7}\,+\,A_{3}\,+\,A_{1}$}{$[4], \,[2], \,[1]$}
\elldata{2114}{$16$}{$D_{5}\,+\,A_{7}\,+\,2\,A_{2}$}{$[1]$}
\elldata{2115}{$16$}{$D_{5}\,+\,A_{7}\,+\,A_{2}\,+\,2\,A_{1}$}{$[2], \,[1]$}
\elldata{2116}{$16$}{$D_{5}\,+\,A_{7}\,+\,4\,A_{1}$}{$[2]$}
\elldata{2117}{$16$}{$D_{5}\,+\,A_{6}\,+\,A_{5}$}{$[1]$}
\elldata{2118}{$16$}{$D_{5}\,+\,A_{6}\,+\,A_{4}\,+\,A_{1}$}{$[1]$}
\elldata{2119}{$16$}{$D_{5}\,+\,A_{6}\,+\,A_{3}\,+\,A_{2}$}{$[1]$}
\elldata{2120}{$16$}{$D_{5}\,+\,A_{6}\,+\,A_{3}\,+\,2\,A_{1}$}{$[1]$}
\elldata{2121}{$16$}{$D_{5}\,+\,A_{6}\,+\,2\,A_{2}\,+\,A_{1}$}{$[1]$}
\elldata{2122}{$16$}{$D_{5}\,+\,A_{6}\,+\,A_{2}\,+\,3\,A_{1}$}{$[1]$}
\elldata{2123}{$16$}{$D_{5}\,+\,2\,A_{5}\,+\,A_{1}$}{$[2], \,[1]$}
\elldata{2124}{$16$}{$D_{5}\,+\,A_{5}\,+\,A_{4}\,+\,A_{2}$}{$[1]$}
\elldata{2125}{$16$}{$D_{5}\,+\,A_{5}\,+\,A_{4}\,+\,2\,A_{1}$}{$[1]$}
\elldata{2126}{$16$}{$D_{5}\,+\,A_{5}\,+\,2\,A_{3}$}{$[1]$}
\elldata{2127}{$16$}{$D_{5}\,+\,A_{5}\,+\,A_{3}\,+\,A_{2}\,+\,A_{1}$}{$[2], \,[1]$}
\elldata{2128}{$16$}{$D_{5}\,+\,A_{5}\,+\,A_{3}\,+\,3\,A_{1}$}{$[2]$}
\elldata{2129}{$16$}{$D_{5}\,+\,A_{5}\,+\,3\,A_{2}$}{$[1]$}
\elldata{2130}{$16$}{$D_{5}\,+\,A_{5}\,+\,2\,A_{2}\,+\,2\,A_{1}$}{$[1]$}
\elldata{2131}{$16$}{$D_{5}\,+\,A_{5}\,+\,A_{2}\,+\,4\,A_{1}$}{$[2]$}
\elldata{2132}{$16$}{$D_{5}\,+\,2\,A_{4}\,+\,A_{3}$}{$[1]$}
\elldata{2133}{$16$}{$D_{5}\,+\,2\,A_{4}\,+\,A_{2}\,+\,A_{1}$}{$[1]$}
\elldata{2134}{$16$}{$D_{5}\,+\,2\,A_{4}\,+\,3\,A_{1}$}{$[1]$}
\elldata{2135}{$16$}{$D_{5}\,+\,A_{4}\,+\,2\,A_{3}\,+\,A_{1}$}{$[1]$}
\elldata{2136}{$16$}{$D_{5}\,+\,A_{4}\,+\,A_{3}\,+\,2\,A_{2}$}{$[1]$}
\elldata{2137}{$16$}{$D_{5}\,+\,A_{4}\,+\,A_{3}\,+\,A_{2}\,+\,2\,A_{1}$}{$[1]$}
\elldata{2138}{$16$}{$D_{5}\,+\,A_{4}\,+\,A_{3}\,+\,4\,A_{1}$}{$[2]$}
\elldata{2139}{$16$}{$D_{5}\,+\,A_{4}\,+\,3\,A_{2}\,+\,A_{1}$}{$[1]$}
\elldata{2140}{$16$}{$D_{5}\,+\,A_{4}\,+\,2\,A_{2}\,+\,3\,A_{1}$}{$[1]$}
\elldata{2141}{$16$}{$D_{5}\,+\,3\,A_{3}\,+\,A_{2}$}{$[2]$}
\elldata{2142}{$16$}{$D_{5}\,+\,3\,A_{3}\,+\,2\,A_{1}$}{$[4], \,[2]$}
\elldata{2143}{$16$}{$D_{5}\,+\,2\,A_{3}\,+\,2\,A_{2}\,+\,A_{1}$}{$[1]$}
\elldata{2144}{$16$}{$D_{5}\,+\,2\,A_{3}\,+\,A_{2}\,+\,3\,A_{1}$}{$[2]$}
\elldata{2145}{$16$}{$D_{5}\,+\,A_{3}\,+\,4\,A_{2}$}{$[1]$}
\elldata{2146}{$16$}{$D_{5}\,+\,A_{3}\,+\,3\,A_{2}\,+\,2\,A_{1}$}{$[1]$}
\elldata{2147}{$16$}{$4\,D_{4}$}{$[2, 2]$}
\elldata{2148}{$16$}{$2\,D_{4}\,+\,A_{8}$}{$[1]$}
\elldata{2149}{$16$}{$2\,D_{4}\,+\,A_{7}\,+\,A_{1}$}{$[2]$}
\elldata{2150}{$16$}{$2\,D_{4}\,+\,A_{5}\,+\,A_{2}\,+\,A_{1}$}{$[2]$}
\elldata{2151}{$16$}{$2\,D_{4}\,+\,2\,A_{4}$}{$[1]$}
\elldata{2152}{$16$}{$2\,D_{4}\,+\,2\,A_{3}\,+\,2\,A_{1}$}{$[2, 2]$}
\elldata{2153}{$16$}{$2\,D_{4}\,+\,4\,A_{2}$}{$[1]$}
\elldata{2154}{$16$}{$D_{4}\,+\,A_{12}$}{$[1]$}
\elldata{2155}{$16$}{$D_{4}\,+\,A_{11}\,+\,A_{1}$}{$[2], \,[1]$}
\elldata{2156}{$16$}{$D_{4}\,+\,A_{10}\,+\,A_{2}$}{$[1]$}
\elldata{2157}{$16$}{$D_{4}\,+\,A_{10}\,+\,2\,A_{1}$}{$[1]$}
\elldata{2158}{$16$}{$D_{4}\,+\,A_{9}\,+\,A_{3}$}{$[1]$}
\elldata{2159}{$16$}{$D_{4}\,+\,A_{9}\,+\,A_{2}\,+\,A_{1}$}{$[2], \,[1]$}
\elldata{2160}{$16$}{$D_{4}\,+\,A_{9}\,+\,3\,A_{1}$}{$[2]$}
\elldata{2161}{$16$}{$D_{4}\,+\,A_{8}\,+\,A_{4}$}{$[1]$}
\elldata{2162}{$16$}{$D_{4}\,+\,A_{8}\,+\,A_{3}\,+\,A_{1}$}{$[1]$}
\elldata{2163}{$16$}{$D_{4}\,+\,A_{8}\,+\,2\,A_{2}$}{$[1]$}
\elldata{2164}{$16$}{$D_{4}\,+\,A_{8}\,+\,A_{2}\,+\,2\,A_{1}$}{$[1]$}
\elldata{2165}{$16$}{$D_{4}\,+\,A_{7}\,+\,A_{5}$}{$[1]$}
\elldata{2166}{$16$}{$D_{4}\,+\,A_{7}\,+\,A_{4}\,+\,A_{1}$}{$[1]$}
\elldata{2167}{$16$}{$D_{4}\,+\,A_{7}\,+\,A_{3}\,+\,A_{2}$}{$[2]$}
\elldata{2168}{$16$}{$D_{4}\,+\,A_{7}\,+\,A_{3}\,+\,2\,A_{1}$}{$[2]$}
\elldata{2169}{$16$}{$D_{4}\,+\,A_{7}\,+\,2\,A_{2}\,+\,A_{1}$}{$[1]$}
\elldata{2170}{$16$}{$D_{4}\,+\,A_{7}\,+\,A_{2}\,+\,3\,A_{1}$}{$[2]$}
\elldata{2171}{$16$}{$D_{4}\,+\,2\,A_{6}$}{$[1]$}
\elldata{2172}{$16$}{$D_{4}\,+\,A_{6}\,+\,A_{5}\,+\,A_{1}$}{$[1]$}
\elldata{2173}{$16$}{$D_{4}\,+\,A_{6}\,+\,A_{4}\,+\,A_{2}$}{$[1]$}
\elldata{2174}{$16$}{$D_{4}\,+\,A_{6}\,+\,A_{4}\,+\,2\,A_{1}$}{$[1]$}
\elldata{2175}{$16$}{$D_{4}\,+\,A_{6}\,+\,A_{3}\,+\,A_{2}\,+\,A_{1}$}{$[1]$}
\elldata{2176}{$16$}{$D_{4}\,+\,A_{6}\,+\,3\,A_{2}$}{$[1]$}
\elldata{2177}{$16$}{$D_{4}\,+\,A_{6}\,+\,2\,A_{2}\,+\,2\,A_{1}$}{$[1]$}
\elldata{2178}{$16$}{$D_{4}\,+\,2\,A_{5}\,+\,A_{2}$}{$[2], \,[1]$}
\elldata{2179}{$16$}{$D_{4}\,+\,2\,A_{5}\,+\,2\,A_{1}$}{$[2]$}
\elldata{2180}{$16$}{$D_{4}\,+\,A_{5}\,+\,A_{4}\,+\,A_{3}$}{$[1]$}
\elldata{2181}{$16$}{$D_{4}\,+\,A_{5}\,+\,A_{4}\,+\,A_{2}\,+\,A_{1}$}{$[1]$}
\elldata{2182}{$16$}{$D_{4}\,+\,A_{5}\,+\,A_{4}\,+\,3\,A_{1}$}{$[2]$}
\elldata{2183}{$16$}{$D_{4}\,+\,A_{5}\,+\,2\,A_{3}\,+\,A_{1}$}{$[2]$}
\elldata{2184}{$16$}{$D_{4}\,+\,A_{5}\,+\,A_{3}\,+\,2\,A_{2}$}{$[1]$}
\elldata{2185}{$16$}{$D_{4}\,+\,A_{5}\,+\,A_{3}\,+\,A_{2}\,+\,2\,A_{1}$}{$[2]$}
\elldata{2186}{$16$}{$D_{4}\,+\,A_{5}\,+\,A_{3}\,+\,4\,A_{1}$}{$[2, 2]$}
\elldata{2187}{$16$}{$D_{4}\,+\,A_{5}\,+\,2\,A_{2}\,+\,3\,A_{1}$}{$[2]$}
\elldata{2188}{$16$}{$D_{4}\,+\,3\,A_{4}$}{$[1]$}
\elldata{2189}{$16$}{$D_{4}\,+\,2\,A_{4}\,+\,A_{3}\,+\,A_{1}$}{$[1]$}
\elldata{2190}{$16$}{$D_{4}\,+\,2\,A_{4}\,+\,2\,A_{2}$}{$[1]$}
\elldata{2191}{$16$}{$D_{4}\,+\,2\,A_{4}\,+\,A_{2}\,+\,2\,A_{1}$}{$[1]$}
\elldata{2192}{$16$}{$D_{4}\,+\,A_{4}\,+\,2\,A_{3}\,+\,2\,A_{1}$}{$[2]$}
\elldata{2193}{$16$}{$D_{4}\,+\,A_{4}\,+\,A_{3}\,+\,2\,A_{2}\,+\,A_{1}$}{$[1]$}
\elldata{2194}{$16$}{$D_{4}\,+\,A_{4}\,+\,3\,A_{2}\,+\,2\,A_{1}$}{$[1]$}
\elldata{2195}{$16$}{$D_{4}\,+\,3\,A_{3}\,+\,A_{2}\,+\,A_{1}$}{$[2]$}
\elldata{2196}{$16$}{$D_{4}\,+\,3\,A_{3}\,+\,3\,A_{1}$}{$[2, 2]$}
\elldata{2197}{$16$}{$D_{4}\,+\,2\,A_{3}\,+\,3\,A_{2}$}{$[1]$}
\elldata{2198}{$16$}{$D_{4}\,+\,2\,A_{3}\,+\,2\,A_{2}\,+\,2\,A_{1}$}{$[2]$}
\elldata{2199}{$16$}{$A_{16}$}{$[1]$}
\elldata{2200}{$16$}{$A_{15}\,+\,A_{1}$}{$[2], \,[1]$}
\elldata{2201}{$16$}{$A_{14}\,+\,A_{2}$}{$[3], \,[1]$}
\elldata{2202}{$16$}{$A_{14}\,+\,2\,A_{1}$}{$[1]$}
\elldata{2203}{$16$}{$A_{13}\,+\,A_{3}$}{$[1]$}
\elldata{2204}{$16$}{$A_{13}\,+\,A_{2}\,+\,A_{1}$}{$[2], \,[1]$}
\elldata{2205}{$16$}{$A_{13}\,+\,3\,A_{1}$}{$[2], \,[1]$}
\elldata{2206}{$16$}{$A_{12}\,+\,A_{4}$}{$[1]$}
\elldata{2207}{$16$}{$A_{12}\,+\,A_{3}\,+\,A_{1}$}{$[1]$}
\elldata{2208}{$16$}{$A_{12}\,+\,2\,A_{2}$}{$[1]$}
\elldata{2209}{$16$}{$A_{12}\,+\,A_{2}\,+\,2\,A_{1}$}{$[1]$}
\elldata{2210}{$16$}{$A_{12}\,+\,4\,A_{1}$}{$[1]$}
\elldata{2211}{$16$}{$A_{11}\,+\,A_{5}$}{$[3], \,[1]$}
\elldata{2212}{$16$}{$A_{11}\,+\,A_{4}\,+\,A_{1}$}{$[1]$}
\elldata{2213}{$16$}{$A_{11}\,+\,A_{3}\,+\,A_{2}$}{$[2], \,[1]$}
\elldata{2214}{$16$}{$A_{11}\,+\,A_{3}\,+\,2\,A_{1}$}{$[4], \,[2], \,[1]$}
\elldata{2215}{$16$}{$A_{11}\,+\,2\,A_{2}\,+\,A_{1}$}{$[3], \,[1]$}
\elldata{2216}{$16$}{$A_{11}\,+\,A_{2}\,+\,3\,A_{1}$}{$[2], \,[1]$}
\elldata{2217}{$16$}{$A_{11}\,+\,5\,A_{1}$}{$[2]$}
\elldata{2218}{$16$}{$A_{10}\,+\,A_{6}$}{$[1]$}
\elldata{2219}{$16$}{$A_{10}\,+\,A_{5}\,+\,A_{1}$}{$[1]$}
\elldata{2220}{$16$}{$A_{10}\,+\,A_{4}\,+\,A_{2}$}{$[1]$}
\elldata{2221}{$16$}{$A_{10}\,+\,A_{4}\,+\,2\,A_{1}$}{$[1]$}
\elldata{2222}{$16$}{$A_{10}\,+\,2\,A_{3}$}{$[1]$}
\elldata{2223}{$16$}{$A_{10}\,+\,A_{3}\,+\,A_{2}\,+\,A_{1}$}{$[1]$}
\elldata{2224}{$16$}{$A_{10}\,+\,A_{3}\,+\,3\,A_{1}$}{$[1]$}
\elldata{2225}{$16$}{$A_{10}\,+\,3\,A_{2}$}{$[1]$}
\elldata{2226}{$16$}{$A_{10}\,+\,2\,A_{2}\,+\,2\,A_{1}$}{$[1]$}
\elldata{2227}{$16$}{$A_{10}\,+\,A_{2}\,+\,4\,A_{1}$}{$[1]$}
\elldata{2228}{$16$}{$A_{9}\,+\,A_{7}$}{$[1]$}
\elldata{2229}{$16$}{$A_{9}\,+\,A_{6}\,+\,A_{1}$}{$[1]$}
\elldata{2230}{$16$}{$A_{9}\,+\,A_{5}\,+\,A_{2}$}{$[2], \,[1]$}
\elldata{2231}{$16$}{$A_{9}\,+\,A_{5}\,+\,2\,A_{1}$}{$[2], \,[1]$}
\elldata{2232}{$16$}{$A_{9}\,+\,A_{4}\,+\,A_{3}$}{$[1]$}
\elldata{2233}{$16$}{$A_{9}\,+\,A_{4}\,+\,A_{2}\,+\,A_{1}$}{$[1]$}
\elldata{2234}{$16$}{$A_{9}\,+\,A_{4}\,+\,3\,A_{1}$}{$[2], \,[1]$}
\elldata{2235}{$16$}{$A_{9}\,+\,2\,A_{3}\,+\,A_{1}$}{$[2], \,[1]$}
\elldata{2236}{$16$}{$A_{9}\,+\,A_{3}\,+\,2\,A_{2}$}{$[1]$}
\elldata{2237}{$16$}{$A_{9}\,+\,A_{3}\,+\,A_{2}\,+\,2\,A_{1}$}{$[2], \,[1]$}
\elldata{2238}{$16$}{$A_{9}\,+\,A_{3}\,+\,4\,A_{1}$}{$[2]$}
\elldata{2239}{$16$}{$A_{9}\,+\,3\,A_{2}\,+\,A_{1}$}{$[1]$}
\elldata{2240}{$16$}{$A_{9}\,+\,2\,A_{2}\,+\,3\,A_{1}$}{$[2], \,[1]$}
\elldata{2241}{$16$}{$A_{9}\,+\,A_{2}\,+\,5\,A_{1}$}{$[2]$}
\elldata{2242}{$16$}{$2\,A_{8}$}{$[3], \,[1]$}
\elldata{2243}{$16$}{$A_{8}\,+\,A_{7}\,+\,A_{1}$}{$[1]$}
\elldata{2244}{$16$}{$A_{8}\,+\,A_{6}\,+\,A_{2}$}{$[1]$}
\elldata{2245}{$16$}{$A_{8}\,+\,A_{6}\,+\,2\,A_{1}$}{$[1]$}
\elldata{2246}{$16$}{$A_{8}\,+\,A_{5}\,+\,A_{3}$}{$[1]$}
\elldata{2247}{$16$}{$A_{8}\,+\,A_{5}\,+\,A_{2}\,+\,A_{1}$}{$[3], \,[1]$}
\elldata{2248}{$16$}{$A_{8}\,+\,A_{5}\,+\,3\,A_{1}$}{$[1]$}
\elldata{2249}{$16$}{$A_{8}\,+\,2\,A_{4}$}{$[1]$}
\elldata{2250}{$16$}{$A_{8}\,+\,A_{4}\,+\,A_{3}\,+\,A_{1}$}{$[1]$}
\elldata{2251}{$16$}{$A_{8}\,+\,A_{4}\,+\,2\,A_{2}$}{$[1]$}
\elldata{2252}{$16$}{$A_{8}\,+\,A_{4}\,+\,A_{2}\,+\,2\,A_{1}$}{$[1]$}
\elldata{2253}{$16$}{$A_{8}\,+\,A_{4}\,+\,4\,A_{1}$}{$[1]$}
\elldata{2254}{$16$}{$A_{8}\,+\,2\,A_{3}\,+\,A_{2}$}{$[1]$}
\elldata{2255}{$16$}{$A_{8}\,+\,2\,A_{3}\,+\,2\,A_{1}$}{$[1]$}
\elldata{2256}{$16$}{$A_{8}\,+\,A_{3}\,+\,2\,A_{2}\,+\,A_{1}$}{$[1]$}
\elldata{2257}{$16$}{$A_{8}\,+\,A_{3}\,+\,A_{2}\,+\,3\,A_{1}$}{$[1]$}
\elldata{2258}{$16$}{$A_{8}\,+\,4\,A_{2}$}{$[3]$}
\elldata{2259}{$16$}{$A_{8}\,+\,3\,A_{2}\,+\,2\,A_{1}$}{$[3], \,[1]$}
\elldata{2260}{$16$}{$A_{8}\,+\,2\,A_{2}\,+\,4\,A_{1}$}{$[1]$}
\elldata{2261}{$16$}{$2\,A_{7}\,+\,A_{2}$}{$[2], \,[1]$}
\elldata{2262}{$16$}{$2\,A_{7}\,+\,2\,A_{1}$}{$[4], \,[2], \,[1]$}
\elldata{2263}{$16$}{$A_{7}\,+\,A_{6}\,+\,A_{3}$}{$[1]$}
\elldata{2264}{$16$}{$A_{7}\,+\,A_{6}\,+\,A_{2}\,+\,A_{1}$}{$[1]$}
\elldata{2265}{$16$}{$A_{7}\,+\,A_{6}\,+\,3\,A_{1}$}{$[1]$}
\elldata{2266}{$16$}{$A_{7}\,+\,A_{5}\,+\,A_{4}$}{$[1]$}
\elldata{2267}{$16$}{$A_{7}\,+\,A_{5}\,+\,A_{3}\,+\,A_{1}$}{$[2], \,[1]$}
\elldata{2268}{$16$}{$A_{7}\,+\,A_{5}\,+\,2\,A_{2}$}{$[1]$}
\elldata{2269}{$16$}{$A_{7}\,+\,A_{5}\,+\,A_{2}\,+\,2\,A_{1}$}{$[2], \,[1]$}
\elldata{2270}{$16$}{$A_{7}\,+\,A_{5}\,+\,4\,A_{1}$}{$[2]$}
\elldata{2271}{$16$}{$A_{7}\,+\,2\,A_{4}\,+\,A_{1}$}{$[1]$}
\elldata{2272}{$16$}{$A_{7}\,+\,A_{4}\,+\,A_{3}\,+\,A_{2}$}{$[1]$}
\elldata{2273}{$16$}{$A_{7}\,+\,A_{4}\,+\,A_{3}\,+\,2\,A_{1}$}{$[2], \,[1]$}
\elldata{2274}{$16$}{$A_{7}\,+\,A_{4}\,+\,2\,A_{2}\,+\,A_{1}$}{$[1]$}
\elldata{2275}{$16$}{$A_{7}\,+\,A_{4}\,+\,A_{2}\,+\,3\,A_{1}$}{$[1]$}
\elldata{2276}{$16$}{$A_{7}\,+\,A_{4}\,+\,5\,A_{1}$}{$[2]$}
\elldata{2277}{$16$}{$A_{7}\,+\,3\,A_{3}$}{$[4]$}
\elldata{2278}{$16$}{$A_{7}\,+\,2\,A_{3}\,+\,A_{2}\,+\,A_{1}$}{$[2], \,[1]$}
\elldata{2279}{$16$}{$A_{7}\,+\,2\,A_{3}\,+\,3\,A_{1}$}{$[4], \,[2]$}
\elldata{2280}{$16$}{$A_{7}\,+\,A_{3}\,+\,3\,A_{2}$}{$[1]$}
\elldata{2281}{$16$}{$A_{7}\,+\,A_{3}\,+\,2\,A_{2}\,+\,2\,A_{1}$}{$[2], \,[1]$}
\elldata{2282}{$16$}{$A_{7}\,+\,A_{3}\,+\,A_{2}\,+\,4\,A_{1}$}{$[2]$}
\elldata{2283}{$16$}{$A_{7}\,+\,A_{3}\,+\,6\,A_{1}$}{$[2, 2]$}
\elldata{2284}{$16$}{$A_{7}\,+\,4\,A_{2}\,+\,A_{1}$}{$[1]$}
\elldata{2285}{$16$}{$A_{7}\,+\,3\,A_{2}\,+\,3\,A_{1}$}{$[1]$}
\elldata{2286}{$16$}{$A_{7}\,+\,2\,A_{2}\,+\,5\,A_{1}$}{$[2]$}
\elldata{2287}{$16$}{$2\,A_{6}\,+\,A_{4}$}{$[1]$}
\elldata{2288}{$16$}{$2\,A_{6}\,+\,A_{3}\,+\,A_{1}$}{$[1]$}
\elldata{2289}{$16$}{$2\,A_{6}\,+\,2\,A_{2}$}{$[1]$}
\elldata{2290}{$16$}{$2\,A_{6}\,+\,A_{2}\,+\,2\,A_{1}$}{$[1]$}
\elldata{2291}{$16$}{$2\,A_{6}\,+\,4\,A_{1}$}{$[1]$}
\elldata{2292}{$16$}{$A_{6}\,+\,2\,A_{5}$}{$[1]$}
\elldata{2293}{$16$}{$A_{6}\,+\,A_{5}\,+\,A_{4}\,+\,A_{1}$}{$[1]$}
\elldata{2294}{$16$}{$A_{6}\,+\,A_{5}\,+\,A_{3}\,+\,A_{2}$}{$[1]$}
\elldata{2295}{$16$}{$A_{6}\,+\,A_{5}\,+\,A_{3}\,+\,2\,A_{1}$}{$[1]$}
\elldata{2296}{$16$}{$A_{6}\,+\,A_{5}\,+\,2\,A_{2}\,+\,A_{1}$}{$[1]$}
\elldata{2297}{$16$}{$A_{6}\,+\,A_{5}\,+\,A_{2}\,+\,3\,A_{1}$}{$[1]$}
\elldata{2298}{$16$}{$A_{6}\,+\,A_{5}\,+\,5\,A_{1}$}{$[2]$}
\elldata{2299}{$16$}{$A_{6}\,+\,2\,A_{4}\,+\,A_{2}$}{$[1]$}
\elldata{2300}{$16$}{$A_{6}\,+\,2\,A_{4}\,+\,2\,A_{1}$}{$[1]$}
\elldata{2301}{$16$}{$A_{6}\,+\,A_{4}\,+\,2\,A_{3}$}{$[1]$}
\elldata{2302}{$16$}{$A_{6}\,+\,A_{4}\,+\,A_{3}\,+\,A_{2}\,+\,A_{1}$}{$[1]$}
\elldata{2303}{$16$}{$A_{6}\,+\,A_{4}\,+\,A_{3}\,+\,3\,A_{1}$}{$[1]$}
\elldata{2304}{$16$}{$A_{6}\,+\,A_{4}\,+\,3\,A_{2}$}{$[1]$}
\elldata{2305}{$16$}{$A_{6}\,+\,A_{4}\,+\,2\,A_{2}\,+\,2\,A_{1}$}{$[1]$}
\elldata{2306}{$16$}{$A_{6}\,+\,A_{4}\,+\,A_{2}\,+\,4\,A_{1}$}{$[1]$}
\elldata{2307}{$16$}{$A_{6}\,+\,3\,A_{3}\,+\,A_{1}$}{$[1]$}
\elldata{2308}{$16$}{$A_{6}\,+\,2\,A_{3}\,+\,2\,A_{2}$}{$[1]$}
\elldata{2309}{$16$}{$A_{6}\,+\,2\,A_{3}\,+\,A_{2}\,+\,2\,A_{1}$}{$[1]$}
\elldata{2310}{$16$}{$A_{6}\,+\,2\,A_{3}\,+\,4\,A_{1}$}{$[2]$}
\elldata{2311}{$16$}{$A_{6}\,+\,A_{3}\,+\,3\,A_{2}\,+\,A_{1}$}{$[1]$}
\elldata{2312}{$16$}{$A_{6}\,+\,A_{3}\,+\,2\,A_{2}\,+\,3\,A_{1}$}{$[1]$}
\elldata{2313}{$16$}{$A_{6}\,+\,4\,A_{2}\,+\,2\,A_{1}$}{$[1]$}
\elldata{2314}{$16$}{$A_{6}\,+\,3\,A_{2}\,+\,4\,A_{1}$}{$[1]$}
\elldata{2315}{$16$}{$3\,A_{5}\,+\,A_{1}$}{$[3], \,[1]$}
\elldata{2316}{$16$}{$2\,A_{5}\,+\,A_{4}\,+\,A_{2}$}{$[1]$}
\elldata{2317}{$16$}{$2\,A_{5}\,+\,A_{4}\,+\,2\,A_{1}$}{$[2], \,[1]$}
\elldata{2318}{$16$}{$2\,A_{5}\,+\,2\,A_{3}$}{$[2], \,[1]$}
\elldata{2319}{$16$}{$2\,A_{5}\,+\,A_{3}\,+\,A_{2}\,+\,A_{1}$}{$[2], \,[1]$}
\elldata{2320}{$16$}{$2\,A_{5}\,+\,A_{3}\,+\,3\,A_{1}$}{$[2]$}
\elldata{2321}{$16$}{$2\,A_{5}\,+\,3\,A_{2}$}{$[3]$}
\elldata{2322}{$16$}{$2\,A_{5}\,+\,2\,A_{2}\,+\,2\,A_{1}$}{$[6], \,[3], \,[2], \,[1]$}
\elldata{2323}{$16$}{$2\,A_{5}\,+\,A_{2}\,+\,4\,A_{1}$}{$[2]$}
\elldata{2324}{$16$}{$2\,A_{5}\,+\,6\,A_{1}$}{$[2, 2]$}
\elldata{2325}{$16$}{$A_{5}\,+\,2\,A_{4}\,+\,A_{3}$}{$[1]$}
\elldata{2326}{$16$}{$A_{5}\,+\,2\,A_{4}\,+\,A_{2}\,+\,A_{1}$}{$[1]$}
\elldata{2327}{$16$}{$A_{5}\,+\,2\,A_{4}\,+\,3\,A_{1}$}{$[1]$}
\elldata{2328}{$16$}{$A_{5}\,+\,A_{4}\,+\,2\,A_{3}\,+\,A_{1}$}{$[2], \,[1]$}
\elldata{2329}{$16$}{$A_{5}\,+\,A_{4}\,+\,A_{3}\,+\,2\,A_{2}$}{$[1]$}
\elldata{2330}{$16$}{$A_{5}\,+\,A_{4}\,+\,A_{3}\,+\,A_{2}\,+\,2\,A_{1}$}{$[1]$}
\elldata{2331}{$16$}{$A_{5}\,+\,A_{4}\,+\,A_{3}\,+\,4\,A_{1}$}{$[2]$}
\elldata{2332}{$16$}{$A_{5}\,+\,A_{4}\,+\,3\,A_{2}\,+\,A_{1}$}{$[1]$}
\elldata{2333}{$16$}{$A_{5}\,+\,A_{4}\,+\,2\,A_{2}\,+\,3\,A_{1}$}{$[1]$}
\elldata{2334}{$16$}{$A_{5}\,+\,A_{4}\,+\,A_{2}\,+\,5\,A_{1}$}{$[2]$}
\elldata{2335}{$16$}{$A_{5}\,+\,3\,A_{3}\,+\,A_{2}$}{$[1]$}
\elldata{2336}{$16$}{$A_{5}\,+\,3\,A_{3}\,+\,2\,A_{1}$}{$[2]$}
\elldata{2337}{$16$}{$A_{5}\,+\,2\,A_{3}\,+\,2\,A_{2}\,+\,A_{1}$}{$[2], \,[1]$}
\elldata{2338}{$16$}{$A_{5}\,+\,2\,A_{3}\,+\,A_{2}\,+\,3\,A_{1}$}{$[2]$}
\elldata{2339}{$16$}{$A_{5}\,+\,2\,A_{3}\,+\,5\,A_{1}$}{$[2, 2]$}
\elldata{2340}{$16$}{$A_{5}\,+\,A_{3}\,+\,4\,A_{2}$}{$[3]$}
\elldata{2341}{$16$}{$A_{5}\,+\,A_{3}\,+\,3\,A_{2}\,+\,2\,A_{1}$}{$[1]$}
\elldata{2342}{$16$}{$A_{5}\,+\,A_{3}\,+\,2\,A_{2}\,+\,4\,A_{1}$}{$[2]$}
\elldata{2343}{$16$}{$A_{5}\,+\,5\,A_{2}\,+\,A_{1}$}{$[3]$}
\elldata{2344}{$16$}{$A_{5}\,+\,4\,A_{2}\,+\,3\,A_{1}$}{$[3]$}
\elldata{2345}{$16$}{$4\,A_{4}$}{$[5], \,[1]$}
\elldata{2346}{$16$}{$3\,A_{4}\,+\,A_{3}\,+\,A_{1}$}{$[1]$}
\elldata{2347}{$16$}{$3\,A_{4}\,+\,2\,A_{2}$}{$[1]$}
\elldata{2348}{$16$}{$3\,A_{4}\,+\,A_{2}\,+\,2\,A_{1}$}{$[1]$}
\elldata{2349}{$16$}{$3\,A_{4}\,+\,4\,A_{1}$}{$[1]$}
\elldata{2350}{$16$}{$2\,A_{4}\,+\,2\,A_{3}\,+\,A_{2}$}{$[1]$}
\elldata{2351}{$16$}{$2\,A_{4}\,+\,2\,A_{3}\,+\,2\,A_{1}$}{$[1]$}
\elldata{2352}{$16$}{$2\,A_{4}\,+\,A_{3}\,+\,2\,A_{2}\,+\,A_{1}$}{$[1]$}
\elldata{2353}{$16$}{$2\,A_{4}\,+\,A_{3}\,+\,A_{2}\,+\,3\,A_{1}$}{$[1]$}
\elldata{2354}{$16$}{$2\,A_{4}\,+\,4\,A_{2}$}{$[1]$}
\elldata{2355}{$16$}{$2\,A_{4}\,+\,3\,A_{2}\,+\,2\,A_{1}$}{$[1]$}
\elldata{2356}{$16$}{$2\,A_{4}\,+\,2\,A_{2}\,+\,4\,A_{1}$}{$[1]$}
\elldata{2357}{$16$}{$A_{4}\,+\,3\,A_{3}\,+\,A_{2}\,+\,A_{1}$}{$[1]$}
\elldata{2358}{$16$}{$A_{4}\,+\,3\,A_{3}\,+\,3\,A_{1}$}{$[2]$}
\elldata{2359}{$16$}{$A_{4}\,+\,2\,A_{3}\,+\,3\,A_{2}$}{$[1]$}
\elldata{2360}{$16$}{$A_{4}\,+\,2\,A_{3}\,+\,2\,A_{2}\,+\,2\,A_{1}$}{$[1]$}
\elldata{2361}{$16$}{$A_{4}\,+\,2\,A_{3}\,+\,A_{2}\,+\,4\,A_{1}$}{$[2]$}
\elldata{2362}{$16$}{$A_{4}\,+\,A_{3}\,+\,4\,A_{2}\,+\,A_{1}$}{$[1]$}
\elldata{2363}{$16$}{$A_{4}\,+\,A_{3}\,+\,3\,A_{2}\,+\,3\,A_{1}$}{$[1]$}
\elldata{2364}{$16$}{$A_{4}\,+\,6\,A_{2}$}{$[3]$}
\elldata{2365}{$16$}{$5\,A_{3}\,+\,A_{1}$}{$[4]$}
\elldata{2366}{$16$}{$4\,A_{3}\,+\,2\,A_{2}$}{$[2], \,[1]$}
\elldata{2367}{$16$}{$4\,A_{3}\,+\,A_{2}\,+\,2\,A_{1}$}{$[4], \,[2]$}
\elldata{2368}{$16$}{$4\,A_{3}\,+\,4\,A_{1}$}{$[4, 2], \,[2, 2]$}
\elldata{2369}{$16$}{$3\,A_{3}\,+\,3\,A_{2}\,+\,A_{1}$}{$[1]$}
\elldata{2370}{$16$}{$3\,A_{3}\,+\,2\,A_{2}\,+\,3\,A_{1}$}{$[2]$}
\elldata{2371}{$16$}{$2\,A_{3}\,+\,4\,A_{2}\,+\,2\,A_{1}$}{$[1]$}
\elldata{2372}{$16$}{$A_{3}\,+\,6\,A_{2}\,+\,A_{1}$}{$[3]$}
\elldata{2373}{$16$}{$8\,A_{2}$}{$[3, 3]$}

\vsr \elldata{No.}{rank}{$ADE$-type}{$G$}

\vsrs \elldata{2374}{$17$}{$2\,E_{8}\,+\,A_{1}$}{$[1]$}
\elldata{2375}{$17$}{$E_{8}\,+\,E_{7}\,+\,A_{2}$}{$[1]$}
\elldata{2376}{$17$}{$E_{8}\,+\,E_{7}\,+\,2\,A_{1}$}{$[1]$}
\elldata{2377}{$17$}{$E_{8}\,+\,E_{6}\,+\,A_{3}$}{$[1]$}
\elldata{2378}{$17$}{$E_{8}\,+\,E_{6}\,+\,A_{2}\,+\,A_{1}$}{$[1]$}
\elldata{2379}{$17$}{$E_{8}\,+\,E_{6}\,+\,3\,A_{1}$}{$[1]$}
\elldata{2380}{$17$}{$E_{8}\,+\,D_{9}$}{$[1]$}
\elldata{2381}{$17$}{$E_{8}\,+\,D_{8}\,+\,A_{1}$}{$[1]$}
\elldata{2382}{$17$}{$E_{8}\,+\,D_{7}\,+\,A_{2}$}{$[1]$}
\elldata{2383}{$17$}{$E_{8}\,+\,D_{7}\,+\,2\,A_{1}$}{$[1]$}
\elldata{2384}{$17$}{$E_{8}\,+\,D_{6}\,+\,A_{3}$}{$[1]$}
\elldata{2385}{$17$}{$E_{8}\,+\,D_{6}\,+\,A_{2}\,+\,A_{1}$}{$[1]$}
\elldata{2386}{$17$}{$E_{8}\,+\,D_{5}\,+\,D_{4}$}{$[1]$}
\elldata{2387}{$17$}{$E_{8}\,+\,D_{5}\,+\,A_{4}$}{$[1]$}
\elldata{2388}{$17$}{$E_{8}\,+\,D_{5}\,+\,A_{3}\,+\,A_{1}$}{$[1]$}
\elldata{2389}{$17$}{$E_{8}\,+\,D_{5}\,+\,2\,A_{2}$}{$[1]$}
\elldata{2390}{$17$}{$E_{8}\,+\,D_{5}\,+\,A_{2}\,+\,2\,A_{1}$}{$[1]$}
\elldata{2391}{$17$}{$E_{8}\,+\,D_{4}\,+\,A_{5}$}{$[1]$}
\elldata{2392}{$17$}{$E_{8}\,+\,D_{4}\,+\,A_{4}\,+\,A_{1}$}{$[1]$}
\elldata{2393}{$17$}{$E_{8}\,+\,D_{4}\,+\,2\,A_{2}\,+\,A_{1}$}{$[1]$}
\elldata{2394}{$17$}{$E_{8}\,+\,A_{9}$}{$[1]$}
\elldata{2395}{$17$}{$E_{8}\,+\,A_{8}\,+\,A_{1}$}{$[1]$}
\elldata{2396}{$17$}{$E_{8}\,+\,A_{7}\,+\,A_{2}$}{$[1]$}
\elldata{2397}{$17$}{$E_{8}\,+\,A_{7}\,+\,2\,A_{1}$}{$[1]$}
\elldata{2398}{$17$}{$E_{8}\,+\,A_{6}\,+\,A_{3}$}{$[1]$}
\elldata{2399}{$17$}{$E_{8}\,+\,A_{6}\,+\,A_{2}\,+\,A_{1}$}{$[1]$}
\elldata{2400}{$17$}{$E_{8}\,+\,A_{6}\,+\,3\,A_{1}$}{$[1]$}
\elldata{2401}{$17$}{$E_{8}\,+\,A_{5}\,+\,A_{4}$}{$[1]$}
\elldata{2402}{$17$}{$E_{8}\,+\,A_{5}\,+\,A_{3}\,+\,A_{1}$}{$[1]$}
\elldata{2403}{$17$}{$E_{8}\,+\,A_{5}\,+\,2\,A_{2}$}{$[1]$}
\elldata{2404}{$17$}{$E_{8}\,+\,A_{5}\,+\,A_{2}\,+\,2\,A_{1}$}{$[1]$}
\elldata{2405}{$17$}{$E_{8}\,+\,2\,A_{4}\,+\,A_{1}$}{$[1]$}
\elldata{2406}{$17$}{$E_{8}\,+\,A_{4}\,+\,A_{3}\,+\,A_{2}$}{$[1]$}
\elldata{2407}{$17$}{$E_{8}\,+\,A_{4}\,+\,A_{3}\,+\,2\,A_{1}$}{$[1]$}
\elldata{2408}{$17$}{$E_{8}\,+\,A_{4}\,+\,2\,A_{2}\,+\,A_{1}$}{$[1]$}
\elldata{2409}{$17$}{$E_{8}\,+\,A_{4}\,+\,A_{2}\,+\,3\,A_{1}$}{$[1]$}
\elldata{2410}{$17$}{$E_{8}\,+\,2\,A_{3}\,+\,A_{2}\,+\,A_{1}$}{$[1]$}
\elldata{2411}{$17$}{$E_{8}\,+\,A_{3}\,+\,3\,A_{2}$}{$[1]$}
\elldata{2412}{$17$}{$E_{8}\,+\,A_{3}\,+\,2\,A_{2}\,+\,2\,A_{1}$}{$[1]$}
\elldata{2413}{$17$}{$2\,E_{7}\,+\,A_{3}$}{$[2], \,[1]$}
\elldata{2414}{$17$}{$2\,E_{7}\,+\,A_{2}\,+\,A_{1}$}{$[1]$}
\elldata{2415}{$17$}{$2\,E_{7}\,+\,3\,A_{1}$}{$[2]$}
\elldata{2416}{$17$}{$E_{7}\,+\,E_{6}\,+\,D_{4}$}{$[1]$}
\elldata{2417}{$17$}{$E_{7}\,+\,E_{6}\,+\,A_{4}$}{$[1]$}
\elldata{2418}{$17$}{$E_{7}\,+\,E_{6}\,+\,A_{3}\,+\,A_{1}$}{$[1]$}
\elldata{2419}{$17$}{$E_{7}\,+\,E_{6}\,+\,2\,A_{2}$}{$[1]$}
\elldata{2420}{$17$}{$E_{7}\,+\,E_{6}\,+\,A_{2}\,+\,2\,A_{1}$}{$[1]$}
\elldata{2421}{$17$}{$E_{7}\,+\,D_{10}$}{$[2], \,[1]$}
\elldata{2422}{$17$}{$E_{7}\,+\,D_{9}\,+\,A_{1}$}{$[1]$}
\elldata{2423}{$17$}{$E_{7}\,+\,D_{8}\,+\,A_{2}$}{$[1]$}
\elldata{2424}{$17$}{$E_{7}\,+\,D_{8}\,+\,2\,A_{1}$}{$[2]$}
\elldata{2425}{$17$}{$E_{7}\,+\,D_{7}\,+\,A_{3}$}{$[1]$}
\elldata{2426}{$17$}{$E_{7}\,+\,D_{7}\,+\,A_{2}\,+\,A_{1}$}{$[1]$}
\elldata{2427}{$17$}{$E_{7}\,+\,D_{7}\,+\,3\,A_{1}$}{$[2]$}
\elldata{2428}{$17$}{$E_{7}\,+\,D_{6}\,+\,D_{4}$}{$[2]$}
\elldata{2429}{$17$}{$E_{7}\,+\,D_{6}\,+\,A_{4}$}{$[1]$}
\elldata{2430}{$17$}{$E_{7}\,+\,D_{6}\,+\,A_{3}\,+\,A_{1}$}{$[2]$}
\elldata{2431}{$17$}{$E_{7}\,+\,D_{6}\,+\,2\,A_{2}$}{$[1]$}
\elldata{2432}{$17$}{$E_{7}\,+\,D_{6}\,+\,A_{2}\,+\,2\,A_{1}$}{$[2]$}
\elldata{2433}{$17$}{$E_{7}\,+\,2\,D_{5}$}{$[1]$}
\elldata{2434}{$17$}{$E_{7}\,+\,D_{5}\,+\,D_{4}\,+\,A_{1}$}{$[2]$}
\elldata{2435}{$17$}{$E_{7}\,+\,D_{5}\,+\,A_{5}$}{$[2], \,[1]$}
\elldata{2436}{$17$}{$E_{7}\,+\,D_{5}\,+\,A_{4}\,+\,A_{1}$}{$[1]$}
\elldata{2437}{$17$}{$E_{7}\,+\,D_{5}\,+\,A_{3}\,+\,A_{2}$}{$[1]$}
\elldata{2438}{$17$}{$E_{7}\,+\,D_{5}\,+\,A_{3}\,+\,2\,A_{1}$}{$[2]$}
\elldata{2439}{$17$}{$E_{7}\,+\,D_{5}\,+\,2\,A_{2}\,+\,A_{1}$}{$[1]$}
\elldata{2440}{$17$}{$E_{7}\,+\,D_{4}\,+\,A_{6}$}{$[1]$}
\elldata{2441}{$17$}{$E_{7}\,+\,D_{4}\,+\,A_{5}\,+\,A_{1}$}{$[2]$}
\elldata{2442}{$17$}{$E_{7}\,+\,D_{4}\,+\,A_{4}\,+\,A_{2}$}{$[1]$}
\elldata{2443}{$17$}{$E_{7}\,+\,D_{4}\,+\,A_{3}\,+\,A_{2}\,+\,A_{1}$}{$[2]$}
\elldata{2444}{$17$}{$E_{7}\,+\,A_{10}$}{$[1]$}
\elldata{2445}{$17$}{$E_{7}\,+\,A_{9}\,+\,A_{1}$}{$[2], \,[1]$}
\elldata{2446}{$17$}{$E_{7}\,+\,A_{8}\,+\,A_{2}$}{$[1]$}
\elldata{2447}{$17$}{$E_{7}\,+\,A_{8}\,+\,2\,A_{1}$}{$[1]$}
\elldata{2448}{$17$}{$E_{7}\,+\,A_{7}\,+\,A_{3}$}{$[1]$}
\elldata{2449}{$17$}{$E_{7}\,+\,A_{7}\,+\,A_{2}\,+\,A_{1}$}{$[2], \,[1]$}
\elldata{2450}{$17$}{$E_{7}\,+\,A_{7}\,+\,3\,A_{1}$}{$[2]$}
\elldata{2451}{$17$}{$E_{7}\,+\,A_{6}\,+\,A_{4}$}{$[1]$}
\elldata{2452}{$17$}{$E_{7}\,+\,A_{6}\,+\,A_{3}\,+\,A_{1}$}{$[1]$}
\elldata{2453}{$17$}{$E_{7}\,+\,A_{6}\,+\,2\,A_{2}$}{$[1]$}
\elldata{2454}{$17$}{$E_{7}\,+\,A_{6}\,+\,A_{2}\,+\,2\,A_{1}$}{$[1]$}
\elldata{2455}{$17$}{$E_{7}\,+\,2\,A_{5}$}{$[1]$}
\elldata{2456}{$17$}{$E_{7}\,+\,A_{5}\,+\,A_{4}\,+\,A_{1}$}{$[1]$}
\elldata{2457}{$17$}{$E_{7}\,+\,A_{5}\,+\,A_{3}\,+\,A_{2}$}{$[2], \,[1]$}
\elldata{2458}{$17$}{$E_{7}\,+\,A_{5}\,+\,A_{3}\,+\,2\,A_{1}$}{$[2]$}
\elldata{2459}{$17$}{$E_{7}\,+\,A_{5}\,+\,2\,A_{2}\,+\,A_{1}$}{$[1]$}
\elldata{2460}{$17$}{$E_{7}\,+\,A_{5}\,+\,A_{2}\,+\,3\,A_{1}$}{$[2]$}
\elldata{2461}{$17$}{$E_{7}\,+\,2\,A_{4}\,+\,A_{2}$}{$[1]$}
\elldata{2462}{$17$}{$E_{7}\,+\,2\,A_{4}\,+\,2\,A_{1}$}{$[1]$}
\elldata{2463}{$17$}{$E_{7}\,+\,A_{4}\,+\,2\,A_{3}$}{$[1]$}
\elldata{2464}{$17$}{$E_{7}\,+\,A_{4}\,+\,A_{3}\,+\,A_{2}\,+\,A_{1}$}{$[1]$}
\elldata{2465}{$17$}{$E_{7}\,+\,A_{4}\,+\,A_{3}\,+\,3\,A_{1}$}{$[2]$}
\elldata{2466}{$17$}{$E_{7}\,+\,A_{4}\,+\,3\,A_{2}$}{$[1]$}
\elldata{2467}{$17$}{$E_{7}\,+\,A_{4}\,+\,2\,A_{2}\,+\,2\,A_{1}$}{$[1]$}
\elldata{2468}{$17$}{$E_{7}\,+\,3\,A_{3}\,+\,A_{1}$}{$[2]$}
\elldata{2469}{$17$}{$E_{7}\,+\,2\,A_{3}\,+\,2\,A_{2}$}{$[1]$}
\elldata{2470}{$17$}{$E_{7}\,+\,2\,A_{3}\,+\,A_{2}\,+\,2\,A_{1}$}{$[2]$}
\elldata{2471}{$17$}{$E_{7}\,+\,A_{3}\,+\,3\,A_{2}\,+\,A_{1}$}{$[1]$}
\elldata{2472}{$17$}{$2\,E_{6}\,+\,D_{5}$}{$[1]$}
\elldata{2473}{$17$}{$2\,E_{6}\,+\,D_{4}\,+\,A_{1}$}{$[1]$}
\elldata{2474}{$17$}{$2\,E_{6}\,+\,A_{5}$}{$[3], \,[1]$}
\elldata{2475}{$17$}{$2\,E_{6}\,+\,A_{4}\,+\,A_{1}$}{$[1]$}
\elldata{2476}{$17$}{$2\,E_{6}\,+\,A_{3}\,+\,A_{2}$}{$[1]$}
\elldata{2477}{$17$}{$2\,E_{6}\,+\,A_{3}\,+\,2\,A_{1}$}{$[1]$}
\elldata{2478}{$17$}{$2\,E_{6}\,+\,2\,A_{2}\,+\,A_{1}$}{$[3]$}
\elldata{2479}{$17$}{$E_{6}\,+\,D_{11}$}{$[1]$}
\elldata{2480}{$17$}{$E_{6}\,+\,D_{10}\,+\,A_{1}$}{$[1]$}
\elldata{2481}{$17$}{$E_{6}\,+\,D_{9}\,+\,A_{2}$}{$[1]$}
\elldata{2482}{$17$}{$E_{6}\,+\,D_{9}\,+\,2\,A_{1}$}{$[1]$}
\elldata{2483}{$17$}{$E_{6}\,+\,D_{8}\,+\,A_{3}$}{$[1]$}
\elldata{2484}{$17$}{$E_{6}\,+\,D_{8}\,+\,A_{2}\,+\,A_{1}$}{$[1]$}
\elldata{2485}{$17$}{$E_{6}\,+\,D_{7}\,+\,D_{4}$}{$[1]$}
\elldata{2486}{$17$}{$E_{6}\,+\,D_{7}\,+\,A_{4}$}{$[1]$}
\elldata{2487}{$17$}{$E_{6}\,+\,D_{7}\,+\,A_{3}\,+\,A_{1}$}{$[1]$}
\elldata{2488}{$17$}{$E_{6}\,+\,D_{7}\,+\,A_{2}\,+\,2\,A_{1}$}{$[1]$}
\elldata{2489}{$17$}{$E_{6}\,+\,D_{6}\,+\,D_{5}$}{$[1]$}
\elldata{2490}{$17$}{$E_{6}\,+\,D_{6}\,+\,A_{5}$}{$[1]$}
\elldata{2491}{$17$}{$E_{6}\,+\,D_{6}\,+\,A_{4}\,+\,A_{1}$}{$[1]$}
\elldata{2492}{$17$}{$E_{6}\,+\,D_{6}\,+\,A_{3}\,+\,A_{2}$}{$[1]$}
\elldata{2493}{$17$}{$E_{6}\,+\,D_{6}\,+\,2\,A_{2}\,+\,A_{1}$}{$[1]$}
\elldata{2494}{$17$}{$E_{6}\,+\,2\,D_{5}\,+\,A_{1}$}{$[1]$}
\elldata{2495}{$17$}{$E_{6}\,+\,D_{5}\,+\,D_{4}\,+\,A_{2}$}{$[1]$}
\elldata{2496}{$17$}{$E_{6}\,+\,D_{5}\,+\,A_{6}$}{$[1]$}
\elldata{2497}{$17$}{$E_{6}\,+\,D_{5}\,+\,A_{5}\,+\,A_{1}$}{$[1]$}
\elldata{2498}{$17$}{$E_{6}\,+\,D_{5}\,+\,A_{4}\,+\,A_{2}$}{$[1]$}
\elldata{2499}{$17$}{$E_{6}\,+\,D_{5}\,+\,A_{4}\,+\,2\,A_{1}$}{$[1]$}
\elldata{2500}{$17$}{$E_{6}\,+\,D_{5}\,+\,2\,A_{3}$}{$[1]$}
\elldata{2501}{$17$}{$E_{6}\,+\,D_{5}\,+\,A_{3}\,+\,A_{2}\,+\,A_{1}$}{$[1]$}
\elldata{2502}{$17$}{$E_{6}\,+\,D_{4}\,+\,A_{7}$}{$[1]$}
\elldata{2503}{$17$}{$E_{6}\,+\,D_{4}\,+\,A_{6}\,+\,A_{1}$}{$[1]$}
\elldata{2504}{$17$}{$E_{6}\,+\,D_{4}\,+\,A_{4}\,+\,A_{3}$}{$[1]$}
\elldata{2505}{$17$}{$E_{6}\,+\,D_{4}\,+\,A_{4}\,+\,A_{2}\,+\,A_{1}$}{$[1]$}
\elldata{2506}{$17$}{$E_{6}\,+\,A_{11}$}{$[3], \,[1]$}
\elldata{2507}{$17$}{$E_{6}\,+\,A_{10}\,+\,A_{1}$}{$[1]$}
\elldata{2508}{$17$}{$E_{6}\,+\,A_{9}\,+\,A_{2}$}{$[1]$}
\elldata{2509}{$17$}{$E_{6}\,+\,A_{9}\,+\,2\,A_{1}$}{$[1]$}
\elldata{2510}{$17$}{$E_{6}\,+\,A_{8}\,+\,A_{3}$}{$[1]$}
\elldata{2511}{$17$}{$E_{6}\,+\,A_{8}\,+\,A_{2}\,+\,A_{1}$}{$[3], \,[1]$}
\elldata{2512}{$17$}{$E_{6}\,+\,A_{8}\,+\,3\,A_{1}$}{$[1]$}
\elldata{2513}{$17$}{$E_{6}\,+\,A_{7}\,+\,A_{4}$}{$[1]$}
\elldata{2514}{$17$}{$E_{6}\,+\,A_{7}\,+\,A_{3}\,+\,A_{1}$}{$[1]$}
\elldata{2515}{$17$}{$E_{6}\,+\,A_{7}\,+\,2\,A_{2}$}{$[1]$}
\elldata{2516}{$17$}{$E_{6}\,+\,A_{7}\,+\,A_{2}\,+\,2\,A_{1}$}{$[1]$}
\elldata{2517}{$17$}{$E_{6}\,+\,A_{6}\,+\,A_{5}$}{$[1]$}
\elldata{2518}{$17$}{$E_{6}\,+\,A_{6}\,+\,A_{4}\,+\,A_{1}$}{$[1]$}
\elldata{2519}{$17$}{$E_{6}\,+\,A_{6}\,+\,A_{3}\,+\,A_{2}$}{$[1]$}
\elldata{2520}{$17$}{$E_{6}\,+\,A_{6}\,+\,A_{3}\,+\,2\,A_{1}$}{$[1]$}
\elldata{2521}{$17$}{$E_{6}\,+\,A_{6}\,+\,2\,A_{2}\,+\,A_{1}$}{$[1]$}
\elldata{2522}{$17$}{$E_{6}\,+\,A_{6}\,+\,A_{2}\,+\,3\,A_{1}$}{$[1]$}
\elldata{2523}{$17$}{$E_{6}\,+\,2\,A_{5}\,+\,A_{1}$}{$[3], \,[1]$}
\elldata{2524}{$17$}{$E_{6}\,+\,A_{5}\,+\,A_{4}\,+\,A_{2}$}{$[1]$}
\elldata{2525}{$17$}{$E_{6}\,+\,A_{5}\,+\,A_{4}\,+\,2\,A_{1}$}{$[1]$}
\elldata{2526}{$17$}{$E_{6}\,+\,A_{5}\,+\,2\,A_{3}$}{$[1]$}
\elldata{2527}{$17$}{$E_{6}\,+\,A_{5}\,+\,A_{3}\,+\,A_{2}\,+\,A_{1}$}{$[1]$}
\elldata{2528}{$17$}{$E_{6}\,+\,A_{5}\,+\,3\,A_{2}$}{$[3]$}
\elldata{2529}{$17$}{$E_{6}\,+\,A_{5}\,+\,2\,A_{2}\,+\,2\,A_{1}$}{$[3]$}
\elldata{2530}{$17$}{$E_{6}\,+\,2\,A_{4}\,+\,A_{3}$}{$[1]$}
\elldata{2531}{$17$}{$E_{6}\,+\,2\,A_{4}\,+\,A_{2}\,+\,A_{1}$}{$[1]$}
\elldata{2532}{$17$}{$E_{6}\,+\,2\,A_{4}\,+\,3\,A_{1}$}{$[1]$}
\elldata{2533}{$17$}{$E_{6}\,+\,A_{4}\,+\,2\,A_{3}\,+\,A_{1}$}{$[1]$}
\elldata{2534}{$17$}{$E_{6}\,+\,A_{4}\,+\,A_{3}\,+\,2\,A_{2}$}{$[1]$}
\elldata{2535}{$17$}{$E_{6}\,+\,A_{4}\,+\,A_{3}\,+\,A_{2}\,+\,2\,A_{1}$}{$[1]$}
\elldata{2536}{$17$}{$E_{6}\,+\,2\,A_{3}\,+\,2\,A_{2}\,+\,A_{1}$}{$[1]$}
\elldata{2537}{$17$}{$E_{6}\,+\,A_{3}\,+\,4\,A_{2}$}{$[3]$}
\elldata{2538}{$17$}{$D_{17}$}{$[1]$}
\elldata{2539}{$17$}{$D_{16}\,+\,A_{1}$}{$[2], \,[1]$}
\elldata{2540}{$17$}{$D_{15}\,+\,A_{2}$}{$[1]$}
\elldata{2541}{$17$}{$D_{15}\,+\,2\,A_{1}$}{$[1]$}
\elldata{2542}{$17$}{$D_{14}\,+\,A_{3}$}{$[1]$}
\elldata{2543}{$17$}{$D_{14}\,+\,A_{2}\,+\,A_{1}$}{$[2], \,[1]$}
\elldata{2544}{$17$}{$D_{14}\,+\,3\,A_{1}$}{$[2]$}
\elldata{2545}{$17$}{$D_{13}\,+\,D_{4}$}{$[1]$}
\elldata{2546}{$17$}{$D_{13}\,+\,A_{4}$}{$[1]$}
\elldata{2547}{$17$}{$D_{13}\,+\,A_{3}\,+\,A_{1}$}{$[1]$}
\elldata{2548}{$17$}{$D_{13}\,+\,2\,A_{2}$}{$[1]$}
\elldata{2549}{$17$}{$D_{13}\,+\,A_{2}\,+\,2\,A_{1}$}{$[1]$}
\elldata{2550}{$17$}{$D_{12}\,+\,D_{5}$}{$[2], \,[1]$}
\elldata{2551}{$17$}{$D_{12}\,+\,D_{4}\,+\,A_{1}$}{$[2]$}
\elldata{2552}{$17$}{$D_{12}\,+\,A_{5}$}{$[1]$}
\elldata{2553}{$17$}{$D_{12}\,+\,A_{4}\,+\,A_{1}$}{$[1]$}
\elldata{2554}{$17$}{$D_{12}\,+\,A_{3}\,+\,A_{2}$}{$[2]$}
\elldata{2555}{$17$}{$D_{12}\,+\,A_{3}\,+\,2\,A_{1}$}{$[2]$}
\elldata{2556}{$17$}{$D_{12}\,+\,2\,A_{2}\,+\,A_{1}$}{$[1]$}
\elldata{2557}{$17$}{$D_{12}\,+\,A_{2}\,+\,3\,A_{1}$}{$[2]$}
\elldata{2558}{$17$}{$D_{11}\,+\,D_{6}$}{$[1]$}
\elldata{2559}{$17$}{$D_{11}\,+\,D_{5}\,+\,A_{1}$}{$[1]$}
\elldata{2560}{$17$}{$D_{11}\,+\,A_{6}$}{$[1]$}
\elldata{2561}{$17$}{$D_{11}\,+\,A_{5}\,+\,A_{1}$}{$[1]$}
\elldata{2562}{$17$}{$D_{11}\,+\,A_{4}\,+\,A_{2}$}{$[1]$}
\elldata{2563}{$17$}{$D_{11}\,+\,A_{4}\,+\,2\,A_{1}$}{$[1]$}
\elldata{2564}{$17$}{$D_{11}\,+\,A_{3}\,+\,A_{2}\,+\,A_{1}$}{$[1]$}
\elldata{2565}{$17$}{$D_{11}\,+\,3\,A_{2}$}{$[1]$}
\elldata{2566}{$17$}{$D_{11}\,+\,2\,A_{2}\,+\,2\,A_{1}$}{$[1]$}
\elldata{2567}{$17$}{$D_{10}\,+\,D_{7}$}{$[1]$}
\elldata{2568}{$17$}{$D_{10}\,+\,D_{6}\,+\,A_{1}$}{$[2]$}
\elldata{2569}{$17$}{$D_{10}\,+\,D_{5}\,+\,A_{2}$}{$[1]$}
\elldata{2570}{$17$}{$D_{10}\,+\,D_{5}\,+\,2\,A_{1}$}{$[2]$}
\elldata{2571}{$17$}{$D_{10}\,+\,D_{4}\,+\,A_{2}\,+\,A_{1}$}{$[2]$}
\elldata{2572}{$17$}{$D_{10}\,+\,A_{7}$}{$[1]$}
\elldata{2573}{$17$}{$D_{10}\,+\,A_{6}\,+\,A_{1}$}{$[1]$}
\elldata{2574}{$17$}{$D_{10}\,+\,A_{5}\,+\,A_{2}$}{$[2], \,[1]$}
\elldata{2575}{$17$}{$D_{10}\,+\,A_{5}\,+\,2\,A_{1}$}{$[2]$}
\elldata{2576}{$17$}{$D_{10}\,+\,A_{4}\,+\,A_{3}$}{$[1]$}
\elldata{2577}{$17$}{$D_{10}\,+\,A_{4}\,+\,A_{2}\,+\,A_{1}$}{$[1]$}
\elldata{2578}{$17$}{$D_{10}\,+\,A_{4}\,+\,3\,A_{1}$}{$[2]$}
\elldata{2579}{$17$}{$D_{10}\,+\,2\,A_{3}\,+\,A_{1}$}{$[2]$}
\elldata{2580}{$17$}{$D_{10}\,+\,A_{3}\,+\,2\,A_{2}$}{$[1]$}
\elldata{2581}{$17$}{$D_{10}\,+\,A_{3}\,+\,A_{2}\,+\,2\,A_{1}$}{$[2]$}
\elldata{2582}{$17$}{$D_{10}\,+\,A_{3}\,+\,4\,A_{1}$}{$[2, 2]$}
\elldata{2583}{$17$}{$D_{10}\,+\,2\,A_{2}\,+\,3\,A_{1}$}{$[2]$}
\elldata{2584}{$17$}{$D_{9}\,+\,D_{8}$}{$[1]$}
\elldata{2585}{$17$}{$D_{9}\,+\,D_{7}\,+\,A_{1}$}{$[1]$}
\elldata{2586}{$17$}{$D_{9}\,+\,D_{6}\,+\,A_{2}$}{$[1]$}
\elldata{2587}{$17$}{$D_{9}\,+\,D_{5}\,+\,A_{3}$}{$[1]$}
\elldata{2588}{$17$}{$D_{9}\,+\,D_{5}\,+\,A_{2}\,+\,A_{1}$}{$[1]$}
\elldata{2589}{$17$}{$D_{9}\,+\,D_{4}\,+\,A_{4}$}{$[1]$}
\elldata{2590}{$17$}{$D_{9}\,+\,A_{8}$}{$[1]$}
\elldata{2591}{$17$}{$D_{9}\,+\,A_{7}\,+\,A_{1}$}{$[1]$}
\elldata{2592}{$17$}{$D_{9}\,+\,A_{6}\,+\,A_{2}$}{$[1]$}
\elldata{2593}{$17$}{$D_{9}\,+\,A_{6}\,+\,2\,A_{1}$}{$[1]$}
\elldata{2594}{$17$}{$D_{9}\,+\,A_{5}\,+\,A_{3}$}{$[1]$}
\elldata{2595}{$17$}{$D_{9}\,+\,A_{5}\,+\,A_{2}\,+\,A_{1}$}{$[1]$}
\elldata{2596}{$17$}{$D_{9}\,+\,A_{5}\,+\,3\,A_{1}$}{$[2]$}
\elldata{2597}{$17$}{$D_{9}\,+\,2\,A_{4}$}{$[1]$}
\elldata{2598}{$17$}{$D_{9}\,+\,A_{4}\,+\,A_{3}\,+\,A_{1}$}{$[1]$}
\elldata{2599}{$17$}{$D_{9}\,+\,A_{4}\,+\,2\,A_{2}$}{$[1]$}
\elldata{2600}{$17$}{$D_{9}\,+\,A_{4}\,+\,A_{2}\,+\,2\,A_{1}$}{$[1]$}
\elldata{2601}{$17$}{$D_{9}\,+\,2\,A_{3}\,+\,2\,A_{1}$}{$[2]$}
\elldata{2602}{$17$}{$D_{9}\,+\,A_{3}\,+\,2\,A_{2}\,+\,A_{1}$}{$[1]$}
\elldata{2603}{$17$}{$D_{9}\,+\,3\,A_{2}\,+\,2\,A_{1}$}{$[1]$}
\elldata{2604}{$17$}{$2\,D_{8}\,+\,A_{1}$}{$[2]$}
\elldata{2605}{$17$}{$D_{8}\,+\,D_{7}\,+\,2\,A_{1}$}{$[2]$}
\elldata{2606}{$17$}{$D_{8}\,+\,D_{6}\,+\,A_{3}$}{$[2]$}
\elldata{2607}{$17$}{$D_{8}\,+\,D_{6}\,+\,A_{2}\,+\,A_{1}$}{$[2]$}
\elldata{2608}{$17$}{$D_{8}\,+\,D_{6}\,+\,3\,A_{1}$}{$[2, 2]$}
\elldata{2609}{$17$}{$D_{8}\,+\,D_{5}\,+\,D_{4}$}{$[2]$}
\elldata{2610}{$17$}{$D_{8}\,+\,D_{5}\,+\,A_{4}$}{$[1]$}
\elldata{2611}{$17$}{$D_{8}\,+\,D_{5}\,+\,A_{3}\,+\,A_{1}$}{$[2]$}
\elldata{2612}{$17$}{$D_{8}\,+\,D_{5}\,+\,A_{2}\,+\,2\,A_{1}$}{$[2]$}
\elldata{2613}{$17$}{$D_{8}\,+\,D_{4}\,+\,A_{3}\,+\,2\,A_{1}$}{$[2, 2]$}
\elldata{2614}{$17$}{$D_{8}\,+\,A_{9}$}{$[1]$}
\elldata{2615}{$17$}{$D_{8}\,+\,A_{8}\,+\,A_{1}$}{$[1]$}
\elldata{2616}{$17$}{$D_{8}\,+\,A_{7}\,+\,A_{2}$}{$[2]$}
\elldata{2617}{$17$}{$D_{8}\,+\,A_{7}\,+\,2\,A_{1}$}{$[2]$}
\elldata{2618}{$17$}{$D_{8}\,+\,A_{6}\,+\,A_{2}\,+\,A_{1}$}{$[1]$}
\elldata{2619}{$17$}{$D_{8}\,+\,A_{5}\,+\,A_{4}$}{$[1]$}
\elldata{2620}{$17$}{$D_{8}\,+\,A_{5}\,+\,A_{3}\,+\,A_{1}$}{$[2]$}
\elldata{2621}{$17$}{$D_{8}\,+\,A_{5}\,+\,2\,A_{2}$}{$[1]$}
\elldata{2622}{$17$}{$D_{8}\,+\,A_{5}\,+\,A_{2}\,+\,2\,A_{1}$}{$[2]$}
\elldata{2623}{$17$}{$D_{8}\,+\,A_{5}\,+\,4\,A_{1}$}{$[2, 2]$}
\elldata{2624}{$17$}{$D_{8}\,+\,2\,A_{4}\,+\,A_{1}$}{$[1]$}
\elldata{2625}{$17$}{$D_{8}\,+\,A_{4}\,+\,A_{3}\,+\,2\,A_{1}$}{$[2]$}
\elldata{2626}{$17$}{$D_{8}\,+\,A_{4}\,+\,2\,A_{2}\,+\,A_{1}$}{$[1]$}
\elldata{2627}{$17$}{$D_{8}\,+\,2\,A_{3}\,+\,A_{2}\,+\,A_{1}$}{$[2]$}
\elldata{2628}{$17$}{$D_{8}\,+\,2\,A_{3}\,+\,3\,A_{1}$}{$[2, 2]$}
\elldata{2629}{$17$}{$D_{8}\,+\,A_{3}\,+\,3\,A_{2}$}{$[1]$}
\elldata{2630}{$17$}{$D_{8}\,+\,A_{3}\,+\,2\,A_{2}\,+\,2\,A_{1}$}{$[2]$}
\elldata{2631}{$17$}{$2\,D_{7}\,+\,A_{2}\,+\,A_{1}$}{$[1]$}
\elldata{2632}{$17$}{$D_{7}\,+\,D_{6}\,+\,A_{4}$}{$[1]$}
\elldata{2633}{$17$}{$D_{7}\,+\,D_{6}\,+\,A_{3}\,+\,A_{1}$}{$[2]$}
\elldata{2634}{$17$}{$D_{7}\,+\,D_{6}\,+\,2\,A_{2}$}{$[1]$}
\elldata{2635}{$17$}{$D_{7}\,+\,2\,D_{5}$}{$[1]$}
\elldata{2636}{$17$}{$D_{7}\,+\,D_{5}\,+\,A_{5}$}{$[1]$}
\elldata{2637}{$17$}{$D_{7}\,+\,D_{5}\,+\,A_{4}\,+\,A_{1}$}{$[1]$}
\elldata{2638}{$17$}{$D_{7}\,+\,D_{5}\,+\,A_{3}\,+\,2\,A_{1}$}{$[2]$}
\elldata{2639}{$17$}{$D_{7}\,+\,D_{5}\,+\,2\,A_{2}\,+\,A_{1}$}{$[1]$}
\elldata{2640}{$17$}{$D_{7}\,+\,D_{4}\,+\,A_{5}\,+\,A_{1}$}{$[2]$}
\elldata{2641}{$17$}{$D_{7}\,+\,D_{4}\,+\,3\,A_{2}$}{$[1]$}
\elldata{2642}{$17$}{$D_{7}\,+\,A_{10}$}{$[1]$}
\elldata{2643}{$17$}{$D_{7}\,+\,A_{9}\,+\,A_{1}$}{$[2], \,[1]$}
\elldata{2644}{$17$}{$D_{7}\,+\,A_{8}\,+\,A_{2}$}{$[1]$}
\elldata{2645}{$17$}{$D_{7}\,+\,A_{8}\,+\,2\,A_{1}$}{$[1]$}
\elldata{2646}{$17$}{$D_{7}\,+\,A_{7}\,+\,A_{3}$}{$[4]$}
\elldata{2647}{$17$}{$D_{7}\,+\,A_{7}\,+\,A_{2}\,+\,A_{1}$}{$[1]$}
\elldata{2648}{$17$}{$D_{7}\,+\,A_{7}\,+\,3\,A_{1}$}{$[2]$}
\elldata{2649}{$17$}{$D_{7}\,+\,A_{6}\,+\,A_{4}$}{$[1]$}
\elldata{2650}{$17$}{$D_{7}\,+\,A_{6}\,+\,A_{3}\,+\,A_{1}$}{$[1]$}
\elldata{2651}{$17$}{$D_{7}\,+\,A_{6}\,+\,2\,A_{2}$}{$[1]$}
\elldata{2652}{$17$}{$D_{7}\,+\,A_{6}\,+\,A_{2}\,+\,2\,A_{1}$}{$[1]$}
\elldata{2653}{$17$}{$D_{7}\,+\,2\,A_{5}$}{$[2], \,[1]$}
\elldata{2654}{$17$}{$D_{7}\,+\,A_{5}\,+\,A_{4}\,+\,A_{1}$}{$[1]$}
\elldata{2655}{$17$}{$D_{7}\,+\,A_{5}\,+\,A_{3}\,+\,A_{2}$}{$[1]$}
\elldata{2656}{$17$}{$D_{7}\,+\,A_{5}\,+\,A_{3}\,+\,2\,A_{1}$}{$[2]$}
\elldata{2657}{$17$}{$D_{7}\,+\,A_{5}\,+\,A_{2}\,+\,3\,A_{1}$}{$[2]$}
\elldata{2658}{$17$}{$D_{7}\,+\,2\,A_{4}\,+\,A_{2}$}{$[1]$}
\elldata{2659}{$17$}{$D_{7}\,+\,2\,A_{4}\,+\,2\,A_{1}$}{$[1]$}
\elldata{2660}{$17$}{$D_{7}\,+\,A_{4}\,+\,A_{3}\,+\,A_{2}\,+\,A_{1}$}{$[1]$}
\elldata{2661}{$17$}{$D_{7}\,+\,A_{4}\,+\,2\,A_{2}\,+\,2\,A_{1}$}{$[1]$}
\elldata{2662}{$17$}{$D_{7}\,+\,3\,A_{3}\,+\,A_{1}$}{$[4]$}
\elldata{2663}{$17$}{$D_{7}\,+\,2\,A_{3}\,+\,2\,A_{2}$}{$[1]$}
\elldata{2664}{$17$}{$D_{7}\,+\,2\,A_{3}\,+\,A_{2}\,+\,2\,A_{1}$}{$[2]$}
\elldata{2665}{$17$}{$2\,D_{6}\,+\,D_{5}$}{$[2]$}
\elldata{2666}{$17$}{$2\,D_{6}\,+\,D_{4}\,+\,A_{1}$}{$[2, 2]$}
\elldata{2667}{$17$}{$2\,D_{6}\,+\,A_{5}$}{$[2]$}
\elldata{2668}{$17$}{$2\,D_{6}\,+\,A_{3}\,+\,A_{2}$}{$[2]$}
\elldata{2669}{$17$}{$2\,D_{6}\,+\,A_{3}\,+\,2\,A_{1}$}{$[2, 2]$}
\elldata{2670}{$17$}{$D_{6}\,+\,2\,D_{5}\,+\,A_{1}$}{$[2]$}
\elldata{2671}{$17$}{$D_{6}\,+\,D_{5}\,+\,A_{6}$}{$[1]$}
\elldata{2672}{$17$}{$D_{6}\,+\,D_{5}\,+\,A_{5}\,+\,A_{1}$}{$[2]$}
\elldata{2673}{$17$}{$D_{6}\,+\,D_{5}\,+\,A_{4}\,+\,A_{2}$}{$[1]$}
\elldata{2674}{$17$}{$D_{6}\,+\,D_{5}\,+\,2\,A_{3}$}{$[2]$}
\elldata{2675}{$17$}{$D_{6}\,+\,D_{5}\,+\,A_{3}\,+\,A_{2}\,+\,A_{1}$}{$[2]$}
\elldata{2676}{$17$}{$D_{6}\,+\,D_{5}\,+\,3\,A_{2}$}{$[1]$}
\elldata{2677}{$17$}{$D_{6}\,+\,D_{4}\,+\,A_{7}$}{$[2]$}
\elldata{2678}{$17$}{$D_{6}\,+\,D_{4}\,+\,A_{5}\,+\,A_{2}$}{$[2]$}
\elldata{2679}{$17$}{$D_{6}\,+\,D_{4}\,+\,2\,A_{3}\,+\,A_{1}$}{$[2, 2]$}
\elldata{2680}{$17$}{$D_{6}\,+\,A_{11}$}{$[2], \,[1]$}
\elldata{2681}{$17$}{$D_{6}\,+\,A_{10}\,+\,A_{1}$}{$[1]$}
\elldata{2682}{$17$}{$D_{6}\,+\,A_{9}\,+\,A_{2}$}{$[2], \,[1]$}
\elldata{2683}{$17$}{$D_{6}\,+\,A_{9}\,+\,2\,A_{1}$}{$[2]$}
\elldata{2684}{$17$}{$D_{6}\,+\,A_{8}\,+\,A_{3}$}{$[1]$}
\elldata{2685}{$17$}{$D_{6}\,+\,A_{8}\,+\,A_{2}\,+\,A_{1}$}{$[1]$}
\elldata{2686}{$17$}{$D_{6}\,+\,A_{7}\,+\,A_{4}$}{$[1]$}
\elldata{2687}{$17$}{$D_{6}\,+\,A_{7}\,+\,A_{3}\,+\,A_{1}$}{$[2]$}
\elldata{2688}{$17$}{$D_{6}\,+\,A_{7}\,+\,2\,A_{2}$}{$[1]$}
\elldata{2689}{$17$}{$D_{6}\,+\,A_{7}\,+\,A_{2}\,+\,2\,A_{1}$}{$[2]$}
\elldata{2690}{$17$}{$D_{6}\,+\,A_{6}\,+\,A_{5}$}{$[1]$}
\elldata{2691}{$17$}{$D_{6}\,+\,A_{6}\,+\,A_{4}\,+\,A_{1}$}{$[1]$}
\elldata{2692}{$17$}{$D_{6}\,+\,A_{6}\,+\,A_{3}\,+\,A_{2}$}{$[1]$}
\elldata{2693}{$17$}{$D_{6}\,+\,A_{6}\,+\,2\,A_{2}\,+\,A_{1}$}{$[1]$}
\elldata{2694}{$17$}{$D_{6}\,+\,2\,A_{5}\,+\,A_{1}$}{$[2]$}
\elldata{2695}{$17$}{$D_{6}\,+\,A_{5}\,+\,A_{4}\,+\,A_{2}$}{$[1]$}
\elldata{2696}{$17$}{$D_{6}\,+\,A_{5}\,+\,A_{4}\,+\,2\,A_{1}$}{$[2]$}
\elldata{2697}{$17$}{$D_{6}\,+\,A_{5}\,+\,2\,A_{3}$}{$[2]$}
\elldata{2698}{$17$}{$D_{6}\,+\,A_{5}\,+\,A_{3}\,+\,A_{2}\,+\,A_{1}$}{$[2]$}
\elldata{2699}{$17$}{$D_{6}\,+\,A_{5}\,+\,A_{3}\,+\,3\,A_{1}$}{$[2, 2]$}
\elldata{2700}{$17$}{$D_{6}\,+\,A_{5}\,+\,2\,A_{2}\,+\,2\,A_{1}$}{$[2]$}
\elldata{2701}{$17$}{$D_{6}\,+\,2\,A_{4}\,+\,A_{3}$}{$[1]$}
\elldata{2702}{$17$}{$D_{6}\,+\,2\,A_{4}\,+\,A_{2}\,+\,A_{1}$}{$[1]$}
\elldata{2703}{$17$}{$D_{6}\,+\,A_{4}\,+\,2\,A_{3}\,+\,A_{1}$}{$[2]$}
\elldata{2704}{$17$}{$D_{6}\,+\,A_{4}\,+\,A_{3}\,+\,2\,A_{2}$}{$[1]$}
\elldata{2705}{$17$}{$D_{6}\,+\,A_{4}\,+\,3\,A_{2}\,+\,A_{1}$}{$[1]$}
\elldata{2706}{$17$}{$D_{6}\,+\,3\,A_{3}\,+\,A_{2}$}{$[2]$}
\elldata{2707}{$17$}{$D_{6}\,+\,3\,A_{3}\,+\,2\,A_{1}$}{$[2, 2]$}
\elldata{2708}{$17$}{$D_{6}\,+\,2\,A_{3}\,+\,2\,A_{2}\,+\,A_{1}$}{$[2]$}
\elldata{2709}{$17$}{$3\,D_{5}\,+\,A_{2}$}{$[1]$}
\elldata{2710}{$17$}{$2\,D_{5}\,+\,D_{4}\,+\,A_{3}$}{$[2]$}
\elldata{2711}{$17$}{$2\,D_{5}\,+\,A_{7}$}{$[4], \,[2], \,[1]$}
\elldata{2712}{$17$}{$2\,D_{5}\,+\,A_{6}\,+\,A_{1}$}{$[1]$}
\elldata{2713}{$17$}{$2\,D_{5}\,+\,A_{5}\,+\,A_{2}$}{$[1]$}
\elldata{2714}{$17$}{$2\,D_{5}\,+\,A_{5}\,+\,2\,A_{1}$}{$[2]$}
\elldata{2715}{$17$}{$2\,D_{5}\,+\,A_{4}\,+\,A_{3}$}{$[1]$}
\elldata{2716}{$17$}{$2\,D_{5}\,+\,A_{4}\,+\,A_{2}\,+\,A_{1}$}{$[1]$}
\elldata{2717}{$17$}{$2\,D_{5}\,+\,2\,A_{3}\,+\,A_{1}$}{$[4]$}
\elldata{2718}{$17$}{$D_{5}\,+\,D_{4}\,+\,A_{8}$}{$[1]$}
\elldata{2719}{$17$}{$D_{5}\,+\,D_{4}\,+\,A_{7}\,+\,A_{1}$}{$[2]$}
\elldata{2720}{$17$}{$D_{5}\,+\,D_{4}\,+\,A_{5}\,+\,A_{2}\,+\,A_{1}$}{$[2]$}
\elldata{2721}{$17$}{$D_{5}\,+\,D_{4}\,+\,2\,A_{4}$}{$[1]$}
\elldata{2722}{$17$}{$D_{5}\,+\,A_{12}$}{$[1]$}
\elldata{2723}{$17$}{$D_{5}\,+\,A_{11}\,+\,A_{1}$}{$[4], \,[2], \,[1]$}
\elldata{2724}{$17$}{$D_{5}\,+\,A_{10}\,+\,A_{2}$}{$[1]$}
\elldata{2725}{$17$}{$D_{5}\,+\,A_{10}\,+\,2\,A_{1}$}{$[1]$}
\elldata{2726}{$17$}{$D_{5}\,+\,A_{9}\,+\,A_{3}$}{$[1]$}
\elldata{2727}{$17$}{$D_{5}\,+\,A_{9}\,+\,A_{2}\,+\,A_{1}$}{$[2], \,[1]$}
\elldata{2728}{$17$}{$D_{5}\,+\,A_{9}\,+\,3\,A_{1}$}{$[2]$}
\elldata{2729}{$17$}{$D_{5}\,+\,A_{8}\,+\,A_{4}$}{$[1]$}
\elldata{2730}{$17$}{$D_{5}\,+\,A_{8}\,+\,A_{3}\,+\,A_{1}$}{$[1]$}
\elldata{2731}{$17$}{$D_{5}\,+\,A_{8}\,+\,2\,A_{2}$}{$[1]$}
\elldata{2732}{$17$}{$D_{5}\,+\,A_{8}\,+\,A_{2}\,+\,2\,A_{1}$}{$[1]$}
\elldata{2733}{$17$}{$D_{5}\,+\,A_{7}\,+\,A_{5}$}{$[1]$}
\elldata{2734}{$17$}{$D_{5}\,+\,A_{7}\,+\,A_{4}\,+\,A_{1}$}{$[1]$}
\elldata{2735}{$17$}{$D_{5}\,+\,A_{7}\,+\,A_{3}\,+\,A_{2}$}{$[2]$}
\elldata{2736}{$17$}{$D_{5}\,+\,A_{7}\,+\,A_{3}\,+\,2\,A_{1}$}{$[4], \,[2]$}
\elldata{2737}{$17$}{$D_{5}\,+\,A_{7}\,+\,2\,A_{2}\,+\,A_{1}$}{$[1]$}
\elldata{2738}{$17$}{$D_{5}\,+\,A_{7}\,+\,A_{2}\,+\,3\,A_{1}$}{$[2]$}
\elldata{2739}{$17$}{$D_{5}\,+\,2\,A_{6}$}{$[1]$}
\elldata{2740}{$17$}{$D_{5}\,+\,A_{6}\,+\,A_{5}\,+\,A_{1}$}{$[1]$}
\elldata{2741}{$17$}{$D_{5}\,+\,A_{6}\,+\,A_{4}\,+\,A_{2}$}{$[1]$}
\elldata{2742}{$17$}{$D_{5}\,+\,A_{6}\,+\,A_{4}\,+\,2\,A_{1}$}{$[1]$}
\elldata{2743}{$17$}{$D_{5}\,+\,A_{6}\,+\,A_{3}\,+\,A_{2}\,+\,A_{1}$}{$[1]$}
\elldata{2744}{$17$}{$D_{5}\,+\,A_{6}\,+\,3\,A_{2}$}{$[1]$}
\elldata{2745}{$17$}{$D_{5}\,+\,A_{6}\,+\,2\,A_{2}\,+\,2\,A_{1}$}{$[1]$}
\elldata{2746}{$17$}{$D_{5}\,+\,2\,A_{5}\,+\,A_{2}$}{$[2], \,[1]$}
\elldata{2747}{$17$}{$D_{5}\,+\,2\,A_{5}\,+\,2\,A_{1}$}{$[2]$}
\elldata{2748}{$17$}{$D_{5}\,+\,A_{5}\,+\,A_{4}\,+\,A_{3}$}{$[1]$}
\elldata{2749}{$17$}{$D_{5}\,+\,A_{5}\,+\,A_{4}\,+\,A_{2}\,+\,A_{1}$}{$[1]$}
\elldata{2750}{$17$}{$D_{5}\,+\,A_{5}\,+\,A_{4}\,+\,3\,A_{1}$}{$[2]$}
\elldata{2751}{$17$}{$D_{5}\,+\,A_{5}\,+\,2\,A_{3}\,+\,A_{1}$}{$[2]$}
\elldata{2752}{$17$}{$D_{5}\,+\,A_{5}\,+\,A_{3}\,+\,2\,A_{2}$}{$[1]$}
\elldata{2753}{$17$}{$D_{5}\,+\,A_{5}\,+\,A_{3}\,+\,A_{2}\,+\,2\,A_{1}$}{$[2]$}
\elldata{2754}{$17$}{$D_{5}\,+\,3\,A_{4}$}{$[1]$}
\elldata{2755}{$17$}{$D_{5}\,+\,2\,A_{4}\,+\,A_{3}\,+\,A_{1}$}{$[1]$}
\elldata{2756}{$17$}{$D_{5}\,+\,2\,A_{4}\,+\,2\,A_{2}$}{$[1]$}
\elldata{2757}{$17$}{$D_{5}\,+\,2\,A_{4}\,+\,A_{2}\,+\,2\,A_{1}$}{$[1]$}
\elldata{2758}{$17$}{$D_{5}\,+\,A_{4}\,+\,2\,A_{3}\,+\,2\,A_{1}$}{$[2]$}
\elldata{2759}{$17$}{$D_{5}\,+\,A_{4}\,+\,A_{3}\,+\,2\,A_{2}\,+\,A_{1}$}{$[1]$}
\elldata{2760}{$17$}{$D_{5}\,+\,3\,A_{3}\,+\,A_{2}\,+\,A_{1}$}{$[4]$}
\elldata{2761}{$17$}{$D_{5}\,+\,2\,A_{3}\,+\,3\,A_{2}$}{$[1]$}
\elldata{2762}{$17$}{$D_{4}\,+\,A_{13}$}{$[1]$}
\elldata{2763}{$17$}{$D_{4}\,+\,A_{12}\,+\,A_{1}$}{$[1]$}
\elldata{2764}{$17$}{$D_{4}\,+\,A_{11}\,+\,A_{2}$}{$[2]$}
\elldata{2765}{$17$}{$D_{4}\,+\,A_{11}\,+\,2\,A_{1}$}{$[2]$}
\elldata{2766}{$17$}{$D_{4}\,+\,A_{10}\,+\,A_{2}\,+\,A_{1}$}{$[1]$}
\elldata{2767}{$17$}{$D_{4}\,+\,A_{9}\,+\,A_{4}$}{$[1]$}
\elldata{2768}{$17$}{$D_{4}\,+\,A_{9}\,+\,A_{3}\,+\,A_{1}$}{$[2]$}
\elldata{2769}{$17$}{$D_{4}\,+\,A_{9}\,+\,2\,A_{2}$}{$[1]$}
\elldata{2770}{$17$}{$D_{4}\,+\,A_{9}\,+\,A_{2}\,+\,2\,A_{1}$}{$[2]$}
\elldata{2771}{$17$}{$D_{4}\,+\,A_{8}\,+\,A_{5}$}{$[1]$}
\elldata{2772}{$17$}{$D_{4}\,+\,A_{8}\,+\,A_{4}\,+\,A_{1}$}{$[1]$}
\elldata{2773}{$17$}{$D_{4}\,+\,A_{7}\,+\,A_{4}\,+\,2\,A_{1}$}{$[2]$}
\elldata{2774}{$17$}{$D_{4}\,+\,A_{7}\,+\,A_{3}\,+\,A_{2}\,+\,A_{1}$}{$[2]$}
\elldata{2775}{$17$}{$D_{4}\,+\,A_{7}\,+\,2\,A_{2}\,+\,2\,A_{1}$}{$[2]$}
\elldata{2776}{$17$}{$D_{4}\,+\,2\,A_{6}\,+\,A_{1}$}{$[1]$}
\elldata{2777}{$17$}{$D_{4}\,+\,A_{6}\,+\,A_{5}\,+\,A_{2}$}{$[1]$}
\elldata{2778}{$17$}{$D_{4}\,+\,A_{6}\,+\,A_{4}\,+\,A_{2}\,+\,A_{1}$}{$[1]$}
\elldata{2779}{$17$}{$D_{4}\,+\,A_{6}\,+\,A_{3}\,+\,2\,A_{2}$}{$[1]$}
\elldata{2780}{$17$}{$D_{4}\,+\,2\,A_{5}\,+\,A_{3}$}{$[2]$}
\elldata{2781}{$17$}{$D_{4}\,+\,2\,A_{5}\,+\,3\,A_{1}$}{$[2, 2]$}
\elldata{2782}{$17$}{$D_{4}\,+\,A_{5}\,+\,2\,A_{4}$}{$[1]$}
\elldata{2783}{$17$}{$D_{4}\,+\,A_{5}\,+\,A_{4}\,+\,A_{3}\,+\,A_{1}$}{$[2]$}
\elldata{2784}{$17$}{$D_{4}\,+\,A_{5}\,+\,2\,A_{3}\,+\,2\,A_{1}$}{$[2, 2]$}
\elldata{2785}{$17$}{$D_{4}\,+\,2\,A_{4}\,+\,2\,A_{2}\,+\,A_{1}$}{$[1]$}
\elldata{2786}{$17$}{$D_{4}\,+\,3\,A_{3}\,+\,2\,A_{2}$}{$[2]$}
\elldata{2787}{$17$}{$A_{17}$}{$[3], \,[1]$}
\elldata{2788}{$17$}{$A_{16}\,+\,A_{1}$}{$[1]$}
\elldata{2789}{$17$}{$A_{15}\,+\,A_{2}$}{$[2], \,[1]$}
\elldata{2790}{$17$}{$A_{15}\,+\,2\,A_{1}$}{$[4], \,[2], \,[1]$}
\elldata{2791}{$17$}{$A_{14}\,+\,A_{3}$}{$[1]$}
\elldata{2792}{$17$}{$A_{14}\,+\,A_{2}\,+\,A_{1}$}{$[3], \,[1]$}
\elldata{2793}{$17$}{$A_{14}\,+\,3\,A_{1}$}{$[1]$}
\elldata{2794}{$17$}{$A_{13}\,+\,A_{4}$}{$[1]$}
\elldata{2795}{$17$}{$A_{13}\,+\,A_{3}\,+\,A_{1}$}{$[2], \,[1]$}
\elldata{2796}{$17$}{$A_{13}\,+\,2\,A_{2}$}{$[1]$}
\elldata{2797}{$17$}{$A_{13}\,+\,A_{2}\,+\,2\,A_{1}$}{$[2], \,[1]$}
\elldata{2798}{$17$}{$A_{13}\,+\,4\,A_{1}$}{$[2]$}
\elldata{2799}{$17$}{$A_{12}\,+\,A_{5}$}{$[1]$}
\elldata{2800}{$17$}{$A_{12}\,+\,A_{4}\,+\,A_{1}$}{$[1]$}
\elldata{2801}{$17$}{$A_{12}\,+\,A_{3}\,+\,A_{2}$}{$[1]$}
\elldata{2802}{$17$}{$A_{12}\,+\,A_{3}\,+\,2\,A_{1}$}{$[1]$}
\elldata{2803}{$17$}{$A_{12}\,+\,2\,A_{2}\,+\,A_{1}$}{$[1]$}
\elldata{2804}{$17$}{$A_{12}\,+\,A_{2}\,+\,3\,A_{1}$}{$[1]$}
\elldata{2805}{$17$}{$A_{11}\,+\,A_{6}$}{$[1]$}
\elldata{2806}{$17$}{$A_{11}\,+\,A_{5}\,+\,A_{1}$}{$[3], \,[1]$}
\elldata{2807}{$17$}{$A_{11}\,+\,A_{4}\,+\,A_{2}$}{$[1]$}
\elldata{2808}{$17$}{$A_{11}\,+\,A_{4}\,+\,2\,A_{1}$}{$[2], \,[1]$}
\elldata{2809}{$17$}{$A_{11}\,+\,2\,A_{3}$}{$[4]$}
\elldata{2810}{$17$}{$A_{11}\,+\,A_{3}\,+\,A_{2}\,+\,A_{1}$}{$[2], \,[1]$}
\elldata{2811}{$17$}{$A_{11}\,+\,A_{3}\,+\,3\,A_{1}$}{$[4], \,[2]$}
\elldata{2812}{$17$}{$A_{11}\,+\,3\,A_{2}$}{$[3]$}
\elldata{2813}{$17$}{$A_{11}\,+\,2\,A_{2}\,+\,2\,A_{1}$}{$[6], \,[3], \,[2], \,[1]$}
\elldata{2814}{$17$}{$A_{11}\,+\,A_{2}\,+\,4\,A_{1}$}{$[2]$}
\elldata{2815}{$17$}{$A_{10}\,+\,A_{7}$}{$[1]$}
\elldata{2816}{$17$}{$A_{10}\,+\,A_{6}\,+\,A_{1}$}{$[1]$}
\elldata{2817}{$17$}{$A_{10}\,+\,A_{5}\,+\,A_{2}$}{$[1]$}
\elldata{2818}{$17$}{$A_{10}\,+\,A_{5}\,+\,2\,A_{1}$}{$[1]$}
\elldata{2819}{$17$}{$A_{10}\,+\,A_{4}\,+\,A_{3}$}{$[1]$}
\elldata{2820}{$17$}{$A_{10}\,+\,A_{4}\,+\,A_{2}\,+\,A_{1}$}{$[1]$}
\elldata{2821}{$17$}{$A_{10}\,+\,A_{4}\,+\,3\,A_{1}$}{$[1]$}
\elldata{2822}{$17$}{$A_{10}\,+\,2\,A_{3}\,+\,A_{1}$}{$[1]$}
\elldata{2823}{$17$}{$A_{10}\,+\,A_{3}\,+\,2\,A_{2}$}{$[1]$}
\elldata{2824}{$17$}{$A_{10}\,+\,A_{3}\,+\,A_{2}\,+\,2\,A_{1}$}{$[1]$}
\elldata{2825}{$17$}{$A_{10}\,+\,3\,A_{2}\,+\,A_{1}$}{$[1]$}
\elldata{2826}{$17$}{$A_{10}\,+\,2\,A_{2}\,+\,3\,A_{1}$}{$[1]$}
\elldata{2827}{$17$}{$A_{9}\,+\,A_{8}$}{$[1]$}
\elldata{2828}{$17$}{$A_{9}\,+\,A_{7}\,+\,A_{1}$}{$[1]$}
\elldata{2829}{$17$}{$A_{9}\,+\,A_{6}\,+\,A_{2}$}{$[1]$}
\elldata{2830}{$17$}{$A_{9}\,+\,A_{6}\,+\,2\,A_{1}$}{$[1]$}
\elldata{2831}{$17$}{$A_{9}\,+\,A_{5}\,+\,A_{3}$}{$[2], \,[1]$}
\elldata{2832}{$17$}{$A_{9}\,+\,A_{5}\,+\,A_{2}\,+\,A_{1}$}{$[2], \,[1]$}
\elldata{2833}{$17$}{$A_{9}\,+\,A_{5}\,+\,3\,A_{1}$}{$[2]$}
\elldata{2834}{$17$}{$A_{9}\,+\,2\,A_{4}$}{$[5], \,[1]$}
\elldata{2835}{$17$}{$A_{9}\,+\,A_{4}\,+\,A_{3}\,+\,A_{1}$}{$[2], \,[1]$}
\elldata{2836}{$17$}{$A_{9}\,+\,A_{4}\,+\,2\,A_{2}$}{$[1]$}
\elldata{2837}{$17$}{$A_{9}\,+\,A_{4}\,+\,A_{2}\,+\,2\,A_{1}$}{$[1]$}
\elldata{2838}{$17$}{$A_{9}\,+\,A_{4}\,+\,4\,A_{1}$}{$[2]$}
\elldata{2839}{$17$}{$A_{9}\,+\,2\,A_{3}\,+\,A_{2}$}{$[1]$}
\elldata{2840}{$17$}{$A_{9}\,+\,2\,A_{3}\,+\,2\,A_{1}$}{$[2]$}
\elldata{2841}{$17$}{$A_{9}\,+\,A_{3}\,+\,2\,A_{2}\,+\,A_{1}$}{$[2], \,[1]$}
\elldata{2842}{$17$}{$A_{9}\,+\,A_{3}\,+\,A_{2}\,+\,3\,A_{1}$}{$[2]$}
\elldata{2843}{$17$}{$A_{9}\,+\,3\,A_{2}\,+\,2\,A_{1}$}{$[1]$}
\elldata{2844}{$17$}{$A_{9}\,+\,2\,A_{2}\,+\,4\,A_{1}$}{$[2]$}
\elldata{2845}{$17$}{$2\,A_{8}\,+\,A_{1}$}{$[3], \,[1]$}
\elldata{2846}{$17$}{$A_{8}\,+\,A_{7}\,+\,A_{2}$}{$[1]$}
\elldata{2847}{$17$}{$A_{8}\,+\,A_{7}\,+\,2\,A_{1}$}{$[1]$}
\elldata{2848}{$17$}{$A_{8}\,+\,A_{6}\,+\,A_{3}$}{$[1]$}
\elldata{2849}{$17$}{$A_{8}\,+\,A_{6}\,+\,A_{2}\,+\,A_{1}$}{$[1]$}
\elldata{2850}{$17$}{$A_{8}\,+\,A_{6}\,+\,3\,A_{1}$}{$[1]$}
\elldata{2851}{$17$}{$A_{8}\,+\,A_{5}\,+\,A_{4}$}{$[1]$}
\elldata{2852}{$17$}{$A_{8}\,+\,A_{5}\,+\,A_{3}\,+\,A_{1}$}{$[1]$}
\elldata{2853}{$17$}{$A_{8}\,+\,A_{5}\,+\,2\,A_{2}$}{$[3]$}
\elldata{2854}{$17$}{$A_{8}\,+\,A_{5}\,+\,A_{2}\,+\,2\,A_{1}$}{$[3], \,[1]$}
\elldata{2855}{$17$}{$A_{8}\,+\,2\,A_{4}\,+\,A_{1}$}{$[1]$}
\elldata{2856}{$17$}{$A_{8}\,+\,A_{4}\,+\,A_{3}\,+\,A_{2}$}{$[1]$}
\elldata{2857}{$17$}{$A_{8}\,+\,A_{4}\,+\,A_{3}\,+\,2\,A_{1}$}{$[1]$}
\elldata{2858}{$17$}{$A_{8}\,+\,A_{4}\,+\,2\,A_{2}\,+\,A_{1}$}{$[1]$}
\elldata{2859}{$17$}{$A_{8}\,+\,A_{4}\,+\,A_{2}\,+\,3\,A_{1}$}{$[1]$}
\elldata{2860}{$17$}{$A_{8}\,+\,2\,A_{3}\,+\,A_{2}\,+\,A_{1}$}{$[1]$}
\elldata{2861}{$17$}{$A_{8}\,+\,A_{3}\,+\,3\,A_{2}$}{$[3]$}
\elldata{2862}{$17$}{$A_{8}\,+\,A_{3}\,+\,2\,A_{2}\,+\,2\,A_{1}$}{$[1]$}
\elldata{2863}{$17$}{$A_{8}\,+\,4\,A_{2}\,+\,A_{1}$}{$[3]$}
\elldata{2864}{$17$}{$A_{8}\,+\,3\,A_{2}\,+\,3\,A_{1}$}{$[3]$}
\elldata{2865}{$17$}{$2\,A_{7}\,+\,A_{3}$}{$[4]$}
\elldata{2866}{$17$}{$2\,A_{7}\,+\,A_{2}\,+\,A_{1}$}{$[2], \,[1]$}
\elldata{2867}{$17$}{$2\,A_{7}\,+\,3\,A_{1}$}{$[4]$}
\elldata{2868}{$17$}{$A_{7}\,+\,A_{6}\,+\,A_{4}$}{$[1]$}
\elldata{2869}{$17$}{$A_{7}\,+\,A_{6}\,+\,A_{3}\,+\,A_{1}$}{$[1]$}
\elldata{2870}{$17$}{$A_{7}\,+\,A_{6}\,+\,2\,A_{2}$}{$[1]$}
\elldata{2871}{$17$}{$A_{7}\,+\,A_{6}\,+\,A_{2}\,+\,2\,A_{1}$}{$[1]$}
\elldata{2872}{$17$}{$A_{7}\,+\,A_{6}\,+\,4\,A_{1}$}{$[2]$}
\elldata{2873}{$17$}{$A_{7}\,+\,2\,A_{5}$}{$[1]$}
\elldata{2874}{$17$}{$A_{7}\,+\,A_{5}\,+\,A_{4}\,+\,A_{1}$}{$[2], \,[1]$}
\elldata{2875}{$17$}{$A_{7}\,+\,A_{5}\,+\,A_{3}\,+\,A_{2}$}{$[1]$}
\elldata{2876}{$17$}{$A_{7}\,+\,A_{5}\,+\,A_{3}\,+\,2\,A_{1}$}{$[2]$}
\elldata{2877}{$17$}{$A_{7}\,+\,A_{5}\,+\,2\,A_{2}\,+\,A_{1}$}{$[2], \,[1]$}
\elldata{2878}{$17$}{$A_{7}\,+\,A_{5}\,+\,A_{2}\,+\,3\,A_{1}$}{$[2]$}
\elldata{2879}{$17$}{$A_{7}\,+\,A_{5}\,+\,5\,A_{1}$}{$[2, 2]$}
\elldata{2880}{$17$}{$A_{7}\,+\,2\,A_{4}\,+\,A_{2}$}{$[1]$}
\elldata{2881}{$17$}{$A_{7}\,+\,2\,A_{4}\,+\,2\,A_{1}$}{$[1]$}
\elldata{2882}{$17$}{$A_{7}\,+\,A_{4}\,+\,A_{3}\,+\,A_{2}\,+\,A_{1}$}{$[1]$}
\elldata{2883}{$17$}{$A_{7}\,+\,A_{4}\,+\,A_{3}\,+\,3\,A_{1}$}{$[2]$}
\elldata{2884}{$17$}{$A_{7}\,+\,A_{4}\,+\,3\,A_{2}$}{$[1]$}
\elldata{2885}{$17$}{$A_{7}\,+\,A_{4}\,+\,2\,A_{2}\,+\,2\,A_{1}$}{$[1]$}
\elldata{2886}{$17$}{$A_{7}\,+\,A_{4}\,+\,A_{2}\,+\,4\,A_{1}$}{$[2]$}
\elldata{2887}{$17$}{$A_{7}\,+\,3\,A_{3}\,+\,A_{1}$}{$[4]$}
\elldata{2888}{$17$}{$A_{7}\,+\,2\,A_{3}\,+\,2\,A_{2}$}{$[2], \,[1]$}
\elldata{2889}{$17$}{$A_{7}\,+\,2\,A_{3}\,+\,A_{2}\,+\,2\,A_{1}$}{$[4], \,[2]$}
\elldata{2890}{$17$}{$A_{7}\,+\,2\,A_{3}\,+\,4\,A_{1}$}{$[4, 2], \,[2, 2]$}
\elldata{2891}{$17$}{$A_{7}\,+\,A_{3}\,+\,3\,A_{2}\,+\,A_{1}$}{$[1]$}
\elldata{2892}{$17$}{$A_{7}\,+\,A_{3}\,+\,2\,A_{2}\,+\,3\,A_{1}$}{$[2]$}
\elldata{2893}{$17$}{$2\,A_{6}\,+\,A_{5}$}{$[1]$}
\elldata{2894}{$17$}{$2\,A_{6}\,+\,A_{4}\,+\,A_{1}$}{$[1]$}
\elldata{2895}{$17$}{$2\,A_{6}\,+\,A_{3}\,+\,A_{2}$}{$[1]$}
\elldata{2896}{$17$}{$2\,A_{6}\,+\,A_{3}\,+\,2\,A_{1}$}{$[1]$}
\elldata{2897}{$17$}{$2\,A_{6}\,+\,2\,A_{2}\,+\,A_{1}$}{$[1]$}
\elldata{2898}{$17$}{$2\,A_{6}\,+\,A_{2}\,+\,3\,A_{1}$}{$[1]$}
\elldata{2899}{$17$}{$A_{6}\,+\,2\,A_{5}\,+\,A_{1}$}{$[1]$}
\elldata{2900}{$17$}{$A_{6}\,+\,A_{5}\,+\,A_{4}\,+\,A_{2}$}{$[1]$}
\elldata{2901}{$17$}{$A_{6}\,+\,A_{5}\,+\,A_{4}\,+\,2\,A_{1}$}{$[1]$}
\elldata{2902}{$17$}{$A_{6}\,+\,A_{5}\,+\,2\,A_{3}$}{$[1]$}
\elldata{2903}{$17$}{$A_{6}\,+\,A_{5}\,+\,A_{3}\,+\,A_{2}\,+\,A_{1}$}{$[1]$}
\elldata{2904}{$17$}{$A_{6}\,+\,A_{5}\,+\,A_{3}\,+\,3\,A_{1}$}{$[2]$}
\elldata{2905}{$17$}{$A_{6}\,+\,A_{5}\,+\,2\,A_{2}\,+\,2\,A_{1}$}{$[1]$}
\elldata{2906}{$17$}{$A_{6}\,+\,2\,A_{4}\,+\,A_{3}$}{$[1]$}
\elldata{2907}{$17$}{$A_{6}\,+\,2\,A_{4}\,+\,A_{2}\,+\,A_{1}$}{$[1]$}
\elldata{2908}{$17$}{$A_{6}\,+\,2\,A_{4}\,+\,3\,A_{1}$}{$[1]$}
\elldata{2909}{$17$}{$A_{6}\,+\,A_{4}\,+\,2\,A_{3}\,+\,A_{1}$}{$[1]$}
\elldata{2910}{$17$}{$A_{6}\,+\,A_{4}\,+\,A_{3}\,+\,2\,A_{2}$}{$[1]$}
\elldata{2911}{$17$}{$A_{6}\,+\,A_{4}\,+\,A_{3}\,+\,A_{2}\,+\,2\,A_{1}$}{$[1]$}
\elldata{2912}{$17$}{$A_{6}\,+\,A_{4}\,+\,3\,A_{2}\,+\,A_{1}$}{$[1]$}
\elldata{2913}{$17$}{$A_{6}\,+\,A_{4}\,+\,2\,A_{2}\,+\,3\,A_{1}$}{$[1]$}
\elldata{2914}{$17$}{$A_{6}\,+\,3\,A_{3}\,+\,A_{2}$}{$[1]$}
\elldata{2915}{$17$}{$A_{6}\,+\,3\,A_{3}\,+\,2\,A_{1}$}{$[2]$}
\elldata{2916}{$17$}{$A_{6}\,+\,2\,A_{3}\,+\,2\,A_{2}\,+\,A_{1}$}{$[1]$}
\elldata{2917}{$17$}{$A_{6}\,+\,A_{3}\,+\,3\,A_{2}\,+\,2\,A_{1}$}{$[1]$}
\elldata{2918}{$17$}{$3\,A_{5}\,+\,A_{2}$}{$[3]$}
\elldata{2919}{$17$}{$3\,A_{5}\,+\,2\,A_{1}$}{$[6], \,[2]$}
\elldata{2920}{$17$}{$2\,A_{5}\,+\,A_{4}\,+\,A_{3}$}{$[2], \,[1]$}
\elldata{2921}{$17$}{$2\,A_{5}\,+\,A_{4}\,+\,A_{2}\,+\,A_{1}$}{$[1]$}
\elldata{2922}{$17$}{$2\,A_{5}\,+\,A_{4}\,+\,3\,A_{1}$}{$[2]$}
\elldata{2923}{$17$}{$2\,A_{5}\,+\,2\,A_{3}\,+\,A_{1}$}{$[2]$}
\elldata{2924}{$17$}{$2\,A_{5}\,+\,A_{3}\,+\,2\,A_{2}$}{$[6], \,[3]$}
\elldata{2925}{$17$}{$2\,A_{5}\,+\,A_{3}\,+\,A_{2}\,+\,2\,A_{1}$}{$[2]$}
\elldata{2926}{$17$}{$2\,A_{5}\,+\,A_{3}\,+\,4\,A_{1}$}{$[2, 2]$}
\elldata{2927}{$17$}{$2\,A_{5}\,+\,3\,A_{2}\,+\,A_{1}$}{$[3]$}
\elldata{2928}{$17$}{$2\,A_{5}\,+\,2\,A_{2}\,+\,3\,A_{1}$}{$[6]$}
\elldata{2929}{$17$}{$A_{5}\,+\,3\,A_{4}$}{$[1]$}
\elldata{2930}{$17$}{$A_{5}\,+\,2\,A_{4}\,+\,A_{3}\,+\,A_{1}$}{$[1]$}
\elldata{2931}{$17$}{$A_{5}\,+\,2\,A_{4}\,+\,2\,A_{2}$}{$[1]$}
\elldata{2932}{$17$}{$A_{5}\,+\,2\,A_{4}\,+\,A_{2}\,+\,2\,A_{1}$}{$[1]$}
\elldata{2933}{$17$}{$A_{5}\,+\,A_{4}\,+\,2\,A_{3}\,+\,A_{2}$}{$[1]$}
\elldata{2934}{$17$}{$A_{5}\,+\,A_{4}\,+\,2\,A_{3}\,+\,2\,A_{1}$}{$[2]$}
\elldata{2935}{$17$}{$A_{5}\,+\,A_{4}\,+\,A_{3}\,+\,2\,A_{2}\,+\,A_{1}$}{$[1]$}
\elldata{2936}{$17$}{$A_{5}\,+\,A_{4}\,+\,A_{3}\,+\,A_{2}\,+\,3\,A_{1}$}{$[2]$}
\elldata{2937}{$17$}{$A_{5}\,+\,A_{4}\,+\,4\,A_{2}$}{$[3]$}
\elldata{2938}{$17$}{$A_{5}\,+\,3\,A_{3}\,+\,A_{2}\,+\,A_{1}$}{$[2]$}
\elldata{2939}{$17$}{$A_{5}\,+\,3\,A_{3}\,+\,3\,A_{1}$}{$[2, 2]$}
\elldata{2940}{$17$}{$A_{5}\,+\,2\,A_{3}\,+\,2\,A_{2}\,+\,2\,A_{1}$}{$[2]$}
\elldata{2941}{$17$}{$A_{5}\,+\,A_{3}\,+\,4\,A_{2}\,+\,A_{1}$}{$[3]$}
\elldata{2942}{$17$}{$A_{5}\,+\,6\,A_{2}$}{$[3, 3]$}
\elldata{2943}{$17$}{$4\,A_{4}\,+\,A_{1}$}{$[5]$}
\elldata{2944}{$17$}{$3\,A_{4}\,+\,A_{3}\,+\,2\,A_{1}$}{$[1]$}
\elldata{2945}{$17$}{$3\,A_{4}\,+\,A_{2}\,+\,3\,A_{1}$}{$[1]$}
\elldata{2946}{$17$}{$2\,A_{4}\,+\,2\,A_{3}\,+\,A_{2}\,+\,A_{1}$}{$[1]$}
\elldata{2947}{$17$}{$2\,A_{4}\,+\,A_{3}\,+\,3\,A_{2}$}{$[1]$}
\elldata{2948}{$17$}{$2\,A_{4}\,+\,A_{3}\,+\,2\,A_{2}\,+\,2\,A_{1}$}{$[1]$}
\elldata{2949}{$17$}{$A_{4}\,+\,3\,A_{3}\,+\,2\,A_{2}$}{$[1]$}
\elldata{2950}{$17$}{$A_{4}\,+\,3\,A_{3}\,+\,A_{2}\,+\,2\,A_{1}$}{$[2]$}
\elldata{2951}{$17$}{$A_{4}\,+\,2\,A_{3}\,+\,3\,A_{2}\,+\,A_{1}$}{$[1]$}
\elldata{2952}{$17$}{$5\,A_{3}\,+\,A_{2}$}{$[4]$}
\elldata{2953}{$17$}{$5\,A_{3}\,+\,2\,A_{1}$}{$[4, 2]$}

\vsr \elldata{No.}{rank}{$ADE$-type}{$G$}

\vsrs \elldata{2954}{$18$}{$2\,E_{8}\,+\,A_{2}$}{$[1]$}
\elldata{2955}{$18$}{$2\,E_{8}\,+\,2\,A_{1}$}{$[1]$}
\elldata{2956}{$18$}{$E_{8}\,+\,E_{7}\,+\,A_{3}$}{$[1]$}
\elldata{2957}{$18$}{$E_{8}\,+\,E_{7}\,+\,A_{2}\,+\,A_{1}$}{$[1]$}
\elldata{2958}{$18$}{$E_{8}\,+\,E_{6}\,+\,D_{4}$}{$[1]$}
\elldata{2959}{$18$}{$E_{8}\,+\,E_{6}\,+\,A_{4}$}{$[1]$}
\elldata{2960}{$18$}{$E_{8}\,+\,E_{6}\,+\,A_{3}\,+\,A_{1}$}{$[1]$}
\elldata{2961}{$18$}{$E_{8}\,+\,D_{10}$}{$[1]$}
\elldata{2962}{$18$}{$E_{8}\,+\,D_{9}\,+\,A_{1}$}{$[1]$}
\elldata{2963}{$18$}{$E_{8}\,+\,D_{7}\,+\,A_{2}\,+\,A_{1}$}{$[1]$}
\elldata{2964}{$18$}{$E_{8}\,+\,D_{6}\,+\,A_{4}$}{$[1]$}
\elldata{2965}{$18$}{$E_{8}\,+\,D_{6}\,+\,2\,A_{2}$}{$[1]$}
\elldata{2966}{$18$}{$E_{8}\,+\,2\,D_{5}$}{$[1]$}
\elldata{2967}{$18$}{$E_{8}\,+\,D_{5}\,+\,A_{5}$}{$[1]$}
\elldata{2968}{$18$}{$E_{8}\,+\,D_{5}\,+\,A_{4}\,+\,A_{1}$}{$[1]$}
\elldata{2969}{$18$}{$E_{8}\,+\,A_{10}$}{$[1]$}
\elldata{2970}{$18$}{$E_{8}\,+\,A_{9}\,+\,A_{1}$}{$[1]$}
\elldata{2971}{$18$}{$E_{8}\,+\,A_{8}\,+\,A_{2}$}{$[1]$}
\elldata{2972}{$18$}{$E_{8}\,+\,A_{8}\,+\,2\,A_{1}$}{$[1]$}
\elldata{2973}{$18$}{$E_{8}\,+\,A_{7}\,+\,A_{2}\,+\,A_{1}$}{$[1]$}
\elldata{2974}{$18$}{$E_{8}\,+\,A_{6}\,+\,A_{4}$}{$[1]$}
\elldata{2975}{$18$}{$E_{8}\,+\,A_{6}\,+\,A_{3}\,+\,A_{1}$}{$[1]$}
\elldata{2976}{$18$}{$E_{8}\,+\,A_{6}\,+\,2\,A_{2}$}{$[1]$}
\elldata{2977}{$18$}{$E_{8}\,+\,A_{6}\,+\,A_{2}\,+\,2\,A_{1}$}{$[1]$}
\elldata{2978}{$18$}{$E_{8}\,+\,2\,A_{5}$}{$[1]$}
\elldata{2979}{$18$}{$E_{8}\,+\,A_{5}\,+\,A_{4}\,+\,A_{1}$}{$[1]$}
\elldata{2980}{$18$}{$E_{8}\,+\,A_{5}\,+\,A_{3}\,+\,A_{2}$}{$[1]$}
\elldata{2981}{$18$}{$E_{8}\,+\,2\,A_{4}\,+\,2\,A_{1}$}{$[1]$}
\elldata{2982}{$18$}{$E_{8}\,+\,A_{4}\,+\,A_{3}\,+\,A_{2}\,+\,A_{1}$}{$[1]$}
\elldata{2983}{$18$}{$E_{8}\,+\,2\,A_{3}\,+\,2\,A_{2}$}{$[1]$}
\elldata{2984}{$18$}{$2\,E_{7}\,+\,D_{4}$}{$[2]$}
\elldata{2985}{$18$}{$2\,E_{7}\,+\,A_{4}$}{$[1]$}
\elldata{2986}{$18$}{$2\,E_{7}\,+\,A_{3}\,+\,A_{1}$}{$[2]$}
\elldata{2987}{$18$}{$2\,E_{7}\,+\,2\,A_{2}$}{$[1]$}
\elldata{2988}{$18$}{$E_{7}\,+\,E_{6}\,+\,D_{5}$}{$[1]$}
\elldata{2989}{$18$}{$E_{7}\,+\,E_{6}\,+\,A_{5}$}{$[1]$}
\elldata{2990}{$18$}{$E_{7}\,+\,E_{6}\,+\,A_{4}\,+\,A_{1}$}{$[1]$}
\elldata{2991}{$18$}{$E_{7}\,+\,E_{6}\,+\,A_{3}\,+\,A_{2}$}{$[1]$}
\elldata{2992}{$18$}{$E_{7}\,+\,D_{11}$}{$[1]$}
\elldata{2993}{$18$}{$E_{7}\,+\,D_{10}\,+\,A_{1}$}{$[2]$}
\elldata{2994}{$18$}{$E_{7}\,+\,D_{9}\,+\,A_{2}$}{$[1]$}
\elldata{2995}{$18$}{$E_{7}\,+\,D_{8}\,+\,A_{2}\,+\,A_{1}$}{$[2]$}
\elldata{2996}{$18$}{$E_{7}\,+\,D_{7}\,+\,A_{4}$}{$[1]$}
\elldata{2997}{$18$}{$E_{7}\,+\,D_{7}\,+\,A_{3}\,+\,A_{1}$}{$[2]$}
\elldata{2998}{$18$}{$E_{7}\,+\,D_{6}\,+\,D_{5}$}{$[2]$}
\elldata{2999}{$18$}{$E_{7}\,+\,D_{6}\,+\,A_{5}$}{$[2]$}
\elldata{3000}{$18$}{$E_{7}\,+\,D_{6}\,+\,A_{3}\,+\,A_{2}$}{$[2]$}
\elldata{3001}{$18$}{$E_{7}\,+\,D_{5}\,+\,A_{6}$}{$[1]$}
\elldata{3002}{$18$}{$E_{7}\,+\,D_{5}\,+\,A_{5}\,+\,A_{1}$}{$[2]$}
\elldata{3003}{$18$}{$E_{7}\,+\,D_{5}\,+\,A_{4}\,+\,A_{2}$}{$[1]$}
\elldata{3004}{$18$}{$E_{7}\,+\,A_{11}$}{$[1]$}
\elldata{3005}{$18$}{$E_{7}\,+\,A_{10}\,+\,A_{1}$}{$[1]$}
\elldata{3006}{$18$}{$E_{7}\,+\,A_{9}\,+\,A_{2}$}{$[2], \,[1]$}
\elldata{3007}{$18$}{$E_{7}\,+\,A_{9}\,+\,2\,A_{1}$}{$[2]$}
\elldata{3008}{$18$}{$E_{7}\,+\,A_{8}\,+\,A_{3}$}{$[1]$}
\elldata{3009}{$18$}{$E_{7}\,+\,A_{8}\,+\,A_{2}\,+\,A_{1}$}{$[1]$}
\elldata{3010}{$18$}{$E_{7}\,+\,A_{7}\,+\,A_{4}$}{$[1]$}
\elldata{3011}{$18$}{$E_{7}\,+\,A_{7}\,+\,A_{3}\,+\,A_{1}$}{$[2]$}
\elldata{3012}{$18$}{$E_{7}\,+\,A_{7}\,+\,2\,A_{2}$}{$[1]$}
\elldata{3013}{$18$}{$E_{7}\,+\,A_{7}\,+\,A_{2}\,+\,2\,A_{1}$}{$[2]$}
\elldata{3014}{$18$}{$E_{7}\,+\,A_{6}\,+\,A_{5}$}{$[1]$}
\elldata{3015}{$18$}{$E_{7}\,+\,A_{6}\,+\,A_{4}\,+\,A_{1}$}{$[1]$}
\elldata{3016}{$18$}{$E_{7}\,+\,A_{6}\,+\,A_{3}\,+\,A_{2}$}{$[1]$}
\elldata{3017}{$18$}{$E_{7}\,+\,A_{6}\,+\,2\,A_{2}\,+\,A_{1}$}{$[1]$}
\elldata{3018}{$18$}{$E_{7}\,+\,A_{5}\,+\,A_{4}\,+\,A_{2}$}{$[1]$}
\elldata{3019}{$18$}{$E_{7}\,+\,A_{5}\,+\,A_{4}\,+\,2\,A_{1}$}{$[2]$}
\elldata{3020}{$18$}{$E_{7}\,+\,A_{5}\,+\,2\,A_{3}$}{$[2]$}
\elldata{3021}{$18$}{$E_{7}\,+\,A_{5}\,+\,A_{3}\,+\,A_{2}\,+\,A_{1}$}{$[2]$}
\elldata{3022}{$18$}{$E_{7}\,+\,A_{4}\,+\,2\,A_{3}\,+\,A_{1}$}{$[2]$}
\elldata{3023}{$18$}{$E_{7}\,+\,A_{4}\,+\,A_{3}\,+\,2\,A_{2}$}{$[1]$}
\elldata{3024}{$18$}{$3\,E_{6}$}{$[3]$}
\elldata{3025}{$18$}{$2\,E_{6}\,+\,D_{6}$}{$[1]$}
\elldata{3026}{$18$}{$2\,E_{6}\,+\,A_{6}$}{$[1]$}
\elldata{3027}{$18$}{$2\,E_{6}\,+\,A_{5}\,+\,A_{1}$}{$[3]$}
\elldata{3028}{$18$}{$2\,E_{6}\,+\,2\,A_{3}$}{$[1]$}
\elldata{3029}{$18$}{$E_{6}\,+\,D_{12}$}{$[1]$}
\elldata{3030}{$18$}{$E_{6}\,+\,D_{11}\,+\,A_{1}$}{$[1]$}
\elldata{3031}{$18$}{$E_{6}\,+\,D_{9}\,+\,A_{3}$}{$[1]$}
\elldata{3032}{$18$}{$E_{6}\,+\,D_{9}\,+\,A_{2}\,+\,A_{1}$}{$[1]$}
\elldata{3033}{$18$}{$E_{6}\,+\,D_{8}\,+\,A_{4}$}{$[1]$}
\elldata{3034}{$18$}{$E_{6}\,+\,D_{7}\,+\,D_{5}$}{$[1]$}
\elldata{3035}{$18$}{$E_{6}\,+\,D_{7}\,+\,A_{4}\,+\,A_{1}$}{$[1]$}
\elldata{3036}{$18$}{$E_{6}\,+\,D_{6}\,+\,A_{6}$}{$[1]$}
\elldata{3037}{$18$}{$E_{6}\,+\,D_{6}\,+\,A_{4}\,+\,A_{2}$}{$[1]$}
\elldata{3038}{$18$}{$E_{6}\,+\,D_{5}\,+\,A_{7}$}{$[1]$}
\elldata{3039}{$18$}{$E_{6}\,+\,D_{5}\,+\,A_{6}\,+\,A_{1}$}{$[1]$}
\elldata{3040}{$18$}{$E_{6}\,+\,D_{5}\,+\,A_{4}\,+\,A_{3}$}{$[1]$}
\elldata{3041}{$18$}{$E_{6}\,+\,A_{12}$}{$[1]$}
\elldata{3042}{$18$}{$E_{6}\,+\,A_{11}\,+\,A_{1}$}{$[3], \,[1]$}
\elldata{3043}{$18$}{$E_{6}\,+\,A_{10}\,+\,A_{2}$}{$[1]$}
\elldata{3044}{$18$}{$E_{6}\,+\,A_{10}\,+\,2\,A_{1}$}{$[1]$}
\elldata{3045}{$18$}{$E_{6}\,+\,A_{9}\,+\,A_{3}$}{$[1]$}
\elldata{3046}{$18$}{$E_{6}\,+\,A_{9}\,+\,A_{2}\,+\,A_{1}$}{$[1]$}
\elldata{3047}{$18$}{$E_{6}\,+\,A_{8}\,+\,A_{4}$}{$[1]$}
\elldata{3048}{$18$}{$E_{6}\,+\,A_{8}\,+\,A_{3}\,+\,A_{1}$}{$[1]$}
\elldata{3049}{$18$}{$E_{6}\,+\,A_{8}\,+\,2\,A_{2}$}{$[3]$}
\elldata{3050}{$18$}{$E_{6}\,+\,A_{8}\,+\,A_{2}\,+\,2\,A_{1}$}{$[3]$}
\elldata{3051}{$18$}{$E_{6}\,+\,A_{7}\,+\,A_{5}$}{$[1]$}
\elldata{3052}{$18$}{$E_{6}\,+\,A_{7}\,+\,A_{4}\,+\,A_{1}$}{$[1]$}
\elldata{3053}{$18$}{$E_{6}\,+\,A_{6}\,+\,A_{5}\,+\,A_{1}$}{$[1]$}
\elldata{3054}{$18$}{$E_{6}\,+\,A_{6}\,+\,A_{4}\,+\,A_{2}$}{$[1]$}
\elldata{3055}{$18$}{$E_{6}\,+\,A_{6}\,+\,A_{4}\,+\,2\,A_{1}$}{$[1]$}
\elldata{3056}{$18$}{$E_{6}\,+\,A_{6}\,+\,A_{3}\,+\,A_{2}\,+\,A_{1}$}{$[1]$}
\elldata{3057}{$18$}{$E_{6}\,+\,2\,A_{5}\,+\,A_{2}$}{$[3]$}
\elldata{3058}{$18$}{$E_{6}\,+\,A_{5}\,+\,A_{4}\,+\,A_{3}$}{$[1]$}
\elldata{3059}{$18$}{$E_{6}\,+\,A_{5}\,+\,A_{3}\,+\,2\,A_{2}$}{$[3]$}
\elldata{3060}{$18$}{$E_{6}\,+\,2\,A_{4}\,+\,A_{3}\,+\,A_{1}$}{$[1]$}
\elldata{3061}{$18$}{$D_{18}$}{$[1]$}
\elldata{3062}{$18$}{$D_{17}\,+\,A_{1}$}{$[1]$}
\elldata{3063}{$18$}{$D_{16}\,+\,A_{2}$}{$[2]$}
\elldata{3064}{$18$}{$D_{16}\,+\,2\,A_{1}$}{$[2]$}
\elldata{3065}{$18$}{$D_{15}\,+\,A_{2}\,+\,A_{1}$}{$[1]$}
\elldata{3066}{$18$}{$D_{14}\,+\,A_{4}$}{$[1]$}
\elldata{3067}{$18$}{$D_{14}\,+\,A_{3}\,+\,A_{1}$}{$[2]$}
\elldata{3068}{$18$}{$D_{14}\,+\,2\,A_{2}$}{$[1]$}
\elldata{3069}{$18$}{$D_{14}\,+\,A_{2}\,+\,2\,A_{1}$}{$[2]$}
\elldata{3070}{$18$}{$D_{13}\,+\,D_{5}$}{$[1]$}
\elldata{3071}{$18$}{$D_{13}\,+\,A_{5}$}{$[1]$}
\elldata{3072}{$18$}{$D_{13}\,+\,A_{4}\,+\,A_{1}$}{$[1]$}
\elldata{3073}{$18$}{$D_{12}\,+\,D_{6}$}{$[2]$}
\elldata{3074}{$18$}{$D_{12}\,+\,D_{5}\,+\,A_{1}$}{$[2]$}
\elldata{3075}{$18$}{$D_{12}\,+\,A_{4}\,+\,2\,A_{1}$}{$[2]$}
\elldata{3076}{$18$}{$D_{12}\,+\,A_{3}\,+\,A_{2}\,+\,A_{1}$}{$[2]$}
\elldata{3077}{$18$}{$D_{12}\,+\,2\,A_{2}\,+\,2\,A_{1}$}{$[2]$}
\elldata{3078}{$18$}{$D_{11}\,+\,A_{6}\,+\,A_{1}$}{$[1]$}
\elldata{3079}{$18$}{$D_{11}\,+\,A_{5}\,+\,A_{2}$}{$[1]$}
\elldata{3080}{$18$}{$D_{11}\,+\,A_{4}\,+\,A_{2}\,+\,A_{1}$}{$[1]$}
\elldata{3081}{$18$}{$D_{11}\,+\,A_{3}\,+\,2\,A_{2}$}{$[1]$}
\elldata{3082}{$18$}{$D_{10}\,+\,D_{7}\,+\,A_{1}$}{$[2]$}
\elldata{3083}{$18$}{$D_{10}\,+\,D_{6}\,+\,A_{2}$}{$[2]$}
\elldata{3084}{$18$}{$D_{10}\,+\,D_{5}\,+\,A_{2}\,+\,A_{1}$}{$[2]$}
\elldata{3085}{$18$}{$D_{10}\,+\,A_{8}$}{$[1]$}
\elldata{3086}{$18$}{$D_{10}\,+\,A_{6}\,+\,A_{2}$}{$[1]$}
\elldata{3087}{$18$}{$D_{10}\,+\,A_{5}\,+\,A_{3}$}{$[2]$}
\elldata{3088}{$18$}{$D_{10}\,+\,A_{5}\,+\,3\,A_{1}$}{$[2, 2]$}
\elldata{3089}{$18$}{$D_{10}\,+\,2\,A_{4}$}{$[1]$}
\elldata{3090}{$18$}{$D_{10}\,+\,A_{4}\,+\,A_{3}\,+\,A_{1}$}{$[2]$}
\elldata{3091}{$18$}{$D_{10}\,+\,2\,A_{3}\,+\,2\,A_{1}$}{$[2, 2]$}
\elldata{3092}{$18$}{$2\,D_{9}$}{$[1]$}
\elldata{3093}{$18$}{$D_{9}\,+\,D_{5}\,+\,A_{4}$}{$[1]$}
\elldata{3094}{$18$}{$D_{9}\,+\,A_{9}$}{$[1]$}
\elldata{3095}{$18$}{$D_{9}\,+\,A_{8}\,+\,A_{1}$}{$[1]$}
\elldata{3096}{$18$}{$D_{9}\,+\,A_{7}\,+\,2\,A_{1}$}{$[2]$}
\elldata{3097}{$18$}{$D_{9}\,+\,A_{6}\,+\,A_{2}\,+\,A_{1}$}{$[1]$}
\elldata{3098}{$18$}{$D_{9}\,+\,A_{5}\,+\,A_{4}$}{$[1]$}
\elldata{3099}{$18$}{$D_{9}\,+\,A_{5}\,+\,A_{3}\,+\,A_{1}$}{$[2]$}
\elldata{3100}{$18$}{$D_{9}\,+\,A_{4}\,+\,2\,A_{2}\,+\,A_{1}$}{$[1]$}
\elldata{3101}{$18$}{$2\,D_{8}\,+\,2\,A_{1}$}{$[2, 2]$}
\elldata{3102}{$18$}{$D_{8}\,+\,D_{6}\,+\,A_{3}\,+\,A_{1}$}{$[2, 2]$}
\elldata{3103}{$18$}{$D_{8}\,+\,2\,D_{5}$}{$[2]$}
\elldata{3104}{$18$}{$D_{8}\,+\,A_{9}\,+\,A_{1}$}{$[2]$}
\elldata{3105}{$18$}{$D_{8}\,+\,A_{7}\,+\,A_{2}\,+\,A_{1}$}{$[2]$}
\elldata{3106}{$18$}{$D_{8}\,+\,A_{6}\,+\,2\,A_{2}$}{$[1]$}
\elldata{3107}{$18$}{$D_{8}\,+\,2\,A_{5}$}{$[2]$}
\elldata{3108}{$18$}{$D_{8}\,+\,A_{5}\,+\,A_{4}\,+\,A_{1}$}{$[2]$}
\elldata{3109}{$18$}{$D_{8}\,+\,A_{5}\,+\,A_{3}\,+\,2\,A_{1}$}{$[2, 2]$}
\elldata{3110}{$18$}{$D_{8}\,+\,2\,A_{3}\,+\,2\,A_{2}$}{$[2]$}
\elldata{3111}{$18$}{$2\,D_{7}\,+\,2\,A_{2}$}{$[1]$}
\elldata{3112}{$18$}{$D_{7}\,+\,D_{6}\,+\,A_{5}$}{$[2]$}
\elldata{3113}{$18$}{$D_{7}\,+\,D_{5}\,+\,A_{5}\,+\,A_{1}$}{$[2]$}
\elldata{3114}{$18$}{$D_{7}\,+\,A_{11}$}{$[4]$}
\elldata{3115}{$18$}{$D_{7}\,+\,A_{10}\,+\,A_{1}$}{$[1]$}
\elldata{3116}{$18$}{$D_{7}\,+\,A_{9}\,+\,A_{2}$}{$[1]$}
\elldata{3117}{$18$}{$D_{7}\,+\,A_{9}\,+\,2\,A_{1}$}{$[2]$}
\elldata{3118}{$18$}{$D_{7}\,+\,A_{7}\,+\,A_{3}\,+\,A_{1}$}{$[4]$}
\elldata{3119}{$18$}{$D_{7}\,+\,A_{7}\,+\,A_{2}\,+\,2\,A_{1}$}{$[2]$}
\elldata{3120}{$18$}{$D_{7}\,+\,A_{6}\,+\,A_{5}$}{$[1]$}
\elldata{3121}{$18$}{$D_{7}\,+\,A_{6}\,+\,A_{4}\,+\,A_{1}$}{$[1]$}
\elldata{3122}{$18$}{$D_{7}\,+\,A_{6}\,+\,A_{3}\,+\,A_{2}$}{$[1]$}
\elldata{3123}{$18$}{$D_{7}\,+\,2\,A_{4}\,+\,A_{2}\,+\,A_{1}$}{$[1]$}
\elldata{3124}{$18$}{$D_{7}\,+\,3\,A_{3}\,+\,A_{2}$}{$[4]$}
\elldata{3125}{$18$}{$3\,D_{6}$}{$[2, 2]$}
\elldata{3126}{$18$}{$2\,D_{6}\,+\,2\,A_{3}$}{$[2, 2]$}
\elldata{3127}{$18$}{$D_{6}\,+\,D_{5}\,+\,A_{7}$}{$[2]$}
\elldata{3128}{$18$}{$D_{6}\,+\,D_{5}\,+\,A_{5}\,+\,A_{2}$}{$[2]$}
\elldata{3129}{$18$}{$D_{6}\,+\,A_{12}$}{$[1]$}
\elldata{3130}{$18$}{$D_{6}\,+\,A_{11}\,+\,A_{1}$}{$[2]$}
\elldata{3131}{$18$}{$D_{6}\,+\,A_{10}\,+\,A_{2}$}{$[1]$}
\elldata{3132}{$18$}{$D_{6}\,+\,A_{9}\,+\,A_{3}$}{$[2]$}
\elldata{3133}{$18$}{$D_{6}\,+\,A_{9}\,+\,A_{2}\,+\,A_{1}$}{$[2]$}
\elldata{3134}{$18$}{$D_{6}\,+\,A_{8}\,+\,A_{4}$}{$[1]$}
\elldata{3135}{$18$}{$D_{6}\,+\,A_{7}\,+\,A_{4}\,+\,A_{1}$}{$[2]$}
\elldata{3136}{$18$}{$D_{6}\,+\,A_{7}\,+\,A_{3}\,+\,A_{2}$}{$[2]$}
\elldata{3137}{$18$}{$D_{6}\,+\,A_{7}\,+\,2\,A_{2}\,+\,A_{1}$}{$[2]$}
\elldata{3138}{$18$}{$D_{6}\,+\,2\,A_{6}$}{$[1]$}
\elldata{3139}{$18$}{$D_{6}\,+\,A_{6}\,+\,A_{4}\,+\,A_{2}$}{$[1]$}
\elldata{3140}{$18$}{$D_{6}\,+\,2\,A_{5}\,+\,2\,A_{1}$}{$[2, 2]$}
\elldata{3141}{$18$}{$D_{6}\,+\,A_{5}\,+\,A_{4}\,+\,A_{3}$}{$[2]$}
\elldata{3142}{$18$}{$D_{6}\,+\,A_{5}\,+\,2\,A_{3}\,+\,A_{1}$}{$[2, 2]$}
\elldata{3143}{$18$}{$D_{6}\,+\,2\,A_{4}\,+\,2\,A_{2}$}{$[1]$}
\elldata{3144}{$18$}{$2\,D_{5}\,+\,A_{8}$}{$[1]$}
\elldata{3145}{$18$}{$2\,D_{5}\,+\,A_{7}\,+\,A_{1}$}{$[4]$}
\elldata{3146}{$18$}{$2\,D_{5}\,+\,2\,A_{4}$}{$[1]$}
\elldata{3147}{$18$}{$D_{5}\,+\,A_{13}$}{$[1]$}
\elldata{3148}{$18$}{$D_{5}\,+\,A_{12}\,+\,A_{1}$}{$[1]$}
\elldata{3149}{$18$}{$D_{5}\,+\,A_{11}\,+\,A_{2}$}{$[2]$}
\elldata{3150}{$18$}{$D_{5}\,+\,A_{11}\,+\,2\,A_{1}$}{$[4]$}
\elldata{3151}{$18$}{$D_{5}\,+\,A_{10}\,+\,A_{2}\,+\,A_{1}$}{$[1]$}
\elldata{3152}{$18$}{$D_{5}\,+\,A_{9}\,+\,A_{4}$}{$[1]$}
\elldata{3153}{$18$}{$D_{5}\,+\,A_{9}\,+\,A_{3}\,+\,A_{1}$}{$[2]$}
\elldata{3154}{$18$}{$D_{5}\,+\,A_{9}\,+\,2\,A_{2}$}{$[1]$}
\elldata{3155}{$18$}{$D_{5}\,+\,A_{9}\,+\,A_{2}\,+\,2\,A_{1}$}{$[2]$}
\elldata{3156}{$18$}{$D_{5}\,+\,A_{8}\,+\,A_{5}$}{$[1]$}
\elldata{3157}{$18$}{$D_{5}\,+\,A_{8}\,+\,A_{4}\,+\,A_{1}$}{$[1]$}
\elldata{3158}{$18$}{$D_{5}\,+\,A_{7}\,+\,A_{4}\,+\,2\,A_{1}$}{$[2]$}
\elldata{3159}{$18$}{$D_{5}\,+\,A_{7}\,+\,A_{3}\,+\,A_{2}\,+\,A_{1}$}{$[4]$}
\elldata{3160}{$18$}{$D_{5}\,+\,2\,A_{6}\,+\,A_{1}$}{$[1]$}
\elldata{3161}{$18$}{$D_{5}\,+\,A_{6}\,+\,A_{5}\,+\,A_{2}$}{$[1]$}
\elldata{3162}{$18$}{$D_{5}\,+\,A_{6}\,+\,A_{4}\,+\,A_{2}\,+\,A_{1}$}{$[1]$}
\elldata{3163}{$18$}{$D_{5}\,+\,A_{6}\,+\,A_{3}\,+\,2\,A_{2}$}{$[1]$}
\elldata{3164}{$18$}{$D_{5}\,+\,2\,A_{5}\,+\,A_{3}$}{$[2]$}
\elldata{3165}{$18$}{$D_{5}\,+\,A_{5}\,+\,2\,A_{4}$}{$[1]$}
\elldata{3166}{$18$}{$D_{5}\,+\,A_{5}\,+\,A_{4}\,+\,A_{3}\,+\,A_{1}$}{$[2]$}
\elldata{3167}{$18$}{$A_{18}$}{$[1]$}
\elldata{3168}{$18$}{$A_{17}\,+\,A_{1}$}{$[3], \,[1]$}
\elldata{3169}{$18$}{$A_{16}\,+\,A_{2}$}{$[1]$}
\elldata{3170}{$18$}{$A_{16}\,+\,2\,A_{1}$}{$[1]$}
\elldata{3171}{$18$}{$A_{15}\,+\,A_{3}$}{$[4]$}
\elldata{3172}{$18$}{$A_{15}\,+\,A_{2}\,+\,A_{1}$}{$[2], \,[1]$}
\elldata{3173}{$18$}{$A_{15}\,+\,3\,A_{1}$}{$[4]$}
\elldata{3174}{$18$}{$A_{14}\,+\,A_{4}$}{$[1]$}
\elldata{3175}{$18$}{$A_{14}\,+\,A_{3}\,+\,A_{1}$}{$[1]$}
\elldata{3176}{$18$}{$A_{14}\,+\,2\,A_{2}$}{$[3]$}
\elldata{3177}{$18$}{$A_{14}\,+\,A_{2}\,+\,2\,A_{1}$}{$[3], \,[1]$}
\elldata{3178}{$18$}{$A_{13}\,+\,A_{5}$}{$[1]$}
\elldata{3179}{$18$}{$A_{13}\,+\,A_{4}\,+\,A_{1}$}{$[2], \,[1]$}
\elldata{3180}{$18$}{$A_{13}\,+\,A_{3}\,+\,A_{2}$}{$[1]$}
\elldata{3181}{$18$}{$A_{13}\,+\,A_{3}\,+\,2\,A_{1}$}{$[2]$}
\elldata{3182}{$18$}{$A_{13}\,+\,2\,A_{2}\,+\,A_{1}$}{$[2], \,[1]$}
\elldata{3183}{$18$}{$A_{13}\,+\,A_{2}\,+\,3\,A_{1}$}{$[2]$}
\elldata{3184}{$18$}{$A_{12}\,+\,A_{6}$}{$[1]$}
\elldata{3185}{$18$}{$A_{12}\,+\,A_{5}\,+\,A_{1}$}{$[1]$}
\elldata{3186}{$18$}{$A_{12}\,+\,A_{4}\,+\,A_{2}$}{$[1]$}
\elldata{3187}{$18$}{$A_{12}\,+\,A_{4}\,+\,2\,A_{1}$}{$[1]$}
\elldata{3188}{$18$}{$A_{12}\,+\,A_{3}\,+\,A_{2}\,+\,A_{1}$}{$[1]$}
\elldata{3189}{$18$}{$A_{12}\,+\,2\,A_{2}\,+\,2\,A_{1}$}{$[1]$}
\elldata{3190}{$18$}{$A_{11}\,+\,A_{6}\,+\,A_{1}$}{$[1]$}
\elldata{3191}{$18$}{$A_{11}\,+\,A_{5}\,+\,A_{2}$}{$[3]$}
\elldata{3192}{$18$}{$A_{11}\,+\,A_{5}\,+\,2\,A_{1}$}{$[6], \,[2]$}
\elldata{3193}{$18$}{$A_{11}\,+\,A_{4}\,+\,A_{2}\,+\,A_{1}$}{$[1]$}
\elldata{3194}{$18$}{$A_{11}\,+\,A_{4}\,+\,3\,A_{1}$}{$[2]$}
\elldata{3195}{$18$}{$A_{11}\,+\,2\,A_{3}\,+\,A_{1}$}{$[4]$}
\elldata{3196}{$18$}{$A_{11}\,+\,A_{3}\,+\,2\,A_{2}$}{$[6], \,[3]$}
\elldata{3197}{$18$}{$A_{11}\,+\,A_{3}\,+\,A_{2}\,+\,2\,A_{1}$}{$[4], \,[2]$}
\elldata{3198}{$18$}{$A_{11}\,+\,3\,A_{2}\,+\,A_{1}$}{$[3]$}
\elldata{3199}{$18$}{$A_{11}\,+\,2\,A_{2}\,+\,3\,A_{1}$}{$[6]$}
\elldata{3200}{$18$}{$A_{10}\,+\,A_{8}$}{$[1]$}
\elldata{3201}{$18$}{$A_{10}\,+\,A_{7}\,+\,A_{1}$}{$[1]$}
\elldata{3202}{$18$}{$A_{10}\,+\,A_{6}\,+\,A_{2}$}{$[1]$}
\elldata{3203}{$18$}{$A_{10}\,+\,A_{6}\,+\,2\,A_{1}$}{$[1]$}
\elldata{3204}{$18$}{$A_{10}\,+\,A_{5}\,+\,A_{3}$}{$[1]$}
\elldata{3205}{$18$}{$A_{10}\,+\,A_{5}\,+\,A_{2}\,+\,A_{1}$}{$[1]$}
\elldata{3206}{$18$}{$A_{10}\,+\,2\,A_{4}$}{$[1]$}
\elldata{3207}{$18$}{$A_{10}\,+\,A_{4}\,+\,A_{3}\,+\,A_{1}$}{$[1]$}
\elldata{3208}{$18$}{$A_{10}\,+\,A_{4}\,+\,2\,A_{2}$}{$[1]$}
\elldata{3209}{$18$}{$A_{10}\,+\,A_{4}\,+\,A_{2}\,+\,2\,A_{1}$}{$[1]$}
\elldata{3210}{$18$}{$A_{10}\,+\,2\,A_{3}\,+\,A_{2}$}{$[1]$}
\elldata{3211}{$18$}{$A_{10}\,+\,A_{3}\,+\,2\,A_{2}\,+\,A_{1}$}{$[1]$}
\elldata{3212}{$18$}{$2\,A_{9}$}{$[5], \,[1]$}
\elldata{3213}{$18$}{$A_{9}\,+\,A_{8}\,+\,A_{1}$}{$[1]$}
\elldata{3214}{$18$}{$A_{9}\,+\,A_{7}\,+\,A_{2}$}{$[1]$}
\elldata{3215}{$18$}{$A_{9}\,+\,A_{6}\,+\,A_{3}$}{$[1]$}
\elldata{3216}{$18$}{$A_{9}\,+\,A_{6}\,+\,A_{2}\,+\,A_{1}$}{$[1]$}
\elldata{3217}{$18$}{$A_{9}\,+\,A_{6}\,+\,3\,A_{1}$}{$[2]$}
\elldata{3218}{$18$}{$A_{9}\,+\,A_{5}\,+\,A_{4}$}{$[2], \,[1]$}
\elldata{3219}{$18$}{$A_{9}\,+\,A_{5}\,+\,A_{3}\,+\,A_{1}$}{$[2]$}
\elldata{3220}{$18$}{$A_{9}\,+\,A_{5}\,+\,A_{2}\,+\,2\,A_{1}$}{$[2]$}
\elldata{3221}{$18$}{$A_{9}\,+\,2\,A_{4}\,+\,A_{1}$}{$[5]$}
\elldata{3222}{$18$}{$A_{9}\,+\,A_{4}\,+\,A_{3}\,+\,2\,A_{1}$}{$[2]$}
\elldata{3223}{$18$}{$A_{9}\,+\,A_{4}\,+\,A_{2}\,+\,3\,A_{1}$}{$[2]$}
\elldata{3224}{$18$}{$A_{9}\,+\,2\,A_{3}\,+\,A_{2}\,+\,A_{1}$}{$[2]$}
\elldata{3225}{$18$}{$A_{9}\,+\,A_{3}\,+\,2\,A_{2}\,+\,2\,A_{1}$}{$[2]$}
\elldata{3226}{$18$}{$2\,A_{8}\,+\,2\,A_{1}$}{$[3], \,[1]$}
\elldata{3227}{$18$}{$A_{8}\,+\,A_{7}\,+\,A_{2}\,+\,A_{1}$}{$[1]$}
\elldata{3228}{$18$}{$A_{8}\,+\,A_{6}\,+\,A_{4}$}{$[1]$}
\elldata{3229}{$18$}{$A_{8}\,+\,A_{6}\,+\,A_{3}\,+\,A_{1}$}{$[1]$}
\elldata{3230}{$18$}{$A_{8}\,+\,A_{6}\,+\,A_{2}\,+\,2\,A_{1}$}{$[1]$}
\elldata{3231}{$18$}{$A_{8}\,+\,A_{5}\,+\,A_{4}\,+\,A_{1}$}{$[1]$}
\elldata{3232}{$18$}{$A_{8}\,+\,A_{5}\,+\,A_{3}\,+\,A_{2}$}{$[3]$}
\elldata{3233}{$18$}{$A_{8}\,+\,A_{5}\,+\,2\,A_{2}\,+\,A_{1}$}{$[3]$}
\elldata{3234}{$18$}{$A_{8}\,+\,2\,A_{4}\,+\,2\,A_{1}$}{$[1]$}
\elldata{3235}{$18$}{$A_{8}\,+\,A_{4}\,+\,A_{3}\,+\,A_{2}\,+\,A_{1}$}{$[1]$}
\elldata{3236}{$18$}{$A_{8}\,+\,A_{4}\,+\,3\,A_{2}$}{$[3]$}
\elldata{3237}{$18$}{$A_{8}\,+\,A_{3}\,+\,3\,A_{2}\,+\,A_{1}$}{$[3]$}
\elldata{3238}{$18$}{$2\,A_{7}\,+\,A_{3}\,+\,A_{1}$}{$[8]$}
\elldata{3239}{$18$}{$2\,A_{7}\,+\,2\,A_{2}$}{$[2], \,[1]$}
\elldata{3240}{$18$}{$2\,A_{7}\,+\,4\,A_{1}$}{$[4, 2]$}
\elldata{3241}{$18$}{$A_{7}\,+\,A_{6}\,+\,A_{5}$}{$[1]$}
\elldata{3242}{$18$}{$A_{7}\,+\,A_{6}\,+\,A_{4}\,+\,A_{1}$}{$[1]$}
\elldata{3243}{$18$}{$A_{7}\,+\,A_{6}\,+\,A_{3}\,+\,A_{2}$}{$[1]$}
\elldata{3244}{$18$}{$A_{7}\,+\,A_{6}\,+\,A_{3}\,+\,2\,A_{1}$}{$[2]$}
\elldata{3245}{$18$}{$A_{7}\,+\,A_{6}\,+\,2\,A_{2}\,+\,A_{1}$}{$[1]$}
\elldata{3246}{$18$}{$A_{7}\,+\,2\,A_{5}\,+\,A_{1}$}{$[2]$}
\elldata{3247}{$18$}{$A_{7}\,+\,A_{5}\,+\,A_{4}\,+\,A_{2}$}{$[1]$}
\elldata{3248}{$18$}{$A_{7}\,+\,A_{5}\,+\,A_{4}\,+\,2\,A_{1}$}{$[2]$}
\elldata{3249}{$18$}{$A_{7}\,+\,A_{5}\,+\,A_{3}\,+\,A_{2}\,+\,A_{1}$}{$[2]$}
\elldata{3250}{$18$}{$A_{7}\,+\,A_{5}\,+\,A_{3}\,+\,3\,A_{1}$}{$[2, 2]$}
\elldata{3251}{$18$}{$A_{7}\,+\,A_{4}\,+\,A_{3}\,+\,2\,A_{2}$}{$[1]$}
\elldata{3252}{$18$}{$A_{7}\,+\,A_{4}\,+\,A_{3}\,+\,A_{2}\,+\,2\,A_{1}$}{$[2]$}
\elldata{3253}{$18$}{$A_{7}\,+\,3\,A_{3}\,+\,A_{2}$}{$[4]$}
\elldata{3254}{$18$}{$A_{7}\,+\,3\,A_{3}\,+\,2\,A_{1}$}{$[4, 2]$}
\elldata{3255}{$18$}{$3\,A_{6}$}{$[7]$}
\elldata{3256}{$18$}{$2\,A_{6}\,+\,A_{4}\,+\,A_{2}$}{$[1]$}
\elldata{3257}{$18$}{$2\,A_{6}\,+\,2\,A_{3}$}{$[1]$}
\elldata{3258}{$18$}{$2\,A_{6}\,+\,2\,A_{2}\,+\,2\,A_{1}$}{$[1]$}
\elldata{3259}{$18$}{$A_{6}\,+\,2\,A_{5}\,+\,2\,A_{1}$}{$[2]$}
\elldata{3260}{$18$}{$A_{6}\,+\,A_{5}\,+\,A_{4}\,+\,A_{3}$}{$[1]$}
\elldata{3261}{$18$}{$A_{6}\,+\,A_{5}\,+\,A_{4}\,+\,A_{2}\,+\,A_{1}$}{$[1]$}
\elldata{3262}{$18$}{$A_{6}\,+\,A_{5}\,+\,2\,A_{3}\,+\,A_{1}$}{$[2]$}
\elldata{3263}{$18$}{$A_{6}\,+\,2\,A_{4}\,+\,A_{3}\,+\,A_{1}$}{$[1]$}
\elldata{3264}{$18$}{$A_{6}\,+\,2\,A_{4}\,+\,A_{2}\,+\,2\,A_{1}$}{$[1]$}
\elldata{3265}{$18$}{$A_{6}\,+\,A_{4}\,+\,2\,A_{3}\,+\,A_{2}$}{$[1]$}
\elldata{3266}{$18$}{$A_{6}\,+\,A_{4}\,+\,A_{3}\,+\,2\,A_{2}\,+\,A_{1}$}{$[1]$}
\elldata{3267}{$18$}{$3\,A_{5}\,+\,A_{3}$}{$[6]$}
\elldata{3268}{$18$}{$3\,A_{5}\,+\,3\,A_{1}$}{$[6, 2]$}
\elldata{3269}{$18$}{$2\,A_{5}\,+\,2\,A_{4}$}{$[1]$}
\elldata{3270}{$18$}{$2\,A_{5}\,+\,A_{4}\,+\,A_{3}\,+\,A_{1}$}{$[2]$}
\elldata{3271}{$18$}{$2\,A_{5}\,+\,A_{4}\,+\,2\,A_{2}$}{$[3]$}
\elldata{3272}{$18$}{$2\,A_{5}\,+\,2\,A_{3}\,+\,2\,A_{1}$}{$[2, 2]$}
\elldata{3273}{$18$}{$2\,A_{5}\,+\,A_{3}\,+\,2\,A_{2}\,+\,A_{1}$}{$[6]$}
\elldata{3274}{$18$}{$2\,A_{5}\,+\,4\,A_{2}$}{$[3, 3]$}
\elldata{3275}{$18$}{$A_{5}\,+\,A_{4}\,+\,2\,A_{3}\,+\,A_{2}\,+\,A_{1}$}{$[2]$}
\elldata{3276}{$18$}{$4\,A_{4}\,+\,2\,A_{1}$}{$[5]$}
\elldata{3277}{$18$}{$2\,A_{4}\,+\,2\,A_{3}\,+\,2\,A_{2}$}{$[1]$}
\elldata{3278}{$18$}{$6\,A_{3}$}{$[4, 4]$}
%%%%%%%%%%%%%%%%%%%%%%%%%%%%%%%
%
% The end of Table 1
%
%%%%%%%%%%%%%%%%%%%%%%%%%%%%%%%%
}

\vfill
\eject

\begin{center}
Table 2.
\end{center}

{\tiny
%%%%%%%%%%%%%%%%%%%%%%%%%%%%%%%
%
% Table 2
%
%%%%%%%%%%%%%%%%%%%%%%%%%%%%%%%%

\medskip 
 \vrule  \hbox{\vbox{\offinterlineskip  \noindent \hrule \Gvsp   
\hbox{ $G=\Z / (  3 )$ 
} 
\Gvsp \hrule \Gvsp 
\vbox{ 
\hbox{  $ 3\,E\sb{6}  $ , 
  $ 2\,E\sb{6}\, +\, A\sb{5}\, +\, A\sb{1}  $ , 
  $ 2\,E\sb{6}\, +\, A\sb{5}  $ , 
  $ 2\,E\sb{6}\, +\, 2\,A\sb{2}\, +\, A\sb{1}  $ , 
  $ 2\,E\sb{6}\, +\, 2\,A\sb{2}  $ , 
  $ E\sb{6}\, +\, A\sb{11}\, +\, A\sb{1}  $ , 
  $ E\sb{6}\, +\, A\sb{11}  $ , 
  $ E\sb{6}\, +\, A\sb{8}\, +\, 2\,A\sb{2}  $ , 
} 
\GGvsp 
\hbox{  $ E\sb{6}\, +\, A\sb{8}\, +\, A\sb{2}\, +\, 2\,A\sb{1}  $ , 
  $ E\sb{6}\, +\, A\sb{8}\, +\, A\sb{2}\, +\, A\sb{1}  $ , 
  $ E\sb{6}\, +\, A\sb{8}\, +\, A\sb{2}  $ , 
  $ E\sb{6}\, +\, 2\,A\sb{5}\, +\, A\sb{2}  $ , 
  $ E\sb{6}\, +\, 2\,A\sb{5}\, +\, A\sb{1}  $ , 
  $ E\sb{6}\, +\, 2\,A\sb{5}  $ , 
} 
\GGvsp 
\hbox{  $ E\sb{6}\, +\, A\sb{5}\, +\, A\sb{3}\, +\, 2\,A\sb{2}  $ , 
  $ E\sb{6}\, +\, A\sb{5}\, +\, 3\,A\sb{2}  $ , 
  $ E\sb{6}\, +\, A\sb{5}\, +\, 2\,A\sb{2}\, +\, 2\,A\sb{1}  $ , 
  $ E\sb{6}\, +\, A\sb{5}\, +\, 2\,A\sb{2}\, +\, A\sb{1}  $ , 
  $ E\sb{6}\, +\, A\sb{5}\, +\, 2\,A\sb{2}  $ , 
} 
\GGvsp 
\hbox{  $ E\sb{6}\, +\, A\sb{3}\, +\, 4\,A\sb{2}  $ , 
  $ E\sb{6}\, +\, 5\,A\sb{2}  $ , 
  $ E\sb{6}\, +\, 4\,A\sb{2}\, +\, 2\,A\sb{1}  $ , 
  $ E\sb{6}\, +\, 4\,A\sb{2}\, +\, A\sb{1}  $ , 
  $ E\sb{6}\, +\, 4\,A\sb{2}  $ , 
  $ A\sb{17}\, +\, A\sb{1}  $ , 
  $ A\sb{17}  $ , 
  $ A\sb{14}\, +\, 2\,A\sb{2}  $ , 
} 
\GGvsp 
\hbox{  $ A\sb{14}\, +\, A\sb{2}\, +\, 2\,A\sb{1}  $ , 
  $ A\sb{14}\, +\, A\sb{2}\, +\, A\sb{1}  $ , 
  $ A\sb{14}\, +\, A\sb{2}  $ , 
  $ A\sb{11}\, +\, A\sb{5}\, +\, A\sb{2}  $ , 
  $ A\sb{11}\, +\, A\sb{5}\, +\, A\sb{1}  $ , 
  $ A\sb{11}\, +\, A\sb{5}  $ , 
  $ A\sb{11}\, +\, A\sb{3}\, +\, 2\,A\sb{2}  $ , 
} 
\GGvsp 
\hbox{  $ A\sb{11}\, +\, 3\,A\sb{2}\, +\, A\sb{1}  $ , 
  $ A\sb{11}\, +\, 3\,A\sb{2}  $ , 
  $ A\sb{11}\, +\, 2\,A\sb{2}\, +\, 2\,A\sb{1}  $ , 
  $ A\sb{11}\, +\, 2\,A\sb{2}\, +\, A\sb{1}  $ , 
  $ A\sb{11}\, +\, 2\,A\sb{2}  $ , 
  $ 2\,A\sb{8}\, +\, 2\,A\sb{1}  $ , 
  $ 2\,A\sb{8}\, +\, A\sb{1}  $ , 
} 
\GGvsp 
\hbox{  $ 2\,A\sb{8}  $ , 
  $ A\sb{8}\, +\, A\sb{5}\, +\, A\sb{3}\, +\, A\sb{2}  $ , 
  $ A\sb{8}\, +\, A\sb{5}\, +\, 2\,A\sb{2}\, +\, A\sb{1}  $ , 
  $ A\sb{8}\, +\, A\sb{5}\, +\, 2\,A\sb{2}  $ , 
  $ A\sb{8}\, +\, A\sb{5}\, +\, A\sb{2}\, +\, 2\,A\sb{1}  $ , 
  $ A\sb{8}\, +\, A\sb{5}\, +\, A\sb{2}\, +\, A\sb{1}  $ , 
} 
\GGvsp 
\hbox{  $ A\sb{8}\, +\, A\sb{5}\, +\, A\sb{2}  $ , 
  $ A\sb{8}\, +\, A\sb{4}\, +\, 3\,A\sb{2}  $ , 
  $ A\sb{8}\, +\, A\sb{3}\, +\, 3\,A\sb{2}\, +\, A\sb{1}  $ , 
  $ A\sb{8}\, +\, A\sb{3}\, +\, 3\,A\sb{2}  $ , 
  $ A\sb{8}\, +\, 4\,A\sb{2}\, +\, A\sb{1}  $ , 
  $ A\sb{8}\, +\, 4\,A\sb{2}  $ , 
} 
\GGvsp 
\hbox{  $ A\sb{8}\, +\, 3\,A\sb{2}\, +\, 3\,A\sb{1}  $ , 
  $ A\sb{8}\, +\, 3\,A\sb{2}\, +\, 2\,A\sb{1}  $ , 
  $ A\sb{8}\, +\, 3\,A\sb{2}\, +\, A\sb{1}  $ , 
  $ A\sb{8}\, +\, 3\,A\sb{2}  $ , 
  $ 3\,A\sb{5}\, +\, A\sb{2}  $ , 
  $ 3\,A\sb{5}\, +\, A\sb{1}  $ , 
  $ 3\,A\sb{5}  $ , 
  $ 2\,A\sb{5}\, +\, A\sb{4}\, +\, 2\,A\sb{2}  $ , 
} 
\GGvsp 
\hbox{  $ 2\,A\sb{5}\, +\, A\sb{3}\, +\, 2\,A\sb{2}  $ , 
  $ 2\,A\sb{5}\, +\, 3\,A\sb{2}\, +\, A\sb{1}  $ , 
  $ 2\,A\sb{5}\, +\, 3\,A\sb{2}  $ , 
  $ 2\,A\sb{5}\, +\, 2\,A\sb{2}\, +\, 2\,A\sb{1}  $ , 
  $ 2\,A\sb{5}\, +\, 2\,A\sb{2}\, +\, A\sb{1}  $ , 
  $ 2\,A\sb{5}\, +\, 2\,A\sb{2}  $ , 
} 
\GGvsp 
\hbox{  $ A\sb{5}\, +\, A\sb{4}\, +\, 4\,A\sb{2}  $ , 
  $ A\sb{5}\, +\, A\sb{3}\, +\, 4\,A\sb{2}\, +\, A\sb{1}  $ , 
  $ A\sb{5}\, +\, A\sb{3}\, +\, 4\,A\sb{2}  $ , 
  $ A\sb{5}\, +\, 5\,A\sb{2}\, +\, A\sb{1}  $ , 
  $ A\sb{5}\, +\, 5\,A\sb{2}  $ , 
  $ A\sb{5}\, +\, 4\,A\sb{2}\, +\, 3\,A\sb{1}  $ , 
} 
\GGvsp 
\hbox{  $ A\sb{5}\, +\, 4\,A\sb{2}\, +\, 2\,A\sb{1}  $ , 
  $ A\sb{5}\, +\, 4\,A\sb{2}\, +\, A\sb{1}  $ , 
  $ A\sb{5}\, +\, 4\,A\sb{2}  $ , 
  $ A\sb{4}\, +\, 6\,A\sb{2}  $ , 
  $ A\sb{3}\, +\, 6\,A\sb{2}\, +\, A\sb{1}  $ , 
  $ A\sb{3}\, +\, 6\,A\sb{2}  $ , 
  $ 7\,A\sb{2}\, +\, A\sb{1}  $ , 
} 
\GGvsp 
\hbox{  $ 7\,A\sb{2}  $ , 
  $ 6\,A\sb{2}\, +\, 3\,A\sb{1}  $ , 
  $ 6\,A\sb{2}\, +\, 2\,A\sb{1}  $ , 
  $ 6\,A\sb{2}\, +\, A\sb{1}  $ , 
  $ 6\,A\sb{2}  $ } 
} 
\Gvsp \hrule }}\vrule 
\par
\medskip
\vrule  \hbox{\vbox{\offinterlineskip  \noindent \hrule \Gvsp   
\hbox{ $G=\Z / (  4 )$ 
} 
\Gvsp \hrule \Gvsp 
\vbox{ 
\hbox{  $ D\sb{7}\, +\, A\sb{11}  $ , 
  $ D\sb{7}\, +\, A\sb{7}\, +\, A\sb{3}\, +\, A\sb{1}  $ , 
  $ D\sb{7}\, +\, A\sb{7}\, +\, A\sb{3}  $ , 
  $ D\sb{7}\, +\, 3\,A\sb{3}\, +\, A\sb{2}  $ , 
  $ D\sb{7}\, +\, 3\,A\sb{3}\, +\, A\sb{1}  $ , 
  $ D\sb{7}\, +\, 3\,A\sb{3}  $ , 
  $ 2\,D\sb{5}\, +\, A\sb{7}\, +\, A\sb{1}  $ , 
} 
\GGvsp 
\hbox{  $ 2\,D\sb{5}\, +\, A\sb{7}  $ , 
  $ 2\,D\sb{5}\, +\, 2\,A\sb{3}\, +\, A\sb{1}  $ , 
  $ 2\,D\sb{5}\, +\, 2\,A\sb{3}  $ , 
  $ D\sb{5}\, +\, A\sb{11}\, +\, 2\,A\sb{1}  $ , 
  $ D\sb{5}\, +\, A\sb{11}\, +\, A\sb{1}  $ , 
  $ D\sb{5}\, +\, A\sb{7}\, +\, A\sb{3}\, +\, A\sb{2}\, +\, A\sb{1}  $ , 
} 
\GGvsp 
\hbox{  $ D\sb{5}\, +\, A\sb{7}\, +\, A\sb{3}\, +\, 2\,A\sb{1}  $ , 
  $ D\sb{5}\, +\, A\sb{7}\, +\, A\sb{3}\, +\, A\sb{1}  $ , 
  $ D\sb{5}\, +\, 3\,A\sb{3}\, +\, A\sb{2}\, +\, A\sb{1}  $ , 
  $ D\sb{5}\, +\, 3\,A\sb{3}\, +\, 2\,A\sb{1}  $ , 
  $ D\sb{5}\, +\, 3\,A\sb{3}\, +\, A\sb{1}  $ , 
} 
\GGvsp 
\hbox{  $ A\sb{15}\, +\, A\sb{3}  $ , 
  $ A\sb{15}\, +\, 3\,A\sb{1}  $ , 
  $ A\sb{15}\, +\, 2\,A\sb{1}  $ , 
  $ A\sb{11}\, +\, 2\,A\sb{3}\, +\, A\sb{1}  $ , 
  $ A\sb{11}\, +\, 2\,A\sb{3}  $ , 
  $ A\sb{11}\, +\, A\sb{3}\, +\, A\sb{2}\, +\, 2\,A\sb{1}  $ , 
  $ A\sb{11}\, +\, A\sb{3}\, +\, 3\,A\sb{1}  $ , 
} 
\GGvsp 
\hbox{  $ A\sb{11}\, +\, A\sb{3}\, +\, 2\,A\sb{1}  $ , 
  $ 2\,A\sb{7}\, +\, A\sb{3}  $ , 
  $ 2\,A\sb{7}\, +\, 3\,A\sb{1}  $ , 
  $ 2\,A\sb{7}\, +\, 2\,A\sb{1}  $ , 
  $ A\sb{7}\, +\, 3\,A\sb{3}\, +\, A\sb{2}  $ , 
  $ A\sb{7}\, +\, 3\,A\sb{3}\, +\, A\sb{1}  $ , 
  $ A\sb{7}\, +\, 3\,A\sb{3}  $ , 
} 
\GGvsp 
\hbox{  $ A\sb{7}\, +\, 2\,A\sb{3}\, +\, A\sb{2}\, +\, 2\,A\sb{1}  $ , 
  $ A\sb{7}\, +\, 2\,A\sb{3}\, +\, 3\,A\sb{1}  $ , 
  $ A\sb{7}\, +\, 2\,A\sb{3}\, +\, 2\,A\sb{1}  $ , 
  $ 5\,A\sb{3}\, +\, A\sb{2}  $ , 
  $ 5\,A\sb{3}\, +\, A\sb{1}  $ , 
  $ 5\,A\sb{3}  $ , 
  $ 4\,A\sb{3}\, +\, A\sb{2}\, +\, 2\,A\sb{1}  $ , 
} 
\GGvsp 
\hbox{  $ 4\,A\sb{3}\, +\, 3\,A\sb{1}  $ , 
  $ 4\,A\sb{3}\, +\, 2\,A\sb{1}  $ } 
} 
\Gvsp \hrule }}\vrule 
\par
\medskip 
\vrule  \hbox{\vbox{\offinterlineskip  \noindent \hrule \Gvsp   
\hbox{ $G=\Z / (  5 )$ 
} 
\Gvsp \hrule \Gvsp 
\vbox{ 
\hbox{  $ 2\,A\sb{9}  $ , 
  $ A\sb{9}\, +\, 2\,A\sb{4}\, +\, A\sb{1}  $ , 
  $ A\sb{9}\, +\, 2\,A\sb{4}  $ , 
  $ 4\,A\sb{4}\, +\, 2\,A\sb{1}  $ , 
  $ 4\,A\sb{4}\, +\, A\sb{1}  $ , 
  $ 4\,A\sb{4}  $ } 
} 
\Gvsp \hrule }}\vrule 
\par
\medskip 
\vrule  \hbox{\vbox{\offinterlineskip  \noindent \hrule \Gvsp   
\hbox{ $G=\Z / (  6 )$ 
} 
\Gvsp \hrule \Gvsp 
\vbox{ 
\hbox{  $ A\sb{11}\, +\, A\sb{5}\, +\, 2\,A\sb{1}  $ , 
  $ A\sb{11}\, +\, A\sb{3}\, +\, 2\,A\sb{2}  $ , 
  $ A\sb{11}\, +\, 2\,A\sb{2}\, +\, 3\,A\sb{1}  $ , 
  $ A\sb{11}\, +\, 2\,A\sb{2}\, +\, 2\,A\sb{1}  $ , 
  $ 3\,A\sb{5}\, +\, A\sb{3}  $ , 
  $ 3\,A\sb{5}\, +\, 2\,A\sb{1}  $ , 
} 
\GGvsp 
\hbox{  $ 2\,A\sb{5}\, +\, A\sb{3}\, +\, 2\,A\sb{2}\, +\, A\sb{1}  $ , 
  $ 2\,A\sb{5}\, +\, A\sb{3}\, +\, 2\,A\sb{2}  $ , 
  $ 2\,A\sb{5}\, +\, 2\,A\sb{2}\, +\, 3\,A\sb{1}  $ , 
  $ 2\,A\sb{5}\, +\, 2\,A\sb{2}\, +\, 2\,A\sb{1}  $ } 
} 
\Gvsp \hrule }}\vrule 
\par
\medskip 
\vrule  \hbox{\vbox{\offinterlineskip  \noindent \hrule \Gvsp   
\hbox{ $G=\Z / (  7 )$ 
} 
\Gvsp \hrule \Gvsp 
\vbox{ 
\hbox{  $ 3\,A\sb{6}  $ } 
} 
\Gvsp \hrule }}\vrule
\par 
\medskip 
\vrule  \hbox{\vbox{\offinterlineskip  \noindent \hrule \Gvsp   
\hbox{ $G=\Z / (  8 )$ 
} 
\Gvsp \hrule \Gvsp 
\vbox{ 
\hbox{  $ 2\,A\sb{7}\, +\, A\sb{3}\, +\, A\sb{1}  $ } 
} 
\Gvsp \hrule }}\vrule 
\par
\medskip 
\vrule  \hbox{\vbox{\offinterlineskip  \noindent \hrule \Gvsp   
\hbox{ $G=\Z / (  2 ) \times  \Z / (  2 )$ 
} 
\Gvsp \hrule \Gvsp 
\vbox{ 
\hbox{  $ D\sb{10}\, +\, A\sb{5}\, +\, 3\,A\sb{1}  $ , 
  $ D\sb{10}\, +\, 2\,A\sb{3}\, +\, 2\,A\sb{1}  $ , 
  $ D\sb{10}\, +\, A\sb{3}\, +\, 4\,A\sb{1}  $ , 
  $ D\sb{10}\, +\, 6\,A\sb{1}  $ , 
  $ 2\,D\sb{8}\, +\, 2\,A\sb{1}  $ , 
  $ D\sb{8}\, +\, D\sb{6}\, +\, A\sb{3}\, +\, A\sb{1}  $ , 
} 
\GGvsp 
\hbox{  $ D\sb{8}\, +\, D\sb{6}\, +\, 3\,A\sb{1}  $ , 
  $ D\sb{8}\, +\, D\sb{4}\, +\, A\sb{3}\, +\, 2\,A\sb{1}  $ , 
  $ D\sb{8}\, +\, D\sb{4}\, +\, 4\,A\sb{1}  $ , 
  $ D\sb{8}\, +\, A\sb{5}\, +\, A\sb{3}\, +\, 2\,A\sb{1}  $ , 
  $ D\sb{8}\, +\, A\sb{5}\, +\, 4\,A\sb{1}  $ , 
} 
\GGvsp 
\hbox{  $ D\sb{8}\, +\, 2\,A\sb{3}\, +\, 3\,A\sb{1}  $ , 
  $ D\sb{8}\, +\, A\sb{3}\, +\, 5\,A\sb{1}  $ , 
  $ D\sb{8}\, +\, 7\,A\sb{1}  $ , 
  $ 3\,D\sb{6}  $ , 
  $ 2\,D\sb{6}\, +\, D\sb{4}\, +\, A\sb{1}  $ , 
  $ 2\,D\sb{6}\, +\, 2\,A\sb{3}  $ , 
  $ 2\,D\sb{6}\, +\, A\sb{3}\, +\, 2\,A\sb{1}  $ , 
} 
\GGvsp 
\hbox{  $ 2\,D\sb{6}\, +\, 4\,A\sb{1}  $ , 
  $ D\sb{6}\, +\, 2\,D\sb{4}\, +\, 2\,A\sb{1}  $ , 
  $ D\sb{6}\, +\, D\sb{4}\, +\, 2\,A\sb{3}\, +\, A\sb{1}  $ , 
  $ D\sb{6}\, +\, D\sb{4}\, +\, A\sb{3}\, +\, 3\,A\sb{1}  $ , 
  $ D\sb{6}\, +\, D\sb{4}\, +\, 5\,A\sb{1}  $ , 
  $ D\sb{6}\, +\, 2\,A\sb{5}\, +\, 2\,A\sb{1}  $ , 
} 
\GGvsp 
\hbox{  $ D\sb{6}\, +\, A\sb{5}\, +\, 2\,A\sb{3}\, +\, A\sb{1}  $ , 
  $ D\sb{6}\, +\, A\sb{5}\, +\, A\sb{3}\, +\, 3\,A\sb{1}  $ , 
  $ D\sb{6}\, +\, A\sb{5}\, +\, 5\,A\sb{1}  $ , 
  $ D\sb{6}\, +\, 3\,A\sb{3}\, +\, 2\,A\sb{1}  $ , 
  $ D\sb{6}\, +\, 2\,A\sb{3}\, +\, 4\,A\sb{1}  $ , 
} 
\GGvsp 
\hbox{  $ D\sb{6}\, +\, A\sb{3}\, +\, 6\,A\sb{1}  $ , 
  $ D\sb{6}\, +\, 8\,A\sb{1}  $ , 
  $ 4\,D\sb{4}  $ , 
  $ 3\,D\sb{4}\, +\, 3\,A\sb{1}  $ , 
  $ 2\,D\sb{4}\, +\, 2\,A\sb{3}\, +\, 2\,A\sb{1}  $ , 
  $ 2\,D\sb{4}\, +\, A\sb{3}\, +\, 4\,A\sb{1}  $ , 
  $ 2\,D\sb{4}\, +\, 6\,A\sb{1}  $ , 
} 
\GGvsp 
\hbox{  $ D\sb{4}\, +\, 2\,A\sb{5}\, +\, 3\,A\sb{1}  $ , 
  $ D\sb{4}\, +\, A\sb{5}\, +\, 2\,A\sb{3}\, +\, 2\,A\sb{1}  $ , 
  $ D\sb{4}\, +\, A\sb{5}\, +\, A\sb{3}\, +\, 4\,A\sb{1}  $ , 
  $ D\sb{4}\, +\, A\sb{5}\, +\, 6\,A\sb{1}  $ , 
  $ D\sb{4}\, +\, 3\,A\sb{3}\, +\, 3\,A\sb{1}  $ , 
} 
\GGvsp 
\hbox{  $ D\sb{4}\, +\, 2\,A\sb{3}\, +\, 5\,A\sb{1}  $ , 
  $ D\sb{4}\, +\, A\sb{3}\, +\, 7\,A\sb{1}  $ , 
  $ D\sb{4}\, +\, 9\,A\sb{1}  $ , 
  $ A\sb{7}\, +\, A\sb{5}\, +\, A\sb{3}\, +\, 3\,A\sb{1}  $ , 
  $ A\sb{7}\, +\, A\sb{5}\, +\, 5\,A\sb{1}  $ , 
  $ A\sb{7}\, +\, 2\,A\sb{3}\, +\, 4\,A\sb{1}  $ , 
} 
\GGvsp 
\hbox{  $ A\sb{7}\, +\, A\sb{3}\, +\, 6\,A\sb{1}  $ , 
  $ A\sb{7}\, +\, 8\,A\sb{1}  $ , 
  $ 2\,A\sb{5}\, +\, 2\,A\sb{3}\, +\, 2\,A\sb{1}  $ , 
  $ 2\,A\sb{5}\, +\, A\sb{3}\, +\, 4\,A\sb{1}  $ , 
  $ 2\,A\sb{5}\, +\, 6\,A\sb{1}  $ , 
  $ A\sb{5}\, +\, 3\,A\sb{3}\, +\, 3\,A\sb{1}  $ , 
} 
\GGvsp 
\hbox{  $ A\sb{5}\, +\, 2\,A\sb{3}\, +\, 5\,A\sb{1}  $ , 
  $ A\sb{5}\, +\, A\sb{3}\, +\, 7\,A\sb{1}  $ , 
  $ A\sb{5}\, +\, 9\,A\sb{1}  $ , 
  $ 4\,A\sb{3}\, +\, 4\,A\sb{1}  $ , 
  $ 3\,A\sb{3}\, +\, 6\,A\sb{1}  $ , 
  $ 2\,A\sb{3}\, +\, 8\,A\sb{1}  $ , 
  $ A\sb{3}\, +\, 10\,A\sb{1}  $ , 
} 
\GGvsp 
\hbox{  $ 12\,A\sb{1}  $ } 
} 
\Gvsp \hrule }}\vrule
\par
\medskip 
\vrule  \hbox{\vbox{\offinterlineskip  \noindent \hrule \Gvsp   
\hbox{ $G=\Z / (  4 ) \times  \Z / (  2 )$ 
} 
\Gvsp \hrule \Gvsp 
\vbox{ 
\hbox{  $ 2\,A\sb{7}\, +\, 4\,A\sb{1}  $ , 
  $ A\sb{7}\, +\, 3\,A\sb{3}\, +\, 2\,A\sb{1}  $ , 
  $ A\sb{7}\, +\, 2\,A\sb{3}\, +\, 4\,A\sb{1}  $ , 
  $ 5\,A\sb{3}\, +\, 2\,A\sb{1}  $ , 
  $ 4\,A\sb{3}\, +\, 4\,A\sb{1}  $ } 
} 
\Gvsp \hrule }}\vrule 
\par
\medskip 
\vrule  \hbox{\vbox{\offinterlineskip  \noindent \hrule \Gvsp   
\hbox{ $G=\Z / (  6 ) \times  \Z / (  2 )$ 
} 
\Gvsp \hrule \Gvsp 
\vbox{ 
\hbox{  $ 3\,A\sb{5}\, +\, 3\,A\sb{1}  $ } 
} 
\Gvsp \hrule }}\vrule 
\par
\medskip 
\vrule  \hbox{\vbox{\offinterlineskip  \noindent \hrule \Gvsp   
\hbox{ $G=\Z / (  3 ) \times  \Z / (  3 )$ 
} 
\Gvsp \hrule \Gvsp 
\vbox{ 
\hbox{  $ 2\,A\sb{5}\, +\, 4\,A\sb{2}  $ , 
  $ A\sb{5}\, +\, 6\,A\sb{2}  $ , 
  $ 8\,A\sb{2}  $ } 
} 
\Gvsp \hrule }}\vrule 
\par
\medskip
\vrule  \hbox{\vbox{\offinterlineskip  \noindent \hrule \Gvsp   
\hbox{ $G=\Z / (  4 ) \times  \Z / (  4 )$ 
} 
\Gvsp \hrule \Gvsp 
\vbox{ 
\hbox{  $ 6\,A\sb{3}  $ } 
} 
\Gvsp \hrule }}\vrule 
\medskip 
\vfill\eject 
%%%%%%%%%%%%%%%%%%%%%%%%%%%%%%%
%
% The end of Table 2
%
%%%%%%%%%%%%%%%%%%%%%%%%%%%%%%%%
}

\end{document}